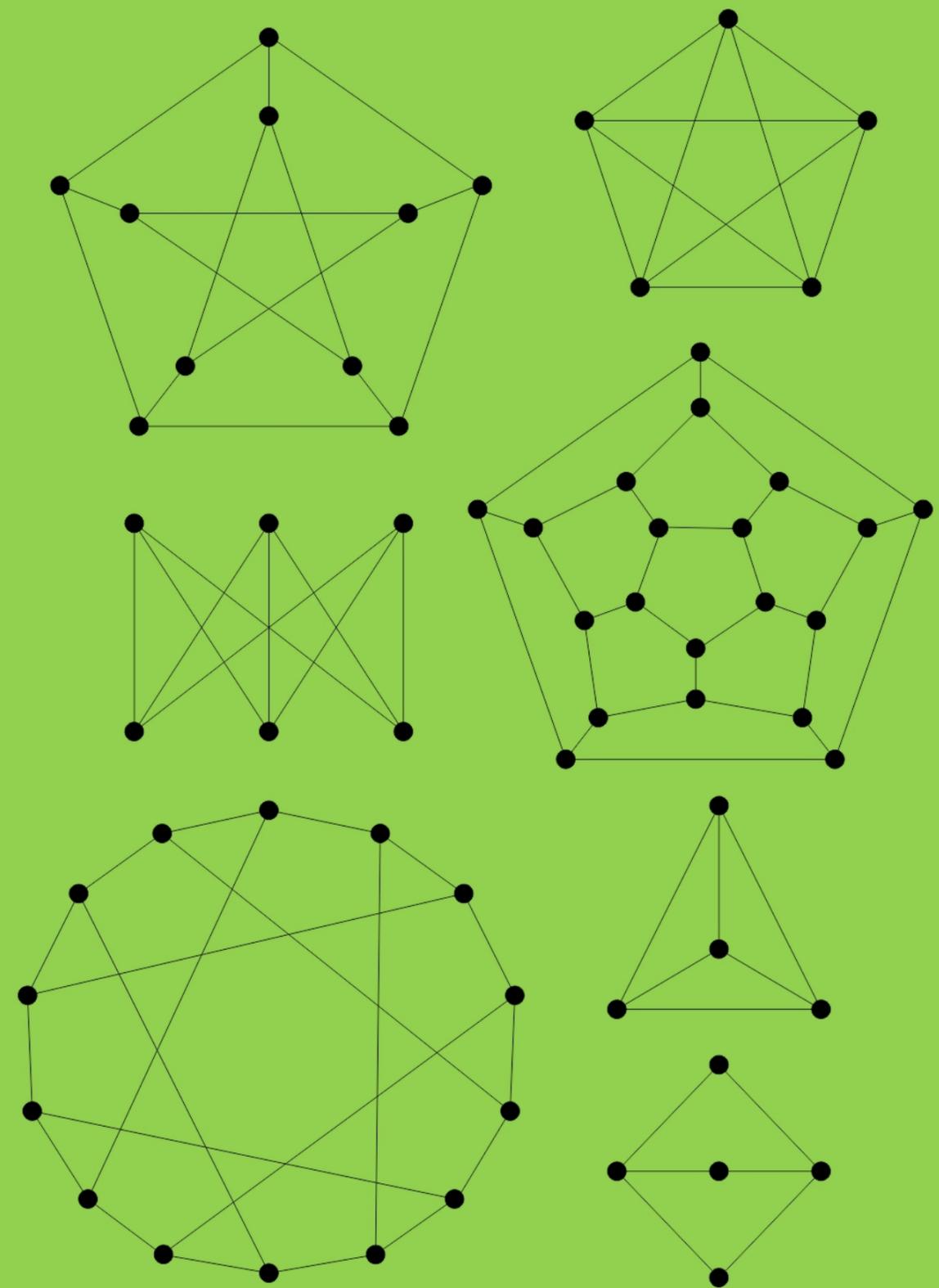

Պ.Ա. Պետրոսյան, Վ.Վ. Մկրտչյան, Ռ.Ռ. Քամալյան

# ԳՐԱՖՆԵՐԻ ՏԵՍՈՒԹՅՈՒՆ

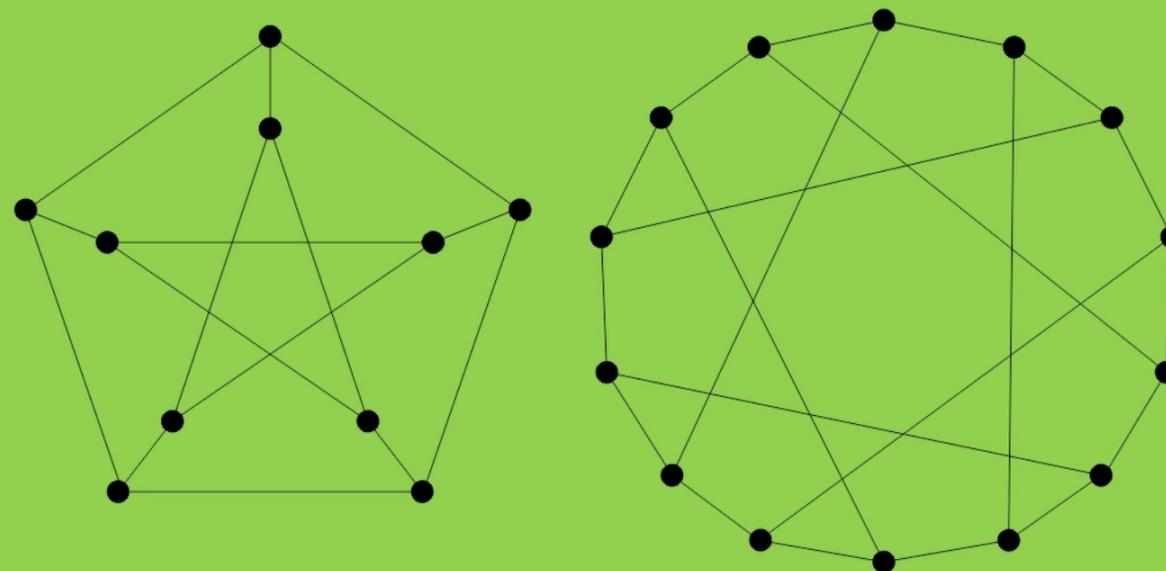

ԵՐԵՎԱՆԻ ՊԵՏԱԿԱՆ ՀԱՄԱԼՍԱՐԱՆ

Պ.Ա. Պետրոսյան, Վ.Վ. Մկրտչյան, Ռ.Ռ. Քամալյան

# ԳՐԱՖՆԵՐԻ ՏԵՍՈՒԹՅՈՒՆ

ՈՒՍՈՒՄՆԱՄԵԹՈԴԱԿԱՆ ՁԵՌՆԱՐԿ

ԵՐԵՎԱՆ
ԵՊՀ ՀՐԱՏԱՐԱԿՉՈՒԹՅՈՒՆ
2015

*Հրատարակության է երաշխավորել ԵՊՀ Ինֆորմատիկայի և կիրառական մաթեմատիկայի ֆակուլտետի խորհուրդը*

Ուսումնամեթոդական ձեռնարկն ընդգրկում է ԵՊՀ Ինֆորմատիկայի և կիրառական մաթեմատիկայի ֆակուլտետի «Գրաֆների տեսություն» դասընթացի նյութը: Ձեռնարկը պարունակում է բազմաթիվ արդյունքներ, որոնք օգտակար կլինեն ինչպես բակալավրիատի ուսանողների, այնպես էլ մագիստրանտների և ասպիրանտների համար: Ձեռնարկի նյութը կարող է ուսումնասիրվել նան սեմինար պարապմունքների ընթացքում:

Նախատեսված է ԵՊՀ Ինֆորմատիկայի և կիրառական մաթեմատիկայի ֆակուլտետի ուսանողների համար:

# ԲՈՎԱՆԴԱԿՈՒԹՅՈՒՆ

# ՆԱԽԱԲԱՆ

Գրաֆների տեսությունը դիսկրետ մաթեմատիկայի հայտնի և արդի ճյուղերից մեկն է: Այն սկիզբ է առել 1736 թ. Լ. Էյլերի կողմից դիտարկված և լուծված «Քյոնիգսբերգյան կամուրջների» հանրահայտ խնդրից: Գրաֆների տեսության հետագա զարգացումը ցույց տվեց, որ այն սերտորեն կապված է մաթեմատիկայի մի շարք բաժինների հետ, որոնցից են՝ խմբերի տեսությունը, ավտոմատների տեսությունը, մատրիցների տեսությունը, կոմբինատոր անալիզը, հավանականությունների տեսությունը, տոպոլոգիան, բարդության տեսությունը և այլն: Մյուս կողմից, բազմաթիվ են նաև այն գիտությունները, որոնցում գրաֆների տեսությունը լայնորեն և արդյունավետ կիրառվում է. ֆիզիկա, քիմիա, կենսաբանություն, գենետիկա, տնտեսագիտություն, հոգեբանություն, լեզվաբանություն և այլն: Այս ամենը կարևոր են դարձնում «Գրաֆների տեսություն» դասընթացի ուսումնասիրությունը:

Ուսումնամեթոդական ձեռնարկը գրված է Երևանի պետական համալսարանի Ինֆորմատիկայի և կիրառական մաթեմատիկայի ֆակուլտետում՝ «Գրաֆների տեսություն» դասընթացի՝ հեղինակների կարդացած դասախոսությունների հիման վրա և ընդգրկում է նշված դասընթացի ուսումնական ծրագրով նախատեսված նյութը: Ձեռնարկը պարունակում է բազմաթիվ արդյունքներ, որոնք օգտակար կլինեն ինչպես բակալավրիատի ուսանողների, այնպես էլ մագիստրանտների և ասպիրանտների համար: Ձեռնարկի նյութը կարող է ուսումնասիրվել նաև սեմինար պարապմունքների ընթացքում: Ձեռնարկը բաղկացած է ութ գլուխից: Այդ գլուխներն ընդգրկում են գրաֆների տեսության դասընթացի երկրորդ կուրսում կարդացվող թեմաները՝ գրաֆներ, գրաֆների տրման եղանակներ, գործողություններ գրաֆների հետ, երկկողմանի գրաֆներ, ծառեր, կապակցվածություն, էյլերյան և համիլտոնյան գրաֆներ, ֆակտորներ, զուգակցումներ, անկախ բազմություններ և ծածկույթներ, հարթ գրաֆներ և գրաֆների ներկումներ, որոնք հիմք են հանդիսանում հետագայում դասավանդվող մի շարք առարկաների համար:

Գրքում օգտագործվող նշանակումների համակարգը հիմնականում վերցված է D. B. West, Introduction to Graph Theory, Prentice-Hall, New Jersey, 2001, J. A. Bondy, U. S. R. Murty, Graph Theory, Springer, 2008 և Ֆ. Харари, Теория графов, Пер. с англ.-М.: Мир, 1973 գրքերից: Որոշ արդյունքների ապացույցներ իրենցից ներկայացնում են հոդվածներում բերված ապացույցների վերամշակված տարբերակներ:

Հեղինակները շնորհակալություն են հայտնում Ռ.Ն. Տոնոյանին, Հ.Յ. Հակոբյանին, ԵՊՀ դիսկրետ մաթեմատիկայի և տեսական ինֆորմատիկայի ամբիոնի աշխատա-կիցներին, ՀՀ ԳԱԱ Ինֆորմատիկայի և ավտոմատացման պրոբլեմների ինստիտուտի «Գրաֆների տեսության էքստրեմալ խնդիրների» լաբորատորիայի աշխատակիցներին՝ ձեռնարկում ներկայացված նյութի բովանդակության և շարադրման եղանակի հետ կապված հարցերում օգտակար առաջարկությունների և դիտողությունների համար: Վերջում հեղինակները շնորհակալություն են հայտնում այն բոլոր ուսանողներին, որոնք իրենց հարցերով նպաստեցին շարադրանքի բարելավմանը: Հեղինակները սիրով կընդունեն բոլոր առաջարկները և դիտողությունները:

# ՀԻՄՆԱԿԱՆ ՆՇԱՆԱԿՈՒՄՆԵՐ

$\mathbb{N}$ - բնական թվերի բազմություն

$\mathbb{Z}_+$ - ամբողջ ոչ բացասական թվերի բազմություն

$\mathbb{R}$ - իրական թվերի բազմություն

$\lfloor a \rfloor$ - ամենամեծ ամբողջ թիվը, որը չի գերազանցում $a$ իրական թվին

$\lceil a \rceil$ - ամենափոքր ամբողջ թիվը, որը փոքր չէ $a$ իրական թվից

$\binom{n}{k}$ - $n$ տարրերից $k$ զուգորդությունների քանակ

$A \subseteq B$ - $A$-ն $B$-ի ենթաբազմություն է

$A \cup B$ - $A$ և $B$ բազմությունների միավորում

$A \cap B$ - $A$ և $B$ բազմությունների հատում

$A \setminus B$ - $A$ և $B$ բազմությունների տարբերություն

$A \times B$ - $A$ և $B$ բազմությունների դեկարտյան արտադրյալ

$|A|$ - $A$ բազմության հզորություն

$V(G)$ - $G$ գրաֆի գագաթների բազմություն

$E(G)$ - $G$ գրաֆի կողերի բազմություն

$uv$ - $u$ և $v$ գագաթները միացնող կող

$d_G(v)$ - $v$ գագաթի աստիճանը $G$ գրաֆում

$\delta(G)$ - $G$ գրաֆի նվազագույն աստիճան

$\Delta(G)$ - $G$ գրաֆի առավելագույն աստիճան

$N_G(v)$ - $G$ գրաֆում $v$ գագաթի շրջակայք

$N_G(S)$ - $G$ գրաֆում $S$-ին հարևան գագաթների բազմություն

$H \subseteq G$ - $H$-ը $G$ գրաֆի ենթագրաֆ է

$G[S]$ - $G$ գրաֆի գագաթների $S$ բազմությամբ ծնված ենթագրաֆ

$G - v$ - $G$ գրաֆից $v$ գագաթի հեռացումից առաջացած գրաֆ

$G - S$ - $G$ գրաֆից $S$ գագաթների բազմության հեռացումից առաջացած գրաֆ

$G - e$ - $G$ գրաֆից $e$ կողի հեռացումից առաջացած գրաֆ

$G + e$ - $G$ գրաֆին $e$ կողի ավելացումից առաջացած գրաֆ

$d_G(u, v)$ - $G$ գրաֆում $u$ և $v$ գագաթների միջև հեռավորություն

$r(G)$ - $G$ գրաֆի շառավիղ

$d(G)$ - $G$ գրաֆի տրամագիծ

$A(G)$ - $G$ գրաֆի հարևանության մատրից

$B(G)$ - $G$ գրաֆի կցության մատրից

$a(G)$ - $G$ գրաֆի անտառների տրոհման թիվ

$c(G)$ - $G$ գրաֆի կապակցվածության բաղադրիչների քանակ

$o(G)$ - $G$ գրաֆի կապակցվածության այն բաղադրիչների քանակը, որոնք պարունակում են կենտ թվով գագաթներ

$cl(G)$ - $G$ գրաֆի համիլտոնյան փակում

$cr(G)$ - $G$ գրաֆի խաչումների թիվ

$cyc(G)$ - $G$ գրաֆի ցիկլոմատիկ թիվ

$\alpha(G)$ - $G$ գրաֆում առավելագույն անկախ բազմության հզորություն

$\alpha'(G)$ - $G$ գրաֆում առավելագույն զուգակցման հզորություն

$\beta(G)$ - $G$ գրաֆում նվազագույն գագաթային ծածկույթի հզորություն

$\beta'(G)$ - $G$ գրաֆում նվազագույն կողային ծածկույթի հզորություն

$\gamma(G)$ - $G$ գրաֆի սեռ

$\varkappa(G)$ - $G$ գրաֆի կապակցվածություն

$\lambda(G)$ - $G$ գրաֆի կողային կապակցվածություն

$\tau(G)$ - $G$ գրաֆի կոշտություն

$\chi(G)$ - $G$ գրաֆի քրոմատիկ թիվ

$\chi'(G)$ - $G$ գրաֆի քրոմատիկ ինդեքս

$\chi''(G)$ - $G$ գրաֆի տոտալ քրոմատիկ թիվ

$\omega(G)$ - $G$ գրաֆի խտություն

$C_n$ - $n$ գագաթ ($n \geq 3$) ունեցող պարզ ցիկլ

$K_n$ - $n$ գագաթ ունեցող լրիվ գրաֆ

$K_{m,n}$ - լրիվ երկկողմանի գրաֆ, որի մի կողմը պարունակում է $m$ գագաթ, իսկ մյուսը՝ $n$ գագաթ

$Q_n$ - $n$-չափանի խորանարդ

$P_n$ - $n$ գագաթ ունեցող պարզ ճանապարհ

$S_n$ - $1, \ldots, n$ թվերի տեղադրությունների խումբ

# Գլուխ 1

# Գրաֆներ: Հիմնական սահմանումներ և պարզագույն հատկություններ

## § 1.1. Գրաֆի սահմանումը, տեսակները և տրման եղանակները

Դիցուք $V = \{v_1, \ldots, v_n\}$-ը ցանկացած ոչ դատարկ վերջավոր բազմություն է, և դիցուք $V^{(2)}$-ը $V$ բազմության տարրերի բոլոր ոչ կարգավոր զույգերի բազմությունն է: Նշենք, որ $|V^{(2)}| = \binom{n}{2}$: Ենթադրենք, որ $E \subseteq V^{(2)}$:

**Սահմանում 1.1.1:** $(V, E)$ կարգավոր զույգին կանվանենք *գրաֆ*, և այն կնշանակենք $G$-ով:

$G = (V, E)$ գրաֆի $V$ բազմության տարրերին կանվանենք գրաֆի *գագաթներ*, իսկ $E$ բազմության տարրերին՝ *կողեր*: Եթե անհրաժեշտ է շեշտել, որ $V$-ն հանդիսանում է $G$ գրաֆի գագաթների բազմություն ($E$-ն հանդիսանում է $G$ գրաֆի կողերի բազմություն), ապա այդ դեպքում մենք կգրենք $V(G)$ ($E(G)$): Եթե $e \in E$ կողը $u, v \in V$ գագաթներից բաղկացած զույգն է, ապա այդ փաստը կգրենք $e = uv$-ով:

Դիցուք $G = (V, E)$ և $G' = (V', E')$ երկու գրաֆներ են:

**Սահմանում 1.1.2:** $G$ և $G'$ գրաֆներր կանվանենք *հավասար* և կգրենք $G = G'$ այն և միայն այն դեպքում, երբ $V = V'$ և $E = E'$:

Նշենք, որ գրաֆները կարելի է դիտարկել որպես հատուկ տիպի համասեռ բինար հարաբերություն, որի հենքային բազմությունը $V$-ն է: Հիշենք, որ $\alpha \subseteq V \times V$ բինար հարաբերությունը կոչվում է

- *ռեֆլեքսիվ*, եթե ցանկացած $v \in V$ համար $v\alpha v$,
- *անտիռեֆլեքսիվ*, եթե ցանկացած $v \in V$ համար $v\bar{\alpha}v$,
- *սիմետրիկ*, եթե ցանկացած $u, v \in V$ համար բավարարվում է հետևյալ պայմանը. եթե $u\alpha v$, ապա $v\alpha u$:

Նկատենք, որ $G = (V, E)$ գրաֆը կարող ենք դիտարկել որպես $\alpha \subseteq V \times V$ անտիռեֆլեքսիվ, սիմետրիկ բինար հարաբերություն, որտեղ ցանկացած $u, v \in V$ համար բավարարվում է հետևյալ պայմանը. $u\alpha v$ այն միայն այն դեպքում, երբ $uv \in E$:

Ստորև կդիտարկենք գրաֆների տրման մի քանի եղանակներ: Նախ նշենք, որ գրաֆը կարելի է տալ, նշելով նրա գագաթների և կողերի բազմությունները: Օրինակ, դիտարկենք $G = (V, E)$ գրաֆը, որտեղ $V = \{v_1, v_2, v_3, v_4, v_5, v_6, v_7\}$ և $E = \{v_1v_2, v_2v_3, v_3v_4, v_1v_4, v_2v_4, v_5v_6, v_6v_7, v_5v_7\}$:

Մեկ այլ եղանակ է գրաֆների տրման երկրաչափական եղանակը, որի էությունը կայանում է հետևյալում. գրաֆի գագաթներին համապատասխանեցնում ենք հարթության կետեր (տարբեր գագաթներին համապատասխանում են տարբեր կետեր), և երկու կետեր միացվում են անընդհատ կորով, որը չի անցնում մեկ այլ գագաթին համապատասխանող կետով` այն և միայն այն դեպքում, երբ նրանց համապատասխանող գագաթները կող են կազմում գրաֆում: Օրինակ, վերը նշված գրաֆը կարելի է պատկերել հետևյալ կերպ.

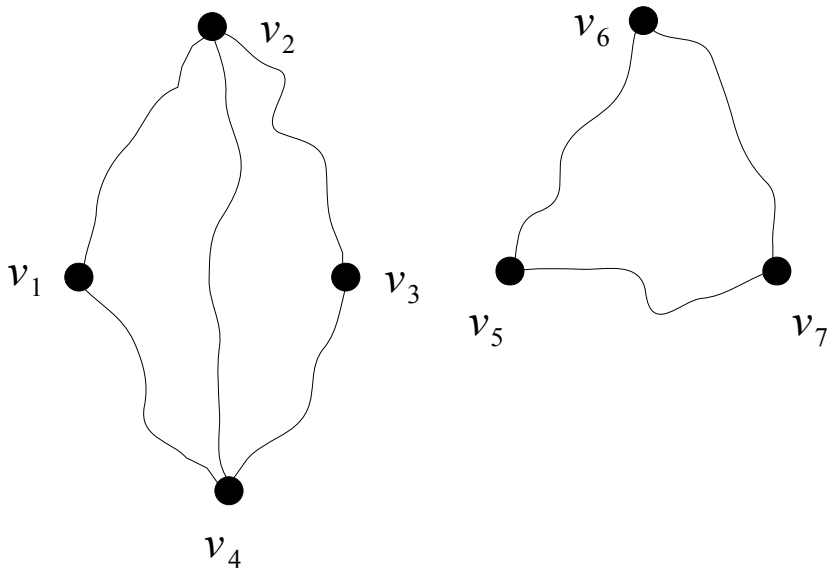

Նկ. 1.1.1

**Սահմանում 1.1.3:** $G$ գրաֆը կանվանենք *նշված* (կամ *համարակալված*), եթե այդ գրաֆի գագաթներին վերագրված են զույգ առ զույգ տարբեր նիշեր:

Գրաֆներ դիտարկելիս մեզ համար կարևոր է, թե որ գագաթներն են միացված կողով, և որոնք` ոչ: Հատուկ նշենք, որ կետերը միացնող կորի ձևը կարևոր չէ: Օրինակ վերը բերված գրաֆը կարելի է պատկերել նաև հետևյալ կերպ.

**G**

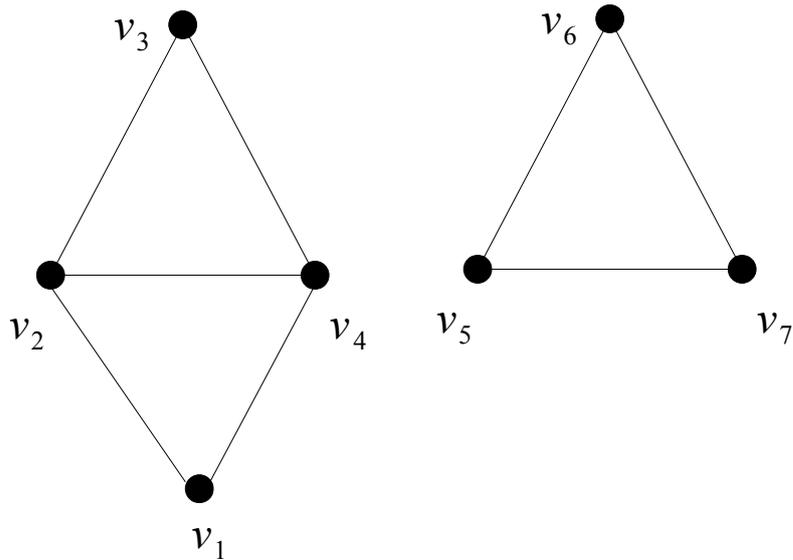

Նկ. 1.1.2

Գրաֆների տրման հաջորդ եղանակները նկարագրելու համար տանք մի քանի սահմանում:

Դիցուք $G = (V, E)$-ն գրաֆ է, $u, v \in V$ և $e, e' \in E$:

**Սահմանում 1.1.4:** $u$ և $v$ գագաթները կանվանենք *հարևան*, եթե $uv \in E$:

**Սահմանում 1.1.5:** $u$ գագաթին և $e$ կողին կանվանենք *կից*, եթե $u \in e$:

**Սահմանում 1.1.6:** $e$ և $e'$ տարբեր կողերը կանվանենք *հարևան*, եթե գոյություն ունի $v \in V$ այնպես, որ $v$-ն կից է $e$-ին և $e'$-ին:

Եթե $G = (V, E)$ գրաֆում $V = \{v_1, \ldots, v_n\}$ և $E = \{e_1, \ldots, e_m\}$, ապա այդ գրաֆին համապատասխանեցնենք $n \times n$ կարգի $A(G) = (a_{ij})_{n \times n}$ մատրիցը հետևյալ կերպ.

$$a_{ij} = \begin{cases} 1, & \text{եթե } v_i \text{ և } v_j \text{ հարևան են,} \\ 0, & \text{հակառակ դեպքում:} \end{cases}$$

$A(G)$ մատրիցը կանվանենք $G$ գրաֆի *հարևանության մատրից*. Նկատենք, որ ցանկացած $i$-ի համար $(1 \leq i \leq n)$ $a_{ii} = 0$, և ցանկացած $i, j$-ի համար $(1 \leq i, j \leq n)$ $a_{ij} = a_{ji}$: Նկ. 1.1.1-ում բերված $G$ գրաֆի հարևանության մատրիցը կլինի

$$A(G) = \begin{pmatrix} 0 & 1 & 0 & 1 & 0 & 0 & 0 \\ 1 & 0 & 1 & 1 & 0 & 0 & 0 \\ 0 & 1 & 0 & 1 & 0 & 0 & 0 \\ 1 & 1 & 1 & 0 & 0 & 0 & 0 \\ 0 & 0 & 0 & 0 & 0 & 1 & 1 \\ 0 & 0 & 0 & 0 & 1 & 0 & 1 \\ 0 & 0 & 0 & 0 & 1 & 1 & 0 \end{pmatrix}$$

Նշենք, որ զրոներից և մեկերից կազմված $n \times n$ կարգի յուրաքանչյուր $A$ մատրից, որը բավարարում է հետևյալ երկու պայմաններին. ցանկացած $i$-ի համար ($1 \leq i \leq n$) $a_{ii} = 0$, և ցանկացած $i,j$-ի համար ($1 \leq i,j \leq n$) $a_{ij} = a_{ji}$, հանդիսանում է համարակալված գագաթներով որևէ գրաֆի հարևանության մատրից:

Եթե $G = (V, E)$ գրաֆում $V = \{v_1, \ldots, v_n\}$ և $E = \{e_1, \ldots, e_m\}$, ապա այդ գրաֆին համապատասխանեցնենք $n \times m$ կարգի $B(G) = (b_{ij})_{n \times m}$ մատրիցը հետևյալ կերպ.

$$b_{ij} = \begin{cases} 1, & \text{եթե } v_i \text{ և } e_j \text{ կից են,} \\ 0, & \text{հակառակ դեպքում:} \end{cases}$$

$B(G)$ մատրիցը կանվանենք $G$ գրաֆի *կցության մատրից*: Նկատենք, որ կցության մատրիցի սյուները զույգ առ զույգ տարբեր են և յուրաքանչյուր սյուն պարունակում է ճիշտ երկու $1$: Նկ. 1.1.1-ում բերված $G$ գրաֆի կցության մատրիցը կլինի

$$B(G) = \begin{pmatrix} 1 & 0 & 0 & 1 & 0 & 0 & 0 & 0 \\ 1 & 1 & 0 & 0 & 1 & 0 & 0 & 0 \\ 0 & 1 & 1 & 0 & 0 & 0 & 0 & 0 \\ 0 & 0 & 1 & 1 & 1 & 0 & 0 & 0 \\ 0 & 0 & 0 & 0 & 0 & 1 & 0 & 1 \\ 0 & 0 & 0 & 0 & 0 & 1 & 1 & 0 \\ 0 & 0 & 0 & 0 & 0 & 0 & 1 & 1 \end{pmatrix},$$

որտեղ ենթադրված է, որ $e_1 = v_1v_2, e_2 = v_2v_3, e_3 = v_3v_4, e_4 = v_1v_4, e_5 = v_2v_4, e_6 = v_5v_6, e_7 = v_6v_7, e_8 = v_5v_7$:

Նշենք, որ զրոներից և մեկերից կազմված $n \times m$ կարգի յուրաքանչյուր $B$ մատրից, որի սյուները զույգ առ զույգ տարբեր են և յուրաքանչյուր սյուն պարունակում է ճիշտ երկու $1$, հանդիսանում է համարակալված գագաթներով և կողերով որևէ գրաֆի կցության մատրից:

Գրաֆների ևս մի ներկայացումը, որը մենք կդիտարկենք, դա գրաֆի *գագաթների հարևանության ցուցակների* միջոցով ներկայացումն է: Եթե $G$ գրաֆում $V(G) = \{v_1, \ldots, v_n\}$, ապա դիտարկենք $n$ երկարությամբ զանգված, որի $i$-րդ բաղադրիչն իրենից ներկայացնում է $v_i$ գագաթին հարևան գագաթների ցուցակը գրված ինչ-որ մի կարգով: Օրինակ, նկ. 1.1.1-ում պատկերված գրաֆի գագաթների հարևանության ցուցակներով ներկայացումը կլինի.

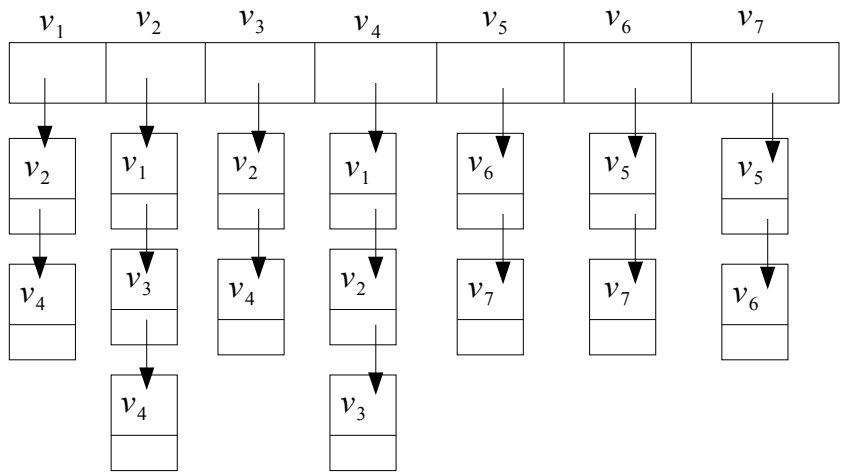

Նկ. 1.1.3

Նշենք, որ գրաֆներից բացի, մենք կդիտարկենք գրաֆների մի քանի ընդհանրացումներ։ Դրանցից առաջինն են *մուլտիգրաֆները*։ Մուլտիգրաֆի դեպքում $(V, E)$ կարգավոր զույգի մեջ $E$-ն $V^{(2)}$ բազմության մուլտիենթաբազմություն է կամ ցուցակ։ Սա նշանակում է, որ $V^{(2)}$ բազմության զույգերը $E$-ում կարող են հանդիպել մեկ անգամից ավելի։ Երկրաչափորեն, մուլտիգրաֆները պատկերվում են գրաֆներին համանման եղանակով, միայն թե գագաթներին համապատասխանող կետերը միացվում են այնքան գծերով, որքան անգամ $E$-ում կրկնվում է համապատասխան զույգը։ Օրինակ, եթե դիտարկենք $G = (V, E)$ մուլտիգրաֆը, որտեղ $V = \{v_1, v_2, v_3, v_4\}$ և $E = \langle v_1v_2, v_1v_2; v_1v_4, v_1v_4; v_2v_4, v_2v_4; v_2v_3; v_3v_4 \rangle$, ապա այն երկրաչափորեն կպատկերենք հետևյալ կերպ։

*G*

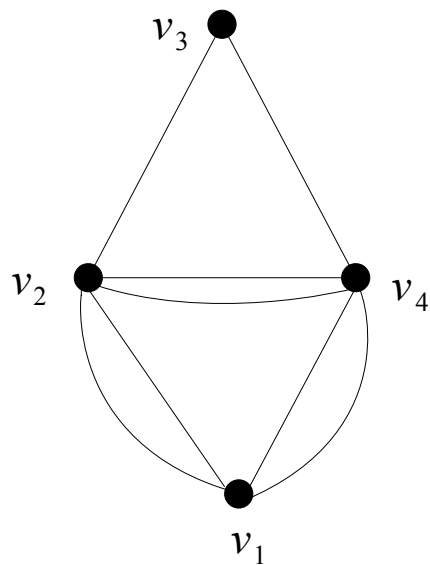

Նկ. 1.1.4

Գրաֆների հաջորդ ընդհանրացումն են *պսևդոգրաֆները*: Պսևդոգրաֆի դեպքում $(V, E)$ կարգավոր զույգի մեջ $E$-ն $V'^{(2)}$ բազմության մուլտիենթաբազմություն է կամ ցուցակ, որտեղ $V'^{(2)} = V^{(2)} \cup \{vv : v \in V\}$: Նշենք, որ $vv$ տիպի զույգերին ընդունված է անվանել $(V, E)$ պսևդոգրաֆի օղակներ: Երկրաչափորեն, պսևդոգրաֆները պատկերվում են մուլտիգրաֆներին համանման եղանակով, միայն թե օղակներին համապատասխանող անընդհատ գծերի սկիզբն ու վերջը համնկնում են: Օրինակ, եթե դիտարկենք $G = (V, E)$ պսևդոգրաֆը, որտեղ

$V = \{v_1, v_2, v_3, v_4\}$ և $E = \langle v_1v_1, v_1v_1; v_1v_2, v_1v_2; v_2v_3; v_3v_4; v_1v_4, v_1v_4; v_2v_4, v_2v_4; v_3v_3 \rangle$

ապա այն երկրաչափորեն կպատկերենք հետևյալ կերպ.

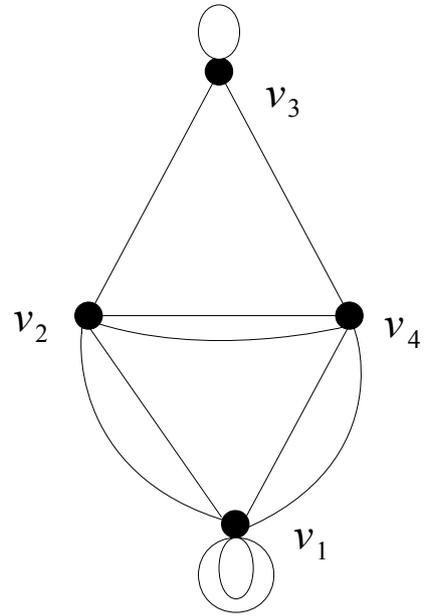

Նկ. 1.1.5

Այնուհետև կդիտարկենք *կողմնորոշված գրաֆները*, որոնք համապատասխանում են այնպիսի $(V, A)$ կարգավոր զույգերին, երբ $A$-ն $V$ բազմության տարրերի կրկնություններով բոլոր կարգավոր զույգերի մուլտիենթաբազմություն է կամ ցուցակ: Սա նշանակում է, որ $A$-ում կարող են հանդիպել $vv$ տիպի զույգեր: $a$-ից և $b$-ից կազմված կարգավոր զույգը կնշանակենք $(a, b)$ սիմվոլով, և այն երկրաչափորեն պատկերելիս, կնկարենք անընդհատ կոր, որը սկսում է $a$-ից և վերջանում $b$-ում, ընդ որում որպեսզի ընդգծենք, որ $b$-ն վերջն է, գծի վերջում կդնենք սլաք: Օրինակ, եթե դիտարկենք $D = (V, A)$ կողմնորոշված գրաֆը, որտեղ

$V = \{v_1, v_2, v_3\}$ և $A = \langle (v_1, v_2), (v_1, v_2); (v_2, v_1); (v_2, v_3); (v_1, v_3); (v_2, v_2) \rangle$, ապա այն երկրաչափորեն կպատկերենք հետևյալ կերպ.

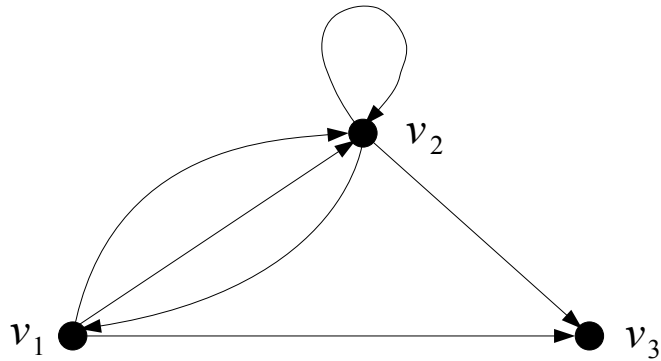

Նկ. 1.1.6

Եվ վերջապես, կդիտարկենք *հիպերգրաֆներ*, որոնք համապատասխանում են այն $(V, E)$ կարգավոր զույգերին, որոնցում $E$-ն $V$-ի ոչ դատարկ ենթաբազմություններից բաղկացած բազմության ենթաբազմություն է: Հիպերգրաֆի օրինակ է $H = (V, E)$ կարգավոր զույգը, որտեղ $V = \{v_1, v_2, v_3, v_4, v_5\}$ և $E = \{\{v_1, v_2\}, \{v_1, v_2, v_3\}, \{v_2, v_3, v_4\}, \{v_3\}, \{v_4, v_5\}\}$ (նկ. 1.1.7):

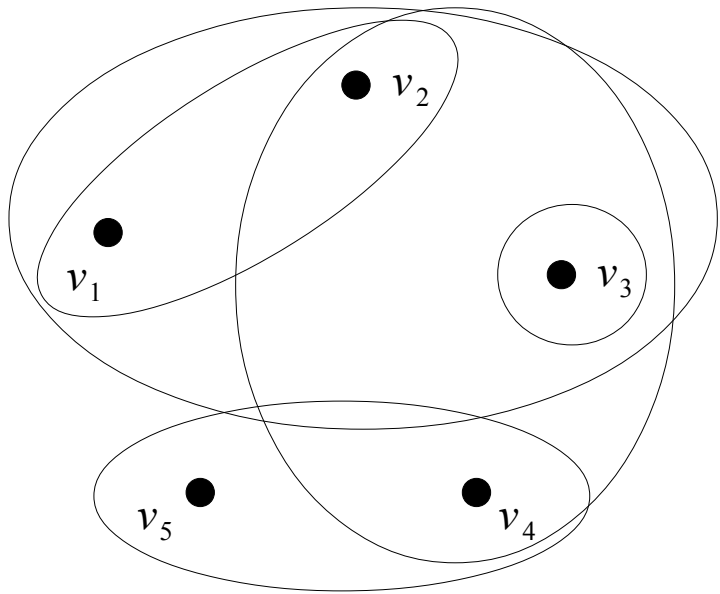

Նկ. 1.1.7

## § 1.2. Աստիճաններ, ենթագրաֆներ և ճանապարհներ

Դիտարկենք $G = (V, E)$ գրաֆը: $G$ գրաֆը կանվանենք $(n, m)$-*գրաֆ*, եթե $|V| = n$ և $|E| = m$: Եթե $S \subseteq V(G)$, ապա կատարենք հետևյալ նշանակումները.

$$N_G(S) = \{u \in V \setminus S: գոյություն ունի\, v \in S որ,\ uv \in E\},$$
$$\partial_G(S) = \{uv \in E: u \in S, v \in V \setminus S\}:$$

$G$ գրաֆում $v \in V$ գագաթի *շրջակայք* ասելով կհասկանանք $N_G(\{v\})$ բազմությունը։ Այն կրճատ կնշանակենք $N_G(v)$-ով։ Ավելին, $v$ գագաթին կից կողերի բազմությունը՝ $\partial_G(\{v\})$-ն կնշանակենք $\partial_G(v)$-ով։

**Սահմանում 1.2.1:** $G$ գրաֆում $v$ գագաթի *աստիճան*, որը կնշանակենք $d_G(v)$-ով կամ $d(v)$-ով, կանվանենք այդ գագաթին կից կողերի քանակը։ Պարզ է, որ $d_G(v) = |\partial_G(v)|$։

Օրինակ, նկ. 1.1.1-ում պատկերված $G$ գրաֆում $v_2$ գագաթի աստիճանը հավասար է երեքի։

$G$ գրաֆում $v$ գագաթը կանվանենք *մեկուսացված*, եթե $d_G(v) = 0$ և կանվանենք *կախված*, եթե $d_G(v) = 1$։ $G$ գրաֆի համար սահմանենք $\delta(G)$ և $\Delta(G)$ թվերը հետևյալ կերպ.

$$\delta(G) = min_{v \in V} d_G(v),\ \Delta(G) = max_{v \in V} d_G(v):$$

$\delta(G)$-ն կանվանենք $G$ գրաֆի *նվազագույն աստիճան*, իսկ $\Delta(G)$-ն՝ *առավելագույն աստիճան*։

Նկատենք, որ ցանկացած $G$ գրաֆում տեղի ունեն հետևյալ անհավասարությունները․

$$0 \leq \delta(G) \leq \Delta(G) \leq |V(G)| - 1:$$

**Թեորեմ 1.2.1 (Լ. Էյլեր):** Կամայական $G = (V, E)$ գրաֆում տեղի ունի

$$\sum_{v \in V(G)} d_G(v) = 2|E(G)|$$

հավասարությունը։

**Ապացույց:** Իրոք, քանի որ ցանկացած կող կից է երկու գագաթի, ապա $\sum_{v \in V(G)} d_G(v)$ գումարում այդ կողը հաշվվում է երկու անգամ, հետևաբար՝

$$\sum_{v \in V(G)} d_G(v) = 2|E(G)|: \quad \blacksquare$$

**Հետևանք 1.2.1:** Կամայական $G = (V, E)$ գրաֆում կենտ աստիճան ունեցող գագաթների քանակը զույգ է։

**Ապացույց:** Իրոք, համաձայն թեորեմ 1.2.1-ի

$$2|E(G)| = \sum_{v \in V(G)} d_G(v) = \sum_{d_G(v)-ն\ զույգ\ է} d_G(v) + \sum_{d_G(v)-ն\ կենտ\ է} d_G(v):$$

Նկատենք, որ հավասարության ձախ մասը, ինչպես նաև աջ մասում գտնվող առաջին գումարելին զույգ թվեր են, հետևաբար զույգ է նաև երկրորդ գումարելին, որտեղից հետևում է, որ զույգ է նաև կենտ աստիճան ունեցող գագաթների քանակը։ ∎

**Դիտողություն 1.2.1:** Նշենք, որ թեորեմ 1.2.1-ը և հետևանք 1.2.1-ը մնում են ճիշտ նաև մուլտիգրաֆների և պսևդոգրաֆների դեպքում, միայն թե պայմանավորվում ենք, որ օղակները պսևդոգրաֆի գագաթի աստիճանն ավելացնում են երկուսով։ Մասնավորապես, սա նշանակում է, որ նկ. 1.1.5-ում պատկերված գրաֆում $v_1$ գագաթի աստիճանը հավասար է ութի։

**Սահմանում 1.2.2:** $G$ գրաֆը կանվանենք *համասեռ* կամ *ռեգուլյար*, եթե $\delta(G) = \Delta(G)$ կամ որ նույնն է, որ եթե նրանում բոլոր գագաթների աստիճանները միևնույն թիվն է։ $G$ գրաֆը կանվանենք $r$-*համասեռ* կամ $r$-*ռեգուլյար*, եթե $\delta(G) = \Delta(G) = r$ ($r \in \mathbb{Z}_+$)։

Թեորեմ 1.2.1-ից անմիջապես հետևում է, որ

**Հետևանք 1.2.2:** Եթե $G$-ն $r$-համասեռ $(n, m)$-գրաֆ է, ապա
$$m = \frac{n \cdot r}{2}:$$

**Սահմանում 1.2.3:** 3-համասեռ գրաֆներին կանվանենք *խորանարդ գրաֆներ*։

Հետևանք 1.2.1-ից բխում է

**Հետևանք 1.2.3:** Խորանարդ գրաֆում գագաթների քանակը զույգ թիվ է։

Ստորև պատկերված են երկու խորանարդ գրաֆներ։

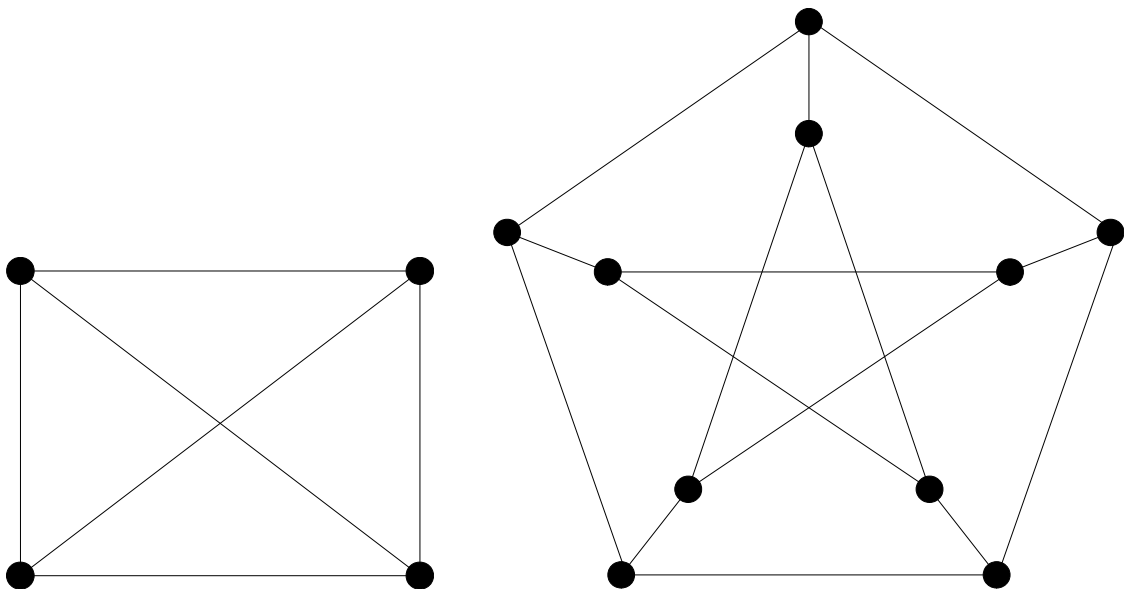

Նկ. 1.2.1

Նկարի աջ մասում պատկերված խորանարդ գրաֆը հայտնի է *Պետերսենի գրաֆ* անունով։

**Սահմանում 1.2.4:** $G$ գրաֆը կոչվում է *լրիվ*, եթե նրանում ցանկացած երկու գագաթ

հարևան են:

$n$ գագաթ ունեցող լրիվ գրաֆը կնշանակենք $K_n$-ով: $K_3$-ը կանվանենք *եռանկյուն*: Նկատենք, որ նկ. 1.2.1-ի ձախ մասում պատկերված է $K_4$-ը:

Դժվար չէ տեսնել, որ $K_n$-ը $(n-1)$-համասեռ գրաֆ է և

$$|E(K_n)| = \binom{n}{2} = \frac{n(n-1)}{2}:$$

**Սահմանում 1.2.5:** $G = (V, E)$ գրաֆը կանվանենք *$r$-կողմանի* ($r \in \mathbb{N}$), եթե $V$ բազմությունը հնարավոր է տրոհել $r$ ենթաբազմությունների այնպես, որ միևնույն ենթաբազմության գագաթները զույգ առ զույգ հարևան չեն: Եթե $r = 2$, ապա $r$-կողմանի գրաֆը կանվանենք *երկկողմանի*: Նկատենք, որ եթե $G = (V, E)$ գրաֆը երկկողմանի է, ապա $V$ բազմությունը հնարավոր է տրոհել երկու ենթաբազմությունների $V_1$ և $V_2$-ի այնպես, որ $G$ գրաֆի ցանկացած կող կից լինի մեկ գագաթի $V_1$-ից և մեկ գագաթի $V_2$-ից:

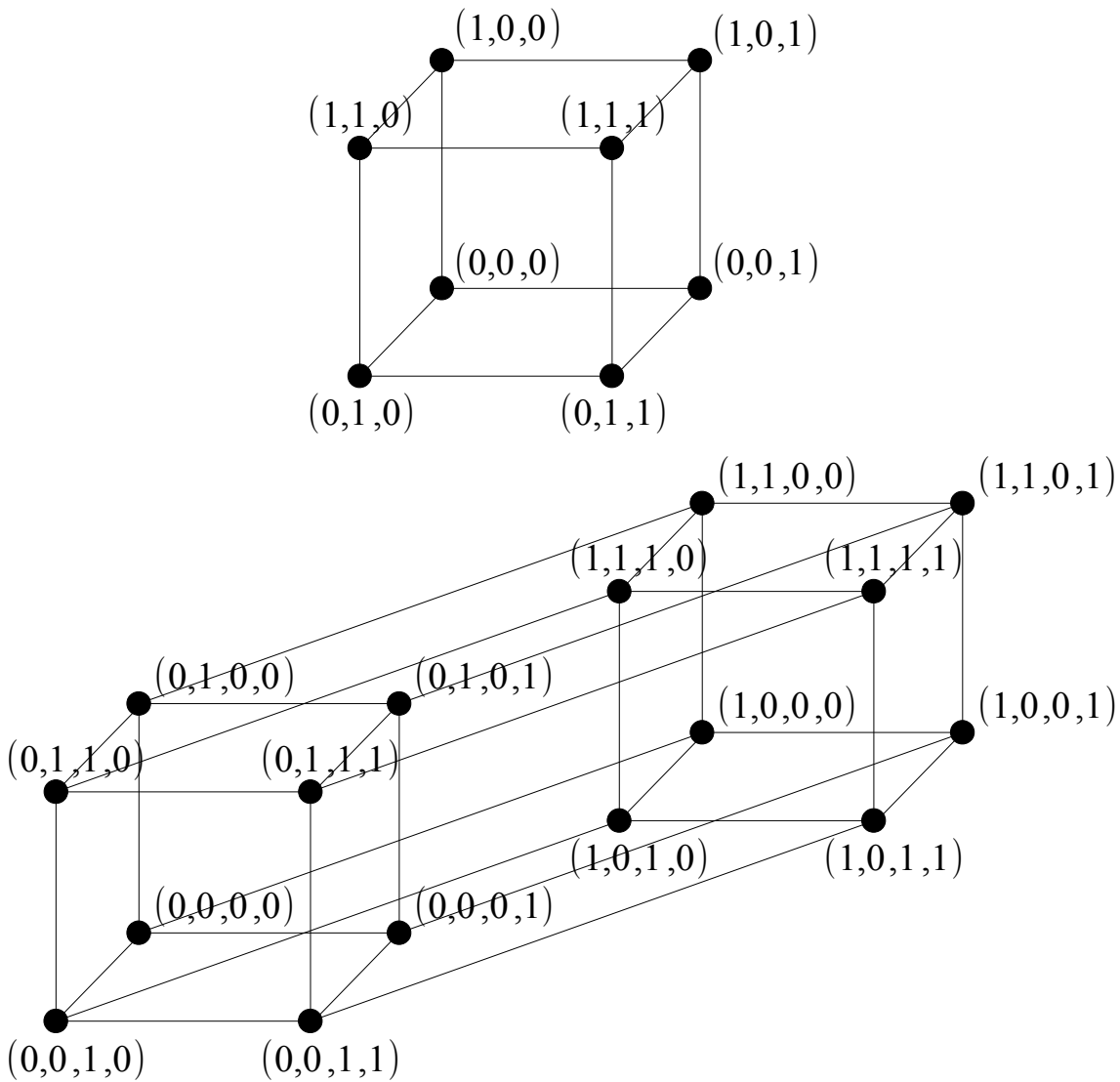

Նկ. 1.2.2

Երկկողմանի գրաֆի օրինակ է *n-չափանի խորանարդը*: n-չափանի խորանարդը նշանակենք $Q_n$-ով: Հիշենք, որ այն սահմանվում էր որպես գրաֆ, որտեղ

$$V(Q_n) = \{\tilde{\alpha}: \tilde{\alpha} = (\alpha_1, \ldots, \alpha_n), \alpha_i \in \{0,1\}, 1 \leq i \leq n\} \text{ և}$$

$$E(Q_n) = \left\{\tilde{\alpha}\tilde{\beta}: \tilde{\alpha}, \tilde{\beta} \in V(Q_n) \text{ և } \sum_{i=1}^{n}|\alpha_i - \beta_i| = 1\right\}:$$

$n$-չափանի խորանարդում որպես $V_1$ վերցնենք այն հավաքածուների բազմությունը, որոնք պարունակում են կենտ թվով մեկեր, իսկ որպես $V_2$՝ այն հավաքածուները, որոնք պարունակում են զույգ թվով մեկեր: Նշենք նաև, որ $n$-չափանի խորանարդը $n$-համասեռ գրաֆ է: Նկ. 1.2.2-ում պատկերված են $Q_3$ և $Q_4$ գրաֆները:

**Սահմանում 1.2.6:** Եթե $G = (V, E)$ երկկողմանի գրաֆում $V_1$ բազմությանը պատկանող յուրաքանչյուր գագաթ միացված է $V_2$ բազմությանը պատկանող յուրաքանչյուր գագաթի, ապա $G$ գրաֆը կանվանենք *լրիվ երկկողմանի գրաֆ*: Եթե այդ դեպքում $|V_1| = m$ և $|V_2| = n$, ապա կգրենք $G = K_{m,n}$:

Նկատենք, որ $|V(K_{m,n})| = m + n$ և $|E(K_{m,n})| = m \cdot n$: Նշենք նաև, որ $K_{1,n}$ գրաֆները կոչվում են *աստղեր*: Ստորև պատկերված են $K_{2,3}, K_{3,3}$ գրաֆները և $K_{1,5}$ աստղը:

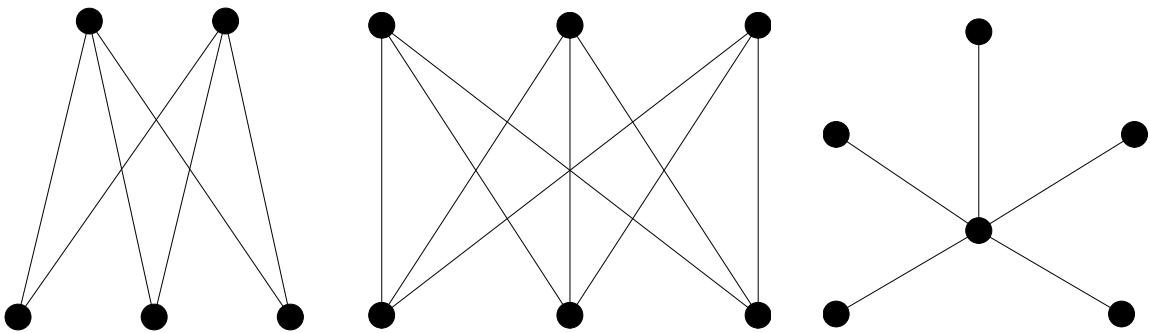

Նկ. 1.2.3

**Թեորեմ 1.2.2:** Այն գրաֆների քանակը, որոնց գագաթների բազմությունը $V = \{v_1, \ldots, v_n\}$-ն է, հավասար է $2^{\binom{n}{2}}$:

**Ապացույց:** Իրոք, քանի որ բոլոր գրաֆներում գագաթների բազմությունը նույնն է, ապա բավական է հաշվել թե քանի իրարից տարբեր հնարավորություն կա կողերի $E$ բազմության համար: Նկատենք, որ $V$ բազմության տարբերից բաղկացած զույգերի քանակը $\binom{n}{2}$-է, և յուրաքանչյուր այդպիսի զույգ կամ մասնակցում է, կամ չի մասնակցում $E$-ի մեջ, և հետևաբար կողերի բազմության համար հնարավոր տարբերակների քանակը $2^{\binom{n}{2}}$-է: ∎

**Թեորեմ 1.2.3:** Այն գրաֆների քանակը, որոնց գագաթների բազմությունը $V = \{v_1, \ldots, v_n\}$-ն է և որոնցում բոլոր գագաթների աստիճանները զույգ թվեր են, հավասար է $2^{\binom{n-1}{2}}$:

**Ապացույց:** Համաձայն թեորեմ 1.2.2-ի, բավական է ապացուցել, որ կարելի է հաստատել փոխմիարժեք արտապատկերում բոլոր այն գրաֆների միջև, որոնց գագաթների բազմությունը $\{v_1, \ldots, v_{n-1}\}$-է, և այն գրաֆների, որոնց գագաթների բազմությունը $\{v_1, \ldots, v_n\}$-ն է և որոնցում բոլոր գագաթների աստիճանները զույգ թվեր են: Վերցնենք ցանկացած $G$ գրաֆ, որի գագաթների բազմությունը $\{v_1, \ldots, v_{n-1}\}$-է: Դիտարկենք $G'$ գրաֆը, որը ստացվում է $G$ գրաֆից հետևյալ կերպ. $G$ գրաֆին ավելացնենք $v_n$ գագաթը, և այն միացնենք կողերով $G$ գրաֆի կենտ աստիճան ունեցող գագաթների հետ: Համաձայն հետևանք 1.2.1-ի, $G'$ գրաֆում բոլոր գագաթների աստիճանները զույգ թվեր են: Ավելին, դժվար չէ տեսնել, որ նկարագրված համապատասխանությունը փոխմիարժեք է: ∎

Դիցուք $G$-ն և $H$-ը գրաֆներ են:

**Սահմանում 1.2.7:** $H$ գրաֆը կոչվում է $G$ գրաֆի *ենթագրաֆ* և կգրենք $H \subseteq G$, եթե $V(H) \subseteq V(G)$ և $E(H) \subseteq E(G)$: Հակառակ դեպքում, կգրենք $H \nsubseteq G$:

**Սահմանում 1.2.8:** $H$ գրաֆը կոչվում է $G$ գրաֆի *կմախքային ենթագրաֆ*, եթե $H \subseteq G$ և $V(H) = V(G)$:

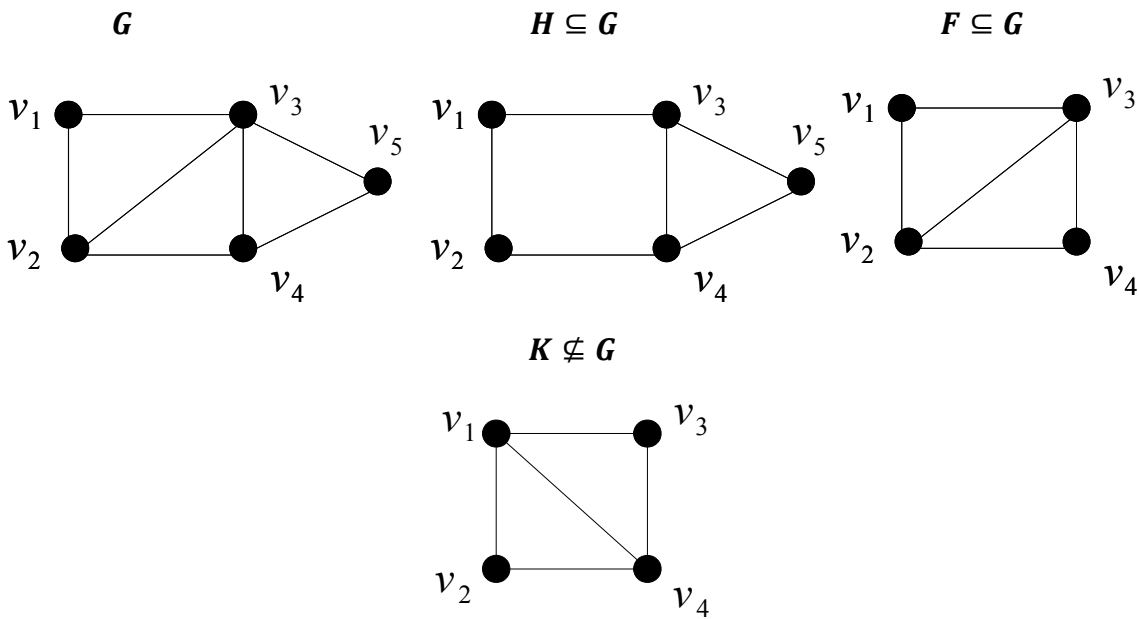

Նկ. 1.2.4

Նկատենք, որ եթե դիտարկենք նկ. 1.2.4-ում պատկերված $G$ գրաֆը, ապա $F$ գրաֆը հանդիսանում է, իսկ $K$ գրաֆը չի հանդիսանում նրա ենթագրաֆ: Մյուս կողմից հեշտ է տեսնել, որ այդ նկարում պատկերված $H$ գրաֆը հանդիսանում է, իսկ $F$ գրաֆը չի հանդիսանում նրա կմախքային ենթագրաֆ:

Դիցուք $G = (V, E)$-ն գրաֆ է և $S \subseteq V(G)$:

**Սահմանում 1.2.9:** $G$ գրաֆի $G[S]$ ենթագրաֆը կոչվում է $S$ *բազմությամբ ծնված ենթագրաֆ* կամ *ծնված ենթագրաֆ*, եթե $V(G[S]) = S$ և $E(G[S]) = \{uv : u, v \in S$ և $uv \in E(G)\}$:

Դիտարկենք հետևյալ օրինակը. եթե $G$ գրաֆը նկ. 1.2.5-ում պատկերված գրաֆն է, ապա այդ նկարում բերված $H$ գրաֆը հանդիսանում է $G$ գրաֆի $\{v_2, v_3, v_4, v_5\}$ բազմությամբ ծնված ենթագրաֆ: Մյուս կողմից հեշտ է տեսնել, որ նույն նկարում պատկերված $F$ և $K$ գրաֆները $G$ գրաֆի ծնված ենթագրաֆներ չեն:

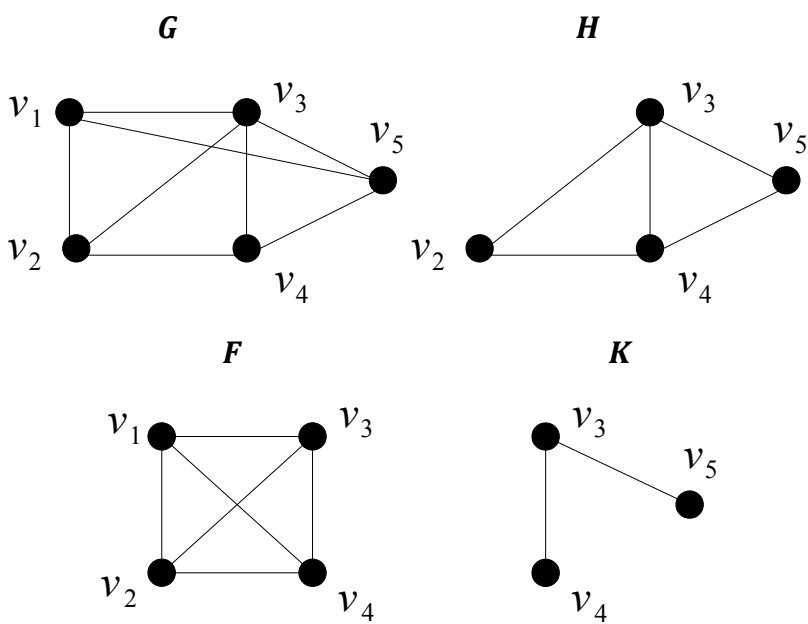

Նկ. 1.2.5

Դիցուք $G = (V, E)$-ն գրաֆ է:

**Սահմանում 1.2.10:** $G$ գրաֆի $u_0, u_1, \ldots, u_{k-1}, u_k$ գագաթներից և $e_1, \ldots, e_k$ կողերից կազմված $u_0, e_1, u_1, \ldots, u_{k-1}, e_k, u_k$ հաջորդականությունը կանվանենք $k$ *երկարությամբ* $u_0$-*ից* $u_k$ *շրջանցում* կամ $k$ *երկարությամբ* $(u_0, u_k)$-*շրջանցում*, եթե $e_i = u_{i-1}u_i$, երբ $1 \leq i \leq k$: Սահմանված $(u_0, u_k)$-շրջանցումը կրճատ կնշանակենք նրա գագաթների $u_0, u_1, \ldots, u_{k-1}, u_k$ հաջորդականությամբ, ենթադրելով, որ յուրաքանչյուր հաջորդ գագաթ հարևան է նախորդին:

Նշենք, որ շրջանցման մեջ գագաթները և կողերը կարող են կրկնվել, իսկ $(u_0, u_k)$-շրջանցման երկարությունը ցույց է տալիս այդ շրջանցման մեջ առկա կողերի քանակը, երբ յուրաքանչյուր կող հաշվվում է այնքան անգամ, որքան այն հանդիպում է շրջանցման մեջ:

Դիտարկենք նկ. 1.2.5-ում պատկերված $G$ գրաֆի գագաթների հետևյալ երկու հաջորդականությունները․ $v_1, v_3, v_2, v_4, v_3, v_2, v_4, v_3, v_5$ և $v_3, v_1, v_4, v_2, v_3, v_5$: Նկատենք, որ նրանցից առաջինը $G$ գրաֆի շրջանցում է, իսկ երկրորդը՝ ոչ:

Նկատենք, որ նկ. 1.2.5-ի $G$ գրաֆի $v_1, v_3, v_2, v_4, v_3, v_2, v_4, v_3, v_5$ շրջանցման երկարությունը հավասար է ութի:

**Սահմանում 1.2.11:** $(u_0, u_k)$-շրջանցումը կանվանենք *փակ*, եթե $u_0 = u_k$:

**Սահմանում 1.2.12:** $G$ գրաֆի $(u_0, u_k)$-շրջանցումը կանվանենք $u_0$-*ից* $u_k$ *ճանապարհի* կամ $(u_0, u_k)$-*ճանապարհի*, եթե $u_0 u_1$-ը,…, $u_{k-1} u_k$-ն $G$ գրաֆի զույգ առ զույգ տարբեր կողեր են: Եթե $P$-ն $G$ գրաֆի ճանապարհի է, ապա $|P|$-ով կնշանակենք այդ ճանապարհի երկարությունը, այսինքն՝ այդ ճանապարհի մեջ առկա կողերի քանակը:

Նկատենք, որ նկ. 1.2.5-ի $G$ գրաֆի $v_1, v_3, v_2, v_4, v_3, v_5$ շրջանցումը ճանապարհի է, իսկ $v_1, v_3, v_2, v_4, v_3, v_2, v_4, v_3, v_5$-ը՝ ոչ:

**Սահմանում 1.2.13:** $G$ գրաֆի $(u_0, u_k)$-ճանապարհը կանվանենք *պարզ* $(u_0, u_k)$-*ճանապարհի*, եթե նրա մեջ մտնող բոլոր գագաթները զույգ առ զույգ տարբեր են:

**Սահմանում 1.2.14:** $G$ գրաֆի $(u_0, u_k)$-ճանապարհը կանվանենք *փակ ճանապարհի* կամ *ցիկլ*, եթե այն փակ շրջանցում է, այսինքն՝ եթե $u_0 = u_k$:

Նկատենք, որ նկ. 1.2.5-ի $G$ գրաֆի $v_1, v_3, v_2, v_4, v_3, v_5$ ճանապարհը պարզ $(v_1, v_5)$-ճանապարհի չէ, իսկ $v_1, v_3, v_4, v_5$-ը նույն գրաֆի պարզ $(v_1, v_5)$-ճանապարհի է:

**Սահմանում 1.2.15:** $G$ գրաֆի ցիկլը կանվանենք *պարզ*, եթե նրանում կրկնվում են միայն առաջին և վերջին գագաթները:

Նկատենք, որ նկ. 1.2.5-ի $G$ գրաֆի $v_3, v_2, v_1, v_3, v_5, v_4, v_3$ ճանապարհը պարզ ցիկլ չէ, իսկ $v_1, v_3, v_4, v_2, v_1$-ը նույն գրաֆի պարզ ցիկլ է:

$n$ գագաթ ունեցող *պարզ ցիկլը* կնշանակենք $C_n$-ով, $n \geq 3$: $n$ գագաթ ունեցող *պարզ ճանապարհը* կնշանակենք $P_n$-ով:

**Սահմանում 1.2.16:** $G$ գրաֆում $u$ և $v$ գագաթների միջև *հեռավորությունը* կսահմանենք որպես կարճագույն $(u, v)$-ճանապարհի երկարություն, եթե $G$ գրաֆում գոյություն ունի առնվազն մեկ $(u, v)$-ճանապարհի, և $+\infty$՝ հակառակ դեպքում: $G$ գրաֆում

$u$ և $v$ գագաթների միջև հեռավորությունը կնշանակենք $d_G(u,v)$-ով կամ $d(u,v)$-ով։

Նկատենք, որ

1. $G$ գրաֆի ցանկացած $u$ և $v$ գագաթների համար $d_G(u,v) \geq 0$, և $d_G(u,v) = 0$ այն և միայն այն դեպքում, երբ $u = v$;
2. $G$ գրաֆի ցանկացած $u$ և $v$ գագաթների համար $d_G(u,v) = d_G(v,u)$;
3. $G$ գրաֆի ցանկացած $u$, $v$ և $w$ գագաթների համար $d_G(u,v) \leq d_G(u,w) + d_G(w,v)$։

Այստեղից հետևում է, որ ցանկացած $G = (V,E)$ գրաֆի համար, որում կամայական երկու գագաթների միջև կա միացնող ճանապարհի, $(V, d_G)$ զույգն իրենից ներկայացնում է մետրիկական տարածություն։

**Լեմմա 1.2.1:** Ենթադրենք, որ $G$ գրաֆում $u$-ն և $v$-ն իրարից տարբեր երկու գագաթներ են։ Այդ դեպքում ցանկացած $(u,v)$-շղթայումից կարելի է առանձնացնել պարզ $(u,v)$-ճանապարհի։

**Ապացույց:** Ապացույցը կկատարենք մակածման եղանակով ըստ $(u,v)$-շղթայման $k$ երկարության։

Ենթադրենք $k = 1$: Այդ դեպքում $(u,v)$-շղթայումը բաղկացած է մեկ կողից, որն էլ կկազմի որոնելի պարզ $(u,v)$-ճանապարհը։ Ենթադրենք, որ պնդումը ճիշտ է բոլոր այն $(u,v)$-շղթայումների համար, որոնց երկարությունը փոքր է $k$-ից, և դիտարկենք $k$ երկարություն ունեցող $u_0, u_1, \ldots, u_k$ շղթայումը, որտեղ $u_0 = u$ և $u_k = v$:

Եթե $u_0, u_1, \ldots, u_k$ շղթայման մեջ բոլոր գագաթները զույգ առ զույգ տարբեր են, ապա այն իրենից ներկայացնում է պարզ $(u,v)$-ճանապարհի, և, հետևաբար, պնդումն ապացուցված է։ Հետևաբար, կարող ենք ենթադրել, որ գոյություն ունեն $i < j$ թվեր այնպես, որ $u_i = u_j$: Դիտարկենք $u_0, u_1, \ldots, u_i, u_{j+1}, \ldots, u_k$ $(u,v)$-շղթայումը։ Նկատենք, որ այս շղթայման երկարությունը փոքր է $k$-ից, և, հետևաբար, համաձայն մակածման ենթադրության, նրանից կարելի է առանձնացնել պարզ $(u,v)$-ճանապարհի։ ∎

**Լեմմա 1.2.2:** $G$ գրաֆի ցանկացած կենտ երկարություն ունեցող փակ շղթայումից կարելի է առանձնացնել կենտ երկարություն ունեցող պարզ ցիկլ։

**Ապացույց:** Ապացույցը կկատարենք մակածման եղանակով ըստ կենտ երկարություն ունեցող փակ շղթայման $k$ երկարության։

Ենթադրենք $k = 3$: Այդ դեպքում կենտ երկարություն ունեցող փակ շղթայումն իրենից ներկայացնում է եռանկյուն, որն էլ կազմում է որոնելի կենտ երկարություն ունեցող պարզ ցիկլը: Ենթադրենք, որ պնդումը ճիշտ է բոլոր այն կենտ երկարություն ունեցող փակ շրջանցումների համար, որոնց երկարությունը փոքր է $k$-ից, և դիտարկենք $k$ երկարություն ունեցող $u_0, u_1, \ldots, u_k$ կենտ, փակ շրջանցումը, որտեղ $u_0 = u_k$:

Եթե $u_0, u_1, \ldots, u_k$ շրջանցման մեջ բոլոր գագաթները զույգ առ զույգ տարբեր են, բացի $u_0 = u_k$-ից, ապա այն իրենից ներկայացնում է կենտ երկարություն ունեցող պարզ ցիկլ, և, հետևաբար, պնդումն ապացուցված է: Հետևաբար, կարող ենք ենթադրել, որ գոյություն ունեն $i < j$ թվեր այնպես, որ $u_i = u_j$: Դիտարկենք $u_0, u_1, \ldots, u_i, u_{j+1}, \ldots, u_k$ և $u_i, u_{i+1}, \ldots, u_j$ փակ շրջանցումները: Նկատենք, որ նրանց երկարությունները փոքր են $k$-ից, ավելին, քանի որ $u_0, u_1, \ldots, u_k$ շրջանցման երկարությունը կենտ է, ապա նշված շրջանցումներից մեկի երկարությունը ևս կենտ է: Համաձայն մակածման ենթադրության նրանից կարելի է առանձնացնել կենտ երկարություն ունեցող պարզ ցիկլ: ∎

## § 1.3. Գործողություններ գրաֆների հետ

Այս պարագրաֆում մենք կդիտարկենք տարբեր գործողություններ գրաֆների հետ: Այդ գործողությունները հնարավորություն են տալիս արդեն գոյություն ունեցող գրաֆների հիման վրա կառուցել նոր գրաֆներ և նաև օգնում են ներկայացնել գրաֆի կառուցվածքը ավելի փոքր և պարզ կառուցվածք ունեցող գրաֆների միջոցով:

1. **Գագաթի հեռացում:** Դիցուք տրված են $G$ գրաֆը ($|V(G)| \geq 2$) և $v \in V(G)$: $G$ գրաֆից $v$ *գագաթի հեռացումը՝* $G - v$ գրաֆը սահմանենք հետևյալ կերպ. $V(G - v) = V(G) \setminus \{v\}$ և $E(G - v) = E(G) \setminus \{e : e - \text{ն կից է } v - \text{ին}\}$:

Նկ. 1.3.1-ում պատկերված է $G$ գրաֆը և այդ գրաֆից $v_1$ գագաթի հեռացումից առաջացած $G - v_1$ գրաֆը:

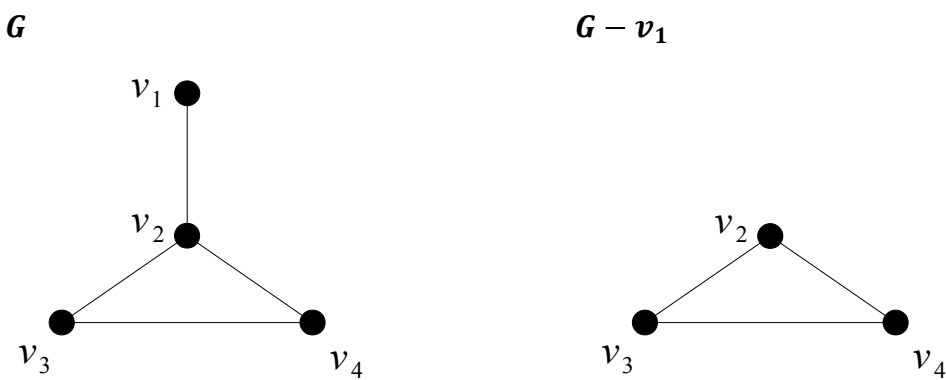

Նկ. 1.3.1

2. **Գագաթի ավելացում:** Դիցուք տրված են $G$ գրաֆը և $v \notin V(G)$: $G$ գրաֆին $v$ *գագաթի ավելացումը*՝ $G + v$ գրաֆը սահմանենք հետևյալ կերպ. $V(G + v) = V(G) \cup \{v\}$ և $E(G + v) = E(G)$:

Նկ. 1.3.2-ում պատկերված է $G$ գրաֆը և այդ գրաֆին $v$ գագաթի ավելացումից առաջացած $G + v$ գրաֆը:

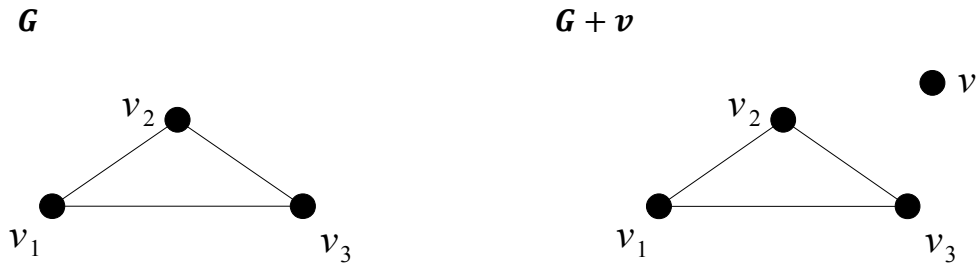

Նկ. 1.3.2

3. **Կողի հեռացում:** Դիցուք տրված են $G$ գրաֆը և $e \in E(G)$: $G$ գրաֆից $e$ *կողի հեռացումը*՝ $G - e$ գրաֆը սահմանենք հետևյալ կերպ. $V(G - e) = V(G)$ և $E(G - e) = E(G) \setminus \{e\}$:

Նկ. 1.3.3-ում պատկերված է $G$ գրաֆը և այդ գրաֆից $v_3v_6$ կողի հեռացումից առաջացած $G - v_3v_6$ գրաֆը:

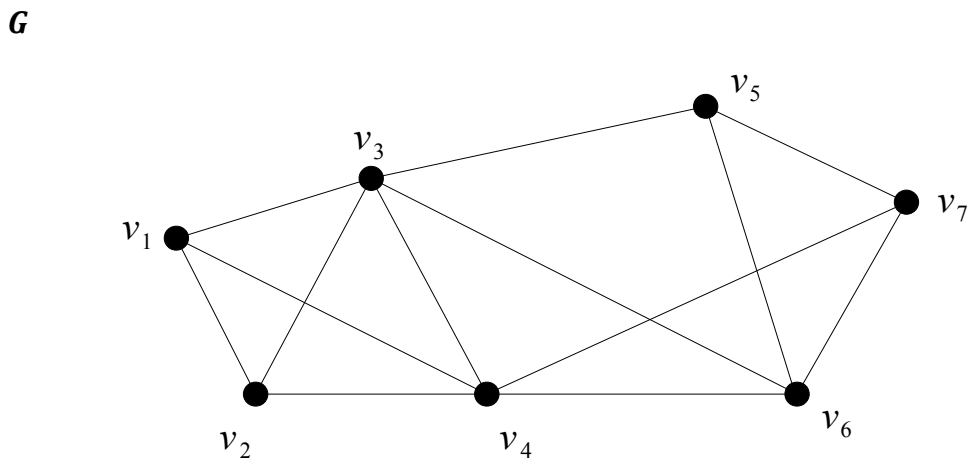

Նկ. 1.3.3

4. **Կողի ավելացում:** Դիցուք տրված են $G$ գրաֆը և $e = uv \notin E(G)$ ($u, v \in V(G)$)։ $G$ գրաֆին $e$ կողի ավելացումը՝ $G + e$ գրաֆը սահմանենք հետևյալ կերպ. $V(G + e) = V(G)$ և $E(G + e) = E(G) \cup \{e\}$։

Նկ. 1.3.4-ում պատկերված է $G$ գրաֆը և այդ գրաֆին $v_4v_5$ կողի ավելացումից առաջացած $G + v_4v_5$ գրաֆը։

$G$

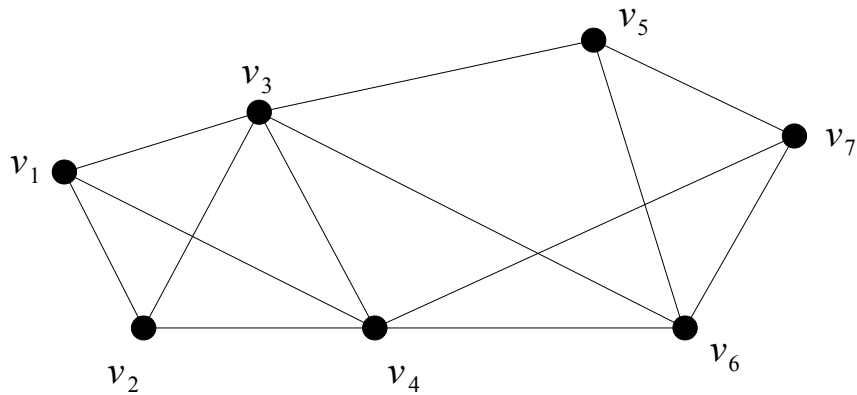

$G + v_4v_5$

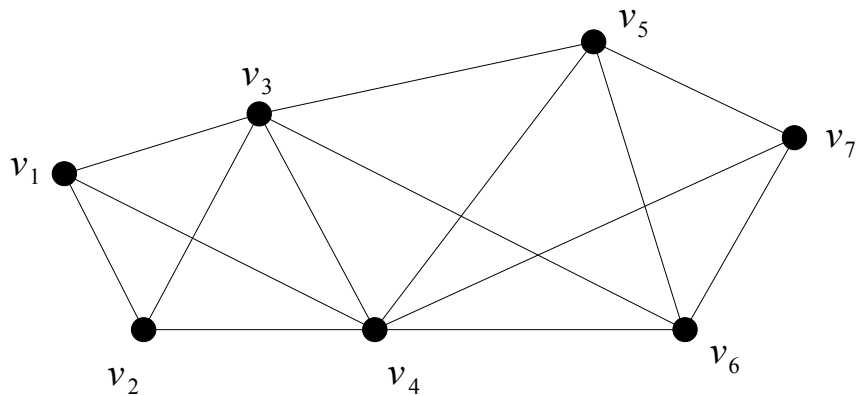

Նկ. 1.3.4

5. **Կողի տրոհում:** Դիցուք տրված են $G$ գրաֆը և $e = uv \in E(G)$։ $G$ գրաֆի $e$ կողի *տրոհումը*՝ $G_e$ գրաֆը սահմանենք հետևյալ կերպ. $V(G_e) = V(G) \cup \{w\}$, $w \notin V(G)$ և $E(G_e) = (E(G)\setminus\{e\}) \cup \{uw, vw\}$։

Նկ. 1.3.5-ում պատկերված է $G$ գրաֆը և այդ գրաֆի $v_4v_7$ կողի տրոհումից առաջացած $G_{v_4v_7}$ գրաֆը։

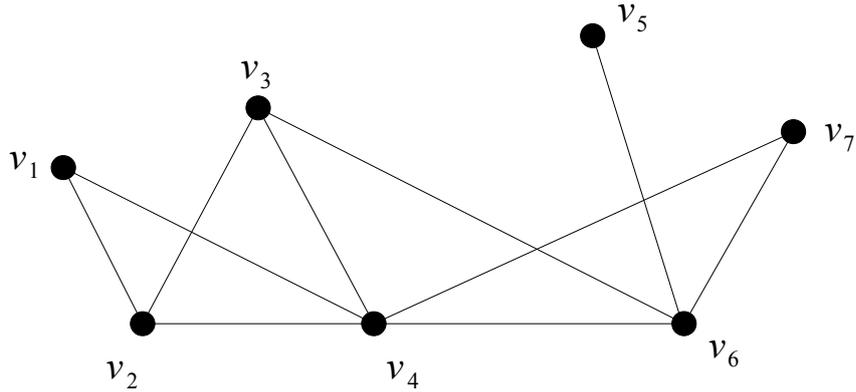

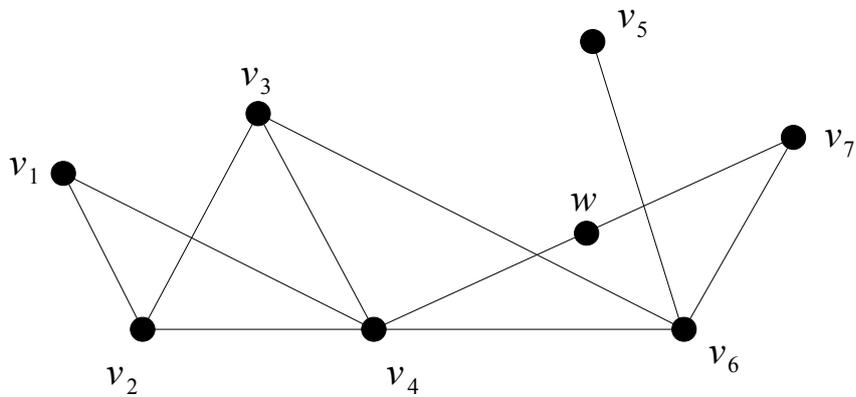

Նկ. 1.3.5

6. **Ենթագրաֆի կծկում**: Դիցուք տրված են $G$ գրաֆը և նրա $H$ ենթագրաֆը: $G$ գրաֆի $H$ *ենթագրաֆի կծկումը*՝ $G/H$ գրաֆը սահմանենք հետևյալ կերպ.

$$V(G/H) = (V(G) \setminus V(H)) \cup \{w\}, \ w \notin V(G) \ \text{և}$$

$$E(G/H) = (E(G) \setminus \{e = uv : u \in V(H) \ \text{կամ} \ v \in V(H)\}) \cup \{e = uw : u \in N_G(V(H)) \setminus V(H)\}:$$

Նկ. 1.3.6-ում պատկերված է $G$ գրաֆը, $H$ ենթագրաֆը և $G$ գրաֆի $H$ ենթագրաֆի կծկումից առաջացած $G/H$ գրաֆի օրինակը:

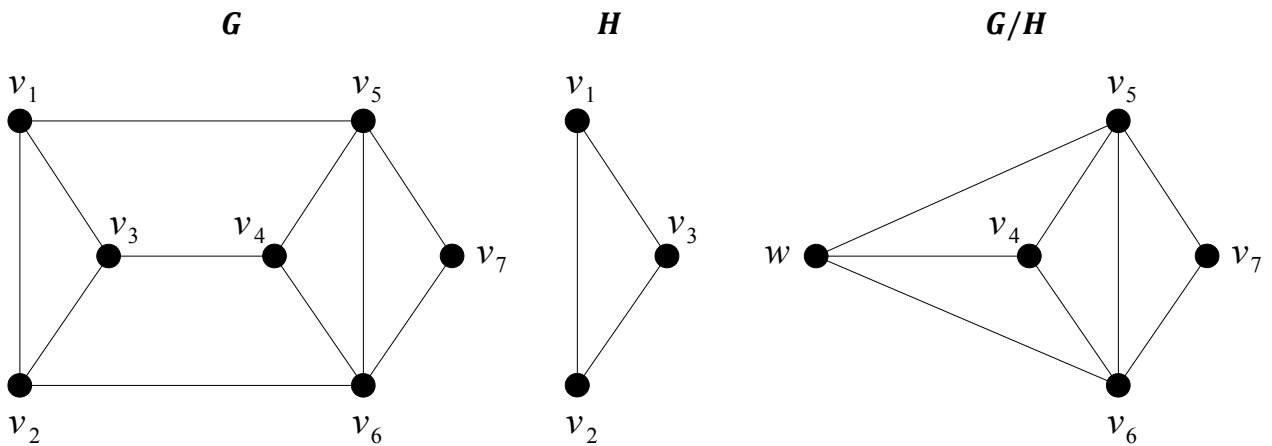

Նկ. 1.3.6

7. **Գրաֆի լրացում:** $G$ գրաֆի *լրացում* կոչվում է $\overline{G}$ գրաֆը, որի համար $V(\overline{G}) = V(G)$ և $E(\overline{G}) = \{uv : u, v \in V(\overline{G})$ և $uv \notin E(G)\}$:

Նկ. 1.3.7-ում պատկերված է $G$ գրաֆը և նրա $\overline{G}$ լրացում գրաֆի օրինակը:

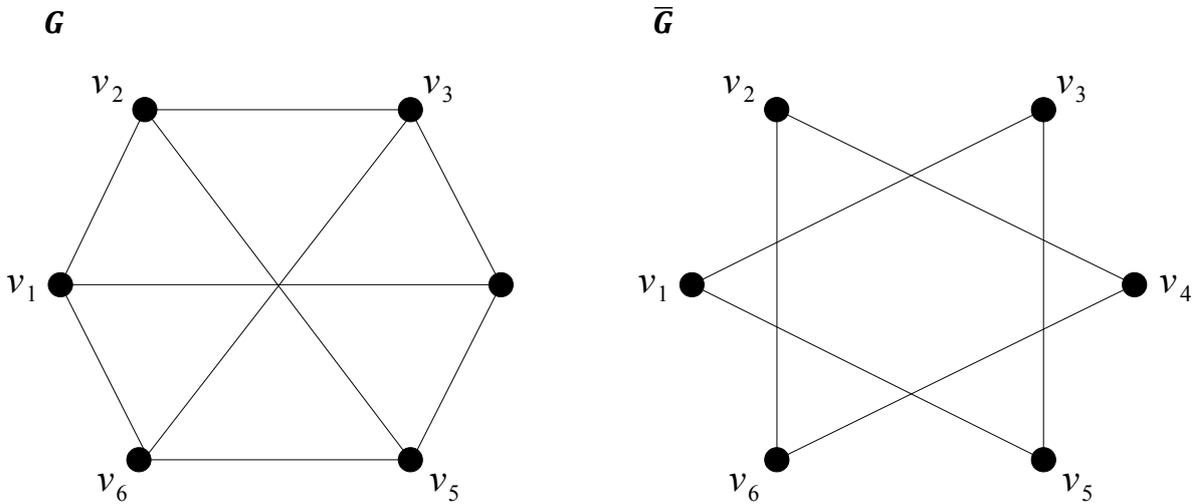

Նկ. 1.3.7

Նկատենք, որ եթե $G$-ն $(n, m)$-գրաֆ է, ապա $\overline{G}$-ը կլինի $\left(n, \binom{n}{2} - m\right)$-գրաֆ:

8. **Գրաֆների միավորում:** Դիցուք տրված են $G$ և $H$ գրաֆները, որոնց համար $V(G) \cap V(H) = \emptyset$: $G$ և $H$ գրաֆների *միավորում* կոչվում է $G \cup H$ գրաֆը, որի համար $V(G \cup H) = V(G) \cup V(H)$ և $E(G \cup H) = E(G) \cup E(H)$:

Նկ. 1.3.8-ում պատկերված են $G$ և $H$ գրաֆները և նրանց միավորում $G \cup H$ գրաֆի օրինակը:

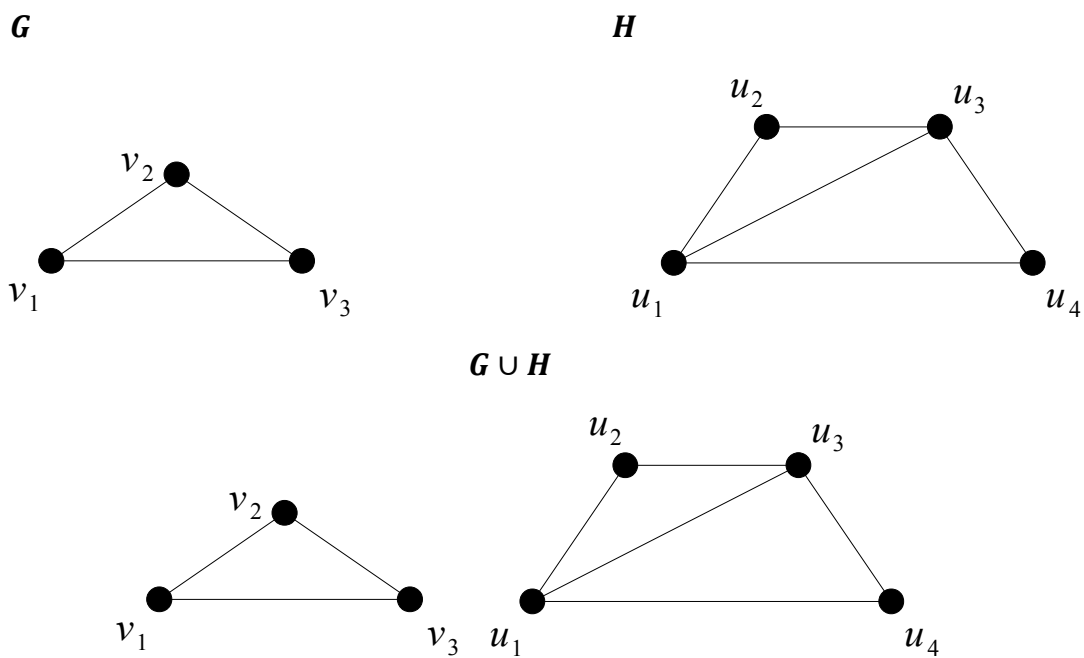

Նկ. 1.3.8

Նկատենք, որ եթե $G$-ն $(n_1, m_1)$-գրաֆ է և $H$-ը $(n_2, m_2)$-գրաֆ է, ապա $G \cup H$-ը կլինի $(n_1 + n_2, m_1 + m_2)$-գրաֆ:

9. **Գրաֆների գումար:** Դիցուք տրված են $G$ և $H$ գրաֆները, որոնց համար $V(G) \cap V(H) = \emptyset$: $G$ և $H$ գրաֆների *գումար* կոչվում է $G + H$ գրաֆը, որի համար $V(G + H) = V(G) \cup V(H)$ և $E(G \cup H) = E(G) \cup E(H) \cup \{uv : u \in V(G)$ և $v \in V(H)\}$:

Նկ. 1.3.9-ում պատկերված են $G$ և $H$ գրաֆները և նրանց գումար $G + H$ գրաֆի օրինակը:

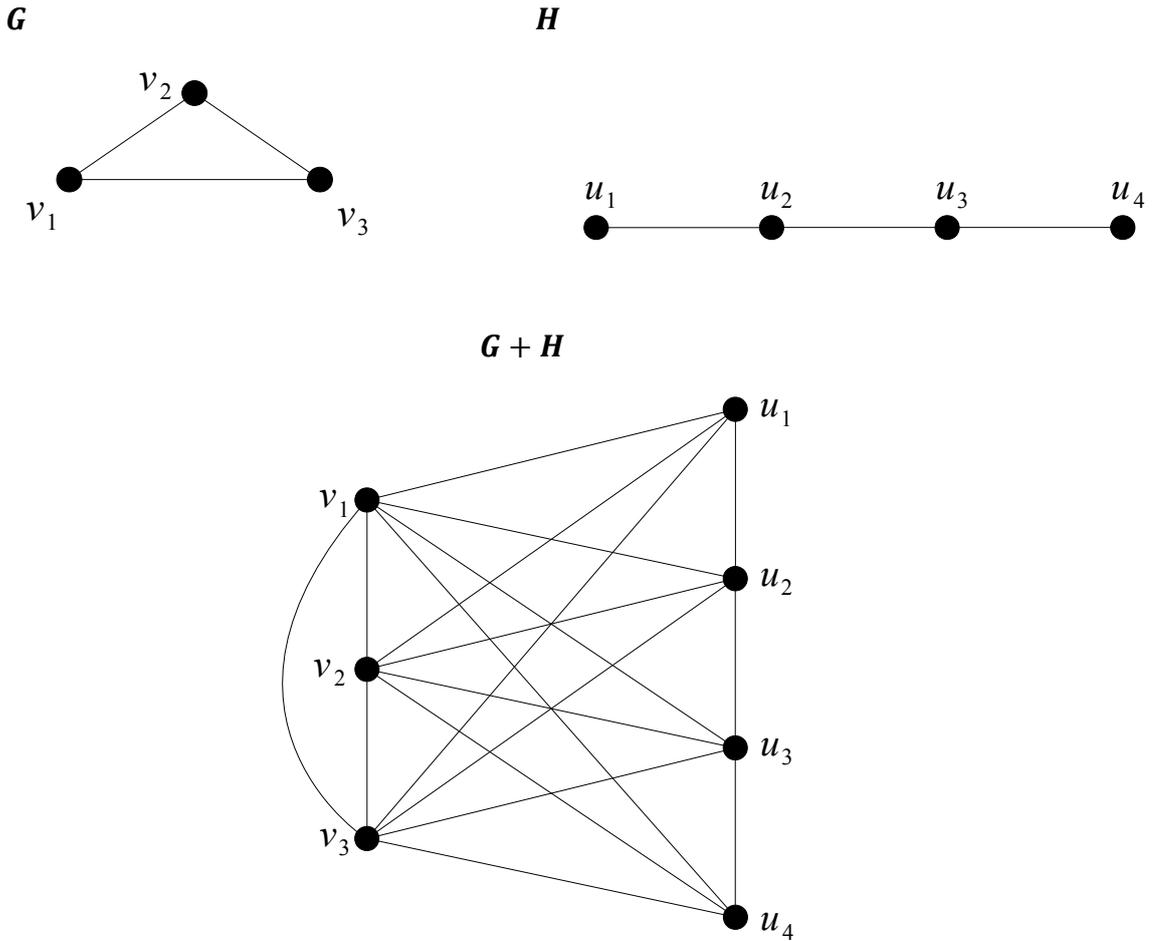

Նկ. 1.3.9

Նկատենք, որ եթե $G$-ն $(n_1, m_1)$-գրաֆ է և $H$-ը $(n_2, m_2)$-գրաֆ է, ապա $G + H$-ը կլինի $(n_1 + n_2, m_1 + m_2 + n_1 \cdot n_2)$-գրաֆ: Նշենք նաև, որ գրաֆների գումար գործողությունը կոմուտատիվ է և ասոցիատիվ:

Հեշտ է տեսնել, որ $K_{m,n}$ լրիվ երկկողմանի գրաֆը կարելի է գրաֆների լրացման և գումարի միջոցով արտահայտել. $K_{m,n} = \overline{K}_m + \overline{K}_n$: Ավելին, գրաֆների լրացումը և գումարը հնարավորություն են տալիս սահմանել $K_{n_1, n_2, \cdots, n_r}$ *լրիվ $r$-կողմանի գրաֆը* հետևյալ կերպ. $K_{n_1, n_2, \cdots, n_r} = \overline{K}_{n_1} + \overline{K}_{n_2} + \cdots + \overline{K}_{n_r}$:

10. **Գրաֆների դեկարտյան արտադրյալ։** Դիցուք տրված են $G$ և $H$ գրաֆները, որոնց համար $V(G) \cap V(H) = \emptyset$։ $G$ և $H$ գրաֆների *դեկարտյան արտադրյալ* կոչվում է $G \square H$ գրաֆը, որի համար

$$V(G \square H) = V(G) \times V(H) \text{ և}$$

$$E(G \square H) = \{(u_1, v_1)(u_2, v_2) : (u_1 = u_2 \text{ և } v_1 v_2 \in E(H)) \text{ կամ } (v_1 = v_2 \text{ և } u_1 u_2 \in E(G))\}:$$

Նկ. 1.3.10-ում պատկերված են $G$ և $H$ գրաֆները և նրանց դեկարտյան արտադրյալ $G \square H$ գրաֆի օրինակը․

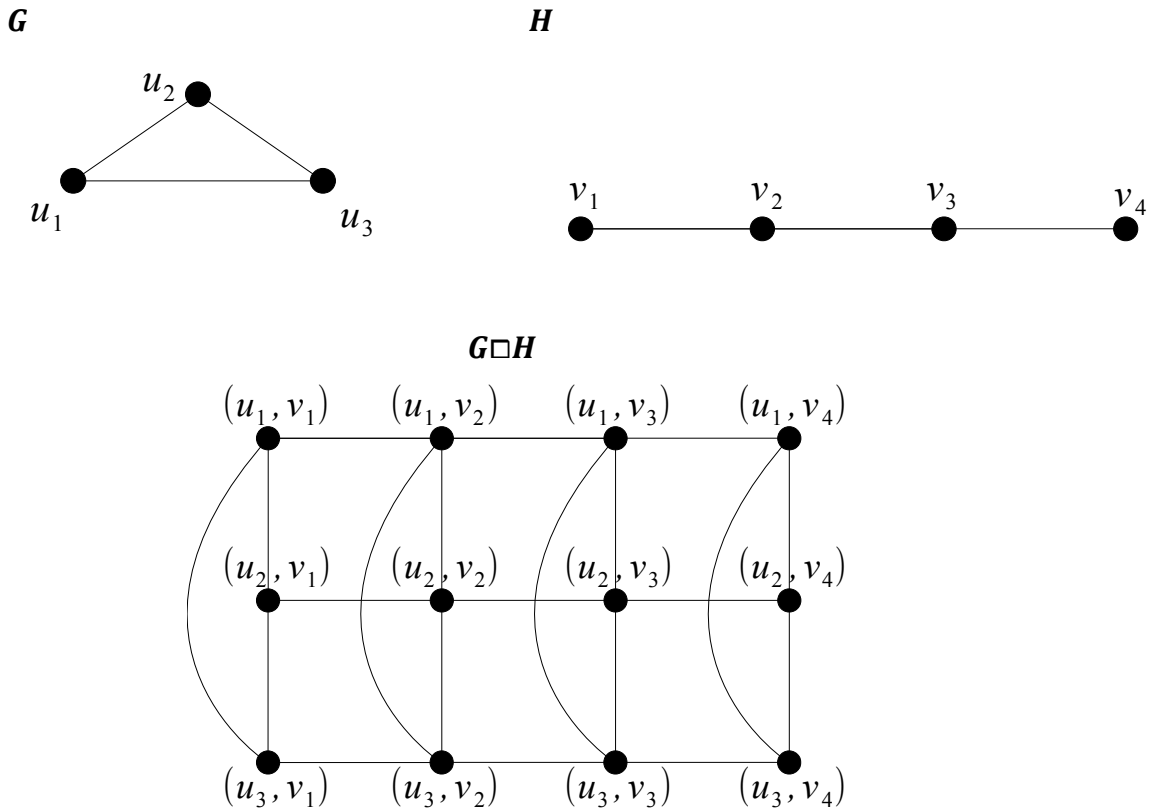

Նկ. 1.3.10

Նկատենք, որ եթե $G$-ն $(n_1, m_1)$-գրաֆ է և $H$-ը $(n_2, m_2)$-գրաֆ է, ապա $G \square H$-ը կլինի $(n_1 \cdot n_2, n_1 \cdot m_2 + m_1 \cdot n_2)$-գրաֆ։ Նշենք նաև, որ գրաֆների դեկարտյան արտադրյալ գործողությունը կոմուտատիվ է և ասոցատիվ։

Հեշտ է տեսնել, որ $Q_n$ $n$-չափանի խորանարդը կարելի է գրաֆների դեկարտյան արտադրյալի միջոցով արտահայտել․ $Q_n = \underbrace{K_2 \square K_2 \square \cdots \square K_2}_{n}$։ Նաև, գրաֆների դեկարտյան արտադրյալը հնարավորություն է տալիս սահմանել $H_{m_1, m_2, \cdots, m_n}$ *Հեմմինգի գրաֆը* հետևյալ կերպ․ $H_{m_1, m_2, \cdots, m_n} = K_{m_1} \square K_{m_2} \square \cdots \square K_{m_n}$։

11. **Գրաֆների թենզորական արտադրյալ:** Դիցուք տրված են $G$ և $H$ գրաֆները, որոնց համար $V(G) \cap V(H) = \emptyset$: $G$ և $H$ գրաֆների *թենզորական արտադրյալ* կոչվում է $G \times H$ գրաֆը, որի համար

$$V(G \times H) = V(G) \times V(H) \text{ և}$$

$$E(G \times H) = \{(u_1, v_1)(u_2, v_2) : u_1 u_2 \in E(G) \text{ և } v_1 v_2 \in E(H)\}:$$

Նկ. 1.3.11-ում պատկերված են $G$ և $H$ գրաֆները և նրանց թենզորական արտադրյալ $G \times H$ գրաֆի օրինակը:

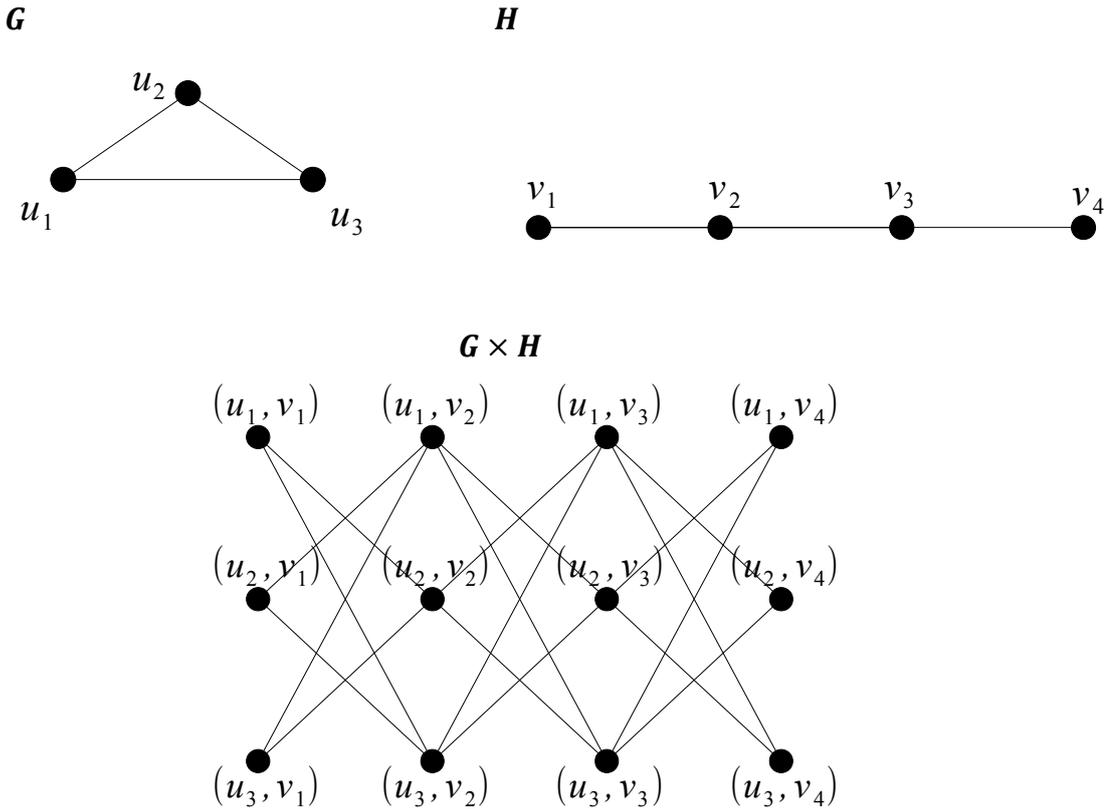

Նկ. 1.3.11

Նկատենք, որ եթե $G$-ն $(n_1, m_1)$-գրաֆ է և $H$-ը $(n_2, m_2)$-գրաֆ է, ապա $G \times H$-ը կլինի $(n_1 \cdot n_2, 2m_1 \cdot m_2)$-գրաֆ: Նշենք նաև, որ գրաֆների թենզորական արտադրյալ գործողությունը կոմուտատիվ է և ասոցիատիվ:

12. **Գրաֆների ուժեղ արտադրյալ:** Դիցուք տրված են $G$ և $H$ գրաֆները, որոնց համար $V(G) \cap V(H) = \emptyset$: $G$ և $H$ գրաֆների *ուժեղ արտադրյալ* կոչվում է $G \boxtimes H$ գրաֆը, որի համար $V(G \boxtimes H) = V(G) \times V(H)$ և $E(G \boxtimes H) = E(G \square H) \cup E(G \times H)$:

Նկ. 1.3.12-ում պատկերված են $G$ և $H$ գրաֆները և նրանց ուժեղ արտադրյալ $G \boxtimes H$ գրաֆի օրինակը:

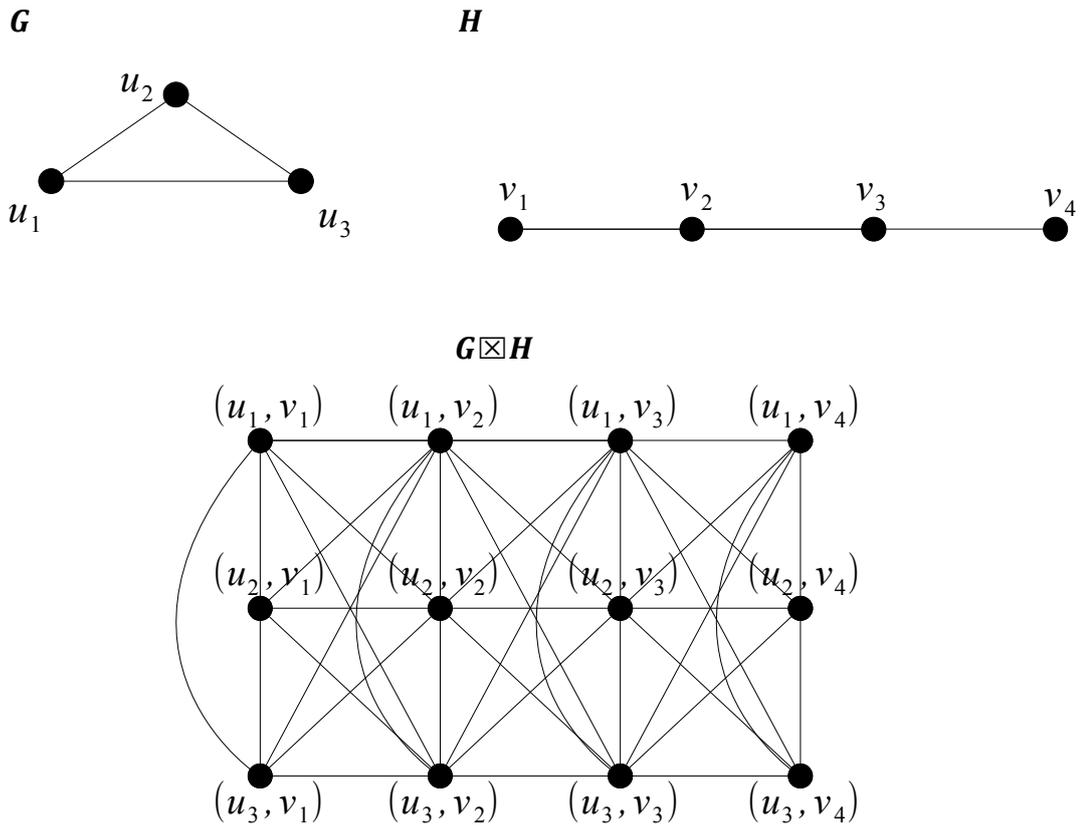

Նկ. 1.3.12

Նկատենք, որ եթե $G$-ն $(n_1, m_1)$-գրաֆ է և $H$-ը $(n_2, m_2)$-գրաֆ է, ապա $G \boxtimes H$-ը կլինի $(n_1 \cdot n_2, n_1 \cdot m_2 + m_1 \cdot n_2 + 2m_1 \cdot m_2)$-գրաֆ: Նշենք նաև, որ գրաֆների ուժեղ արտադրյալ գործողությունը կոմուտատիվ է և ասոցիատիվ:

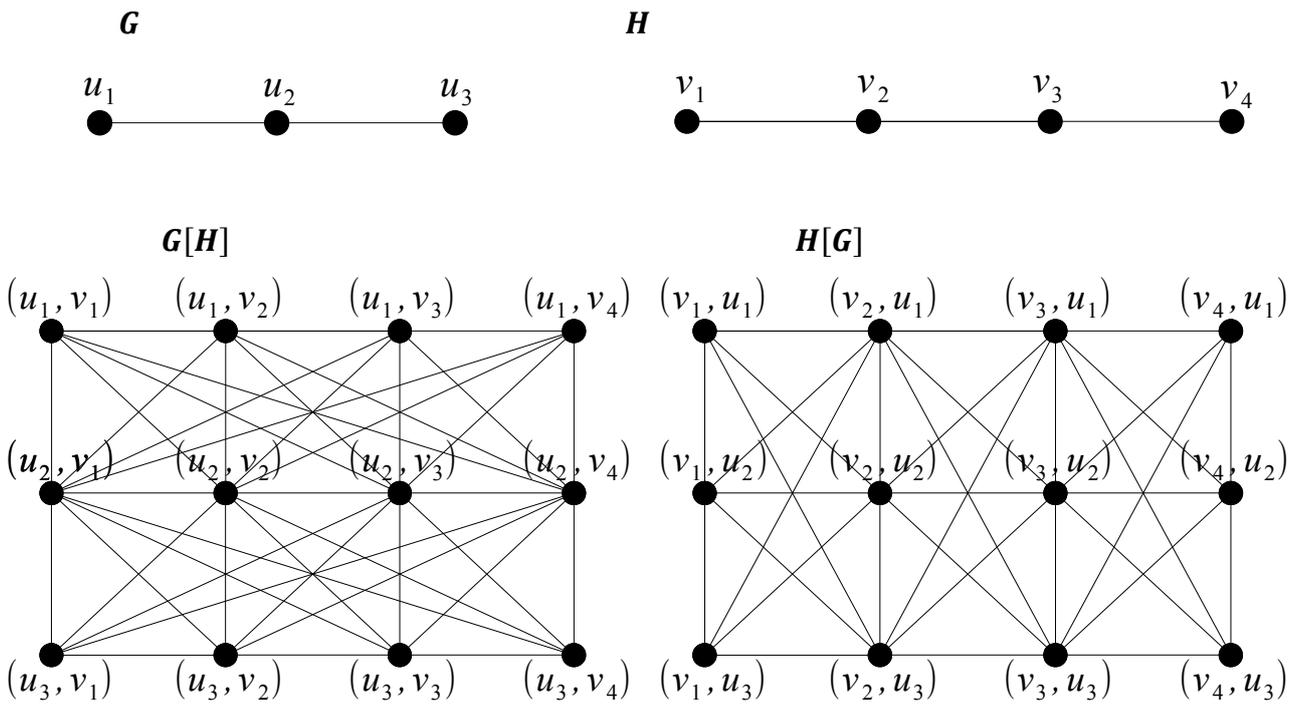

Նկ. 1.3.13

13. **Գրաֆների կոմպոզիցիա:** Դիցուք տրված են $G$ և $H$ գրաֆները, որոնց համար $V(G) \cap V(H) = \emptyset$: $G$ և $H$ գրաֆների *կոմպոզիցիա* կոչվում է $G[H]$ գրաֆը, որի համար

$$V(G[H]) = V(G) \times V(H) \text{ և}$$

$$E(G[H]) = \{(u_1, v_1)(u_2, v_2) : u_1u_2 \in E(G) \text{ կամ } (u_1 = u_2 \text{ և } v_1v_2 \in E(H))\}:$$

Նկ. 1.3.13-ում պատկերված են $G$ և $H$ գրաֆները և նրանց $G[H]$ և $H[G]$ կոմպոզիցիաների օրինակները:

Նկատենք, որ եթե $G$-ն $(n_1, m_1)$-գրաֆ է և $H$-ը $(n_2, m_2)$-գրաֆ է, ապա $G[H]$-ը կլինի $(n_1 \cdot n_2, n_1 \cdot m_2 + m_1 \cdot (n_2)^2)$-գրաֆ: Նկ. 1.3.13-ում պատկերված $G[H]$ և $H[G]$ գրաֆները ցույց են տալիս, որ գրաֆների կոմպոզիցիա գործողությունը կոմուտատիվ չէ: Նշենք նաև, որ գրաֆների կոմպոզիցիա գործողությունը ասոցիատիվ է:

Գրաֆների տարբեր արտադրյալներին ավելի մանրամասն կարելի է ծանոթանալ Համակ, Իմրիի և Կլավզարի գրքում [17]:

# § 1.4. Գրաֆների իզոմորֆիզմ, հոմոմորֆիզմ, ավտոմորֆիզմ և գրաֆի ավտոմորֆիզմների խումբը

Սահմանենք գրաֆների իզոմորֆիզմի գաղափարը:

**Սահմանում 1.4.1:** $G$ և $H$ գրաֆները կոչվում են *իզոմորֆ*, եթե գոյություն ունի $f: V(G) \to V(H)$ փոխմիարժեք համապատասխանություն, որ $uv \in E(G)$ այն և միայն այն դեպքում, երբ $f(u)f(v) \in E(H)$: Եթե $G$ և $H$ գրաֆները իզոմորֆ են կգրենք $G \cong H$:

Դիտարկենք նկ. 1.4.1-ում պատկերված գրաֆները: Համոզվենք, որ այդ նկարում պատկերված գրաֆներից $G, H$ և $K$ գրաֆները զույգ առ զույգ իզոմորֆ են, իսկ $F$-ը՝ իզոմորֆ չէ այդ գրաֆներին: Իրոք, եթե դիտարկենք $f_1: V(G) \to V(H)$, $f_2: V(H) \to V(K)$ և $f_3: V(K) \to V(G)$ արտապատկերումները, որտեղ $f_1(u_1) = v_1$, $f_1(u_2) = v_3$, $f_1(u_3) = v_5$, $f_1(w_1) = v_2$, $f_1(w_2) = v_4$, $f_1(w_3) = v_6$, $f_2(v_i) = y_i$ $(1 \leq i \leq 6)$ և $f_3(y_1) = u_1$, $f_3(y_2) = w_1$, $f_3(y_3) = u_2$, $f_3(y_4) = w_2$, $f_3(y_5) = u_3$, $f_3(y_6) = w_3$, ապա պարզ է,

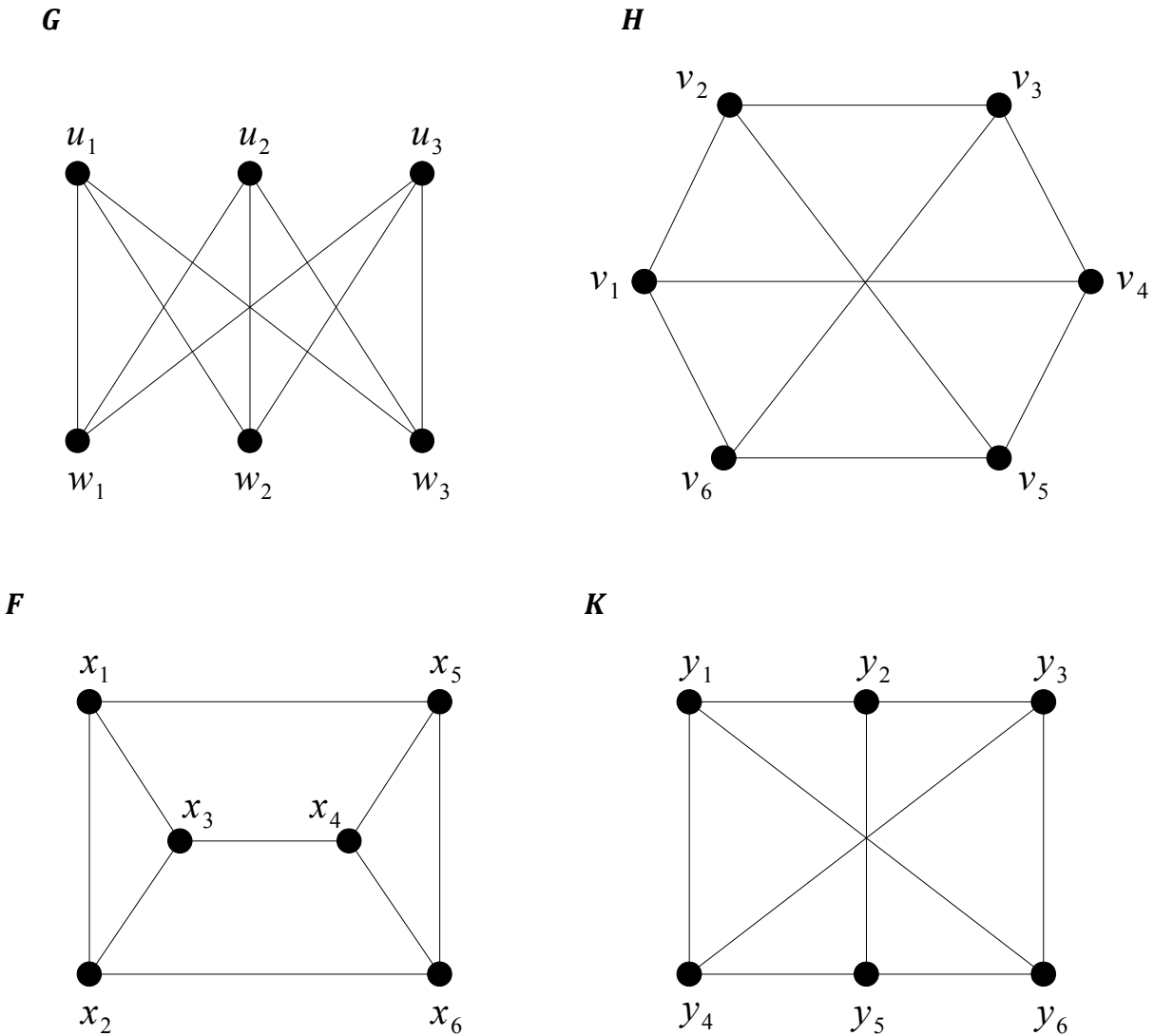

Նկ. 1.4.1

որ այդ արտապատկերումները հանդիսանում են փոխմիարժեք համապատասխանություններ $G, H$ և $K$ գրաֆների գագաթների միջև, որոնք ինչպես և պահանջվում է իզոմորֆիզմի սահմանումը, պահպանում են «կող լինելու» հատկությունը։ Նկատենք, որ $F$ գրաֆը պարունակում է $K_3$-ը որպես ենթագրաֆ, մինչդեռ $G, H$ և $K$ գրաֆները $K_3$ ենթագրաֆ չունեն։ Հետևաբար, $F$-ը իզոմորֆ չէ այդ գրաֆներից և ոչ մեկին։ Նկատենք, որ եթե $G \cong H$, ապա այդ գրաֆների գագաթները կարելի է այնպես վերահամարակալել, որ նրանց հարևանության մատրիցները համընկնեն։ Նշենք նաև, որ ցանկացած $G$ գրաֆի համար $G \cong G$, և, եթե $G \cong H$, ապա $H \cong G$, ինչպես, եթե $G \cong H$ և $H \cong F$, ապա $G \cong F$։ Այստեղից հետևում է, որ գրաֆների իզոմորֆիզմը գրաֆների բազմության վրա որոշված համարժեքության հարաբերություն է։ Հետագայում, որպես կանոն, իզոմորֆ գրաֆները միմյանցից չենք տարբերի։

Գրաֆների իզոմորֆիզմի հետ է կապված գրաֆների տեսության հայտնի և բարդ հիպոթեզներից մեկը, որը ձևակերպել են Կելլին և Ուլամը:

**Հիպոթեզ 1.4.1:** Դիցուք ունենք $G$ և $H$ գրաֆները, որտեղ $V(G) = \{u_1, \ldots, u_n\}$, $V(H) = \{v_1, \ldots, v_n\}$ և $n \geq 3$: Եթե ցանկացած $i$-ի ($1 \leq i \leq n$) համար տեղի ունի $G - u_i \cong H - v_i$ պայմանը, ապա $G \cong H$:

Նշենք, որ $n \geq 3$ պայմանը էական է այս հիպոթեզում, քանի որ $K_2$ և $\overline{K_2}$ գրաֆները իզոմորֆ չեն: Հայտնի է, որ Կելլիի և Ուլամի հիպոթեզը ճիշտ է ցիկլ չպարունակող և միակ ցիկլ պարունակող կապակցված գրաֆների համար [21,36]: Հարարիի կողմից առաջարկվել է այս հիպոթեզի կողային տարբերակը.

**Հիպոթեզ 1.4.2:** Դիցուք ունենք $G$ և $H$ գրաֆները, որտեղ $E(G) = \{e_1, \ldots, e_m\}$, $E(H) = \{f_1, \ldots, f_m\}$ և $m \geq 4$. Եթե ցանկացած $i$-ի ($1 \leq i \leq m$) համար տեղի ունի $G - e_i \cong H - f_i$ պայմանը, ապա $G \cong H$:

Նշենք, որ $m \geq 4$ պայմանը այստեղ նույնպես էական է, քանի որ $K_3 + v$ և $K_{1,3}$ գրաֆները իզոմորֆ չեն: Այս հիպոթեզի ապացույցի ուղղությամբ Լ. Լովասը և S. Մյուլլերը հասել են որոշ հաջողությունների [24,29]: Մասնավորապես, Լովասը ցույց է տվել, որ եթե $n$ գագաթ ունեցող $G$ գրաֆում $|E(G)| > \frac{1}{2}\binom{n}{2}$-ից, ապա Հարարիի հիպոթեզը ճիշտ է:

**Սահմանում 1.4.2:** $G$ գրաֆի *հոմոմորֆիզմ* $H$ գրաֆի վրա կոչվում է $f: V(G) \to V(H)$ արտապատկերումը, որի դեպքում եթե $uv \in E(G)$, ապա $f(u)f(v) \in E(H)$: Եթե գոյություն ունի $f$ $G$ գրաֆի հոմոմորֆիզմ $H$ գրաֆի վրա, ապա կգրենք $f: G \to H$ կամ $G \to H$:

Դիտարկենք նկ. 1.4.2-ում պատկերված գրաֆները: Համոզվենք, որ նկ. 1.4.2-ում պատկերված գրաֆների համար $G \to H$ և $F \to K$: Իրոք, եթե դիտարկենք $f_1: V(G) \to V(H)$ և $f_2: V(F) \to V(K)$ արտապատկերումները, որտեղ $f_1(v_1) = u_1$, $f_1(v_2) = u_2$, $f_1(v_3) = u_1$, $f_1(v_4) = u_3$, $f_1(v_5) = u_1$, $f_1(v_6) = u_2$, $f_1(v_7) = u_1$, $f_1(v_8) = u_3$ և $f_2(x_1) = y_1$, $f_2(x_2) = y_2$, $f_2(x_3) = y_1$, $f_2(x_4) = y_2$, ապա պարզ է, որ $f_1$-ը $G$ գրաֆի հոմոմորֆիզմ է $H$ գրաֆի վրա, իսկ $f_2$-ը՝ $F$ գրաֆի հոմոմորֆիզմ է $K$ գրաֆի վրա:

Նկատենք, որ եթե $f: V(G) \to V(H)$-ը փոխմիարժեք համապատասխանություն է, $f: G \to H$ և $f^{-1}: H \to G$, ապա $G \cong H$: Նշենք նաև, որ ցանկացած $G$ գրաֆի համար $G \to G$, և եթե $G \to H$ և $H \to F$, ապա $G \to F$: Գրաֆների հոմոմորֆիզմի մասին ավելի մանրամասն կարելի է ծանոթանալ Հելլի և Նեշետրիլի գրքում [19]:

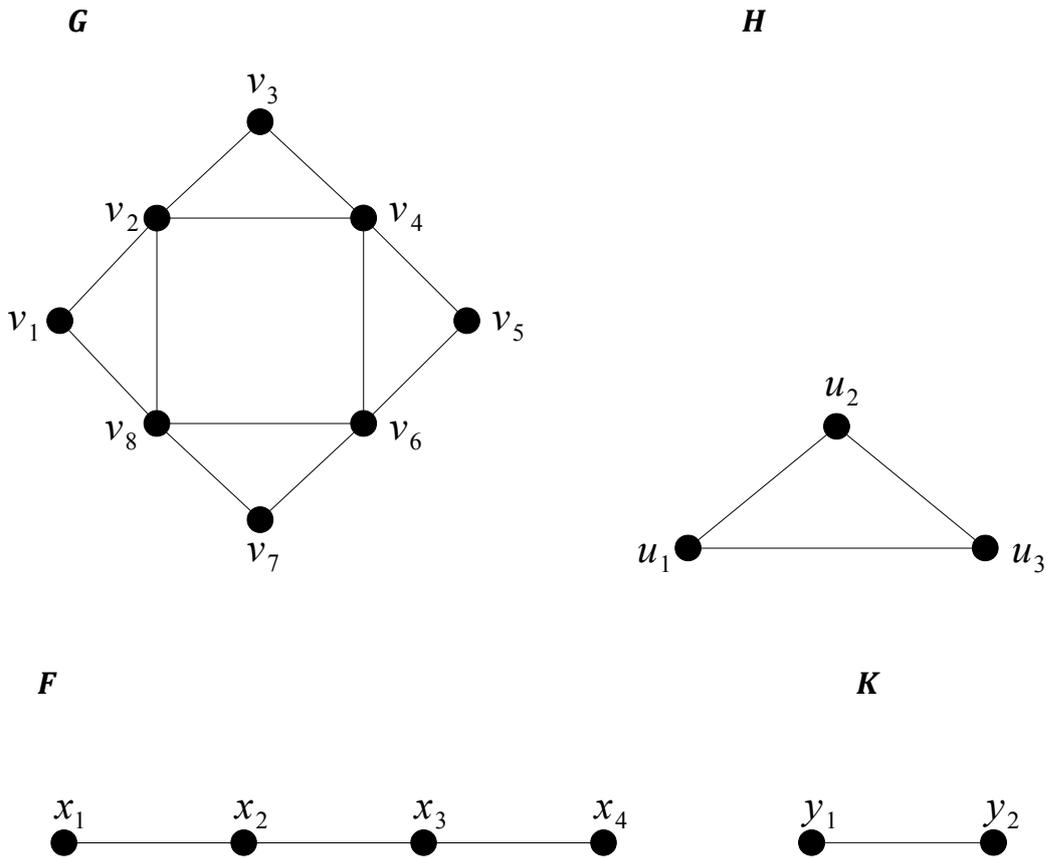

Նկ. 1.4.2

Այժմ սահմանենք գրաֆների ավտոմորֆիզմի գաղափարը:

**Սահմանում 1.4.3:** $G$ գրաֆի *ավտոմորֆիզմ* կոչվում է $G$ գրաֆի իզոմորֆիզմը իր վրա:

Պարզ է, որ $G$ գրաֆի ավտոմորֆիզմն այդ գրաֆի գագաթների այնպիսի տեղադրություն է, որի դեպքում պահպանվում է գագաթների հարևանությունը: $G$ գրաֆի բոլոր ավտոմորֆիզմների բազմությունը նշանակենք $Aut(G)$-ով: Նկատենք, որ ցանկացած $G$ գրաֆի համար $Aut(G) \neq \emptyset$, քանի որ $G$-ն ունի նույնական ավտոմորֆիզմ $f: V(G) \to V(G)$, որի դեպքում $f(v) = v$ $G$ գրաֆի ցանկացած $v$ գագաթի համար: Պարզ է, որ եթե $f$-ը $G$-ի ավտոմորֆիզմ է, ապա $f^{-1}$-ը նաև $G$-ի ավտոմորֆիզմ է, և, եթե $f$-ը և $g$-ն $G$-ի գագաթների հարևանությունը պահպանող տեղադրություններ են (ավտոմորֆիզմներ են), ապա $f \cdot g$ տեղադրությունը նաև գագաթների հարևանությունը պահպանող տեղադրություն է (ավտոմորֆիզմ է): Այստեղից հետևում է, որ ցանկացած $G$ գրաֆի համար $Aut(G)$-ն խումբ է տեղադրությունների բազմապատկման նկատմամբ: Հեշտ է տեսնել, որ $Aut(K_n) = S_n$: Պարզվում է, որ տեղի ունի ավելի ընդհանուր փաստ.

**Թեորեմ 1.4.1 (Ռ. Ֆրուխտ):** Ցանկացած վերջավոր $F$ խմբի համար գոյություն ունի $G$ գրաֆ, որի $Aut(G)$-ն իզոմորֆ է $F$-ին:

**Սահմանում 1.4.4:** $G$ գրաֆի $u$ և $v$ գագաթները կոչվում են *նման*, եթե գոյություն ունի $G$ գրաֆի այնպիսի $f$ ավտոմորֆիզմ, որի դեպքում $f(u) = v$:

**Սահմանում 1.4.5:** $G$ գրաֆի $e = uv$ և $e' = xy$ կողերը կոչվում են *նման*, եթե գոյություն ունի $G$ գրաֆի այնպիսի $f$ ավտոմորֆիզմ, որի դեպքում $f(u)f(v) = xy$:

**Սահմանում 1.4.6:** $G$ գրաֆը կոչվում է *գագաթային սիմետրիկ գրաֆ*, եթե նրա ցանկացած երկու գագաթներ նման են:

**Սահմանում 1.4.7:** $G$ գրաֆը կոչվում է *կողային սիմետրիկ գրաֆ*, եթե նրա ցանկացած երկու կողեր նման են:

Դիտարկենք նկ. 1.4.3-ում պատկերված գրաֆները:

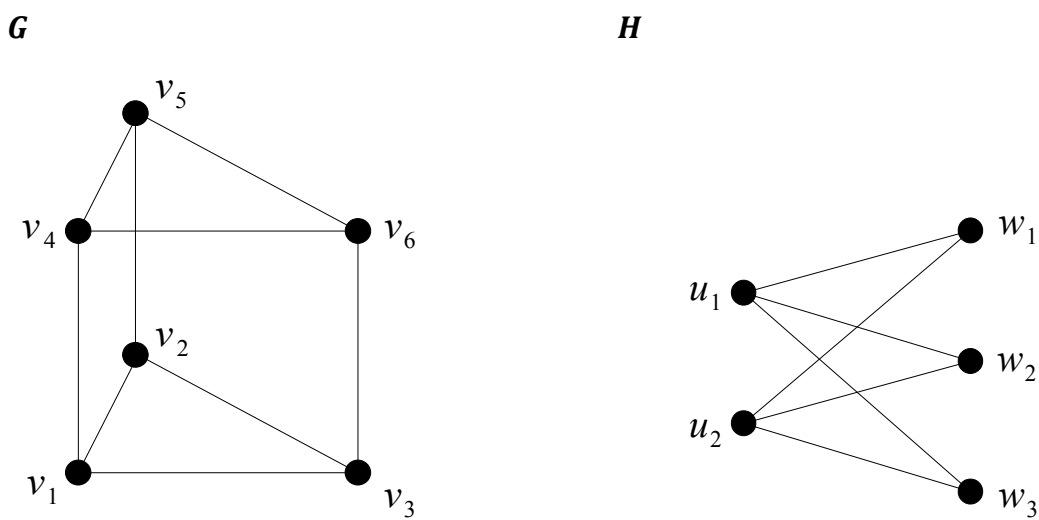

նկ. 1.4.3

Հեշտ է տեսնել, որ նկ. 1.4.3-ում պատկերված գրաֆներից $G$-ն գագաթային սիմետրիկ գրաֆ է, բայց կողային սիմետրիկ գրաֆ չէ, իսկ $H$-ը՝ կողային սիմետրիկ է, բայց գագաթային սիմետրիկ չէ:

**Սահմանում 1.4.8:** Եթե $G$-ն ունի միայն նույնական ավտոմորֆիզմ, ապա $G$-ն կոչվում է *ասիմետրիկ գրաֆ*:

Դիտարկենք նկ. 1.4.4-ում պատկերված $G$ գրաֆը:

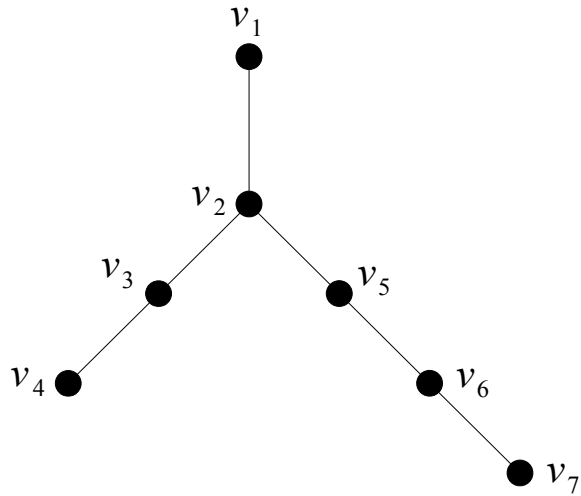

Նկ. 1.4.4

Հեշտ է տեսնել, որ նկ. 1.4.4-ում պատկերված $G$ գրաֆը ասիմետրիկ գրաֆ է։

Եզրափակելով պարագրաֆը, բերենք երկու հայտնի գագաթային և կողային սիմետրիկ գրաֆների պատկերները։

**Պետերսենի գրաֆ**.                **Հիվուդի գրաֆ**.

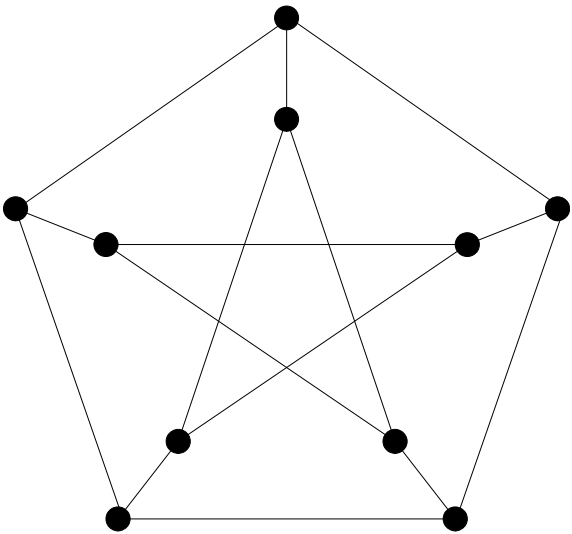
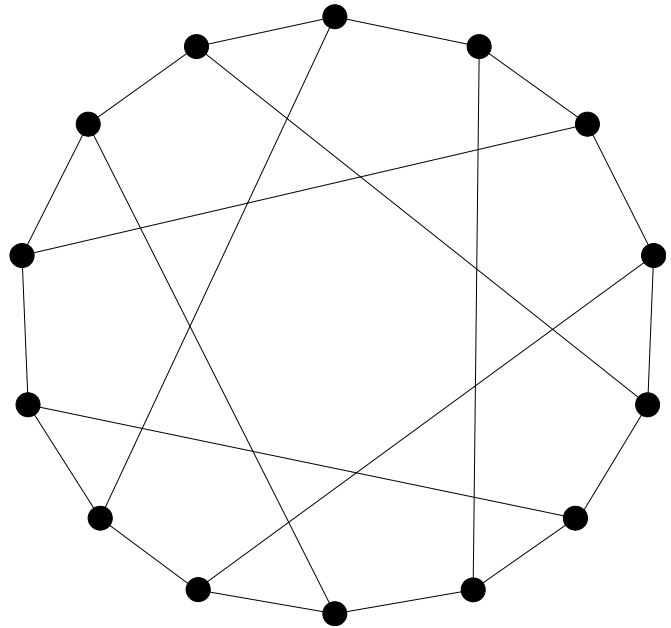

Նկ. 1.4.5

## § 1.5. Էքստրեմալ, կցության և հատումների գրաֆներ

Այս պարագրաֆում մենք կդիտարկենք էքստրեմալ, կցության և հատումների գրաֆներ։ Էքստրեմալ գրաֆների տեսությունում հիմնականում ուսումնասիրվում են ամենամեծ կամ ամենափոքր (*Էքստրեմալ*) գրաֆները, որոնք օժտված են որոշակի հատկությամբ:

Կասենք, որ $G$ գրաֆը չի պարունակում $F$ գրաֆ, եթե $F$ գրաֆը $G$ գրաֆի ենթագրաֆ չէ: Եթե $F = C_3$ կամ $F = K_3$, ապա ասում են, որ $G$ գրաֆը չի պարունակում եռանկյուն:

Էքստրեմալ գրաֆների տեսության առաջին թեորեմ է համարվում Մանթելի կողմից ապացուցված թեորեմը եռանկյուն չպարունակող գրաֆների մասին:

**Թեորեմ 1.5.1:** Եթե $G$ $(n,m)$-գրաֆը չի պարունակում եռանկյուն, ապա $m \leq \left\lfloor\frac{n^2}{4}\right\rfloor$:

**Ապացույց 1:** Դիցուք $V(G) = \{v_1, ..., v_n\}$: Քանի, որ $G$ գրաֆը չի պարունակում եռանկյուն, ուստի ցանկացած $uv \in E(G)$-ի համար տեղի ունի հետևյալը. $N_G(u) \cap N_G(v) = \emptyset$: Այստեղից հետևում է, որ ցանկացած $uv \in E(G)$-ի համար $d_G(u) + d_G(v) \leq n$: Գումարելով այս անհավասարությունը ըստ բոլոր $uv$ կողերի, կստանանք.

$$\sum_{uv \in E(G)} \big(d_G(u) + d_G(v)\big) = \sum_{v \in V(G)} \big(d_G(v)\big)^2 \leq n \cdot m:$$

Այժմ դիտարկենք $\overline{d} = \big(d_G(v_1), ..., d_G(v_n)\big)$ և $\overline{1} = (1, ..., 1)$ վեկտորները $\mathbb{R}^n$-ից: Ըստ թեորեմ 1.2.1-ի և Կոշու-Բունյակովսկու անհավասարության, կստանանք.

$$(2m)^2 = \left(\sum_{v \in V(G)} d_G(v)\right)^2 = (\overline{d}, \overline{1})^2 \leq (\overline{d}, \overline{d}) \cdot (\overline{1}, \overline{1}) = n \cdot \sum_{v \in V(G)} \big(d_G(v)\big)^2:$$

Այստեղից, հաշվի առնելով $\sum_{v \in V(G)} \big(d_G(v)\big)^2 \leq n \cdot m$ անհավասարությունը, ստանում ենք հետևյալը.

$$(2m)^2 \leq n \sum_{v \in V(G)} \big(d_G(v)\big)^2 \leq n^2 \cdot m,$$

ուստի $m \leq \frac{n^2}{4}$: ∎

**Ապացույց 2:** Ապացույցը կատարենք մաթմաթման եղանակով ըստ $n$-ի: Մենք կապացուցենք թեորեմը զույգ $n$-ի համար, կենտ $n$-ի համար ապացույցը կատարվում է համանման ձևով: Դիցուք $n = 2k$: Հեշտ է տեսնել, որ թեորեմը ճիշտ է $n \leq 4$-ի դեպքում: Ենթադրենք, որ թեորեմը ճիշտ է ցանկացած $G'$ գրաֆի համար, որը չի պարունակում եռանկյուն և որի գագաթների քանակը $n = 2k$-ից մեծ չէ: Դիտարկենք $G$ $(n,m)$-գրաֆը, որը չի պարունակում եռանկյուն և $n = 2k + 2$: Մենք կարող ենք ենթադրել, որ գոյություն ունի $e = uv$ կող $G$ գրաֆում, քանի որ հակառակ դեպքում թեորեմն ակնհայտ է: Դիցուք $G' = G - u - v$: Պարզ է, որ $G'$ գրաֆը չի պարունակում եռանկյուն և $|V(G')| = 2k$: Ըստ մակածման ենթադրությամբ, կստանանք $|E(G')| \leq \left\lfloor \frac{4k^2}{4} \right\rfloor = k^2$: Մյուս կողմից, քանի որ $G$ գրաֆը չի պարունակում եռանկյուն, $u$ և $v$ գագաթների համար տեղի ունի հետևյալը. $N_G(u) \cap N_G(v) = \emptyset$: Այստեղից հետևում է, որ $u$ և $v$ գագաթների համար $d_G(u) + d_G(v) \leq 2k + 2$: Այժմ գնահատենք $G$ գրաֆի կողերի քանակը.

$$m = |E(G)| \leq |E(G')| + d_G(u) + d_G(v) - 1 \leq k^2 + 2k + 1 = \left\lfloor \frac{(2k+2)^2}{4} \right\rfloor,$$

ուստի $m \leq \left\lfloor \frac{n^2}{4} \right\rfloor$: ∎

**Հետևանք 1.5.1:** Կամայական $G\left(n, \left\lfloor \frac{n^2}{4} \right\rfloor + 1\right)$-գրաֆ ($n \geq 3$) պարունակում է եռանկյուն:

Ցույց տանք, որ թեորեմ 1.5.1-ի վերին գնահատականը հասանելի է: Իրոք, եթե մենք դիտարկենք լրիվ երկկողմանի $K_{\left\lfloor \frac{n}{2} \right\rfloor, \left\lceil \frac{n}{2} \right\rceil}$ գրաֆը, ապա հեշտ է տեսնել, որ $\left|V\left(K_{\left\lfloor \frac{n}{2} \right\rfloor, \left\lceil \frac{n}{2} \right\rceil}\right)\right| = n$ և $\left|E\left(K_{\left\lfloor \frac{n}{2} \right\rfloor, \left\lceil \frac{n}{2} \right\rceil}\right)\right| = \left\lfloor \frac{n}{2} \right\rfloor \cdot \left\lceil \frac{n}{2} \right\rceil = \left\lfloor \frac{n^2}{4} \right\rfloor$:

Այժմ ապացուցենք Ռեյմանի կողմից ստացված բավարար պայմանը գրաֆում չորս երկարությամբ ցիկլի գոյության համար:

**Թեորեմ 1.5.2:** Եթե $G$ գրաֆը բավարարում է $\sum_{v \in V(G)} \binom{d_G(v)}{2} > \binom{n}{2}$ պայմանին, ապա $G$-ն պարունակում է չորս երկարությամբ ցիկլ:

**Ապացույց:** Նշանակենք $p_2(v)$-ով երկու երկարություն ունեցող այն ճանապարհների քանակը $G$ գրաֆում, որոնց կենտրոնական գագաթը $v$-ն է: Պարզ է, որ ցանկացած $v$ գագաթի համար $p_2(v)$-ն հավասար է $\binom{d_G(v)}{2}$-ի: Քանի որ երկու երկարություն ունեցող յուրաքանչյուր ճանապարհ $G$ գրաֆում ունի միակ կենտրոնական գագաթ, ուստի $G$ գրաֆի բոլոր երկու երկարություն ունեցող ճանապարհների քանակը կլինի. $\sum_{v \in V(G)} p_2(v) = \sum_{v \in V(G)} \binom{d_G(v)}{2}$:

Մյուս կողմից, ամեն մի այդպիսի ճանապարհի ունի ճիշտ երկու ծայրակետ, հետևաբար, բոլոր երկու երկարություն ունեցող ճանապարհները կարելի է տրոհել $\binom{n}{2}$ դասերի, ըստ այդ ճանապարհների ծայրակետերի: Քանի որ $\sum_{v \in V(G)} \binom{d_G(v)}{2} > \binom{n}{2}$, ուստի այդ դասերից մեկը պարունակում է առնվազն երկու հատ երկու երկարություն ունեցող տարբեր ճանապարհներ միևնույն ծայրակետերով, իսկ դա նշանակում է, որ $G$-ն պարունակում է չորս երկարությամբ ցիկլ: ∎

Ռեյմանի կողմից նաև ապացուցվել է թեորեմ $C_4$ չպարունակող գրաֆների մասին։

**Թեորեմ 1.5.3։** Եթե $G$ $(n,m)$-գրաֆը չի պարունակում $C_4$, ապա $m \leq \frac{n}{4}\left(1 + \sqrt{4n-3}\right)$։

**Ապացույց։** Դիցուք $V(G) = \{u_1, \ldots, u_n\}$։ Նշանակենք $F$-ով $G$ գրաֆի բոլոր կարգավորված հետնյալ եռյակների բազմությունը. $F = \{(u,v,w): uv \in E(G), uw \in E(G)$ և $v \neq w\}$։ Հաշվենք $|F|$-ը։ Պարզ է, որ յուրաքանչյուր $u$ գագաթի ներդրումը $|F|$-ի մեջ կլինի. $d_G(u)(d_G(u)-1)$, ուստի

$$|F| = \sum_{u \in V(G)} d_G(u)(d_G(u) - 1)։$$

Մյուս կողմից, ամեն մի կարգավորված $(v,w)$ զույգ կարող է մասնակցել ամենաշատը մի $(u,v,w)$ կարգավորված եռյակի մեջ, հակառակ դեպքում՝ $G$-ն կպարունակի $C_4$։ Այստեղից հետևում է, որ $|F| \leq n(n-1)$։ Հետևաբար,

$$n(n-1) \geq |F| = \sum_{u \in V(G)} d_G(u)(d_G(u) - 1) = \sum_{u \in V(G)} (d_G(u))^2 - \sum_{u \in V(G)} d_G(u)։$$

Ըստ թեորեմ 1.2.1-ի, կստանանք

$$n(n-1) \geq \sum_{u \in V(G)} (d_G(u))^2 - 2m։$$

Դիտարկենք $\overline{d} = (d_G(u_1), \ldots, d_G(u_n))$ և $\overline{1} = (1, \ldots, 1)$ վեկտորները $\mathbb{R}^n$-ից։ Թեորեմ 1.2.1-ից և Կոշու-Բունյակովսկու անհավասարությունից, կստանանք.

$$(2m)^2 = \left(\sum_{u \in V(G)} d_G(u)\right)^2 = (\overline{d}, \overline{1})^2 \leq (\overline{d}, \overline{d}) \cdot (\overline{1}, \overline{1}) = n \sum_{u \in V(G)} (d_G(u))^2։$$

Այստեղից, հաշվի առնելով $n(n-1) \geq \sum_{u \in V(G)} (d_G(u))^2 - 2m$ անհավասարությունը, ստանում ենք հետնյալը.

$$n(n-1) \geq \sum_{u \in V(G)} (d_G(u))^2 - 2m \geq \frac{(2m)^2}{n} - 2m։$$

Վերջին անհավասարությունից գալիս ենք $4m^2 - 2mn - n^2(n-1) \leq 0$ քառակուսային անհավասարմանը ըստ $m$-ի, որը լուծելիս ստանում ենք $m \leq \frac{n}{4}\left(1 + \sqrt{4n-3}\right)$ անհավասարությունը։ ∎

Հայտնի է, որ թեորեմ 1.5.3-ի վերին գնահատականը հասանելի է վերջավոր պրոյեկտիվ երկրաչափության կցության գրաֆների վրա։

Այժմ սահմանենք *Տուրանի $ex(n,F)$ թիվը* հետնյալ կերպ.

$$ex(n, F) = max_G\{|E(G)|: |V(G)| = n \text{ և } F \nsubseteq G\}:$$

Այլ կերպ ասած, $ex(n, F)$-ը այն կողերի ամենամեծ քանակն է, որը կարող է ունենալ $n$ գագաթ ունեցող և $F$ գրաֆ չպարունակող գրաֆը:

Սահմանենք *Տուրանի $T_{n,r}$ գրաֆը*: Տուրանի $T_{n,r}$ գրաֆը իրենից ներկայացնում է $n$ գագաթ ունեցող լրիվ $r$-կողմանի գրաֆ, որի $n - r\left\lfloor\frac{n}{r}\right\rfloor$ կողմերը պարունակում են $\left\lceil\frac{n}{r}\right\rceil$ հատ գագաթներ, իսկ բոլոր մնացած կողմերը պարունակում են $\left\lfloor\frac{n}{r}\right\rfloor$ հատ գագաթներ: Նկ. 1.5.1-ում պատկերված է Տուրանի $T_{7,3}$ գրաֆը: Նկատենք, որ $|E(T_{n,r})| \leq \left(1 - \frac{1}{r}\right)\frac{n^2}{2}$:

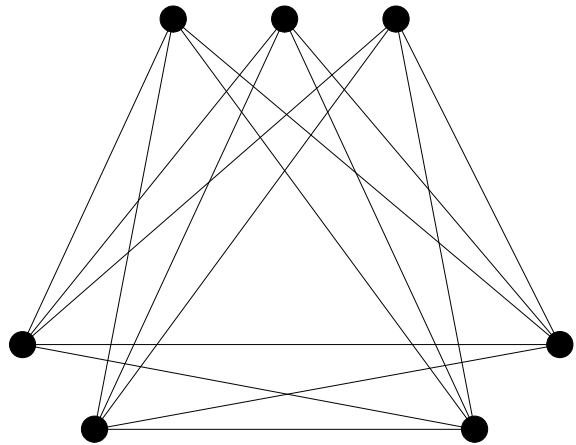

Նկ. 1.5.1

Այժմ ձևակերպենք էքստրեմալ գրաֆների տեսության դասական արդյունքներից մեկը. Տուրանի թեորեմը $K_{r+1}$ չպարունակող գրաֆների մասին:

**Թեորեմ 1.5.4:** $n$ գագաթ ունեցող և $K_{r+1}$ չպարունակող գրաֆներից Տուրանի $T_{n,r}$ գրաֆը ունի ամենաշատ կողեր. $ex(n, K_{r+1}) = |E(T_{n,r})|:$

**Ապացույց:** Նախ նկատենք, որ Տուրանի $T_{n,r}$ գրաֆը չի պարունակում $K_{r+1}$ գրաֆը: Ցույց տանք, որ $n$ գագաթ ունեցող բոլոր $r$-կողմանի գրաֆներից, Տուրանի $T_{n,r}$ գրաֆը ամենաշատ կողերն ունի: Նկատենք, որ $n$ գագաթ ունեցող $r$-կողմանի գրաֆներից ամենաշատ կողերն ունի $n$ գագաթ ունեցող լրիվ $r$-կողմանի գրաֆը, քանի որ եթե $r$-կողմանի գրաֆը լրիվ չէ, ապա մենք կարող ենք այնպես ավելացնել կողեր այդ գրաֆին, որ ստացված գրաֆը ևս լինի $r$-կողմանի: Այստեղից հետևում է, որ բավական է ցույց տալ, որ $n$ գագաթ ունեցող լրիվ $r$-կողմանի գրաֆներից ամենաշատ կողերն ունի Տուրանի $T_{n,r}$ գրաֆը:

Դիցուք $K_{n_1,n_2,\cdots,n_r}$-ը լրիվ $r$-կողմանի գրաֆ է և $\sum_{i=1}^{r} n_i = n$: Ցույց տանք, որ $|n_i - n_j| \leq 1$ ցանկացած $i, j \in \{1, \ldots, r\}$-ի համար: Ենթադրենք հակառակը՝ գոյություն ունեն

$i_0, j_0 \in \{1, \ldots, r\}$-ի, որ $n_{i_0} \geq n_{j_0} + 2$։ Դիտարկենք $K_{n_1, \cdots, n_{i_0}-1, \cdots, n_{j_0}+1, \cdots, n_r}$ գրաֆը։ Պարզ է, որ $K_{n_1, \cdots, n_{i_0}-1, \cdots, n_{j_0}+1, \cdots, n_r}$ գրաֆը ևս լրիվ $r$-կողմանի գրաֆ է և $\sum_{i=1}^{r} n_i = n$, բայց այդ գրաֆի կողերի քանակը հավասար կլինի․

$$\left|E\left(K_{n_1, \cdots, n_{i_0}-1, \cdots, n_{j_0}+1, \cdots, n_r}\right)\right| = \left|E(K_{n_1, n_2, \cdots, n_r})\right| - n_{j_0} + n_{i_0} - 1 > \left|E(K_{n_1, n_2, \cdots, n_r})\right|,$$

որը հակասություն է։

Այժմ ցույց տանք, որ եթե $n$ գագաթ ունեցող $G$ գրաֆը չի պարունակում $K_{r+1}$, ապա գոյություն ունի $H$ լրիվ $r$-կողմանի գրաֆ, որում $V(H) = V(G)$ և $|E(H)| \geq |E(G)|$։ Ապացույցը կատարենք մակածման եղանակով ըստ $r$-ի։ Եթե $r = 1$, ապա $|E(G)| = 0$ և վերցնելով $H = G$-ի, կստանանք անհրաժեշտ $H$ գրաֆի գոյությունը։ Ենթադրենք $r \geq 2$ և պնդումը ճիշտ է ցանկացած գրաֆի համար, որը չի պարունակում $K_r$։ Դիցուք $G$-ն $n$ գագաթ ունեցող գրաֆ է, որը չի պարունակում $K_{r+1}$ և $x \in V(G)$, $d_G(x) = \Delta(G)$։ Դիտարկենք $G$ գրաֆի ծնված $G' = G[N_G(x)]$ ենթագրաֆը։ Քանի որ $G$ գրաֆի $x$ գագաթը հարևան է $G'$-ի բոլոր գագաթներին, ուստի $G'$-ը չի պարունակում $K_r$։ Ըստ մակածման ենթադրության, կստանանք, որ գոյություն ունի $H'$ լրիվ $(r-1)$-կողմանի գրաֆ, որում $V(H') = V(G') = N_G(x)$ և $|E(H')| \geq |E(G')|$։

Կառուցենք $H$ լրիվ $r$-կողմանի գրաֆ, որում $V(H) = V(G)$ և $|E(H)| \geq |E(G)|$։ Դիցուք $S = V(G) \setminus N_G(x)$։ Սահմանենք $H$ գրաֆը հետևյալ կերպ․ $V(H) = S \cup V(H')$, $E(H) = E(H') \cup \{uv : u \in S, v \in N_G(x)\}$ (նկ. 1.5.2)։ Քանի որ $S$-ի գագաթները զույգ առ զույգ հարևան չեն $H$ գրաֆում, ուստի $H$-ը լրիվ $r$-կողմանի գրաֆ է։ Ապացուցենք, որ $|E(H)| \geq |E(G)|$։

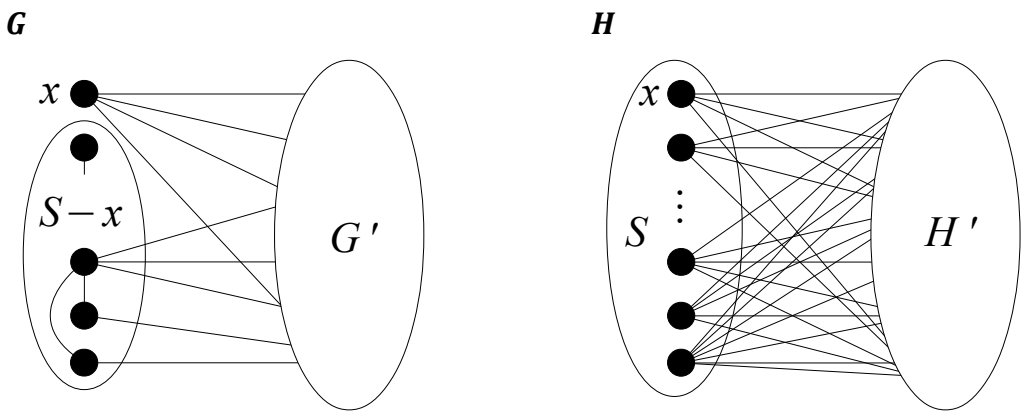

Նկ. 1.5.2

$H$ գրաֆի սահմանումից հետևում է, որ

$$|E(H)| = |E(H')| + |S| \cdot |N_G(x)| = |E(H')| + \Delta(G)(n - \Delta(G))։$$

Մյուս կողմից, հեշտ է տեսնել, որ

$$|E(G)| \leq |E(G')| + \sum_{v \in S} d_G(v) \leq |E(G')| + |S| \cdot \Delta(G) = |E(G')| + \Delta(G)(n - \Delta(G)):$$

Այստեղից, հաշվի առնելով $|E(H')| \geq |E(G')|$ անհավասարությունը, ստանում ենք հետևյալը.

$$|E(H)| = |E(H')| + \Delta(G)(n - \Delta(G)) \geq |E(G')| + \Delta(G)(n - \Delta(G)) \geq |E(G)|: \blacksquare$$

**Հետևանք 1.5.2:** Ցանկացած $G\ (n, |E(T_{n,r})| + 1)$-գրաֆ $(n \geq r + 1)$ պարունակում է $K_{r+1}:$

Նշենք, որ Տուրանի թեորեմը ընդհանրացնում է Մանթելի թեորեմը, քանի որ $ex(n, K_3) = ex(n, C_3) = \left\lfloor \frac{n^2}{4} \right\rfloor:$ Ռեյմանի թեորեմից ստանում ենք, որ $ex(n, C_4) \leq \frac{n}{4}(1 + \sqrt{4n-3}):$ Հայտնի է նաև, որ $ex(n, C_{2k+1}) = O(n^2)$, իսկ $ex(n, C_{2k}) = O\left(n^{1+\frac{1}{k}}\right):$ Էքստրեմալ գրաֆների տեսությանն ավելի մանրամասն կարելի է ծանոթանալ Բոլլոբաշի գրքում [7]:

Այժմ անցնենք կցության գրաֆներին: Դիցուք տրված է $S = \{s_1, s_2, \dots, s_n\}$ բազմությունը և այդ բազմության ենթաբազմությունների $\mathfrak{F} = \{F_1, F_2, \dots, F_m\}$ ընտանիքը: Սահմանենք $(S, \mathfrak{F})$ զույգի համար հետևյալ *կցության* $G(S, \mathfrak{F})$ *գրաֆը*.

$$V(G(S, \mathfrak{F})) = S \cup \{F_1, F_2, \dots, F_m\}\ \text{և}$$

$$E(G(S, \mathfrak{F})) = \{s_i F_j : s_i \in S, F_j \in \mathfrak{F}\ \text{և}\ s_i \in F_j, 1 \leq i \leq n, 1 \leq j \leq m\}:$$

Հեշտ է տեսնել, որ $G(S, \mathfrak{F})$ կցության գրաֆը երկկողմանի գրաֆ է:

Դիցուք $S = \{1, 2, 3, 4, 5\}$ և $\mathfrak{F} = \{\{1, 2, 3\}, \{1, 3, 4\}, \{2, 3, 5\}, \{4, 5\}\}:$ Ստորև պատկերված է $G(S, \mathfrak{F})$ կցության գրաֆը:

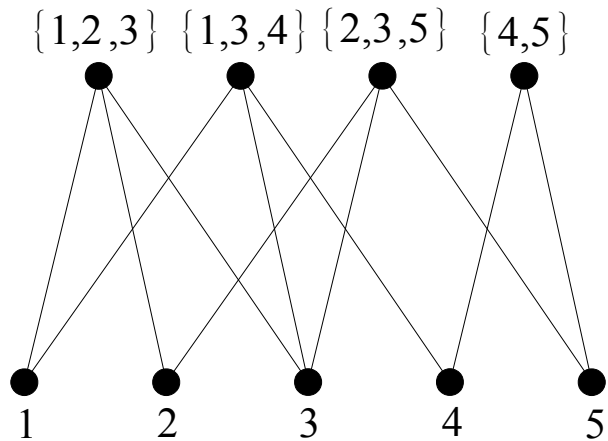

Նկ. 1.5.3

Նկատենք, որ յուրաքանչյուր $H = (V, E)$ հիպերգրաֆին ևս կարելի է համապատասխանեցնել $G(V, E)$ կցության գրաֆ։ Կցության գրաֆների հետաքրքիր օրինակներ են հանդիսանում վերջավոր աֆինական և պրոյեկտիվ երկրաչափության կցության գրաֆները։ Նկ. 1.5.4-ում պատկերված են վերջավոր երկրորդ կարգի աֆինական և պրոյեկտիվ երկրաչափության կցության գրաֆները։

Վերջում անդրադառնանք նաև հատումների գրաֆներին։ Դիցուք տրված է $S$ բազմությունը և այդ բազմության ենթաբազմությունների $\mathfrak{F} = \{F_1, F_2, \ldots, F_m\}$ ընտանիքը։ Սահմանենք $(S, \mathfrak{F})$ զույգի համար հատումների հետևյալ $\Omega(S, \mathfrak{F})$ գրաֆը․

$$V(\Omega(S, \mathfrak{F})) = \mathfrak{F} = \{F_1, F_2, \ldots, F_m\} \text{ և}$$

$$E(\Omega(S, \mathfrak{F})) = \{F_i F_j : F_i, F_j \in \mathfrak{F} \text{ և } F_i \cap F_j \neq \emptyset, 1 \leq i \neq j \leq m\}։$$

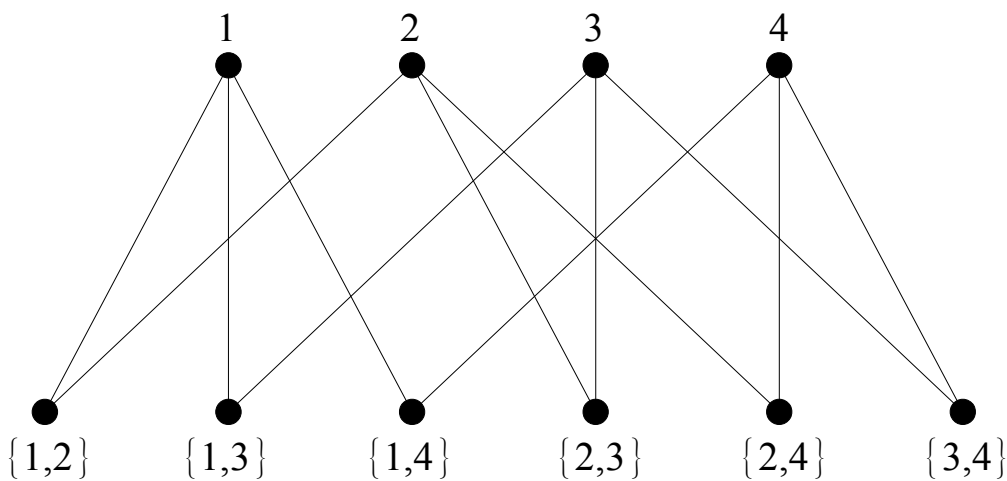

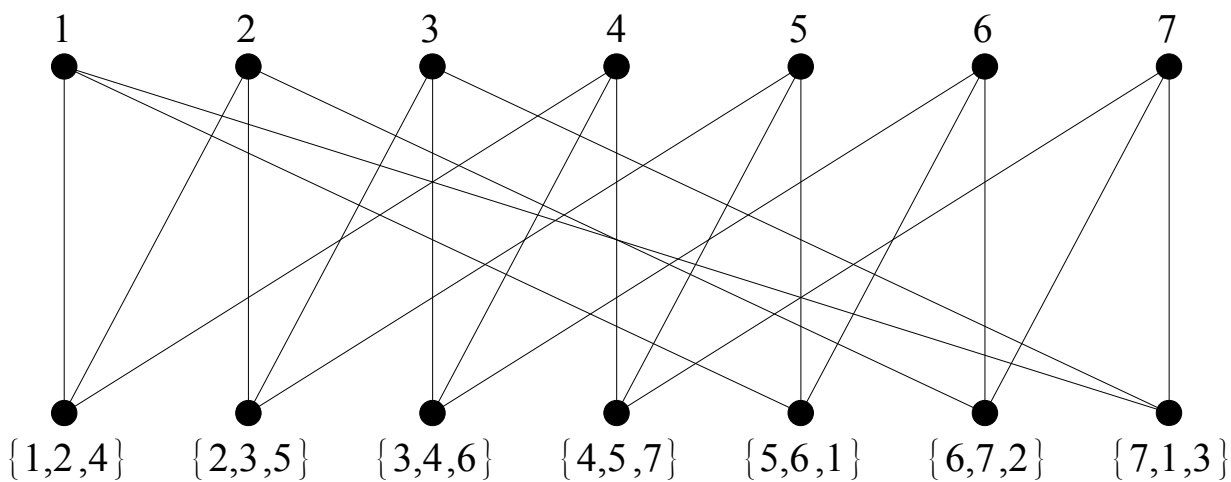

Նկ. 1.5.4.

**Սահմանում 1.5.1:** $G$ գրաֆը կոչվում է *հատումների գրաֆ*, եթե գոյություն ունի $S = \{s_1, s_2, \ldots, s_n\}$ բազմություն և այդ բազմության ենթաբազմությունների $\mathfrak{F} = \{F_1, F_2, \ldots, F_m\}$ ընտանիք, որ $G \cong \Omega(S, \mathfrak{F})$։

Դիցուք $S = \{1, 2, 3, 4, 5\}$ և $\mathfrak{F} = \{\{1, 2, 3\}, \{1, 3, 4\}, \{2, 3, 5\}, \{4, 5\}\}$։ Ստորև պատկերված

է $\Omega(S, \mathfrak{F})$ հատումների գրաֆը:

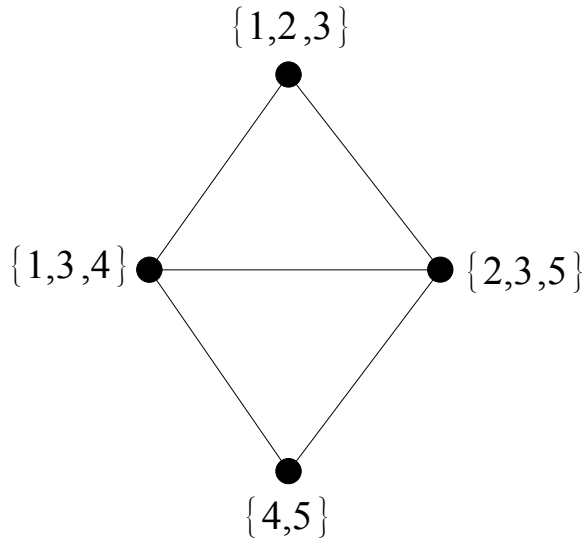

Նկ. 1.5.5.

Հատումների գրաֆներին վերաբերվող առաջին արդյունքը ստացվել է Մարչևսկու կողմից:

**Թեորեմ 1.5.5:** Կամայական գրաֆ հանդիսանում է հատումների գրաֆ:

**Ապացույց:** Դիցուք $G = (V, E)$ գրաֆ է, որտեղ $V(G) = \{v_1, v_2, \ldots, v_p\}$: Նկատենք, որ թեորեմը ապացուցելու համար բավական է կառուցել $S$ բազմությունը և այդ բազմության ենթաբազմությունների այնպիսի $\mathfrak{F}$ ընտանիք, որ բավարարվի $G \cong \Omega(S, \mathfrak{F})$ պայմանը: Սահմանենք $S$ բազմությունը և այդ բազմության ենթաբազմությունների $\mathfrak{F}$ ընտանիքը հետևյալ կերպ.

$$\mathfrak{F} = \{F_1, F_2, \ldots, F_p\}, \text{որտեղ } F_i = \{v_i\} \cup \partial_G(v_i) \ (1 \leq i \leq p) \text{ և}$$
$$S = V(G) \cup E(G):$$

Անմիջականորեն ստուգվում է, որ $G \cong \Omega(S, \mathfrak{F})$: ∎

Հատումների գրաֆների հետաքրքիր օրինակներ են հանդիսանում կողային և միջակայքների գրաֆները: $G = (V, E)$ գրաֆի համար դիտարկենք $\Omega(V, E)$ հատումների գրաֆը. այն անվանում են $G$ գրաֆի *կողային գրաֆ* և նշանակում են $L(G)$-ով: Հեշտ է տեսնել, որ $L(G)$ գրաֆում գագաթներին համապատասխանում են $G$ գրաֆի կողերը և $L(G)$ գրաֆի երկու գագաթները հարևան են, եթե նրանց համապատասխանող կողերը $G$ գրաֆում հարևան են: Նկ. 1.5.6-ում պատկերված է $G$ գրաֆը և նրա կողային $L(G)$ գրաֆի օրինակը:

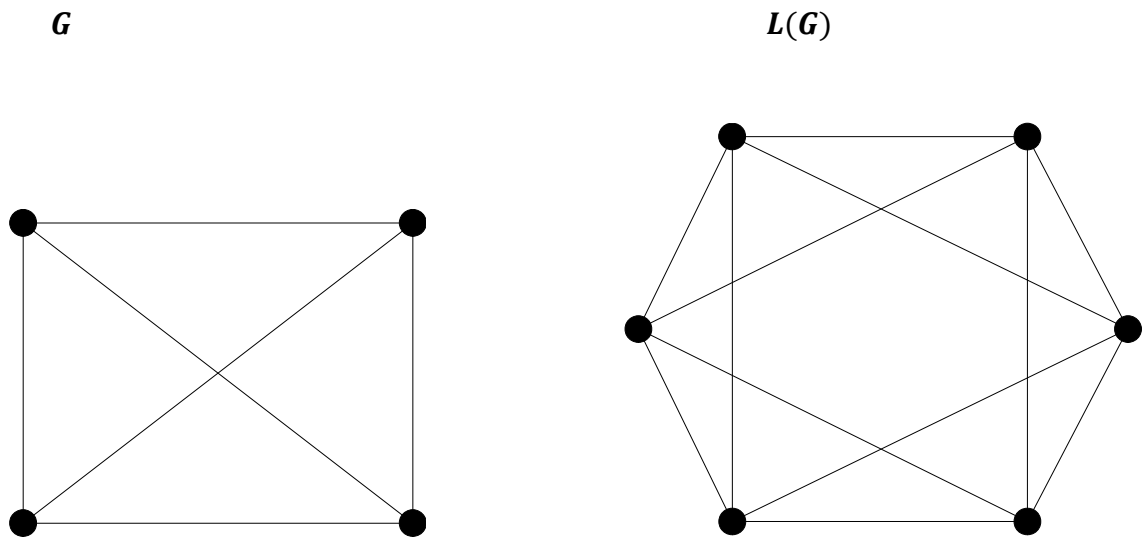

Նկ. 1.5.6

Դիցուք տրված է $\mathbb{R}$ իրական թվերի բազմությունը և $\mathfrak{F} = \{I_1, I_2, \ldots, I_m\}$ փակ հատվածների ընտանիքը։ Այդ դեպքում հատումների $\Omega(\mathbb{R}, \mathfrak{F})$ գրաֆը անվանում են *միջակայքերի գրաֆ*։ Նկ. 1.5.7-ում պատկերված է միջակայքերի գրաֆի օրինակ $\mathfrak{F} = \{[1,5], [1,3], [2,4], [2,6], [3,5], [4,6]\}$ հատվածների ընտանիքի դեպքում:

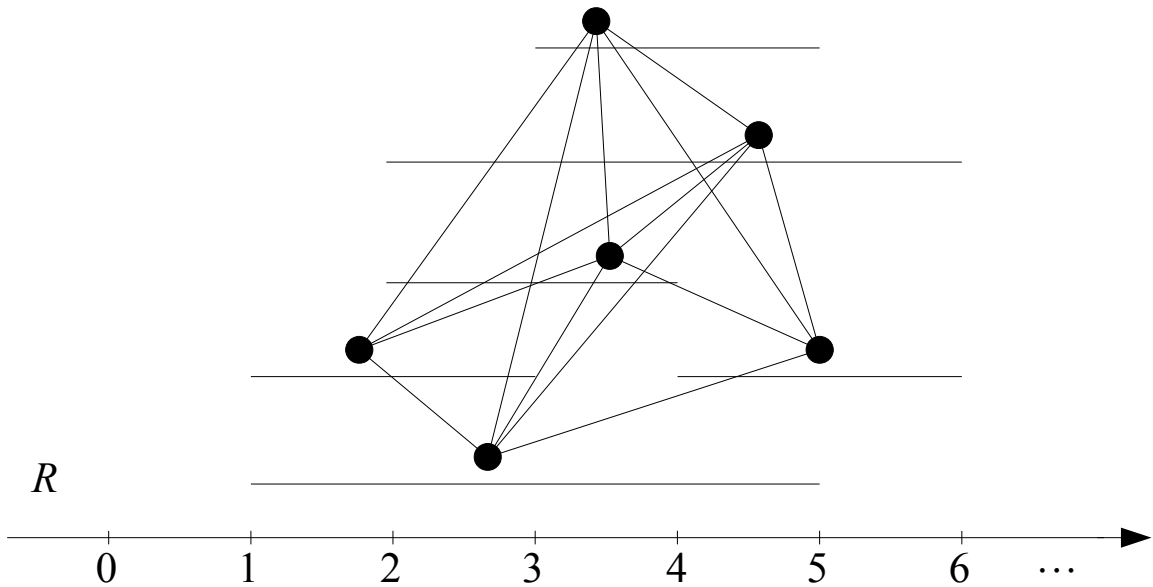

Նկ. 1.5.7

# Գլուխ 2

## Գրաֆների դասեր

### § 2.1. Կապակցվածության բաղադրիչներ և կապակցված գրաֆներ

Դիցուք $G = (V, E)$-ն գրաֆ է:

**Սահմանում 2.1.1:** $G$ գրաֆը կանվանենք *կապակցված*, եթե նրա ցանկացած երկու $u$ և $v$ գագաթների համար $G$ գրաֆում գոյություն ունի $(u, v)$-ճանապարհի:

Նշենք, որ կապակցված գրաֆի օրինակներ են հանդիսանում լրիվ և լրիվ երկկողմանի գրաֆները: Կապակցված են նաև ստորև պատկերված երեք գրաֆները:

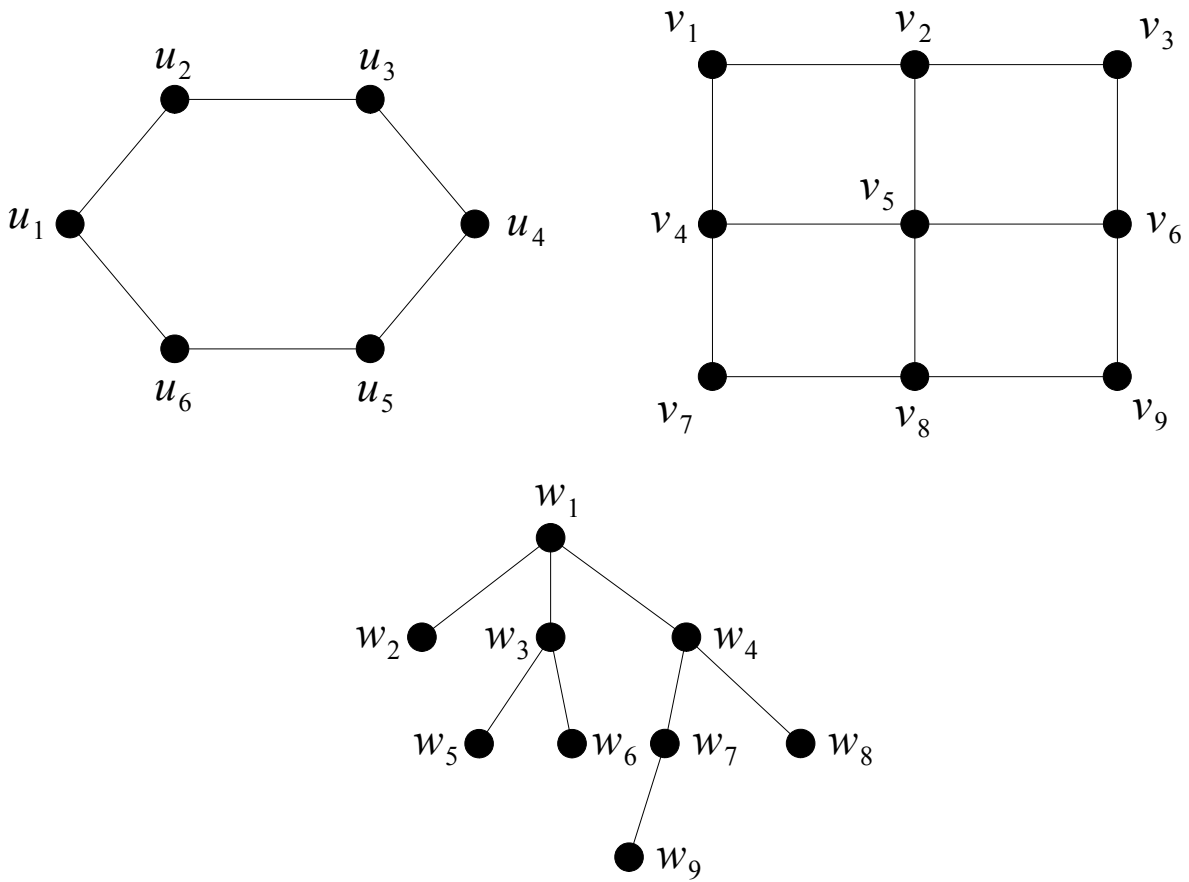

Նկ. 2.1.1

Եթե $G = (V, E)$-ն ցանկացած գրաֆ է, ապա դիտարկենք $V$ բազմության վրա սահմանված $\alpha$ բինար հարաբերությունը, որտեղ ցանկացած $u, v \in V$ համար

$u\alpha v$ այն և միայն այն դեպքում, երբ $G$ գրաֆում գոյություն ունի $(u, v)$-ճանապարհ:

Նկատենք, որ $\alpha$ բինար հարաբերությունը բավարարում է հետևյալ երեք պայմաններին.

- *ռեֆլեքսիվություն*, այսինքն՝ ցանկացած $v \in V$ համար $v\alpha v$,
- *սիմետրիկություն*, այսինքն՝ ցանկացած $u, v \in V$ համար, եթե $u\alpha v$, ապա $v\alpha u$,
- *տրանզիտիվություն*, այսինքն՝ ցանկացած $u, v, w \in V$ համար, եթե $u\alpha v$ և $v\alpha w$, ապա $u\alpha w$:

Հետևաբար, $\alpha$ բինար հարաբերությունն իրենից ներկայացնում է համարժեքության հարաբերություն, որտեղից հետևում է, որ $V$ բազմությունը կարելի է տրոհել $V_1, \ldots, V_p$ ենթաբազմությունների այնպես, որ

- $V = V_1 \cup \ldots \cup V_p$, $V_i \cap V_j = \emptyset$ երբ $1 \le i \ne j \le p$,
- ցանկացած $u, v \in V$ համար $u\alpha v$ այն և միայն այն դեպքում, երբ գոյություն ունի $k, 1 \le k \le p$, այնպես, որ $u, v \in V_k$:

Դիտարկենք $G$ գրաֆի $G_j = G[V_j]$ ենթագրաֆները, $1 \le j \le p$: $G$ գրաֆի $G_1, \ldots, G_p$ ենթագրաֆներն ընդունված է անվանել $G$ գրաֆի *կապակցվածության* կամ *կապակցված բաղադրիչներ*: Նկատենք, որ գրաֆի կապակցվածության բաղադրիչները կապակցված գրաֆներ են: Ավելին, նկատենք, որ § 1.3-ում սահմանված գրաֆների միավորում գործողությունը թույլ է տալիս ստանալ $G$ գրաֆի հետևյալ ներկայացումը

$$G = G_1 \cup \ldots \cup G_p:$$

**Դիտողություն 2.1.1**: $G$ գրաֆը կապակցված է այն և միայն այն դեպքում, երբ այն ունի կապակցվածության մեկ բաղադրիչ:

Նկ. 2.1.2-ում պատկերված գրաֆը կապակցված չէ և այն ունի կապակցվածության չորս բաղադրիչ:

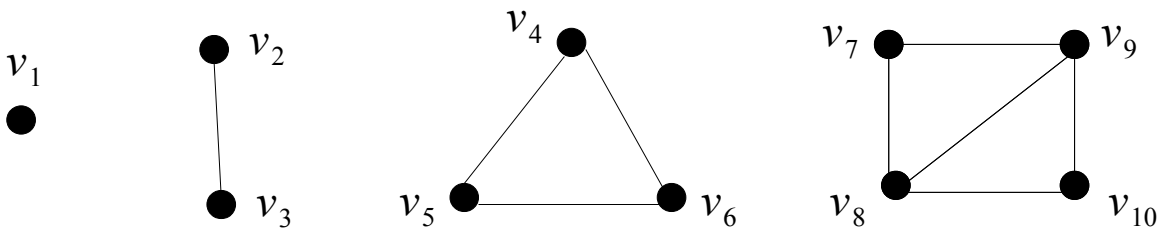

Նկ. 2.1.2

**Դիտողություն 2.1.2**: Կամայական $G$ գրաֆում գոյություն ունի ամենաերկար ճանապարհի:

Եթե $G = (V, E)$-ն կապակցված $(n, m)$-գրաֆ է, ապա $cyc(G) = m - n + 1$ թիվը կանվանենք $G$ գրաֆի *ցիկլոմատիկ թիվ*: Ստորև կապացուցենք կապակցված գրաֆների ցիկլոմատիկ թվին առնչվող մեկ թեորեմ:

**Թեորեմ 2.1.1:** Կապակցված $G = (V, E)$ գրաֆի համար $cyc(G) \geq 0$: Ավելին,

1. $cyc(G) = 0$ այն և միայն այն դեպքում, երբ $G$ գրաֆում ցիկլ չկա,
2. $cyc(G) = 1$ այն և միայն այն դեպքում, երբ $G$ գրաֆում կա ճիշտ մեկ ցիկլ:

**Ապացույց:** Դիտարկենք $G$ կապակցված գրաֆը, և դիտարկենք այն կառուցելու հետևյալ եղանակը.

**Քայլ 1.** Ընտրենք որևէ $e = uv$ կող և նրան կից $u$ և $v$ գագաթները: $G$ գրաֆի ստացված $G_1$ ենթագրաֆը պարունակում է մեկ կող և երկու գագաթ, և, հետևաբար, նրա համար $cyc(G_1) = 0$:

**Քայլ 2.** $G$ գրաֆի ստացված ենթագրաֆին հերթականորեն ավելացնենք $E$ բազմության կողեր և նրանց կից գագաթներ այնպես, որ յուրաքանչյուր քայլում ավելացվող կողն ունենա գոնե մեկ ընդհանուր գագաթ $G$ գրաֆի արդեն կառուցված ենթագրաֆի հետ: Նկատենք, որ դա հնարավոր է, քանի որ ըստ ենթադրության $G$ գրաֆը կապակցված է: Ավելին, ամեն անգամ $G$ գրաֆի ստացվող ենթագրաֆները կապակցած են:

Ակնհայտ է, որ քայլ 2-ն ամեն անգամ կատարելիս, մենք ավելացնում ենք մեկ նոր կող և ամենաշատը մեկ նոր գագաթ: Այստեղից հետևում է, որ կողերի և գագաթների քանակների տարբերությունը մնում է ոչ բացասական: Քանի որ $G$ գրաֆի սկզբնական $G_1$ ենթագրաֆի համար $cyc(G_1) = 0$, և $G_1$-ից հնարավոր է ստանալ $G$-ն, ապա պարզ է, որ $cyc(G) \geq 0$: Ավելին, նկատենք, որ $cyc(G) = 0$ այն և միայն այն դեպքում, երբ յուրաքանչյուր անգամ մեկ կող ավելացնելիս ավելացվում է ճիշտ մեկ գագաթ, և, հետևաբար, այդպիսի գրաֆը ցիկլ չի կարող պարունակել:

Նկատենք, որ $cyc(G) = 1$ այն և միայն այն դեպքում, երբ վերը նկարագրված եղանակով $G$ գրաֆը կառուցելիս ճիշտ մեկ անգամ է հանդիպում $h = xy$ կող, որը կից է արդեն կառուցած ենթագրաֆի $x$ և $y$ գագաթներին: Առանց ընդհանրությունը խախտելու, մենք կարող ենք ենթադրել, որ $h$ կողն ավելացվում է վերջին քայլում: Քանի որ նախորդ քայլերում ավելացված գագաթներն ու կողերը ցիկլ չեն առաջացնում, ապա դժվար չէ տեսնել, որ $h$ կողն ավելացնելիս առաջանում է ճիշտ մեկ ցիկլ: ∎

## § 2.2. Երկկողմանի գրաֆներ

Հիշենք, որ § 1.2-ում $G = (V, E)$ գրաֆն անվանեցինք երկկողմանի, եթե $V$ բազմությունը հնարավոր է տրոհել երկու ենթաբազմությունների $V_1$ և $V_2$-ի այնպես, որ $G$ գրաֆի ցանկացած կող կից է մեկ գագաթի $V_1$-ից և մեկ գագաթի $V_2$-ից: Նշենք, որ նկ. 2.1.1-ում պատկերված երեք գրաֆներն էլ երկկողմանի են: Իրոք, դիտարկենք նրանցից առաջինի գագաթների բազմության հետևյալ տրոհումը. $U_1 = \{u_1, u_3, u_5\}$ և $U_2 = \{u_2, u_4, u_6\}$, և նկատենք, որ այդ գրաֆի ցանկացած կող միացնում է մեկական գագաթ $U_1$-ից և $U_2$-ից: Երկրորդ գրաֆի երկկողմանիության մեջ համոզվելու համար դիտարկենք նրա գագաթների հետևյալ տրոհումը. $V_1 = \{v_1, v_3, v_5, v_7, v_9\}$ և $V_2 = \{v_2, v_4, v_6, v_8\}$, և նկատենք, որ այդ գրաֆի ցանկացած կող միացնում է մեկական գագաթ $V_1$-ից և $V_2$-ից: Եվ, վերջապես, երրորդ գրաֆի երկկողմանիության մեջ համոզվելու համար դիտարկենք նրա գագաթների հետևյալ տրոհումը. $W_1 = \{w_1, w_5, w_6, w_7, w_8\}$ և $W_2 = \{w_2, w_3, w_4, w_9\}$, և նկատենք, որ այդ գրաֆի ցանկացած կող միացնում է մեկական գագաթ $W_1$-ից և $W_2$-ից:

Ստորև կապացուցենք Քյոնիգի թեորեմը, որը նկարագրում է երկկողմանի գրաֆները:

**Թեորեմ 2.2.1 (Դ. Քյոնիգ):** Որպեսզի $G = (V, E)$ գրաֆը լինի երկկողմանի, անհրաժեշտ է և բավարար, որ այն չպարունակի կենտ երկարություն ունեցող պարզ ցիկլ:

**Ապացույց:** Նախ նկատենք, որ կենտ երկարություն ունեցող պարզ ցիկլը երկկողմանի չէ, հետևաբար, ցանկացած գրաֆ, որը պարունակում է կենտ երկարություն ունեցող պարզ ցիկլ չի կարող լինել երկկողմանի: Սա նշանակում է, որ եթե գրաֆը երկկողմանի է, ապա նրա բոլոր պարզ ցիկլերն ունեն զույգ երկարություն:

Հիմա ենթադրենք, որ $G$ գրաֆը չունի կենտ երկարություն ունեցող պարզ ցիկլ և ցույց տանք, որ այն երկկողմանի է:

Նախ պնդումն ապացուցենք այն մասնավոր դեպքում, երբ $G$ գրաֆը կապակցված է: Վերցնենք $G$ գրաֆի ցանկացած $u$ գագաթ: $V_1$-ով նշանակենք $G$ գրաֆի այն գագաթների բազմությունը, որոնք $u$ գագաթից գտնվում են զույգ հեռավորության վրա, իսկ $V_2$-ով այն գագաթների բազմությունը, որոնք $u$ գագաթից գտնվում են կենտ հեռավորության վրա: Նկատենք, որ $V_1$ և $V_2$ բազմությունները կազմում են $V$ բազմության տրոհում: Ավելին, $u \in V_1$:

Ցույց տանք, որ $G$ գրաֆի ցանկացած կող միացնում է մեկ գագաթ $V_1$-ից և մեկ գագաթ $V_2$-ից: Ենթադրենք հակառակը, դիցուք $G$ գրաֆի $e = vw$ կողը միացնում է $v$ և $w$ գագաթները, որոնք միաժամանակ պատկանում են $V_1$-ին կամ $V_2$-ին: Սա նշանակում է, որ գոյություն ունեն $u$ գագաթը $v$ և $w$ գագաթներին միացնող $P_v$ և $P_w$ ճանապարհներ այնպես, որ այդ ճանապարհների երկարություններն ունեն նույն զույգությունը: Դիտարկենք $G$ գրաֆի հետևյալ շրջանցումը. $u$ գագաթից $P_v$ ճանապարհի երկայնքով շարժվենք մինչև $v$ գագաթ, այնուհետև $e = vw$ կողով շարժվենք դեպի $w$ գագաթ, որից հետո $P_w$ ճանապարհով $w$ գագաթից վերադառնանք $u$ գագաթ: Նկատենք, որ նկարագրված շրջանցումն իրենից ներկայացնում է կենտ երկարություն ունեցող փակ շրջանցում: Համաձայն լեմմա 1.2.2-ի նրանից կարելի է անջատել $G$ գրաֆի կենտ երկարություն ունեցող պարզ ցիկլ, ինչը հակասում է թեորեմի պայմաններին: Հետևաբար, թեորեմի պնդումը ճիշտ է այն մասնավոր դեպքում, երբ $G$ գրաֆը կապակցված է:

Հիմա դիտարկենք ցանկացած $G$ գրաֆը, և դիցուք $G_1, \ldots, G_p$-ն նրա կապակցվածության բաղադրիչներն են: Նկատենք, որ քանի որ $G$ գրաֆը չի պարունակում կենտ երկարություն ունեցող պարզ ցիկլ, ապա նրա կապակցվածության բաղադրիչները ևս չեն պարունակի այդպիսի ցիկլ: Համաձայն վերը ապացուցվածի, $G_1, \ldots, G_p$ կապակցվածության բաղադրիչները հանդիսանում են երկկողմանի գրաֆներ, և հետևաբար $j = 1, \ldots, p$-ի համար $V(G_j)$ բազմությունը կարելի է տրոհել $V_1^{(j)}$ և $V_2^{(j)}$ բազմությունների այնպես, որ $G_j$ գրաֆի ցանկացած կող միացնում է մեկական գագաթ $V_1^{(j)}$-ից և $V_2^{(j)}$-ից: Նշանակենք՝

$$V^{(1)} = V_1^{(1)} \cup \ldots \cup V_1^{(p)} \text{ և } V^{(2)} = V_2^{(1)} \cup \ldots \cup V_2^{(p)}:$$

Նկատենք, որ $G$ գրաֆի ցանկացած կող միացնում է մեկական գագաթ $V^{(1)}$-ից և $V^{(2)}$-ից, և, հետևաբար, $G$ գրաֆը նույնպես երկկողմանի է: ∎

Դիցուք $G$-ն գրաֆ է: $G$ գրաֆի կցության $B(G)$ մատրիցը կանվանենք *տոտալ ունիմոդուլյար մատրից*, եթե այդ մատրիցի յուրաքանչյուր քառակուսային ենթամատրիցի որոշիչը հավասար է $0, 1$ կամ $-1$-ի: Այժմ տանք երկկողմանի գրաֆների մեկ այլ նկարագրում:

**Թեորեմ 2.2.2:** Որպեսզի $G = (V, E)$ գրաֆը լինի երկկողմանի, անհրաժեշտ է և բավարար, որ նրա կցության $B(G)$ մատրիցը լինի տոտալ ունիմոդուլյար:

**Ապացույց։** Նախ ցույց տանք, որ եթե $G$ գրաֆի կցության $B(G)$ մատրիցը տոտալ ունիմոդուլյար է, ապա $G$-ն երկկողմանի գրաֆ է։ Ենթադրենք հակառակը՝ $G$ գրաֆի կցության $B(G)$ մատրիցը տոտալ ունիմոդուլյար է, բայց $G$-ն երկկողմանի գրաֆ չէ։ Ըստ թեորեմ 2.2.1-ի $G$-ն պարունակում է կենտ երկարություն ունեցող պարզ ցիկլ։ Դիցուք այդ պարզ ցիկլի երկարությունը $2l + 1$ է։ Դիտարկենք այդ ցիկլի գագաթներին և կողերին համապատասխանող $B(G)$ մատրիցի ենթամատրիցը։ Դիցուք այդ մատրիցը $B'$-ն է։ Պարզ է, որ $B'$-ը $(2l + 1) \times (2l + 1)$ կարգի քառակուսային մատրից է։ Այժմ դիտարկենք $B''$ մատրիցը, որը ստացվում է $B'$-ից որոշ տողեր և սյուներ տեղափոխելով այնպես, որ $B''$ մատրիցը ընդունի հետևյալ տեսքը.

$$B'' = \begin{pmatrix} 1 & 0 & 0 & \cdots & 0 & 1 \\ 1 & 1 & 0 & \cdots & 0 & 0 \\ 0 & 1 & 1 & \cdots & 0 & 0 \\ \vdots & \vdots & \vdots & \ddots & \vdots & \vdots \\ 0 & 0 & 0 & \cdots & 1 & 1 \end{pmatrix}:$$

Հեշտ է տեսնել, որ $det(B'') = 1 + (-1)^{2l} = 2$։ Մյուս կողմից պարզ է, որ $B'$ մատրիցի որոշիչը կարող է տարբերվել $det(B'')$-ից միայն նշանով, իսկ դա հակասում է $B(G)$ մատրիցի տոտալ ունիմոդուլյար լինելուն։

Այժմ ցույց տանք, որ եթե $G$-ն երկկողմանի գրաֆ է, ապա $G$ գրաֆի կցության $B(G)$ մատրիցը տոտալ ունիմոդուլյար է։ Դիտարկենք $B(G)$ մատրիցի ցանկացած $Q$ $k \times k$ կարգի քառակուսային ենթամատրիցը։ Ապացույցը կատարենք մակածման եղանակով ըստ $k$-ի։ Եթե $k = 1$-ի, ապա ակնհայտ է, որ $det(Q) = 0$ կամ $det(Q) = 1$։ Ենթադրենք, որ պնդումը ճիշտ է $B(G)$ մատրիցի $Q'$ $k' \times k'$ կարգի ցանկացած քառակուսային ենթամատրիցի համար, որտեղ $k' < k$։ Դիտարկենք $Q$ $k \times k$ կարգի քառակուսային ենթամատրիցը։ Եթե $Q$-ն պարունակում է սյուն, որի բոլոր տարրերը զրոներ են, ապա պարզ է, որ $det(Q) = 0$։ Եթե $Q$-ն պարունակում է սյուն, որի ճիշտ մեկ տարրն է $1$, ապա վերլուծելով $det(Q)$-ն ըստ այդ սյանը մենք ըստ մակածման ենթադրության կստանանք, որ $det(Q) = 0, 1$ կամ $-1$-ի։ Այստեղից հետևում է, որ մենք կարող ենք ենթադրել, որ $Q$ մատրիցի յուրաքանչյուր սյուն պարունակում է ճիշտ երկու հատ $1$։ Քանի որ $G$-ն երկկողմանի գրաֆ է, ուստի այդ մեկերից մեկը կպատկանի $G$-ի մի կողմին, իսկ մյուսը՝ մյուս կողմին։ Պարզ է, որ մենք կարող ենք ենթադրել, որ $Q$ մատրիցի առաջին $r$ տողերին համապատասխանում է $G$ երկկողմանի գրաֆի մի կողմը, իսկ մյուս $k - r$ տողերին՝ այդ գրաֆի մյուս կողմը։ Քանի որ $G$-ն երկկողմանի գրաֆ է, ուստի $Q$ մատրիցի յուրաքանչյուր սյուն կպարունակի մեկ հատ **1** առաջին $r$ տողերից և ճիշտ մեկ հատ **1** մյուս $k-r$ տողերից։ Այստեղից հետևում է, որ $Q$ մատրիցի առաջին $r$ տողերի գումարը հավասար է այդ մատրիցի մյուս $k-r$ տողերի գումարին, ուստի $Q$ մատրիցի տողերը գծորեն կախված են և $det(Q) = 0$։ ∎

Վերջում ապացուցենք մի թեորեմ, որը ցույց է տալիս, որ ցանկացած գրաֆ պարունակում է կողերով հարուստ կմախքային երկկողմանի ենթագրաֆ։

**Թեորեմ 2.2.3 (Պ. Էրդյոշ):** Կամայական $G$ գրաֆ պարունակում է կմախքային երկկողմանի $H$ ենթագրաֆ, որում $|E(H)| \geq \frac{|E(G)|}{2}$։

**Ապացույց:** Դիտարկենք $G$ գրաֆի գագաթների $V(G)$ բազմության բոլոր հնարավոր տրոհումները երկու ենթաբազմությունների և ընտրենք այն մեկը, որի դեպքում այդ երկու ենթաբազմությունների միջև կողերի քանակը առավելագույնն է։ Դիցուք այդ տրոհումը $V(G) = U \cup W$-ն է։ Պարզ է, որ այդ տրոհումը ծնում է կմախքային երկկողմանի $H$ ենթագրաֆ։ Ցույց տանք, որ ցանկացած $v \in V(G)$-ի համար $d_H(v) \geq \frac{d_G(v)}{2}$։ Ենթադրենք հակառակը՝ գոյություն ունի $u \in V(G)$-ին, որ $d_H(u) < \frac{d_G(u)}{2}$։ Առանց ընդհանրությունը խախտելու կարող ենք ենթադրել, որ $u \in U$։ Քանի որ $d_H(u) < \frac{d_G(u)}{2}$, ուստի $u$ գագաթը $U$-ում ավելի շատ հարևան գագաթներ ունի քան $W$-ում, իսկ դա նշանակում է, որ եթե մենք $u$ գագաթը $U$-ից տեղափոխենք $W$, ապա կստանանք $V(G) = (U\setminus\{u\}) \cup (W \cup \{u\})$ նոր տրոհումը, որի $(U\setminus\{u\})$ և $(W \cup \{u\})$ ենթաբազմությունների միջև կողերի քանակը ավելի մեծ է, որը հակասում է սկզբնական տրոհման ընտրությանը։ Այստեղից հետևում է, որ

$$|E(H)| = \frac{1}{2}\sum_{v \in V(G)} d_H(v) \geq \frac{1}{2}\sum_{v \in V(G)} \frac{d_G(v)}{2} = \frac{|E(G)|}{2}։ \blacksquare$$

Երկկողմանի գրաֆներին ավելի մանրամասն կարելի է ծանոթանալ Հասրաթյանի, Դենլիի և Հազվիսթի գրքում [4]։

## § 2.3. Ծառեր

**Սահմանում 2.3.1:** Ցիկլ չպարունակող կապակցված գրաֆը կանվանենք *ծառ*։
**Սահմանում 2.3.2:** Ցիկլ չպարունակող գրաֆը կանվանենք *անտառ*։

Նկատենք, որ անտառն այնպիսի գրաֆ է, որի կապակցվածության բոլոր բաղադրիչներն իրենցից ներկայացնում են ծառեր: Ավելին, կամայական անտառ երկկողմանի գրաֆ է:

Նշենք, որ նկ. 2.1.2-ում պատկերված գրաֆը անտառ չէ: Ավելին, նկ. 2.1.1-ում պատկերված գրաֆներից առաջին երկուսը ծառեր չեն, իսկ երրորդը ծառ է: Վերջապես, ծառի օրինակ է նաև 2.3.1-ում պատկերված գրաֆը:

Ստորև կապացուցենք մի թեորեմ, որն առաջարկում է ծառի մի քանի նկարագրություններ:

**Թեորեմ 2.3.1:** $G = (V, E)$ $(n, m)$-գրաֆի համար հետևյալ պայմանները իրար համարժեք են.

(1) $G$-ն ծառ է,

(2) $G$ գրաֆում ցանկացած երկու գագաթ միացված են ճիշտ մեկ ճանապարհով,

(3) $G$-ն կապակցված է և $m = n - 1$,

(4) $G$-ն չունի ցիկլ և $m = n - 1$,

(5) $G$-ն չունի ցիկլ և $G$-ի ցանկացած երկու ոչ հարևան $u$ և $v$ գագաթների համար $G + uv$ գրաֆն ունի ճիշտ մեկ ցիկլ:

**Ապացույց:** Նախ ցույց տանք, որ (1)-ից հետևում է (2)-ը: Իրոք, դիցուք $G$-ն ծառ է: Այդ դեպքում, քանի որ $G$-ն կապակցված է, նրանում ցանկացած երկու գագաթ միացված են առնվազն մեկ ճանապարհով: Ցույց տանք, որ ցանկացած երկու գագաթ միացված են ճիշտ մեկ ճանապարհով:

Ենթադրենք, $G$-ում գոյություն ունեն երկու $u$ և $v$ գագաթներ, որոնք միացված են իրարից տարբեր $P_1$ և $P_2$ ճանապարհներով: $u$ գագաթից $P_1$ ճանապարհի երկայնքով շարժվենք դեպի $v$ գագաթ: Քանի որ $P_1 \neq P_2$, ապա գոյություն կունենա այնպիսի $w$ գագաթ, որը պատկանում է $P_1$-ին և $P_2$-ին, որի հաջորդը $P_1$ ճանապարհի վրա չի պատկանում $P_2$-ին: Քանի որ $v$ գագաթը գտնվում է $P_1$-ի և $P_2$-ի վրա, ապա $w$ գագաթից $P_1$ ճանապարհի երկայնքով շարժվելուց, մենք կհանդիպենք $w'$ գագաթի, որը պատկանում է $P_1$-ին և $P_2$-ին, բայց որի նախորդները, որոնք ընկած են մինչև $w$ գագաթ չեն պատկանում $P_1$-ին և $P_2$-ին միաժամանակ: Նկատենք, որ $P_1$-ի և $P_2$-ի $w$ և $w'$ գագաթները միացնող ենթաճանապարհները միասին կազմում են ցիկլ, ինչը հակասում է $G$-ի ցիկլ չունենալու պայմանին: Հետևաբար $G$-ում ցանկացած երկու գագաթ միացված են ճիշտ մեկ ճանապարհով:

Հիմա ցույց տանք, որ (2)-ից հետևում է (3)-ը: Ենթադրենք, որ $G$-ում ցանկացած երկու գագաթ միացված են ճիշտ մեկ ճանապարհով։ Նկատենք, որ այս պայմանից հետևում է, որ $G$-ն կապակցված է։ Մակածման եղանակով ըստ $n$-ի ցույց տանք, որ $m = n - 1$։

Նկատենք, որ $m = n - 1$ հավասարությունը ակնհայտ է $n = 1, 2$ դեպքերում։ Ենթադրենք, որ այն ճիշտ է (2) պայմանին բավարարող բոլոր գրաֆների համար, որոնց գագաթների քանակը փոքր է $n$-ից, և դիտարկենք (2) պայմանին բավարարող $G$ գրաֆը։ Վերցնենք $G$ գրաֆի ցանկացած $e$ կող, և դիտարկենք $G - e$ գրաֆը։ Քանի որ $G$ գրաֆում ցանկացած երկու գագաթ միացված են ճիշտ մեկ ճանապարհով, ապա $G - e$ գրաֆը պետք է ունենա կապակցվածության ճիշտ երկու բաղադրիչ։ Դիցուք այդ բաղադրիչները $G_1$-ը և $G_2$-ն են։ Նկատենք, որ $G_1$-ում և $G_2$-ում ցանկացած երկու գագաթ միացված են ճիշտ մեկ ճանապարհով և նրանցում գագաթների քանակը փոքր է $n$-ից։ Համաձայն մակածման ենթադրության

$$|E(G_1)| = |V(G_1)| - 1 \text{ և } |E(G_2)| = |V(G_2)| - 1,$$

որտեղից հետևում է, որ

$$m = 1 + |E(G_1)| + |E(G_2)| = 1 + |V(G_1)| + |V(G_2)| - 2 = n - 1։$$

Ցույց տանք, որ (3)-ից հետևում է (4)-ը։ Դիցուք $G$-ն կապակցված է և $m = n - 1$։ Նկատենք, որ բավական է ցույց տալ, որ $G$-ն չունի ցիկլ։ Ենթադրենք հակառակը, այսինքն ենթադրենք, որ $G$-ն ունի $p$ երկարությամբ ցիկլ։ Նկատենք, որ այս ցիկլի վրա գտնվում են $G$ գրաֆի $p$ գագաթ և $p$ կող։ $G$ գրաֆի մնացած $n - p$ գագաթներին համապատասխանեցնենք կողեր հետևյալ կերպ. դիցուք $v$-ն ցիկլի վրա չգտնվող գագաթ է, դիտարկենք $v$-ից ցիկլի գագաթներ տանող կարճագույն ճանապարհները և ընտրենք նրանցից ամենակարճը։ $v$-ին համապատասխանեցնենք այդ ամենակարճ ճանապարհի վրա գտնվող և $v$-ին կից կողը։ Նկատենք, որ ցիկլի վրա չգտնվող $n - p$ գագաթներին կհամապատասխանեն $n - p$ իրարից տարբեր կողեր, և, հետևաբար, $G$ գրաֆում կողերի քանակը կլինի առնվազն

$$m \geq p + n - p = n,$$

ինչը հակասում է $m = n - 1$ հավասարությանը։

Ցույց տանք, որ (4)-ից հետևում է (5)-ը։ Ենթադրենք $G$-ն չունի ցիկլ և $m = n - 1$։ Նկատենք, որ բավական է ցույց տալ, որ $G$-ի ցանկացած երկու ոչ հարևան $u$ և $v$ գագաթների համար $G + uv$ գրաֆն ունի ճիշտ մեկ ցիկլ։ Քանի որ $G$-ն չունի ցիկլ, $G$-ի կապակցվածության բաղադրիչները ծառեր են, և քանի որ արդեն ցույց ենք տվել, որ (1)-ից հետևում է (3)-ը, ապա կունենանք, որ $m = n - k$, որտեղ $k$-ն $G$-ի կապակցվածության բաղադրիչների քանակն է։ Մյուս կողմից, քանի որ $m = n - 1$, ապա կունենանք, որ $k = 1$, այսինքն $G$-ն ծառ է, և, հետևաբար, համաձայն (2)-ի, $G$-ում ցանկացած երկու ոչ հարևան $u$ և $v$ գագաթներ միացված են ճիշտ մեկ ճանապարհով։ Այստեղից հետևում է, որ $G + uv$ գրաֆում կլինի ճիշտ մեկ ցիկլ։

Ցույց տանք, որ (5)-ից հետևում է (1)-ը։ Ենթադրենք, որ $G$-ն չունի ցիկլ և $G$-ի ցանկացած երկու ոչ հարևան $u$ և $v$ գագաթների համար $G + uv$ գրաֆն ունի ճիշտ մեկ ցիկլ։ Նկատենք, որ բավական է ցույց տալ, որ $G$-ն կապակցված է։ Դիտարկենք ցանկացած երկու $u$ և $v$ գագաթներ։ Եթե նրանք հարևան են, ապա պարզ է, որ միացված են ճանապարհով։ Ենթադրենք, որ $u$ և $v$ գագաթները հարևան չեն։ Դիտարկենք $G + uv$ գրաֆը։ Համաձայն մեր ենթադրության, $G + uv$ գրաֆն ունի ճիշտ մեկ ցիկլ։ Այստեղից հետևում է, որ $u$ և $v$ գագաթները $G$ գրաֆում միացված են ճանապարհով, հետևաբար, $G$-ն կապակցված գրաֆ է։ ∎

Նկատենք, որ ապացուցված թեորեմից և թեորեմ 2.1.1-ից հետևում է, որ ծառը կարելի է սահմանել որպես այնպիսի կապակցված գրաֆ, որի ցիկլոմատիկ թիվը հավասար է զրոյի։

Հիշենք, որ § 1.2-ում $G$ գրաֆի $v$ գագաթը անվանեցինք կախված, եթե $d_G(v) = 1$։

**Հետևանք 2.3.1։** Եթե $T = (V, E)$-ն ծառ է, որում $|V| \geq 2$, ապա $T$-ն պարունակում է առնվազն երկու կախված գագաթ։

**Ապացույց 1։** Համաձայն դիտողություն 2.1.2-ի, $T$ ծառում գոյություն ունի ամենաերկար ճանապարհ։ Նկատենք, որ քանի որ $|V| \geq 2$, ապա այդ ճանապարհի ծայրակետերը երկուսն են, որոնք, ինչը դժվար չէ տեսնել, $T$ ծառի կախված գագաթներ են։ ∎

**Ապացույց 2։** Քանի որ $T$-ն կապակցված է և $|V| \geq 2$, ապա նրանում ցանկացած գագաթի աստիճանն առնվազն մեկ է։ Օգտվելով թեորեմ 1.2.1-ից և թեորեմ 2.3.1-ից, կստանանք,

$$\sum_{v \in V} d(v) = 2|E| = 2|V| - 2,$$

որտեղից հետևում է, որ առնվազն երկու գագաթի աստիճան պետք է լինի մեկ։ ∎

**Թեորեմ 2.3.2։** Դիցուք $T$-ն ծառ է, որում $|E(T)| = k$, և $G = (V, E)$-ն գրաֆ է, որում $\delta(G) \geq k$։ Այդ դեպքում $G$ գրաֆը պարունակում է $T$ ծառին իզոմորֆ ենթագրաֆ։

**Ապացույց։** Ապացույցը կատարենք մակածման եղանակով ըստ $k$-ի։ Ենթադրենք

$k = 0$: Այդ դեպքում $T$-ն բաղկացած է մեկ գագաթից, և պնդումն ակնհայտ է: Ենթադրենք, որ թեորեմի պնդումը ճիշտ է բոլոր այն ծառերի համար, որոնց կողերի քանակը փոքր է $k$-ից, և դիտարկենք $k > 0$ կող պարունակող $T$ ծառը: Համաձայն հետևանք 2.3.1-ի, $T$ ծառը պարունակում է $u$ կախված գագաթ: $v$-ով նշանակենք $u$ գագաթի միակ հարևան գագաթը $T$-ում, և դիցուք $T' = T - u$: Նկատենք, որ $T'$-ը նույնպես ծառ է, որը պարունակում է $k - 1$ կող: Քանի որ $\delta(G) \geq k > k - 1$, ապա, համաձայն մակածման ենթադրության, $G$ գրաֆը պարունակում է $T'$ ծառին իզոմորֆ ենթագրաֆ: Դիցուք այդ ենթագրաֆում $v$ գագաթին համապատասխան գագաթը $x$-ն է: Քանի որ $\delta(G) \geq k$, ապա $d_G(x) \geq k$, որտեղից հետևում է, որ $x$ գագաթը հարևան է այնպիսի $y$ գագաթի, որին համապատասխան գագաթ չկա $G$ գրաֆի $T'$ ծառին իզոմորֆ ենթագրաֆում: Ավելացնելով $y$ գագաթը և $xy$ կողն այդ ենթագրաֆին, մենք կստանանք $G$ գրաֆի $T$ ծառին իզոմորֆ ենթագրաֆ: ∎

Ստորև կդիտարկենք տրված $n$ նշված գագաթներով իրարից տարբեր ծառերի հաշվման խնդիրը:

**Թեորեմ 2.3.3 (Կելլի):** $\{1, 2, \ldots, n\}$ բազմությունը որպես գագաթների բազմություն ունեցող ծառերի քանակը հավասար է $n^{n-2}$:

**Ապացույց:** Ապացույցը կկատարենք Պրյուֆերի կոդավորման եղանակով: Նախ նկատենք, որ թեորեմի պնդումը ճիշտ է $n = 1$ և $n = 2$ դեպքերում, հետևաբար կարող ենք ենթադրել, որ $n \geq 3$: Այս դեպքում պնդումն ապացուցելու համար մենք ցույց կտանք, որ կարելի է հաստատել փոխմիարժեք արտապատկերում $\{1, 2, \ldots, n\}$ բազմությունը որպես գագաթների բազմություն ունեցող ծառերի և $\{1, 2, \ldots, n\}$ բազմության տարրերի $n - 2$ երկարությամբ բոլոր հաջորդականությունների միջև: Քանի որ վերջիններիս քանակը $n^{n-2}$ է, ապա այստեղից էլ կստացվի թեորեմի ապացույցը:

Դիցուք $T$-ն $1, 2, \ldots, n$ գագաթներով որևէ ծառ է: Համաձայն հետևանք 2.3.1-ի $T$ ծառում գոյություն ունի կախված գագաթ: Դիտարկենք նրանցից ամենափոքրը: Դիցուք այն $i_1$-ն է, ավելին դիցուք $i_1$-ի միակ հարևանը $T$-ում $j_1$-ն է: Դիտարկենք $T - i_1$ ծառը: Կրկին համաձայն հետևանք 2.3.1-ի նրանում գոյություն ունի ամենափոքր $i_2$ կախված գագաթ, որի միակ հարևանը $T - i_1$-ում դիցուք $j_2$-ն է: Դիտարկենք $T - i_1 - i_2$ ծառը: Նկարագրված պրոցեսը կատարենք $n - 2$ անգամ: Արդյունքում $T - i_1 - i_2 - \cdots - i_{n-2}$ ծառը բաղկացած կլինի մեկ կողից:

$T$ ծառին համապատասխանեցնենք $n - 2$ երկարություն ունեցող $l(T)$ հաջորդականությունը, որտեղ $l(T)$-ն որոշվում է հետևյալ կերպ.

$$l(T) = j_1 j_2 \ldots j_{n-2}:$$

Նկատենք, որ նկ. 2.3.1-ում պատկերված ծառին կհամապատասխանի **10, 7, 5, 8, 2, 3, 8, 2** հաջորդականությունը:

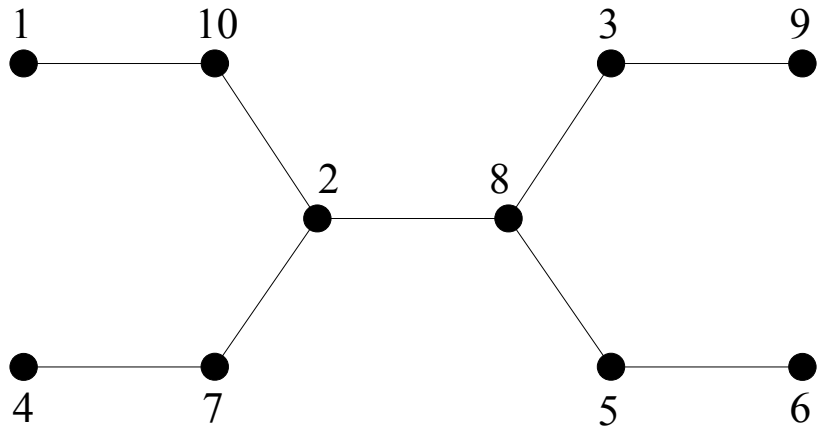

Նկ. 2.3.1

Ավելին, նկատենք, որ եթե $T_1$-ը և $T_2$-ը իրարից տարբեր ծառեր են, ապա $l(T_1) \neq l(T_2)$: Սա նշանակում է, որ թեորեմի ապացույցն ավարտելու համար բավական է ցույց տալ, որ $\{1, 2, \ldots, n\}$ բազմության տարրերի $n-2$ երկարությամբ ցանկացած հաջորդականության համապատասխանում է ծառ:

Դիցուք $l_1 l_2 \ldots l_{n-2}$-ը $\{1, 2, \ldots, n\}$ բազմության տարրերի $n-2$ երկարությամբ ցանկացած հաջորդականություն է: Նրան համապատասխանեցնենք գրաֆ հետևյալ կերպ. գտնենք $\{1, 2, \ldots, n\}$ բազմության ամենափոքր տարրը, որը չի մասնակցում $l_1 l_2 \ldots l_{n-2}$ հաջորդականության մեջ: Դիցուք այն $k_1$-ն է: $k_1$-ը միացնենք կողով $l_1$-ի հետ, որից հետո $\{1, 2, \ldots, n\}$ բազմությունից հեռացնենք $k_1$-ը, իսկ $l_1 l_2 \ldots l_{n-2}$ հաջորդականությունից հեռացնենք $l_1$-ը: Այնուհետև վարվենք նույն կերպ, այսինքն գտնենք $\{1, 2, \ldots, n\} \setminus \{k_1\}$ բազմության ամենափոքր տարրը, որը չի մասնակցում $l_2 \ldots l_{n-2}$ հաջորդականության մեջ: Դիցուք այն $k_2$-ն է: $k_2$-ը միացնենք կողով $l_2$-ի հետ, որից հետո $\{1, 2, \ldots, n\} \setminus \{k_1\}$ բազմությունից հեռացնենք $k_2$-ը, իսկ $l_2 \ldots l_{n-2}$ հաջորդականությունից հեռացնենք $l_2$-ը: Կատարենք նշված քայլերը $n-2$ անգամ, որի արդյունքում $\{1, 2, \ldots, n\}$ բազմության մեջ կմնա երկու տարր, որոնք էլ միացնենք կողով և դրանով ավարտենք $l_1 l_2 \ldots l_{n-2}$ հաջորդականությանը համապատասխանող գրաֆի կառուցումը:

Նկատենք, որ այս եղանակով **10, 7, 5, 8, 2, 3, 8, 2** հաջորդականությանը կհամապատասխանի նկ. 2.3.1-ում պատկերված ծառը:

Մակածման եղանակով ցույց տանք, որ վերը նկարագրված եղանակով ցանկացած

$l_1 l_2 \ldots l_{n-2}$ հաջորդականության համապատասխանեցրել ենք ծառ։ Նկատենք, որ այս պնդումն ակնհայտ է, երբ $n \leq 3$։

Ենթադրենք, որ $n-3$ երկարությամբ ցանկացած հաջորդականության համապատասխանում է $n-1$ գագաթանի ծառ, և դիտարկենք $n-2$ երկարությամբ $l_1 l_2 \ldots l_{n-2}$ հաջորդականությունը։ Նկատենք, որ $l_1 l_2 \ldots l_{n-2}$ հաջորդականությանը համապատասխանող գրաֆում $l_1$-ը միացված է կողով $k_1$-ի հետ, ընդ որում $k_1$-ի աստիճանը հավասար է մեկի։ Եթե այդ գրաֆից հեռացնենք $k_1$ գագաթը, ապա ստացված գրաֆը կհամապատասխանի $n-3$ երկարություն ունեցող $l_2 \ldots l_{n-2}$ հաջորդականությանը, որը, համաձայն մակածման ենթադրության, ծառ է։ Նկատենք, որ այդ դեպքում ծառ կլինի նաև $l_1 l_2 \ldots l_{n-2}$ հաջորդականությանը համապատասխանող գրաֆը, քանի որ այն ստացվում է $l_1$ գագաթը $k_1$-ի հետ կողով միացնելով։ ∎

Հիշենք, որ § 1.2-ում ներմուծեցինք գրաֆի կմախքային ենթագրաֆի գաղափարը։ $H$ գրաֆը կոչվում էր $G$ գրաֆի կմախքային ենթագրաֆ, եթե $V(H) = V(G)$ և $E(H) \subseteq E(G)$։

**Թեորեմ 2.3.4:** Եթե $G$-ն կապակցված գրաֆ է, ապա այն պարունակում է *կմախքային ծառ* (ծառ հանդիսացող կմախքային ենթագրաֆ)։

**Ապացույց:** Ապացույցը կատարենք մակածման եղանակով ըստ $G$ գրաֆի կողերի բազմության հզորության։ Նկատենք, որ պնդումն ակնհայտ է $|E(G)| = 0, 1$ դեպքերում։ Ենթադրենք, որ պնդումը ճիշտ է բոլոր այն կապակցված գրաֆների համար, որոնց կողերի քանակը փոքր է $|E(G)|$-ից, և դիտարկենք $G$ կապակցված գրաֆը։ Եթե $G$ գրաֆը չի պարունակում ցիկլ, ապա $G$-ն ծառ է, և թեորեմն ապացուցված է։ Հետևաբար կարող ենք ենթադրել, որ $G$-ն պարունակում է ցիկլ։ Դիցուք $e$-ն այդ ցիկլի որևէ կող է։ Դիտարկենք $G - e$ գրաֆը։ Նկատենք, որ այն կապակցված է և պարունակում է $|E(G)|$-ից քիչ կող։ Համաձայն մակածման ենթադրության, $G - e$ գրաֆը պարունակում է կմախքային ծառ, որն էլ հանդիսանում է $G$ գրաֆի կմախքային ծառ։ ∎

Նկատենք, որ կմախքային ենթագրաֆների լեզվով թեորեմ 2.3.3-ը կարելի է ձևակերպել հետևյալ կերպ. $n$ գագաթ պարունակող լրիվ գրաֆի կմախքային ծառերի քանակը հավասար է $n^{n-2}$։ Հաշվի առնելով թեորեմ 2.3.4-ը, նշված արդյունքը կարելի է փորձել ընդհանրացնել՝ պարզելով, թե ինչի է հավասար ցանկացած կապակցված գրաֆի կմախքային ծառերի քանակը։

Նշենք, որ այս հարցի պատասխանը տալիս է Կիրխհոֆի թեորեմը, որը ձևակերպելու համար կատարենք որոշ նշանակումներ.

Դիցուք $G$-ն կապակցված գրաֆ է և $V(G) = \{v_1, \ldots, v_n\}$։ Ենթադրենք, որ $A(G)$-ն $G$ գրաֆի § 1.1-ում սահմանված հարևանության մատրիցն է, և դիտարկենք $n \times n$ կարգի $D(G) = (d_{ij})_{n \times n}$ մատրիցը, որը սահմանվում է հետևյալ կերպ.

$$d_{ij} = \begin{cases} d(v_i), & \text{եթե } i = j, \\ 0, & \text{հակառակ դեպքում։} \end{cases}$$

$G$ գրաֆի *լապլասյան* կանվանենք հետևյալ ձևով սահմանված մատրիցը.

$$Lap(G) = D(G) - A(G)։$$

**Թեորեմ 2.3.5 (Կիրխոֆ):** Ցանկացած կապակցված $G$ գրաֆի համար նրա լապլասյանի բոլոր հանրահաշվական լուծումները իրար հավասար են և նրանց ընդհանուր արժեքը հավասար է $G$ գրաֆի կմախքային ծառերի քանակին։

Հիշենք, որ կապակցված $G$ գրաֆում ցանկացած երկու $u$ և $v$ գագաթների միջև $d(u, v)$ հեռավորությունը սահմանվում է որպես $u$ և $v$ գագաթները միացնող ամենակարճ ճանապարհի երկարությունը (§ 1.2)։ Կապակցված $G$ գրաֆի և նրա ցանկացած $u$ գագաթի համար կատարենք հետևյալ նշանակումները.

$$\varepsilon(u) = \max_{v \in V(G)} d(u, v), \quad r(G) = \min_{u \in V(G)} \varepsilon(u), \quad d(G) = \max_{u \in V(G)} \varepsilon(u)$$

$\varepsilon(u)$-ն կոչվում է $u$ գագաթի *էրցենտրիսիտետ*, $r(G)$-ն կոչվում է $G$ գրաֆի *շառավիղ*, իսկ $d(G)$-ն կոչվում է $G$ գրաֆի *տրամագիծ*։ $G$ գրաֆի $u$ գագաթը կոչվում է *կենտրոնական*, եթե $\varepsilon(u) = r(G)$, իսկ $G$ գրաֆի բոլոր կենտրոնական գագաթների բազմությունը կոչվում է $G$ գրաֆի *կենտրոն*։

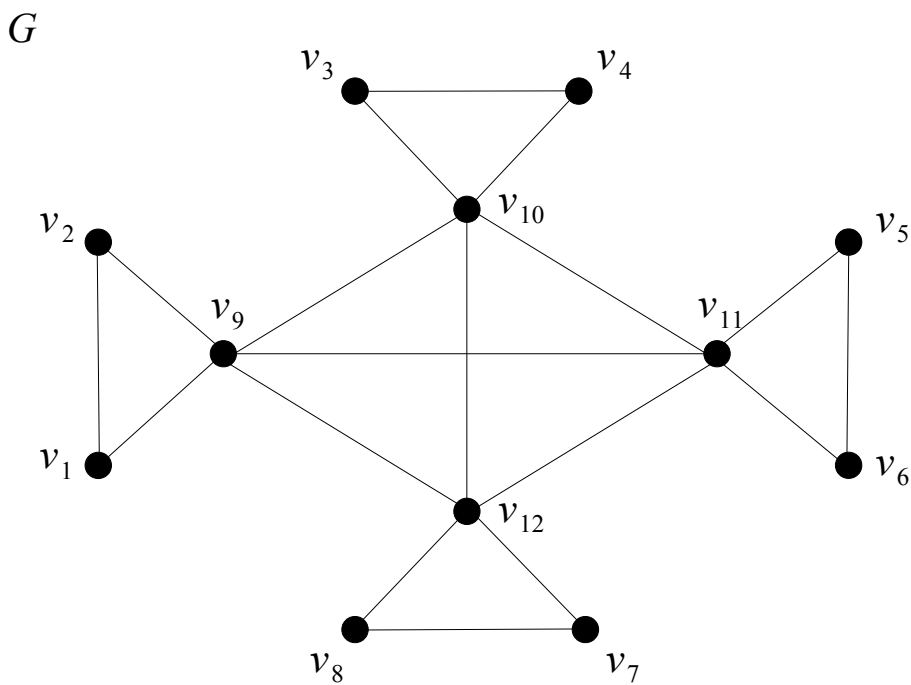

Նկ. 2.3.2

Նկ. 2.3.2-ում պատկերված $G$ գրաֆում $\varepsilon(v_1) = \varepsilon(v_2) = \varepsilon(v_3) = \varepsilon(v_4) = \varepsilon(v_5) = \varepsilon(v_6) = \varepsilon(v_7) = \varepsilon(v_8) = 3$, $\varepsilon(v_9) = \varepsilon(v_{10}) = \varepsilon(v_{11}) = \varepsilon(v_{12}) = 2$, հետևաբար $r(G) = 2$ և $d(G) = 3$: Գրաֆի կենտրոնական գագաթներն են $v_9$-ը, $v_{10}$-ը, $v_{11}$-ը և $v_{12}$-ը, որոնք էլ կազմում են գրաֆի կենտրոնը: Նկատենք, որ այս գրաֆում կենտրոնը բաղկացած է չորս գագաթից: Ծառերի կենտրոնի վերաբերյալ հայտնի է Ժորդանի թեորեմը:

**Թեորեմ 2.3.6 (Ժորդան):** Ցանկացած ծառի կենտրոնը բաղկացած է ոչ ավելի քան երկու գագաթից:

**Ապացույց:** Ապացույցը կատարենք մակածման եղանակով ըստ ծառի գագաթների քանակի: Նախ նկատենք, որ պնդումն ակնհայտ է մեկ և երկու գագաթ պարունակող ծառերի համար: Ենթադրենք, որ պնդումը ճիշտ է բոլոր այն ծառերի համար, որոնց գագաթների քանակը փոքր է $n$-ից, և դիտարկենք $n \geq 3$ գագաթ պարունակող $T$ ծառը և նրա որևէ $u$ գագաթ: Նկատենք, որ $T$ ծառում $u$ գագաթից ամենահեռու գտնվող գագաթը միշտ կախված է, հետևաբար, եթե դիտարկենք $T'$ ծառը, որը ստացվում է $T$ ծառից հեռացնելով նրա բոլոր կախված գագաթները, ապա $T'$ ծառում բոլոր գագաթների էքսցենտրիսիտետները ստացվում են $T$ ծառում նրանց ունեցած էքսցենտրիսիտետներից հանելով մեկ: Այստեղից հետևում է, որ $T$ ծառի և $T'$ ծառի կենտրոնները համընկնում են: Նկատենք, որ $|V(T')| < |V(T)| = n$, հետևաբար, $T'$ ծառի կենտրոնը բաղկացած է ոչ ավելի, քան երկու գագաթից: Այստեղից հետևում է, որ $T$ ծառի կենտրոնը ևս բաղկացած է ոչ ավելի, քան երկու գագաթից: ∎

# Գլուխ 3

# Կապակցվածություն

## § 3.1. Միակցման կետեր և կամուրջներ

Դիցուք $G = (V, E)$-ն գրաֆ է: $c(G)$-ով նշանակենք $G$ գրաֆի կապակցված բաղադրիչների քանակը:

**Սահմանում 3.1.1:** $G$ գրաֆի $v$ գագաթը կոչվում է *միակցման կետ*, եթե $c(G - v) > c(G)$:

Նկատենք, որ եթե $v$-ն $G$ կապակցված գրաֆի միակցման կետ է, ապա $G - v$ գրաֆը կապակցված չէ:

**Սահմանում 3.1.2:** $G$ գրաֆի $e$ կողը կոչվում է *կամուրջ*, եթե $c(G - e) > c(G)$:

Նկատենք, որ եթե $e$-ն $G$ կապակցված գրաֆի կամուրջ է, ապա $G - e$ գրաֆը կապակցված չէ:

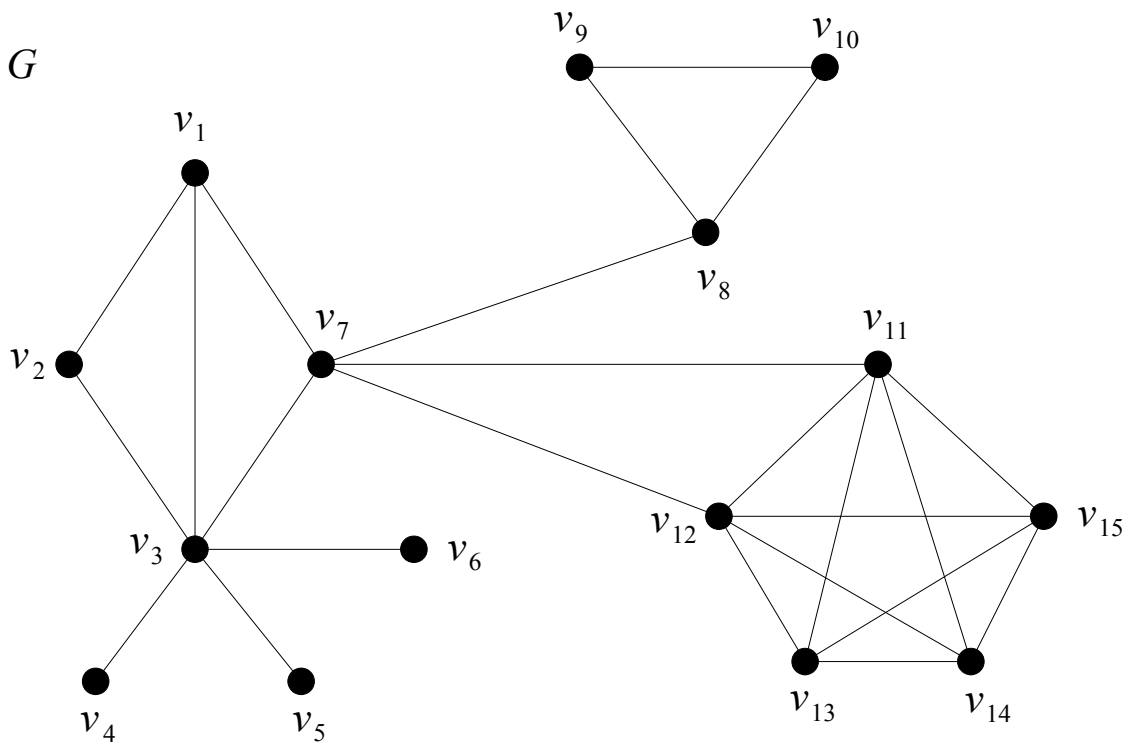

Նկ. 3.1.1

Դիտարկենք նկ. 3.1.1-ում պատկերված $G$ գրաֆը: Հեշտ է տեսնել, որ $v_3, v_7, v_8$

գագաթները հանդիսանում են $G$ գրաֆի միակցման կետեր, իսկ $v_3v_4, v_3v_5, v_3v_6, v_7v_8$ կողերը՝ կամուրջներ։

Ստորև կձևակերպենք և կապացուցենք թեորեմ միակցման կետերի մասին։

**Թեորեմ 3.1.1:** Եթե $v$-ն $G$ կապակցված գրաֆի գագաթ է, ապա հետևյալ երեք պնդումներն իրար համարժեք են.

(1) $v$-ն $G$ գրաֆի միակցման կետ է,

(2) $G$ գրաֆում գոյություն ունեն $v$-ից տարբեր $u$ և $w$ գագաթներ, որ $v$-ն պատկանում է ցանկացած պարզ $(u, w)$-ճանապարհին,

(3) գոյություն ունի $V(G)\setminus\{v\}$ գագաթների բազմության տրոհում $U$ և $W$ ենթաբազմությունների, որ ցանկացած $u \in U$ և $w \in W$ գագաթների համար $v$-ն պատկանում է յուրաքանչյուր պարզ $(u, w)$-ճանապարհին։

**Ապացույց:** Նախ ցույց տանք, որ (1)-ից հետևում է (3)-ը։ Իրոք, քանի որ $v$-ն $G$ գրաֆի միակցման կետ է, ուստի $G − v$ գրաֆը կապակցված չէ և, հետևաբար, պարունակում է առնվազն երկու կապակցված բաղադրիչ։ Դիտարկենք $V(G)\setminus\{v\}$ գագաթների բազմության տրոհում $U$ և $W$ ենթաբազմությունների, որտեղ $U$-ն $G − v$ գրաֆի որևէ կապակցված բաղադրիչի գագաթների բազմությունն է, իսկ $W$-ն՝ $V(G)\setminus(U \cup \{v\})$ գագաթների բազմությունն է։ Պարզ է, որ այս դեպքում ցանկացած $u \in U$ և $w \in W$ գագաթների համար, $u$ և $w$ գագաթները պատկանում են $G − v$ գրաֆի կապակցվածության տարբեր բաղադրիչներին, իսկ դա նշանակում է, որ $G$ գրաֆում յուրաքանչյուր պարզ $(u, w)$-ճանապարհի պարունակում է $v$ գագաթը։

Նկատենք, որ (3)-ից հետևում է (2)-ը, քանի որ (2)-ը (3)-ի մասնավոր դեպքն է։

Ցույց տանք, որ (2)-ից հետևում է (1)-ը։ Իրոք, եթե գոյություն ունեն $v$-ից տարբեր այնպիսի $u$ և $w$ գագաթներ, որ $v$-ն պատկանում է ցանկացած պարզ $(u, w)$-ճանապարհին, ապա $u$ և $w$ գագաթները պատկանում են $G − v$ գրաֆի տարբեր կապակցվածության բաղադրիչներին, իսկ դա նշանակում է, որ $v$-ն $G$ գրաֆի միակցման կետ է։ ∎

Պարզվում է, որ առնվազն երկու գագաթ ունեցող ցանկացած կապակցված գրաֆում գոյություն ունեն գագաթներ, որոնք միակցման կետեր չեն։

**Թեորեմ 3.1.2:** Եթե $G$-ն առնվազն երկու գագաթ ունեցող կապակցված գրաֆ է, ապա այն պարունակում է առնվազն երկու գագաթ, որոնք միակցման կետեր չեն։

**Ապացույց:** Ապացույցի համար դիտարկենք $G$ գրաֆի $u$ և $v$ գագաթները, որոնց համար $d(u,v) = d(G)$։ Ենթադրենք, $v$-ն միակցման կետ է $G$ գրաֆում։ Այդ դեպքում գոյություն ունի $w$ գագաթ, որը պատկանում է $G - v$ գրաֆի կապակցվածության այն բաղադրիչին, որը չի պարունակում $u$ գագաթը։ Քանի որ $u$ և $w$ գագաթները պատկանում են $G - v$ գրաֆի կապակցվածության տարբեր բաղադրիչներին, ուստի $v$-ն պատկանում է $G$ գրաֆի ցանկացած պարզ $(u,w)$-ճանապարհին և, հետևաբար, $d(u,w) > d(u,v)$, որը հակասում է $d(u,v) = d(G)$ պայմանին։ Այստեղից հետևում է, որ $v$-ն միակցման կետ չէ։ Համանման ձևով կարելի է ցույց տալ, որ $u$-ն միակցման կետ չէ։ ∎

Այժմ անցնենք կամուրջների ուսումնասիրմանը։ Առաջին թեորեմը, որը մենք կապացուցենք, տալիս է կամուրջների նկարագրում։

**Թեորեմ 3.1.3:** $G$ գրաֆի $e$ կողը կամուրջ է այն և միայն այն դեպքում, երբ այն չի պատկանում $G$ գրաֆի ոչ մի պարզ ցիկլին։

**Ապացույց:** Նախ նկատենք, որ $G$ գրաֆի $e = uv$ կողը կամուրջ է այն և միայն այն դեպքում, երբ $u$ և $v$ գագաթները պատկանում են $G - e$ գրաֆի կապակցվածության տարբեր բաղադրիչներին։

Ենթադրենք $e = uv$-ն $G$ գրաֆի կամուրջ է և $C = P, u, v, Q$-ն $G$ գրաֆի պարզ ցիկլ է, որտեղ $P$-ն և $Q$-ն $G$ գրաֆի պարզ ճանապարհներ են։ Հեշտ է տեսնել, որ այս դեպքում $Q, P$-ն պարզ $(v, u)$-ճանապարհ է $G - e$ գրաֆում, ինչը հակասում է $uv$ կողի $G$-ում կամուրջ լինելու ենթադրությանը։

Ենթադրենք $e = uv$-ն չի պատկանում $G$ գրաֆի ոչ մի պարզ ցիկլին և $e$-ն $G$ գրաֆի կամուրջ չէ։ Այստեղից հետևում է, որ $u$ և $v$ գագաթները պատկանում են $G - e$ գրաֆի կապակցվածության միևնույն բաղադրիչին, ուստի $G - e$ գրաֆում գոյություն ունի $P$ պարզ $(u, v)$-ճանապարհ։ Հեշտ է տեսնել, որ այս դեպքում $C = P, v, u$-ն $G$ գրաֆի պարզ ցիկլ է, որին պատկանում է $e$ կողը։ ∎

Նկատենք, որ նոր ապացուցված թեորեմից և սահմանում 2.3.1-ից բխում է.

**Հետևանք 3.1.1:** Կապակցված $G$ գրաֆը ծառ է այն և միայն այն դեպքում, երբ $G$ գրաֆի յուրաքանչյուր կող կամուրջ է։

**Դիտողություն 3.1.1:** Նշենք, որ եթե $G$ գրաֆում կամուրջին կից է գագաթ, որը կախված չէ, ապա այդ գագաթը միակցման կետ է։ Իրոք, ենթադրենք $e = uv$ կողը $G$ գրաֆի կամուրջ է և $d_G(v) > 1$։ Այդ դեպքում գոյություն ունի $u$-ից տարբեր և $v$-ին հարևան $w$ գագաթ $G$ գրաֆում։ Եթե $G$ գրաֆում գոյություն ունի $P$ պարզ $(u, w)$-ճանապարհ, որը չի պարունակում $v$ գագաթը, ապա $C = P, v, u$-ն $G$ գրաֆի պարզ ցիկլ է, ինչը հակասում է $e = uv$ կողի $G$-ում կամուրջ լինելուն։ Հետևաբար, $v$ գագաթը պատկանում է ցանկացած պարզ $(u,w)$-ճանապարհին և, համաձայն թեորեմ 3.1.1-ի (2) պնդման, այն միակցման կետ է։ Պետք է նշել, որ հակառակ պնդումը սխալ է. ոչ բոլոր միակցման կետերը կից են կամուրջների։ Այսպես, օրինակ, նկ. 3.1.2-ում պատկերված գրաֆի $v_4$ գագաթը միակցման կետ է, սակայն $v_2v_4, v_3v_4, v_4v_5, v_4v_6$ կողերից ոչ մեկը կամուրջ չէ պատկերված գրաֆում։

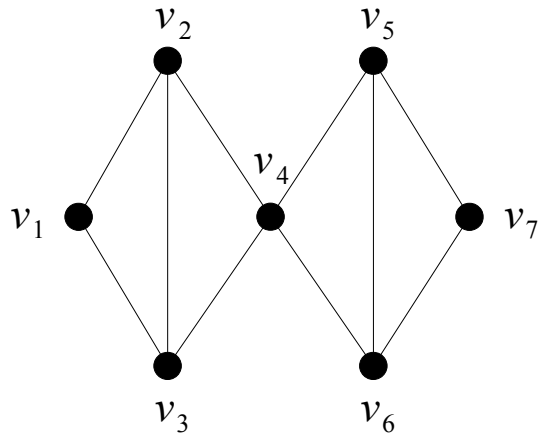

Նկ. 3.1.2

Թեորեմ 3.1.1-ի ապացույցի դատողությունները կրկնելով և հաշվի առնելով թեորեմ 3.1.3-ը, կարելի է ապացուցել հետևյալ թեորեմը։

**Թեորեմ 3.1.4:** Եթե $e$-ն $G$ կապակցված գրաֆի կող է, ապա հետևյալ չորս պնդումներն իրար համարժեք են.

(1) $e$-ն $G$ գրաֆի կամուրջ է,

(2) $e$ կողը չի պատկանում $G$ գրաֆի ոչ մի պարզ ցիկլին,

(3) $G$ գրաֆում գոյություն ունեն այնպիսի $u$ և $v$ գագաթներ, որ $e$-ն պատկանում է ցանկացած պարզ $(u,v)$-ճանապարհին,

(4) գոյություն ունի $V(G)$ գագաթների բազմության տրոհում $U$ և $W$ ենթաբազմությունների, որ ցանկացած $u \in U$ և $w \in W$ գագաթների համար $e$-ն պատկանում է յուրաքանչյուր պարզ $(u,w)$-ճանապարհին։

# § 3.2. Կապակցվածություն և կողային կապակցվածություն

Դիցուք $G = (V, E)$-ն գրաֆ է։

**Սահմանում 3.2.1:** $G$ գրաֆի *կապակցվածություն* (կնշանակենք $\varkappa(G)$-ով) կոչվում է գագաթների նվազագույն քանակը, որոնց հեռացման արդյունքում առաջանում է ոչ կապակցված գրաֆ կամ $K_1$։

Նկատենք, որ այս սահմանումից հետևում է, որ եթե $G$-ն կապակցված չէ, ապա $\varkappa(G) = 0$, իսկ եթե $G$-ն կապակցված է և պարունակում է միակցման կետ, ապա $\varkappa(G) = 1$։ Նշենք նաև, որ $K_n$ լրիվ գրաֆից գագաթներ հեռացնելով հնարավոր չէ ստանալ ոչ կապակցված գրաֆ, իսկ $K_1$ գրաֆը ստացվում է $K_n$-ից $n-1$ հատ գագաթներ հեռացնելով, ուստի $\varkappa(K_n) = n - 1$։ Դժվար չէ ցույց տալ, որ $\varkappa(K_{m,n}) = \min\{m, n\}$։

**Սահմանում 3.2.2:** $G$ գրաֆի *կողային կապակցվածություն* (կնշանակենք $\lambda(G)$-ով) կոչվում է կողերի նվազագույն քանակը, որոնց հեռացման արդյունքում առաջանում է ոչ կապակցված գրաֆ կամ $K_1$։

Նկատենք, որ այս սահմանումից հետևում է, որ $\lambda(K_1) = 0$ և եթե $G$-ն կապակցված չէ, ապա $\lambda(G) = 0$, իսկ եթե $G$-ն կապակցված է և պարունակում է կամուրջ, ապա $\lambda(G) = 1$։ Դժվար չէ ցույց տալ, որ $\lambda(K_n) = n - 1$։

Պարզվում է, որ ցանկացած գրաֆում կապացվածության, կողային կապակցվածության և գրաֆի նվազագույն աստիճանի միջև գոյություն ունի կապ։

**Թեորեմ 3.2.1 (Ուիտնի):** Կամայական $G$ գրաֆում տեղի ունի

$$\varkappa(G) \leq \lambda(G) \leq \delta(G)$$

անհավասարությունը։

**Ապացույց:** Նախ ցույց տանք, որ $\lambda(G) \leq \delta(G)$։ Իրոք, եթե $E(G) = \emptyset$, ապա $\lambda(G) = 0$, հակառակ դեպքում՝ դիտարկենք $G$ գրաֆի նվազագույն աստիճան ունեցող որևէ գագաթ։ Դիցուք այդ գագաթը $v$-ն է։ Հեռացնենք $G$ գրաֆից $\partial_G(v)$-ի կողերը։ Պարզ է, որ ստացված գրաֆը կլինի ոչ կապակցված գրաֆ, ուստի $\lambda(G) \leq |\partial_G(v)| = \delta(G)$։

Այժմ համոզվենք, որ $\varkappa(G) \leq \lambda(G)$։ Դիտարկենք դեպքեր: Եթե $G$-ն կապակցված չէ կամ $G = K_1$, ապա $\varkappa(G) = \lambda(G) = 0$։ Եթե $G$-ն կապակցված է և պարունակում է կամուրջ, ապա $\lambda(G) = 1$։ Այս դեպքում, ըստ դիտողություն 3.1.1-ի, $G$-ն ունի միակցման կետ կամ $G = K_2$, ուստի $\varkappa(G) = 1$։ Ենթադրենք $\lambda(G) \geq 2$։ Պարզ է, որ $G$ գրաֆում գոյություն ունեն $\lambda(G) - 1$ հատ կողեր, որոնց հեռացումը առաջացնում է կամուրջ պարունակող գրաֆ։ Դիցուք այդ կամուրջը $e = uv$-ն է։ Յուրաքանչյուր հեռացված կողի համար ընտրենք այդ կողին կից $u$ և $v$ գագաթներից տարբեր գագաթ։ Այնուհետև հեռացնենք $G$ գրաֆից բոլոր ընտրված գագաթները։ Ակնհայտ է, որ այդ հեռացման արդյունքում ստացված գրաֆում կբացակայեն վերը նշված $\lambda(G) - 1$ հատ կողերը։ Եթե արդյունքում ստացված գրաֆը կապակցված չէ, ապա $\varkappa(G) < \lambda(G)$, հակառակ դեպքում՝ այդ գրաֆը պարունակում է $e = uv$ կամուրջ և, հետևաբար, այդ գրաֆից հեռացնելով $u$ կամ $v$ գագաթը, մենք կստանանք ոչ կապակցված կամ $K_1$ գրաֆ։ Այստեղից հետևում է, որ $\varkappa(G) \leq \lambda(G)$։ ■

Չարտրանդի և Հարարիի կողմից ցույց է տրվել, որ թեորեմ 3.2.1-ում նշված՝ $\varkappa(G) \leq \lambda(G) \leq \delta(G)$ անհավասարությունը հնարավոր չէ ուժեղացնել։

**Թեորեմ 3.2.2 (Չարտրանդ, Հարարի):** Ցանկացած $a, b, c$ բնական թվերի համար, որտեղ $0 < a \leq b \leq c$, գոյություն ունի $G$ գրաֆ, որում $\varkappa(G) = a$, $\lambda(G) = b$ և $\delta(G) = c$։

**Ապացույց:** Կառուցենք $G$ գրաֆը հետևյալ կերպ.

$$V(G) = U \cup W, \text{ որտեղ } U = \{u_1, \ldots, u_c, u_{c+1}\}, W = \{w_1, \ldots, w_c, w_{c+1}\} \text{ և}$$

$$E(G) = \{u_i u_j, w_i w_j : 1 \leq i < j \leq c+1\} \cup \{u_i w_i : 1 \leq i \leq a\} \cup \{u_a w_j : a+1 \leq j \leq b\}.$$

Հաշվի առնելով $\varkappa(K_n) = \lambda(K_n) = n - 1$ հավասարությունը, դժվար չէ համոզվել, որ $G$ գրաֆում $\varkappa(G) = a$, $\lambda(G) = b$ և $\delta(G) = c$։ ■

Նկ. 3.2.1-ում պատկերված է թեորեմ 3.2.2-ի ապացույցում կառուցված $G$ գրաֆը $a = 2, b = 3$ և $c = 4$-ի դեպքում։

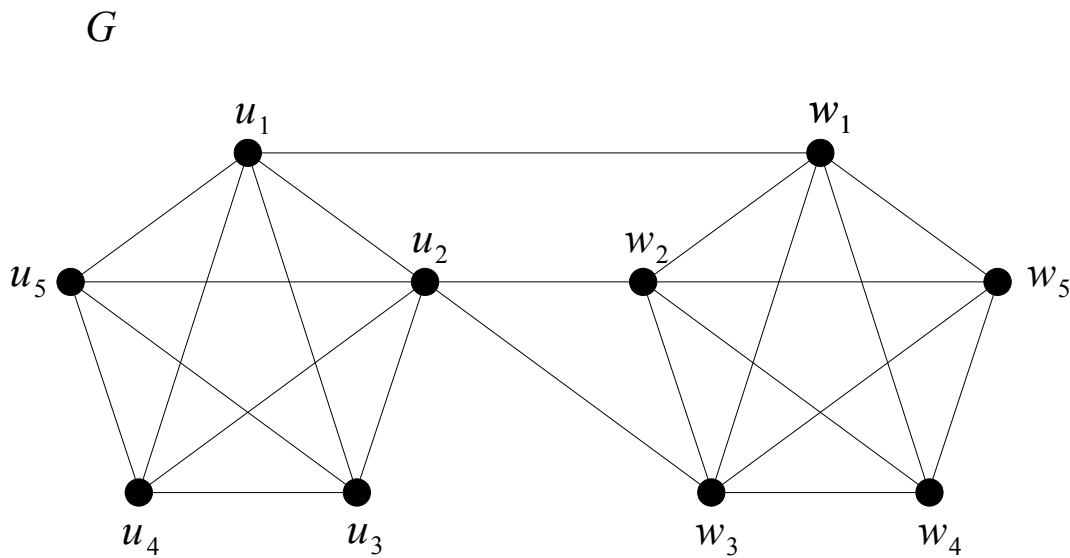

Նկ. 3.2.1

Առանց ապացույցի նշենք ևս մեկ այլ հետաքրքիր փաստ:

**Թեորեմ 3.2.3 (Չարտրանդ):** Եթե $n$ գագաթ ունեցող $G$ գրաֆի համար տեղի ունի $\delta(G) \geq \left\lfloor \frac{n}{2} \right\rfloor$ պայմանը, ապա $\lambda(G) = \delta(G)$։

**Սահմանում 3.2.3:** $G$ գրաֆը կոչվում է **$k$-կապակցված** գրաֆ, եթե $\varkappa(G) \geq k$։

Նկատենք, որ $G$ ($G \neq K_1$) գրաֆը **1**-կապակցված է այն և միայն այն դեպքում, երբ այն կապակցված է։

**Սահմանում 3.2.4:** $G$ գրաֆը կոչվում է **$k$-կողային կապակցված** գրաֆ, եթե $\lambda(G) \geq k$։

Նշենք, որ թեորեմ 3.2.1-ից հետևում է, որ եթե $G$ գրաֆը $k$-կապակցված է, ապա այն նաև $k$-կողային կապակցված է։ Մյուս կողմից նկ. 3.2.1-ում պատկերված $G$ գրաֆը ցույց է տալիս, որ հակառակ պնդումը սխալ է։ Նկատենք նաև, որ եթե $n$ գագաթ ունեցող $G$ գրաֆը $k$-կապակցված է կամ $k$-կողային կապակցված է, ապա թեորեմ 3.2.1-ից հետևում է, որ $\delta(G) \geq k$ և, հետևաբար, այդ գրաֆը պարունակում է առնվազն $\left\lceil \frac{n \cdot k}{2} \right\rceil$ հատ կող։ Պարզվում է, որ այս գնահատականը հասանելի է։ Դրա համար սահմանենք $n$ գագաթ ($n \geq 3$) ունեցող *պարզ ցիկլի $k$-րդ աստիճան* $C_n^k$ գրաֆը ($1 \leq k < \frac{n}{2}$) հետևյալ կերպ. դասավորում ենք շրջանագծի վրա $v_1, \ldots, v_n$ գագաթները, այնուհետև յուրաքանչյուր $v_i$ գագաթը միացնում ենք կողով շրջանագծի աջ և ձախ մասում գտնվող մոտակա $k$ հատ գագաթներին։ Նկ. 3.2.2-ում պատկերված է $C_8^3$ գրաֆը։ Հեշտ է տեսնել, որ $C_n^k$ գրաֆը **$2k$-համասեռ** գրաֆ է և $\left|E(C_n^k)\right| = n \cdot k$։

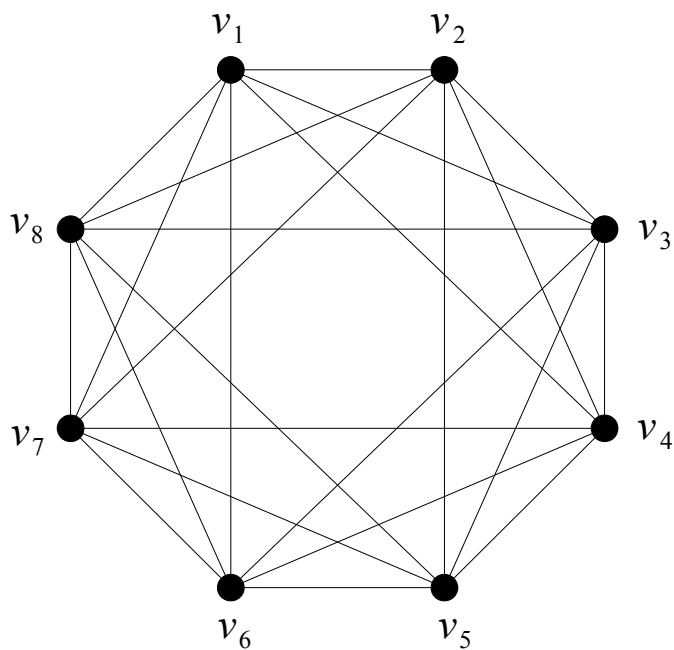

Նկ. 3.2.2

**Թեորեմ 3.2.4 (Հարարի):** Յանկացած $n \geq 3$ և $1 \leq k < \frac{n}{2}$ թվերի համար տեղի ունի $\varkappa(C_n^k) = 2k$ հավասարությունը:

**Ապացույց:** Քանի որ $\delta(C_n^k) = 2k$, ուստի, ըստ թեորեմ 3.2.1-ի, $\varkappa(C_n^k) \leq 2k$: Ցույց տանք, որ $\varkappa(C_n^k) \geq 2k$: Ընդհակառակը, ընտրենք ցանկացած $S \subseteq V(C_n^k)$, որի համար $|S| < 2k$: Համոզվենք, որ $C_n^k - S$ գրաֆը կապակցված է: Ընտրենք ցանկացած $u, v \in V(C_n^k)\setminus S$: Պարզ է, որ $u$ և $v$ գագաթները $C_n^k$ գրաֆից հեռացնելուց հետո շրջանագծի վրա դասավորված գագաթները կտրոհվեն երկու $A$ և $B$ կիսաշրջանագծերի: Այստեղից հետևում է, որ $C_n^k - S$ գրաֆում $u$-ից $v$ գագաթ մենք կարող ենք հասնել շարժվելով $A$ կամ $B$ կիսաշրջանի վրայով: Մյուս կողմից, քանի որ $|S| < 2k$, ուստի $A$ կամ $B$ կիսաշրջանի վրա առկա են $k - 1$-ից ոչ ավելի $S$-ի գագաթներ: Որոշակիության համար ենթադրենք, որ $A$ կիսաշրջանի վրա են առկա $k - 1$-ից ոչ ավելի $S$-ի գագաթներ: Քանի որ, ըստ $C_n^k$ գրաֆի սահմանման, յուրաքանչյուր գագաթ միացված է կողով շրջանագծի աջ և ձախ մասում գտնվող մոտակա $k$ հատ գագաթներին, ուստի $C_n^k - S$ գրաֆի $A$ կիսաշրջանում մենք կարող ենք նշել $(u, v)$-ճանապարհի: ∎

## §3.3. 2-կապակցված և 3-կապակցված գրաֆներ

Դիցուք $G = (V, E)$-ն կապակցված գրաֆ է և $u, v \in V(G)$:

**Սահմանում 3.3.1:** $G$ գրաֆի $P$ և $Q$ պարզ $(u, v)$-ճանապարհները կոչվում են *գագաթներով չհատվող պարզ $(u, v)$-ճանապարհներ*, եթե $V(P) \cap V(Q) = \{u, v\}$:

**Սահմանում 3.3.2:** $G$ գրաֆի $P$ և $Q$ պարզ $(u, v)$-ճանապարհները կոչվում են *կողերով չհատվող պարզ $(u, v)$-ճանապարհներ*, եթե $E(P) \cap E(Q) = \emptyset$:

Նկատենք, որ գագաթներով չհատվող պարզ ճանապարհները նաև կողերով չհատվող պարզ ճանապարհներ են: Հեշտ է տեսնել, որ հակառակ պնդումը սխալ է:

Ստորև կձևակերպենք և կապացուցենք 2-կապակցված գրաֆների նկարագրումը տվող թեորեմ:

**Թեորեմ 3.3.1 (Ուիտնի):** Առնվազն երեք գագաթ պարունակող $G$ գրաֆը 2-կապակցված գրաֆ է այն և միայն այն դեպքում, երբ $G$ գրաֆի ցանկացած $u$ և $v$ գագաթները միացված են առնվազն երկու գագաթներով չհատվող պարզ $(u, v)$-ճանապարհներով:

**Ապացույց:** Նախ նկատենք, որ եթե $G$ գրաֆում ցանկացած $u$ և $v$ գագաթները միացված են առնվազն երկու գագաթներով չհատվող պարզ $(u,v)$-ճանապարհներով, ապա այդ գրաֆից մեկ գագաթ հեռացնելուց հետո ստացված գրաֆը կլինի կապակցված, ուստի $G$-ն միակցման կետեր չի պարունակում։ Այստեղից հետևում է, որ $G$-ն 2-կապակցված գրաֆ է։

Ենթադրենք $G$-ն 2-կապակցված գրաֆ է։ Ցույց տանք, որ այդ դեպքում $G$ գրաֆի ցանկացած $u$ և $v$ գագաթները միացված են առնվազն երկու գագաթներով չհատվող պարզ $(u,v)$-ճանապարհներով։ Ապացույցը կատարենք մակածման եղանակով ըստ $d_G(u,v)$-ի։ Եթե $d_G(u,v) = 1$, ապա $uv \in E(G)$։ Ինչպես արդեն նշել ենք նախորդ պարագրաֆում, եթե $G$ գրաֆը 2-կապակցված է, ապա այն նաև 2-կողային կապակցված է, ուստի $G - uv$ գրաֆը կապակցված է։ Այստեղից հետևում է, որ $G - uv$ գրաֆում գոյություն ունի $(u,v)$-ճանապարհ և, հետևաբար, ըստ լեմմա 1.2.1-ի, $G - uv$ գրաֆում գոյություն ունի $P$ պարզ $(u,v)$-ճանապարհ։ Վերցնելով որպես պարզ $(u,v)$-ճանապարհ $G$ գրաֆի $uv$ կողը, մենք կստանանք, որ այդ գրաֆում գոյություն ունեն գագաթներով չհատվող երկու պարզ $(u,v)$-ճանապարհներ։

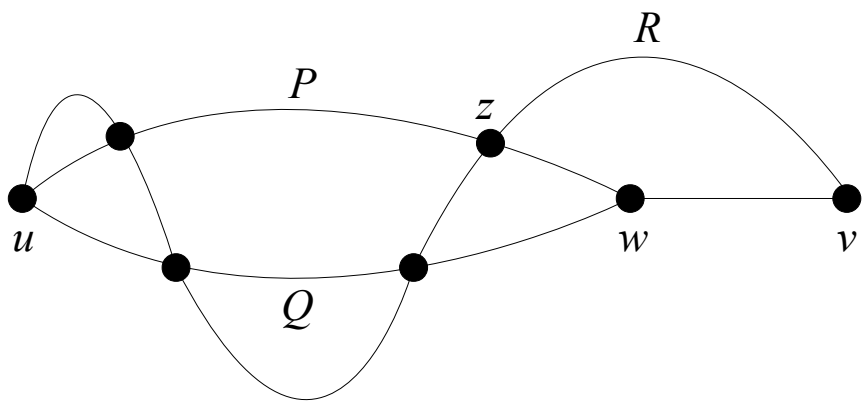

Նկ. 3.3.1

Ենթադրենք $d_G(u,v) = k > 1$ և պնդումը ճիշտ է $G$ գրաֆի ցանկացած $x$ և $y$ գագաթների դեպքում, որոնց համար $1 \leq d_G(x,y) < k$։ Դիցուք $w$ գագաթը անմիջապես նախորդում է $v$-ին կարճագույն $(u,v)$-ճանապարհի վրա։ Այդ դեպքում պարզ է, որ $d_G(u,w) = k - 1$ և, հետևաբար, ըստ մակածման ենթադրության, $G$ գրաֆում գոյություն ունեն $P$ և $Q$ գագաթներով չհատվող պարզ $(u,w)$-ճանապարհներ։ Քանի որ $G - w$ գրաֆը կապակցված է, ուստի այդ գրաֆում գոյություն ունի $(u,v)$-ճանապարհ և, հետևաբար, ըստ լեմմա 1.2.1-ի գոյություն ունի $R$ պարզ $(u,v)$-ճանապարհի։ Եթե $V(P) \cap V(R) = \{u\}$ կամ $V(Q) \cap V(R) = \{u\}$, ապա $R$-ը և $P, v$-ն կամ $R$-ը և $Q, v$-ն անհրաժեշտ գագաթներով չհատվող պարզ $(u,v)$-ճանապարհներն են $G$ գրաֆում։ Այժմ ենթադրենք, որ $\big(V(P) \cap V(R)\big)\setminus\{u\} \neq \emptyset$ և $\big(V(Q) \cap V(R)\big)\setminus\{u\} \neq \emptyset$։ Դիցուք $z$-ը $R$-ի վերջին գագաթն է, որը պատկանում է $V(P) \cup V(Q)$-ին (նկ. 3.3.1)։ Որոշակիության համար ենթադրենք, որ $z \in V(P)$։ Այդ դեպքում $P$ պարզ $(u,w)$-ճանապարհի $(u,z)$-ենթաճանապարհը և $R$ պարզ $(u,v)$-ճանապարհի $(z,v)$-ենթաճանապարհը առաջացնում են $G$ գրաֆի պարզ $(u,v)$-ճանապարհ։ Հեշտ է տեսնել, որ ստացված պարզ $(u,v)$-ճանապարհը և $Q, v$ պարզ $(u,v)$-ճանապարհը անհրաժեշտ գագաթներով չհատվող պարզ $(u,v)$-ճանապարհներն են $G$ գրաֆում։ ∎

Այժմ ապացուցենք ընդլայնման մի լեմմա, որն ունի պարզ ձևակերպում և հաճախ օգտագործվում է տարբեր թեորեմներ ապացուցելու ժամանակ։

**Լեմմա 3.3.1:** Եթե $G$-ն $k$-կապակցված գրաֆ է, իսկ $G'$ գրաֆը ստացվում է $G$-ին նոր $v$ գագաթ ավելացնելով և միացնելով այն $G$ գրաֆի առնվազն $k$ հատ գագաթներին, ապա $G'$-ը ևս կլինի $k$-կապակցված գրաֆ։

**Ապացույց:** Դիցուք $S \subseteq V(G')$ և $c(G'-S) \geq 2$։ Եթե $v \in S$, ապա $c\big(G-(S\setminus\{v\})\big) \geq 2$, որտեղից, քանի որ $G$-ն $k$-կապակցված գրաֆ է, բխում է $|S| \geq k+1$։ Եթե $v \notin S$ և $N_{G'}(v) \subseteq S$, ապա $|S| \geq k$։ Վերջապես, եթե $v \notin S$ և $N_{G'}(v) \not\subseteq S$, ապա քանի որ $G$-ն $k$-կապակցված գրաֆ է, ուստի $|S| \geq k$։ ∎

Ստորև կձևակերպենք և կապացուցենք ավելի ընդհանուր 2-կապակցված գրաֆների նկարագրումը։

**Թեորեմ 3.3.2:** Եթե $G$-ն առնվազն երեք գագաթ պարունակող գրաֆ է, ապա հետևյալ չորս պնդումները իրար համարժեք են.

(1) $G$-ն կապակցված գրաֆ է և այն չի պարունակում միակցման կետեր,

(2) $G$ գրաֆի ցանկացած $u$ և $v$ գագաթները միացված են առնվազն երկու գագաթներով չհատվող պարզ $(u,v)$-ճանապարհներով,

(3) $G$ գրաֆում ցանկացած $u$ և $v$ գագաթների համար գոյություն ունի այնպիսի $C$ պարզ ցիկլ, որ $u,v \in V(C)$,

(4) $G$ գրաֆում ցանկացած $e$ և $e'$ կողերի համար գոյություն ունի այնպիսի $C$ պարզ ցիկլ, որ $e,e' \in E(C)$։

**Ապացույց:** Նախ նկատենք, որ թեորեմ 3.3.1-ից հետևում է, որ (1)-ը և (2)-ը իրար համարժեք են։ Նկատենք նաև, որ եթե $G$ գրաֆի ցանկացած $u$ և $v$ գագաթները միացված են առնվազն երկու գագաթներով չհատվող պարզ $(u,v)$-ճանապարհներով, ապա այդ

Ճանապարհները առաջացնում են $C$ պարզ ցիկլ, որ $u, v \in V(C)$ և հակառակը՝ եթե $G$ գրաֆում ցանկացած $u$ և $v$ գագաթների համար գոյություն ունի $C$ պարզ ցիկլ, որ $u, v \in V(C)$, ապա այդ ցիկլը կարելի է տրոհել երկու գագաթներով չհատվող պարզ $(u, v)$-ճանապարհների։ Այստեղից հետնում է, որ (2)-ը և (3)-ը ևս իրար համարժեք են։

Ցույց տանք, որ (4)-ից հետնում է (3)-ը։ Նախ նկատենք, որ եթե $G$ գրաֆում ցանկացած $e$ և $e'$ կողերի համար գոյություն ունի այնպիսի $C$ պարզ ցիկլ, որ $e, e' \in E(C)$, ապա $\delta(G) \geq 2$։ Վերցնենք $G$ գրաֆում ցանկացած $u$ և $v$ գագաթներ և ընտրենք այդ գագաթներին կից $e$ և $e'$ կողերը։ Այդ դեպքում պարզ է, որ $e$ և $e'$ կողերը պարունակող $C$ պարզ ցիկլը նաև կպարունակի $u$ և $v$ գագաթները, ուստի (4)-ից հետնում է (3)-ը։

Թեորեմը ապացուցելու համար բավական է ցույց տալ, որ (1) և (3)-ից հետնում է (4)-ը։ Ենթադրենք, որ $G$-ն 2-կապակցված գրաֆ է և $uv, xy \in E(G)$։ Ավելացնենք $G$ գրաֆին $w$ գագաթը և միացնենք այն $u$ և $v$ գագաթներին, այնուհետև ավելացնենք $z$ գագաթը և միացնենք այն $x$ և $y$ գագաթներին։ Համաձայն լեմմա 3.3.1-ի, ստացված $G'$ գրաֆը ևս կլինի 2-կապակցված գրաֆ և, հետևաբար, ըստ (3)-ի $G'$ գրաֆում $w$ և $z$ գագաթների համար գոյություն ունի այնպիսի $C'$ պարզ ցիկլ, որ $w, z \in V(C')$։ Քանի որ $d_{G'}(w) = d_{G'}(z) = 2$, ուստի $uw, wv \in E(C')$ և $xz, zy \in E(C')$, բայց $uv \notin E(C')$ և $xy \notin E(C')$։ Այդ ցիկլից հեռացնելով $w$ և $z$ գագաթները և ավելացնելով $uv$ և $xy$ կողերը, մենք կստանանք $G$ գրաֆում այդ կողերը պարունակող պարզ ցիկլ։ ∎

**Լեմմա 3.3.2:** Եթե $G$-ն 2-կապակցված գրաֆ է և $e \in E(G)$, ապա $G_e$ գրաֆը, որը ստացվում է $G$-ից $e$ կողը տրոհելով, ևս կլինի 2-կապակցված գրաֆ։

**Ապացույց:** Դիցուք $e = uv \in E(G)$ և $uw, wv \in E(G_e)$։ Համաձայն թեորեմ 3.3.2-ի (4) պնդման այս թեորեմն ապացուցելու համար բավական է ցույց տալ, որ $G_e$ գրաֆում ցանկացած $f$ և $f'$ կողերի համար գոյություն ունի այնպիսի $C$ պարզ ցիկլ, որ $f, f' \in E(C)$։ Եթե $f, f' \in E(G)$, ապա $G$-ում, քանի որ այն 2-կապակցված գրաֆ է, գոյություն ունի այնպիսի $C$ պարզ ցիկլ, որ $f, f' \in E(C)$։ Եթե $uv \in E(C)$, ապա այդ կողը $C$ ցիկլում մենք կփոխարինենք $uw$ և $wv$ կողերով, հակառակ դեպքում՝ $G$ գրաֆի $C$ պարզ ցիկլը հանդիսանում է նաև $G_e$ գրաֆի պարզ ցիկլ։ Եթե $f \in E(G)$ և $f' \in \{uw, wv\}$, ապա փոխարինելով $e$ կողը $G$ գրաֆի $f$ և $e$ կողերը պարունակող պարզ ցիկլում $uw$ և $wv$ կողերով, մենք կստանանք $G_e$ գրաֆի $f$ և $f'$-ը պարունակող պարզ ցիկլ։ Վերջապես, եթե $f, f' = \{uw, wv\}$, ապա $G$-ն, քանի որ այն 2-կապակցված գրաֆ է և, հետևաբար, նաև 2-կողային կապակցված է, չի պարունակում կամուրջներ։ Ըստ թեորեմ 3.1.3-ի $G$ գրաֆում գոյություն ունի $e$ կողը պարունակող $C'$ պարզ ցիկլ, որում $e$ կողը փոխարինելով $uw$ և $wv$ կողերով մենք կստանանք $G_e$ գրաֆի $f$ և $f'$-ը պարունակող պարզ ցիկլ: ∎

Դիցուք $G = (V, E)$-ն գրաֆ է և $u, v \in V(G)$ (պարտադիր չէ, որ $u \neq v$):

**Սահմանում 3.3.3:** $G$ գրաֆին *ճանապարհի ավելացում* համարում ենք գրաֆին առնվազն $l$ ($l \geq 1$) երկարություն ունեցող պարզ $(u, v)$-ճանապարհի ավելացում, որի դեպքում $G$-ին ավելացվում է $l - 1$ հատ նոր գագաթ: Եթե $u \neq v$, ապա այդ ավելացված ճանապարհը կոչվում է *բաց ականջ* (կամ ուղղակի *ականջ*), հակառակ դեպքում՝ կոչվում է *փակ ականջ*:

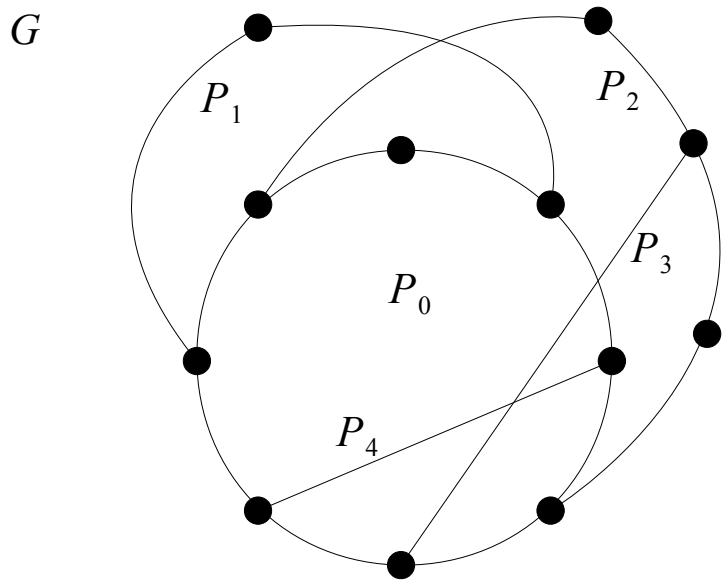

Նկ. 3.3.2

**Սահմանում 3.3.4:** $G$ գրաֆի *ականջային դեկոմպոզիցիա* կոչվում է այդ գրաֆի $G = P_0 \cup P_1 \cup \cdots \cup P_k$ տրոհումը, որի դեպքում $P_0$-ն $G$ գրաֆի պարզ ցիկլ է և ցանկացած $i \geq 1$ համար $P_i$-ն ականջ է արդեն կառուցված $P_0 \cup \cdots \cup P_{i-1}$ գրաֆի համար:

Նկ. 3.3.2-ում պատկերված է $G$ գրաֆը և նրա $G = P_0 \cup P_1 \cup P_2 \cup P_3 \cup P_4$ ականջային դեկոմպոզիցիան:

**Սահմանում 3.3.5:** $G$ գրաֆի *փակ ականջային դեկոմպոզիցիա* կոչվում է այդ գրաֆի $G = P_0 \cup P_1 \cup \cdots \cup P_k$ տրոհումը, որի դեպքում $P_0$-ն $G$ գրաֆի պարզ ցիկլ է և ցանկացած $i \geq 1$ համար $P_i$-ն բաց կամ փակ ականջ է արդեն կառուցված $P_0 \cup \cdots \cup P_{i-1}$ գրաֆի համար:

Նկ. 3.3.3-ում պատկերված է $G$ գրաֆը և նրա $G = P_0 \cup P_1 \cup P_2 \cup P_3 \cup P_4$ փակ ականջային դեկոմպոզիցիան:

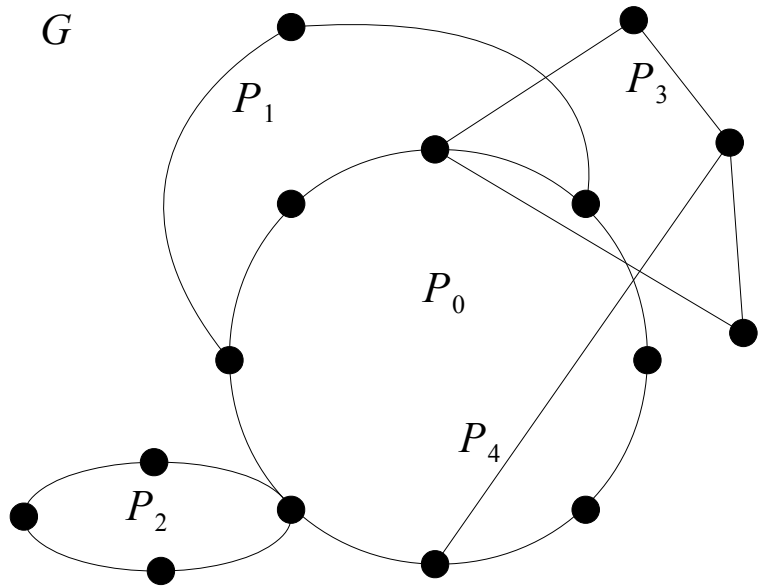

Նկ. 3.3.3

Ստորև մենք կձևակերպենք և կապացուցենք **2-կապակցված** և **2-կողային կապակցված** գրաֆների կոնստրուկտիվ նկարագրումը:

**Թեորեմ 3.3.3 (Ուիտնի):** Առնվազն երեք գագաթ պարունակող գրաֆը **2-կապակցված** գրաֆ է այն և միայն այն դեպքում, երբ այն ունի ականջային դեկոմպոզիցիա: Ավելին, ցանկացած պարզ ցիկլ **2-կապակցված** գրաֆում հանդիսանում է սկզբնական ցիկլ որևէ ականջային դեկոմպոզիցիայի համար:

**Ապացույց:** Նախ ցույց տանք, որ եթե գրաֆն ունի ականջային դեկոմպոզիցիա, ապա այն **2-կապակցված** գրաֆ է: Քանի որ ցանկացած պարզ ցիկլ **2-կապակցված** գրաֆ է, ուստի պնդումը ապացուցելու համար բավական է համոզվել, որ ճանապարհի ավելացումը պահպանում է **2-կապակցված** լինելու հատկությունը: Դիցուք $u$-ն և $v$-ն ($u \neq v$) $P$ ականջի ծայրակետեր են $G$ **2-կապակցված** գրաֆում: Համաձայն լեմմա 3.3.2-ի $G + uv$ գրաֆը ևս **2-կապակցված** գրաֆ է: Հեշտ է տեսնել, որ կողի տրոհման գործողության հաջորդական կիրառումը $G + uv$ գրաֆը դարձնում է $G \cup P$ գրաֆ, որը համաձայն լեմմա 3.3.2-ի ևս կլինի **2-կապակցված** գրաֆ:

Դիցուք տրված է $G$ **2-կապակցված** գրաֆը: Մենք կկառուցենք այդ գրաֆի համար ականջային դեկոմպոզիցիա, վերցնելով, որպես սկզբնական ցիկլ, այդ գրաֆի ցանկացած $C$ պարզ ցիկլ: Դիցուք $G_0 = C$: Ենթադրենք, որ մենք արդեն կառուցել ենք $G_i$ ենթագրաֆը որոշ ականջներ ավելացնելով: Եթե $G_i \neq G$, ապա ընտրենք $uv \in E(G) \setminus E(G_i)$ և $xy \in E(G_i)$: Քանի որ $G$-ն **2-կապակցված** գրաֆ է, ապա, համաձայն թեորեմ 3.3.2-ի (4) պնդման, $G$ գրաֆում գոյություն ունի այնպիսի $C'$ պարզ ցիկլ, որ $uv, xy \in E(C')$: Դիցուք $P$-ն $C'$ պարզ ցիկլի ենթաճանապարհի է, որը պարունակում է *uv* կողը և ճիշտ երկու գագաթներ $G_i$-ից, որոնք այդ ճանապարհի ծայրակետերն են։ Հեշտ է տեսնել, որ այդ *P*-ն հանդիսանում է $G_i$ գրաֆի ականջ և, հետևաբար, մենք կառուցեցինք $G_i$-ից ավելի մեծ $G_{i+1}$ ենթագրաֆ, որտեղ $G_{i+1} = G_i \cup P$։ Այս կառուցման ընթացքը կավարտվի, երբ *G*-ն ամբողջությամբ սպառվի։ ∎

**Թեորեմ 3.3.4 (Ուիտնի):** Առնվազն երեք գագաթ պարունակող գրաֆը **2-կողային կապակցված գրաֆ** է այն և միայն այն դեպքում, երբ այն ունի փակ ականջային դեկոմպոզիցիա։ Ավելին, ցանկացած պարզ ցիկլ **2-կողային կապակցված գրաֆում** հանդիսանում է սկզբնական ցիկլ որևէ փակ ականջային դեկոմպոզիցիայի համար։

**Ապացույց:** Նախ ցույց տանք, որ եթե գրաֆն ունի փակ ականջային դեկոմպոզիցիա, ապա այն **2-կապակցված գրաֆ** է։ Քանի որ ըստ թեորեմ 3.1.3-ի կամուրջները չեն պատկանում պարզ ցիկլերին, ուստի կապակցված գրաֆը **2-կողային կապակցված գրաֆ** է այն և միայն այն դեպքում, երբ նրա յուրաքանչյուր կող պատկանում է որևէ պարզ ցիկլի։ Պարզ է, որ սկզբնական պարզ ցիկլը **2-կողային կապակցված գրաֆ** է։ Նկատենք, որ բաց կամ փակ *P* ականջ ավելացնելու ժամանակ **2-կողային կապակցված** *G* գրաֆին, այդ գրաֆում արդեն գոյություն ունի *P* ականջի ծայրակետերը (որոնք կարող են նաև համընկնել) միացնող պարզ ճանապարհի, ուստի *P* ականջի բոլոր կողերը կպատկանեն պարզ ցիկլի։ Այստեղից հետևում է, որ $G \cup P$ գրաֆը ևս կլինի **2-կողային կապակցված գրաֆ**։

Դիցուք տրված է *G* **2-կողային կապակցված գրաֆը**։ Մենք կկառուցենք այդ գրաֆի համար փակ ականջային դեկոմպոզիցիա, վերցնելով, որպես սկզբնական ցիկլ, այդ գրաֆի ցանկացած *C* պարզ ցիկլ։ Դիցուք $G_0 = C$։ Ենթադրենք մենք արդեն կառուցել ենք $G_i$ ենթագրաֆ որոշ բաց կամ փակ ականջներ ավելացնելով։ Եթե $G_i \neq G$, ապա ընտրենք $uv \in E(G) \setminus E(G_i)$ և, հաշվի առնելով, որ *G*-ն կապակցված գրաֆ է, ենթադրենք, որ $u \in V(G_i)$։ Քանի որ *G*-ն **2-կողային կապակցված գրաֆ** է, ուստի, համաձայն թեորեմ 3.1.3-ի, *G* գրաֆում գոյություն ունի այնպիսի *C′* պարզ ցիկլ, որ $uv \in E(C')$։ *C′* ցիկլի վրայով շարժվելով մինչև $G_i$ գրաֆի գագաթ հանդիպելը, մենք կստենք բաց կամ փակ *P* ականջ և կկառուցենք $G_i$-ից ավելի մեծ $G_{i+1}$ ենթագրաֆ, որտեղ $G_{i+1} = G_i \cup P$։ Այս կառուցման ընթացքը կավարտվի, երբ *G*-ն ամբողջությամբ սպառվի։ ∎

Այս պարագրաֆի վերջում անդրադառնանք նաև **3-կապակցված գրաֆներին**, որոնց նկարագրումը տրվել է Տատտի կողմից։ Այդ նկարագրման մեջ Տատտի կողմից

ներմուծվել են անիվները:

**Սահմանում 3.3.6:** Եթե $n \geq 4$, ապա $n$ գագաթ ունեցող $W_n$ *անիվը* սահմանենք հետևյալ կերպ. $W_n = K_1 + C_{n-1}$:

Ստորև մենք կնշենք առանց ապացույցի 3-կապակցված գրաֆների մասին Տատտի թեորեմը:

**Թեորեմ 3.3.5 (Տատտ):** $G$ գրաֆը 3-կապակցված գրաֆ է այն և միայն այն դեպքում, երբ այն անիվ է կամ ստացվում է անիվից հետևյալ երկու գործողությունների միջոցով.

(1) նոր կող ավելացնելով,

(2) առնվազն չորս աստիճան ունեցող $v$ գագաթի փոխարինումով երկու $v'$ և $v''$ իրար հարևան գագաթներով այնպես, որ այն գագաթը, որը հարևան էր $v$-ին, հարևան լինի $v'$ և $v''$ գագաթներից միայն մեկին, ընդ որում $d(v') \geq 3$ և $d(v'') \geq 3$:

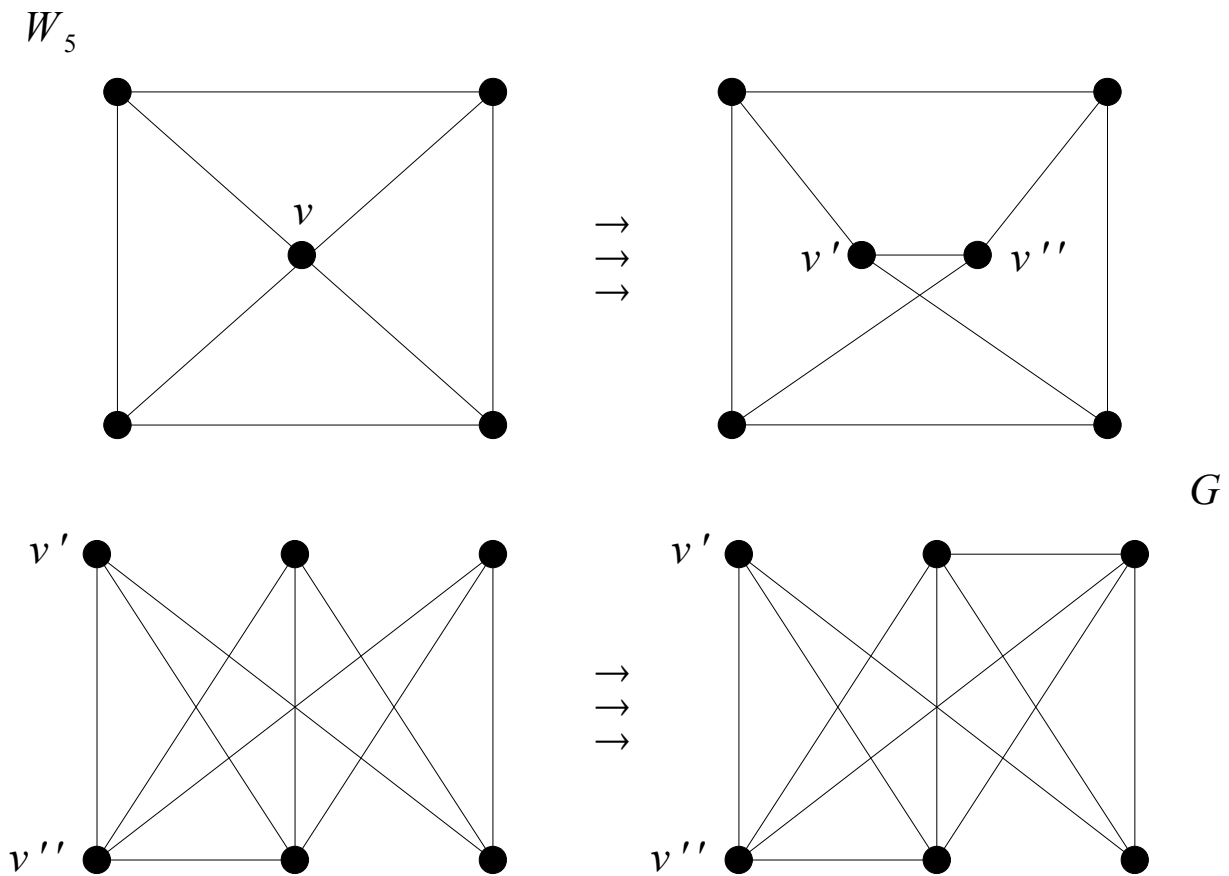

Նկ. 3.3.4

Նկ. 3.3.4-ում պատկերված $G$ գրաֆը 3-կապակցված է, քանի որ այն կարելի է ստանալ $W_5$-ից հաջորդաբար կիրառելով թեորեմ 3.3.5-ում նշված (1) և (2)

գործողությունները։ Այդ նկարում նաև պատկերված է $W_5$-ից $G$ գրաֆի ստանալու ամբողջ ընթացքը։

Վերջում ձևակերպենք և ապացուցենք Տոմասսենի թեորեմը 3-կապակցված գրաֆների մասին։

**Թեորեմ 3.3.6 (Տոմասսեն):** Առնվազն հինգ գագաթ ունեցող $G$ 3-կապակցված գրաֆում գոյություն ունի կող, որի կծկումը բերում է 3-կապակցված գրաֆի։

**Ապացույց:** Ենթադրենք հակառակը, ցանկացած $e \in E(G)$-ի համար $G/e$ գրաֆը 3-կապակցված չէ։ Հետևաբար, ցանկացած $e = xy \in E(G)$-ի համար գոյություն ունի $S \subseteq V(G/e)$, որ $|S| = 2$ և $c(G/e - S) \geq 2$, ընդ որում $S$-ի գագաթներից մեկը $e$ կողի կծկումից առաջացող գագաթն է։ Պարզ է, որ եթե $S$-ի մյուս գագաթը $z$-ն է, ապա $c(G - x - y - z) \geq 2$։ $G$ գրաֆի բոլոր կողերից ընտրենք այն $e = xy$ կողը և նրան համապատասխան մյուս $z$ գագաթը, որի դեպքում $G - x - y - z$ գրաֆի ամենաշատ գագաթներ պարունակող կապակցված բաղադրիչը $H$-ն է։ Դիցուք $H'$-ը $G - x - y - z$ գրաֆի մյուս կապակցված բաղադրիչն է։ Քանի որ $\{x, y, z\}$ բազմությունը ամենաշատ գագաթներ պարունակող բազմություն է, որ $c(G - x - y - z) \geq 2$, ուստի $x, y$ և $z$-ից յուրաքանչյուրը ունի հարևան գագաթ $H$-ում և $H'$-ում (նկ. 3.3.5)։ Դիցուք $u$-ն հարևան է $z$-ին $H'$-ում, իսկ $v$-ն մյուս գագաթն է, որ $c(G - u - v - z) \geq 2$։ Դիցուք $H_{xy} = G[V(H) \cup \{x, y\}]$։ Պարզ է, որ $H_{xy}$ ենթագրաֆը կապակցված է։ Նկատենք, որ եթե $v \in V(H_{xy})$, ապա $H_{xy} - v$ ենթագրաֆը ևս կապակցված է, քանի որ հակառակ դեպքում $c(G - v - z) \geq 2$։ Այստեղից հետևում է, որ $H_{xy} - v$ ենթագրաֆը պատկանում է $G - u - v - z$ գրաֆի այն կապակցված բաղադրիչին, որը պարունակում է ավելի շատ գագաթներ, քան $H$-ը, իսկ դա հակասում է $x, y$ և $z$-ի ընտրությանը։ ∎

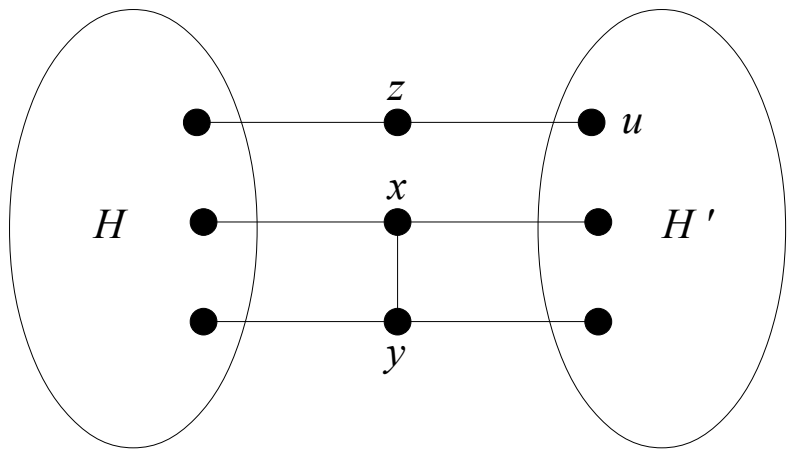

Նկ. 3.3.5

# § 3.4. *k*-կապակցված և *k*-կողային կապակցված գրաֆներ, Մենգերի թեորեմ

Դիցուք $G = (V, E)$-ն գրաֆ է, $S \subseteq V(G)$ և $u, v \in V(G) \setminus S$:

**Սահմանում 3.4.1:** $G$ գրաֆի գագաթների $S$ բազմությունը կոչվում է $(u, v)$-*կտրվածք*, եթե $u$ և $v$ գագաթները պատկանում են $G - S$ գրաֆի կապակցվածության տարբեր բաղադրիչների:

Նկատենք, որ եթե $G$ գրաֆում $uv \in E(G)$, ապա այդ գրաֆում գոյություն չունի $(u, v)$-կտրվածք:

**Սահմանում 3.4.2:** $G$ գրաֆի կողերի $F$ բազմությունը կոչվում է *կողային* $(u, v)$-*կտրվածք*, եթե $u$ և $v$ գագաթները պատկանում են $G - F$ գրաֆի կապակցվածության տարբեր բաղադրիչների:

Հիշենք նաև նախորդ պարագրաֆում բերված գագաթներով չհատվող պարզ $(u, v)$-ճանապարհների և կողերով չհատվող պարզ $(u, v)$-ճանապարհների սահմանումը: $G$ գրաֆի $P$ և $Q$ պարզ $(u, v)$-ճանապարհները կոչվում են գագաթներով (կողերով) չհատվող պարզ $(u, v)$-ճանապարհներ, եթե $V(P) \cap V(Q) = \{u, v\}$ ($E(P) \cap E(Q) = \emptyset$):

Ստորև կձևակերպենք և կապացուցենք գրաֆների կապակցվածությանը վերաբերվող դասական արդյունքներից մեկը:

**Թեորեմ 3.4.1 (Մենգեր):** Եթե $G$ գրաֆում $u$ և $v$ գագաթները հարևան չեն, ապա այդ գրաֆում գագաթների նվազագույն քանակը $(u, v)$-կտրվածքում հավասար է գագաթներով չհատվող պարզ $(u, v)$-ճանապարհների առավելագույն քանակին:

**Ապացույց:** Նախ նկատենք, որ եթե $S$-ը $G$ գրաֆի $(u, v)$-կտրվածք է, ապա ցանկացած պարզ $(u, v)$-ճանապարհ անցնում է $S$-ի որևէ գագաթով: Մյուս կողմից, քանի որ $G$ գրաֆի գագաթներով չհատվող պարզ $(u, v)$-ճանապարհները հատվում են միայն ծայրակետերում, ուստի $S$-ի ոչ մի գագաթ չի կարող միաժամանակ մասնակցել երկու տարբեր գագաթներով չհատվող պարզ $(u, v)$-ճանապարհներում: Այստեղից հետևում է, որ $S$-ը պարունակում է առնվազն այնքան գագաթներ, որքան գագաթներով չհատվող պարզ $(u, v)$-ճանապարհներ կան $G$ գրաֆում:

Պարզ է, որ թեորեմն ապացուցելու համար բավական է ցույց տալ, որ եթե $S$-ը $G$ գրաֆի նվազագույն հզորություն ունեցող $(u, v)$-կտրվածք է, ապա $G$ գրաֆում գոյություն ունեն $|S|$ հատ գագաթներով չհատվող պարզ $(u,v)$-ճանապարհներ։ Այդ պնդման ապացույցը կատարենք մակածման եղանակով ըստ $|V(G)| + |E(G)|$-ի։ Դիտարկենք այդ պնդման պայմանին բավարարող ամենափոքր գրաֆները։ Եթե $|V(G)| = 2$, ապա $V(G) = \{u, v\}$ և $E(G) = \emptyset$։ Հետևաբար, $u$ և $v$ գագաթները արդեն պատկանում են $G$ գրաֆի կապակցվածության տարբեր բաղադրիչների, ուստի $|S| = \emptyset$։ Մյուս կողմից ակնհայտ է, որ այդ գրաֆում գոյություն չունի պարզ $(u,v)$-ճանապարհ։ Եթե $V(G) = \{u, w, v\}$ և $E(G) = \{uw\}$ կամ $E(G) = \{wv\}$, ապա հեշտ է տեսնել, որ նորից $|S| = \emptyset$ և $G$ գրաֆում գոյություն չունի պարզ $(u,v)$-ճանապարհ։ Եթե $V(G) = \{u, w, v\}$ և $E(G) = \{uw, wv\}$, ապա հեշտ է տեսնել, որ $S = \{w\}$ և $G$ գրաֆում գոյություն ունի միակ $P = u, w, v$ պարզ $(u,v)$-ճանապարհ։ Ենթադրենք պնդումը ճիշտ է ցանկացած $G'$ գրաֆի դեպքում, որի համար $|V(G')| + |E(G')| < |V(G)| + |E(G)|$։ Դիտարկենք $G$ գրաֆը։ Նկատենք, որ ըստ մակածման ենթադրության մենք կարող ենք համարել, որ $G$-ն կապակցված գրաֆ է։ Դիցուք $S$-ը $G$ գրաֆի նվազագույն հզորություն ունեցող $(u,v)$-կտրվածք է։ Ցույց տանք, որ $G$ գրաֆում գոյություն ունեն $|S|$ հատ գագաթներով չհատվող պարզ $(u,v)$-ճանապարհներ։ Քննարկենք երեք դեպք։

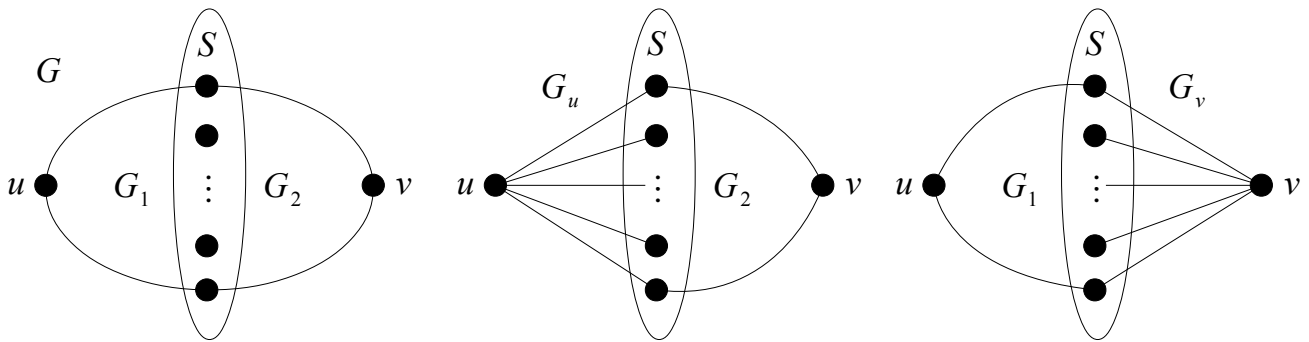

Նկ. 3.4.1

Դեպք 1: $S \not\subseteq N_G(u)$ և $S \not\subseteq N_G(v)$։

Այս դեպքում նկատենք, որ $G - S$ գրաֆը բաղկացած է $G_1$ և $G_2$ գրաֆներից, ընդ որում $u \in V(G_1), v \in V(G_2)$ և $|V(G_1)|, |V(G_2)| \geq 2$։ Դիցուք $G_u = G/G_1$ և $G_v = G/G_2$։ Այլ կերպ ասած, $G_u$ գրաֆը ստացվում է $G$ գրաֆից, կծկելով առնվազն երկու գագաթ պարունակող $G_1$ գրաֆը $u$ գագաթին, իսկ $G_v$-ն՝ կծկելով առնվազն երկու գագաթ պարունակող $G_2$ գրաֆը $v$ գագաթին (նկ. 3.4.1)։ Հեշտ է տեսնել, որ $S$-ը հանդիսանում է $G_u$ և $G_v$ գրաֆների նվազագույն հզորություն ունեցող $(u,v)$-կտրվածք։ Մյուս կողմից, քանի որ $|V(G_u)| < |V(G)|$ և $|V(G_v)| < |V(G)|$, ուստի համաձայն մակածման ենթադրության $G_u$ և

$G_v$ գրաֆներում գոյություն ունեն $|S|$ հատ զագաթներով չհատվող պարզ $(u,v)$-ճանապարհներ։ Միացնելով $G_v$ գրաֆի $u$-ից $S$ յուրաքանչյուր ճանապարհը $G_u$ գրաֆի $S$-ից $v$ համապատասխան ճանապարհի հետ, մենք կստանանք $G$ գրաֆում $|S|$ հատ զագաթներով չհատվող պարզ $(u,v)$-ճանապարհներ։

Դեպք 2: $S \subseteq N_G(u)$ կամ $S \subseteq N_G(v)$ և գոյություն ունի $w \in S$ զագաթ, որ $w \in N_G(u) \cap N_G(v)$։

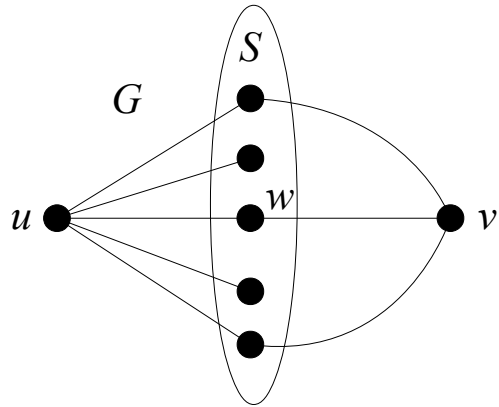

Նկ. 3.4.2

Որոշակիության համար ենթադրենք, որ $S \subseteq N_G(u)$։ Դիտարկենք $G' = G - w$ գրաֆը։ Հեշտ է տեսնել, որ $S \setminus \{w\}$ բազմությունը հանդիսանում է $G'$ գրաֆի նվազագույն հզորություն ունեցող $(u,v)$-կտրվածք (նկ. 3.4.2)։ Համաձայն մակածման ենթադրության, $G'$ գրաֆում գոյություն ունեն $|S \setminus \{w\}| = |S| - 1$ հատ զագաթներով չհատվող պարզ $(u,v)$-ճանապարհներ։ Ավելացնելով այդ ճանապարհներին $G$ գրաֆի $P = u, w, v$ պարզ $(u,v)$-ճանապարհը, մենք կստանանք $G$ գրաֆում $|S|$ հատ զագաթներով չհատվող պարզ $(u,v)$-ճանապարհներ։

Դեպք 3: $S \subseteq N_G(u)$ կամ $S \subseteq N_G(v)$ և գոյություն չունի $w \in S$ զագաթ, որ $w \in N_G(u) \cap N_G(v)$։

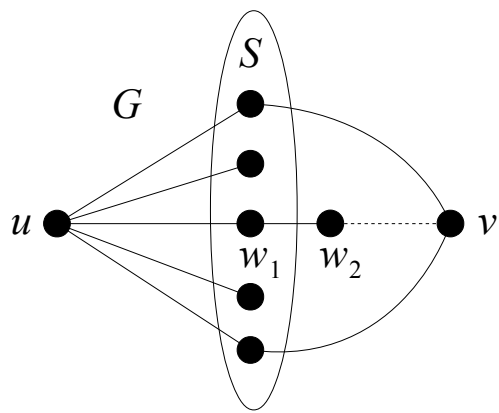

Նկ. 3.4.3

Որոշակիության համար նորից ենթադրենք, որ $S \subseteq N_G(u)$։ Դիտարկենք $G$ գրաֆում կարճագույն $(u,v)$-ճանապարհը։ Դիցուք այդ ճանապարհը $P = u, w_1, \ldots, w_k, v$-ն է։ Քանի որ $P$-ն կարճագույն $(u,v)$-ճանապարհ է և գոյություն չունի $w \in S$ գագաթ, որ $w \in N_G(u) \cap N_G(v)$, ուստի $P$-ն պարզ $(u,v)$-ճանապարհ է, $w_1 \in S$ և $w_2 \neq v$ (նկ. 3.4.3)։ Դիտարկենք $G' = G/w_1w_2$ գրաֆը։ Նկատենք, որ $w_2 \notin S$, քանի որ հակառակ դեպքում $G$ գրաֆի $P' = u, w_2, \ldots, w_k, v$ ճանապարհը կլիներ ավելի կարճ, քան $P$-ն։ Այստեղից հետևում է, որ $S$ բազմությունը $G'$ գրաֆի նվազագույն հզորություն ունեցող $(u,v)$-կտրվածք է։ Համաձայն մակածման ենթադրության, $G'$ գրաֆում գոյություն ունեն $|S|$ հատ գագաթներով չհատվող պարզ $(u,v)$-ճանապարհներ։ Քանի որ $G'$ գրաֆի գագաթներով չհատվող պարզ $(u,v)$-ճանապարհները նաև $G$ գրաֆի գագաթներով չհատվող պարզ $(u,v)$-ճանապարհներ են, ուստի $G$ գրաֆում գոյություն ունեն $|S|$ հատ գագաթներով չհատվող պարզ $(u,v)$-ճանապարհներ։ ∎

Նշենք, որ գոյություն ունեն Մենգերի թեորեմի այլ տարբերակներ։ Նշենք դրանցից մի քանիսը։

**Թեորեմ 3.4.2 (Մենգեր):** Եթե $G$ գրաֆի $u$ և $v$ գագաթները իրարից տարբեր են, ապա այդ գրաֆում կողերի նվազագույն քանակը կողային $(u,v)$-կտրվածքում հավասար է կողերով չհատվող պարզ $(u,v)$-ճանապարհների առավելագույն քանակին։

**Թեորեմ 3.4.3 (Ուիտնի):** Առնվազն $k+1$ գագաթ պարունակող $G$ գրաֆը $k$-կապակցված գրաֆ է այն և միայն այն դեպքում, երբ $G$ գրաֆի ցանկացած $u$ և $v$ գագաթները միացված են առնվազն $k$ հատ գագաթներով չհատվող պարզ $(u,v)$-ճանապարհներով։

**Թեորեմ 3.4.4 (Ուիտնի):** Առնվազն $k+1$ գագաթ պարունակող $G$ գրաֆը $k$-կողային կապակցված գրաֆ է այն և միայն այն դեպքում, երբ $G$ գրաֆի ցանկացած $u$ և $v$ գագաթները միացված են առնվազն $k$ հատ կողերով չհատվող պարզ $(u,v)$-ճանապարհներով։

Դիցուք $G = (V, E)$-ն գրաֆ է և $x \in V(G)$ և $U \subseteq V(G)$։

**Սահմանում 3.4.3:** $G$ գրաֆի պարզ ճանապարհների բազմությունը, որոնց մի ծայրակետը $x$-ն է, իսկ մյուսը $U$-ից է, ընդ որում այդ ճանապարհները հատվում են միայն $x$ գագաթում, կոչվում է $(x, U)$-*հովհար*։ Եթե $x \in U$, ապա $(x, U)$-հովհարը պարունակում է զրո երկարություն ունեցող ճանապարհի։

Ստորև մենք կձևակերպենք և կապացուցենք Դիրակի հովհարների լեմման։

**Թեորեմ 3.4.5 (Դիրակ):** $G$ գրաֆը $k$-կապակցված գրաֆ է այն և միայն այն դեպքում, երբ այն պարունակում է առնվազն $k+1$ գագաթ և ցանկացած $x$-ի և $|U| \geq k$ պայմանին բավարարող $U \subseteq V(G)$-ի համար $G$-ն ունի $(x, U)$-հովհար:

**Ապացույց:** Ենթադրենք, $G$ գրաֆը $k$-կապակցված է, $x \in V(G)$ և $U \subseteq V(G)$ ենթաբազմությունը բավարարում է $|U| \geq k$ պայմանին: Սահմանենք $G'$ գրաֆը հետևյալ կերպ. $V(G') = V(G) \cup \{y\}$, որտեղ $y \notin V(G)$, և $E(G') = E(G) \cup \{uy : u \in U\}$: Համաձայն լեմմա 3.3.1-ի $G'$ գրաֆը ևս $k$-կապակցված է և, հետևաբար, ըստ թեորեմ 3.4.3-ի $G'$ գրաֆում գոյություն ունեն առնվազն $k$ հատ գագաթներով չհատվող պարզ $(x, y)$-ճանապարհներ: Հեռացնելով $y$ գագաթը այդ ճանապարհներից, մենք կստանանք, որ $G$-ն ունի $(x, U)$-հովհար:

Ենթադրենք, $G$ գրաֆը պարունակում է առնվազն $k+1$ գագաթ և ցանկացած $x$-ի և $|U| \geq k$ պայմանին բավարարող $U \subseteq V(G)$-ի համար $G$-ն ունի $(x, U)$-հովհար: Ցույց տանք, որ $G$ գրաֆը $k$-կապակցված է: Ցանկացած $v \in V(G)$-ի և $U = V(G) - v$-ի համար գոյություն ունի $(v, U)$-հովհար, ուստի $\delta(G) \geq k$: Համաձայն թեորեմ 3.4.3-ի, $G$ գրաֆը $k$-կապակցված է այն և միայն այն դեպքում, երբ այդ գրաֆի ցանկացած $x$ և $y$ գագաթները միացված են առնվազն $k$ հատ գագաթներով չհատվող պարզ $(x, y)$-ճանապարհներով: Դիցուք $x, y \in V(G)$ և $U = N_G(y)$: Քանի որ $G$-ն ունի $(x, U)$-հովհար, ուստի, շարունակելով $(x, U)$-հովհարի $k$ հատ ճանապարհները մինչև $y$ գագաթը, մենք կստանանք, որ $x$ և $y$ գագաթները միացված են առնվազն $k$ հատ գագաթներով չհատվող պարզ $(x, y)$-ճանապարհներով: ∎

## Գլուխ 4

## Գրաֆների շրջանցումներ

### § 4.1. Էյլերյան ճանապարհներ և ցիկլեր

Դիցուք $G = (V, E)$-ն գրաֆ է։

**Սահմանում 4.1.1:** Կասենք, որ $G$ գրաֆում գոյություն ունի *Էյլերյան ճանապարհ*, եթե $G$-ում գոյություն ունի ճանապարհ, որը պարունակում է $G$-ի բոլոր գագաթները և կողերը։ Եթե այդ ճանապարհը ցիկլ է, ապա այն կանվանենք *Էյլերյան ցիկլ*։

Հեշտ է տեսնել, որ հարթության վրա էյլերյան ճանապարհ (ցիկլ) պարունակող գրաֆի պատկերը կարելի է նկարել առանց մատիտը թղթից կտրելու, յուրաքանչյուր կողով անցնելով ճիշտ մեկ անգամ (և վերադառնալ սկզբնական կետ):

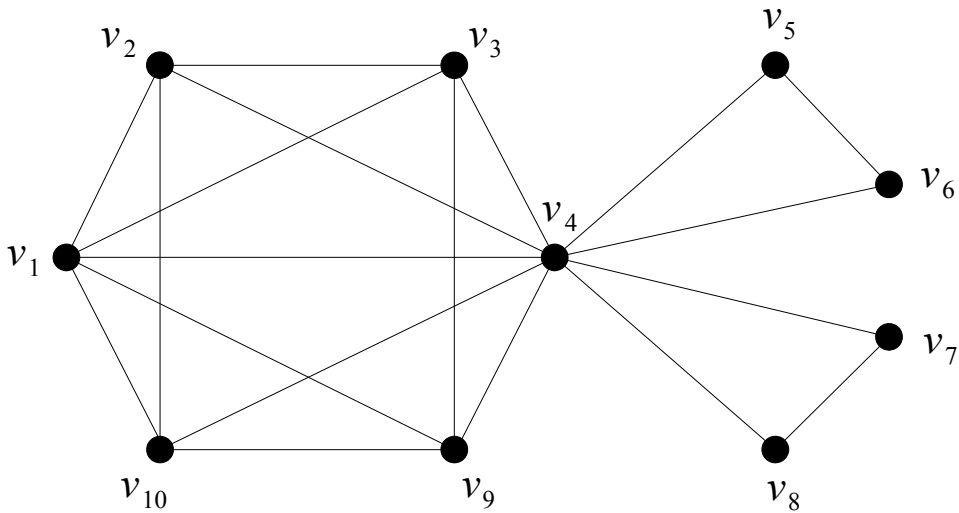

Նկ. 4.1.1

**Սահմանում 4.1.2:** $G$ գրաֆը կանվանենք *Էյլերյան* գրաֆ, եթե այն պարունակում է էյլերյան ցիկլ։

Բնական է դիտարկել հետևյալ խնդիրը. ինչպիսի պայմանների պետք է բավարարի գրաֆը, որպեսզի այն ունենա էյլերյան ճանապարհ կամ ցիկլ (լինի էյլերյան)։ Պարզ է, որ եթե $G$ գրաֆը պարունակում է էյլերյան ճանապարհ կամ ցիկլ, ապա $G$-ն պետք է լինի կապակցված: Սակայն ոչ բոլոր կապակցված գրաֆներն ունեն Էյլերյան ճանապարհ կամ Էյլերյան ցիկլ: Այսպես, օրինակ, նկ. 4.1.1-ում պատկերված կապակցված գրաֆը ունի Էյլերյան ճանապարհ, բայց չունի Էյլերյան ցիկլ: Իրոք, դիտարկենք այդ գրաֆում $(v_1, v_4)$-ճանապարհը՝ $v_1, v_2, v_3, v_4, v_5, v_6, v_4, v_7, v_8, v_4, v_9, v_3, v_1, v_9, v_{10}, v_4, v_2, v_{10}, v_1, v_4$: Հեշտ է տեսնել, որ նշված ճանապարհը Էյլերյան ճանապարհ է: Նկ. 4.1.2-ի ձախ մասում պատկերված է Էյլերյան գրաֆ (համոզվելու համար բավական է դիտարկել $u_1, u_2, u_3, u_4, u_5, u_6, u_3, u_5, u_2, u_6, u_1$ ցիկլը), իսկ աջ մասում՝ գրաֆ, որը չունի Էյլերյան ճանապարհի և հետևաբար նաև չունի Էյլերյան ցիկլ:

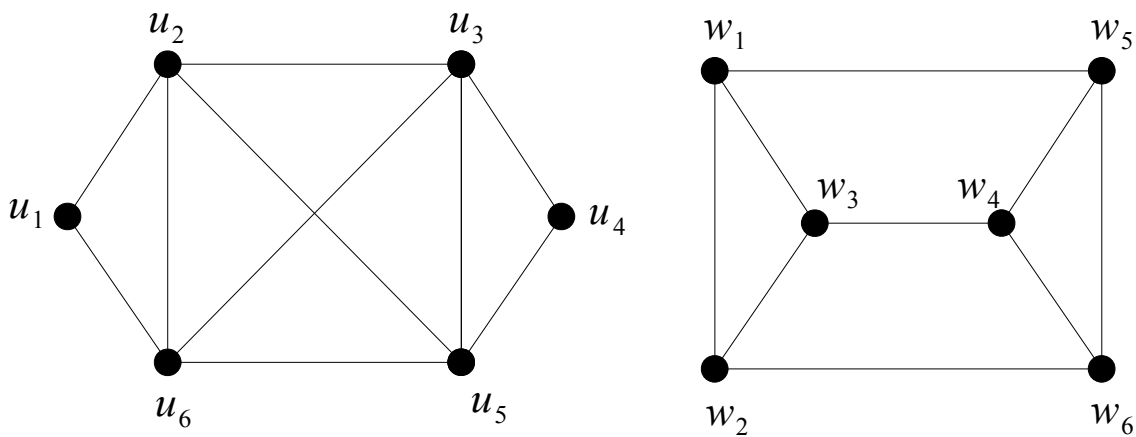

Նկ. 4.1.2

Նախ ապացուցենք մի լեմմա, որն ունի բավականին պարզ ձևակերպում և հաճախ է կիրառվում տարբեր թեորեմներ ապացուցելու ժամանակ:

**Լեմմա 4.1.1:** Եթե $G$ գրաֆի համար տեղի ունի $\delta(G) \geq 2$ պայմանը, ապա $G$-ն պարունակում է պարզ ցիկլ:

**Ապացույց:** Դիտարկենք $G$ գրաֆի առավելագույն երկարություն ունեցող որևէ պարզ ճանապարհ՝ $P = v_1, \ldots, v_s$: Քանի որ $P$-ի երկարությունը առավելագույնն է, ուստի այդ պարզ ճանապարհի $v_1$ և $v_s$ գագաթների բոլոր հարևանները պատկանում են $P$-ին (հակառակ դեպքում հնարավոր կլիներ նշել $P$-ից ավելի երկար պարզ ճանապարհ): Մյուս կողմից, քանի որ $\delta(G) \geq 2$, ուստի $v_1$ գագաթը կունենա հարևան $v_k$ գագաթ, որտեղ $3 \leq k \leq s$: Այստեղից հետևում է, որ $G$ գրաֆը պարունակում է $C = v_1, \ldots, v_k, v_1$ պարզ ցիկլ: ∎

**Սահմանում 4.1.3:** $G$ գրաֆը կանվանենք *զույգ* գրաֆ, եթե նրա բոլոր գագաթների աստիճանները զույգ թվեր են:

Ստորև կձևակերպենք և կապացուցենք Վեբլենի թեորեմը զույգ գրաֆների մասին:

**Թեորեմ 4.1.1:** $G$ գրաֆը զույգ է այն և միայն այն դեպքում, երբ $G$ գրաֆի կողերի բազմությունը կարելի է տրոհել կողերով չհատվող պարզ ցիկլերի:

**Ապացույց:** Դիցուք $G$-ն զույգ գրաֆ է: Ցույց տանք, որ $G$ գրաֆի կողերի բազմությունը կարելի է տրոհել կողերով չհատվող պարզ ցիկլերի: Ապացույցը կատարենք մակածման եղանակով ըստ $|E(G)|$-ի: Եթե $E(G) = \emptyset,$ ապա պարզ է որ պնդումը ճիշտ է, քանի որ այդ դեպքում $G$ գրաֆի կողերով չհատվող պարզ ցիկլերի բազմությունը դատարկ է: Ենթադրենք $E(G) \neq \emptyset$ և պնդումը ճիշտ է ցանկացած $H$ զույգ գրաֆի դեպքում, որի համար $|E(H)| < |E(G)|$: Դիտարկենք $G$ զույգ գրաֆը: Դիցուք $S = \{v : v \in V(G)$ և $d_G(v) > 0\}$: Պարզ է, որ $G$ գրաֆի $G[S]$ ենթագրաֆը ևս զույգ գրաֆ է և $\delta(G[S]) \geq 2$: Ըստ լեմմա 4.1.1-ի $G[S]$ գրաֆը պարունակում է $C$ պարզ ցիկլ: Դիտարկենք $G' = G \setminus E(C)$ գրաֆը: Պարզ է, որ $G'$ գրաֆը ևս զույգ գրաֆ է և, հետևաբար, ըստ մակածման ենթադրության $E(G') = E(C_1) \cup \cdots \cup E(C_l)$, որտեղ $C_1, \ldots, C_l$-ը $G'$ գրաֆի կողերով չհատվող պարզ ցիկլերն են: Այստեղից հետևում է, որ $E(G) = E(C_1) \cup \cdots \cup E(C_l) \cup E(C)$-ն իրենից ներկայացնում է $G$ գրաֆի կողերի բազմության տրոհում կողերով չհատվող պարզ ցիկլերի:

Այժմ ապացուցենք, որ եթե $G$ գրաֆի կողերի բազմությունը կարելի է տրոհել կողերով չհատվող պարզ ցիկլերի, ապա $G$-ն զույգ գրաֆ է: Իրոք, եթե դիտարկենք կամայական $v \in V(G)$ գագաթ և ենթադրենք, որ այդ գագաթը մասնակցում է $G$ գրաֆի կողերի բազմության տրոհման մեջ $C_{i_1}, \ldots, C_{i_k}$ կողերով չհատվող պարզ ցիկլերում, ապա պարզ է, որ $d_G(v) = 2k$: Այստեղից հետևում է, որ $G$-ն զույգ գրաֆ է: ∎

Այժմ կապացուցենք Էյլերի թեորեմը, որը նկարագրում է էյլերյան գրաֆները:

**Թեորեմ 4.1.2 (Լ. Էյլեր):** Որպեսզի $G$ գրաֆը լինի էյլերյան, անհրաժեշտ է և բավարար, որ այն լինի կապակցված և զույգ գրաֆ:

**Ապացույց:** Ինչպես նշել ենք, եթե $G$ գրաֆը էյլերյան է, ապա այն կապակցված է: Ցույց տանք, որ կամայական $v \in V(G)$ գագաթի համար $d_G(v)$-ն զույգ թիվ է: Ենթադրենք, որ $C = v_1, \ldots, v_s, v_1$-ն էյլերյան ցիկլ է և $v$ գագաթը մասնակցում է այդ ցիկլում $k$ անգամ: Այդ դեպքում հեշտ է տեսնել, որ եթե $v \neq v_1$, ապա $d_G(v) = 2k$, հակառակ դեպքում՝ $d_G(v) = 2(k - 1)$ (քանի որ յուրաքանչյուր անգամ $v$ գագաթը այցելելու ժամանակ $C$ ցիկլը մտնում է այդ գագաթ մի կողով և դուրս է գալիս մեկ այլ կողով): Այստեղից հետևում է, որ $G$ գրաֆի յուրաքանչյուր գագաթի աստիճանը զույգ թիվ է:

Ցույց տանք, որ եթե $G$-ն կապակցված և զույգ գրաֆ է, ապա $G$-ն էյլերյան է: Դիտարկենք $G$ գրաֆի առավելագույն երկարություն ունեցող որևէ ճանապարհ. $P =$

$v_1, ..., v_s$: Նախ ապացուցենք, որ $P$-ն ցիկլ է: Ենթադրենք հակառակը՝ $v_1 \neq v_s$: Քանի որ $v_1$ գագաթին կից կենտ թվով կողեր են պատկանում $P$ ճանապարհին, ուստի գոյություն կունենա $uv_1 \in E(G)\setminus E(P)$: Դիտարկենք $G$ գրաֆի $P' = u, v_1, ..., v_s$ ճանապարհը: Պարզ է, որ $|P'| > |P|$, ինչը հակասում է $P$ ճանապարհի ընտրությանը: Հետևաբար, $v_1 = v_s$ և $P = v_1, ..., v_{s-1}, v_1$-ը ցիկլ է:

Թեորեմի ապացույցն ավարտելու համար բավական է ցույց տալ, որ $P = v_1, ..., v_{s-1}, v_1$-ն էյլերյան ցիկլ է: Ենթադրենք հակառակը՝ $P$-ն էյլերյան չէ: Քանի որ $P$-ն էյլերյան չէ, ուստի գոյություն կունենա $uw \in E(G)\setminus E(P)$: Մյուս կողմից, քանի որ $G$-ն կապակցված գրաֆ է, ուստի գոյություն կունենան ճանապարհներ, որոնք $u$ կամ $w$ գագաթը միացնում են $P$-ի գագաթների հետ: Որոշակիության համար դիտարկենք նրանցից կարճագույնը: Դիցուք այդ ճանապարհը $P'' = u_1, ..., u_t$-ն է, որտեղ $u_1 = u$ և $u_t = v_l$ ($t \geq 1$): Այդ դեպքում դիտարկենք $G$ գրաֆի $P''' = w, u = u_1, ..., u_t = v_l, v_{l+1}, ..., v_{s-1}, v_1, ..., v_l$ ճանապարհը: Պարզ է, որ $|P'''| > |P|$, ինչը հակասում է $P$ ճանապարհի ընտրությանը: Հետևաբար, $P$-ն էյլերյան ցիկլ է: ∎

Թեորեմ 4.1.1 և 4.1.2-ից ստանում ենք էյլերյան գրաֆների մեկ այլ նկարագրում:

**Թեորեմ 4.1.3:** $G$ գրաֆը էյլերյան է այն և միայն այն դեպքում, երբ $G$-ն կապակցված գրաֆ է և այդ գրաֆի կողերի բազմությունը կարելի է տրոհել կողերով չհատվող պարզ ցիկլերի:

Նշենք նաև առանց ապացույցի էյլերյան գրաֆների ևս մի նկարագրում:

**Թեորեմ 4.1.4 (Տոյդա, Մակ-Կիի):** $G$ գրաֆը էյլերյան է այն և միայն այն դեպքում, երբ $G$-ն կապակցված գրաֆ է և այդ գրաֆի յուրաքանչյուր կող պատկանում է կենտ թվով պարզ ցիկլերի:

Ստորև նկարագրվում է էյլերյան գրաֆում էյլերյան ցիկլ կառուցելու Ֆլյորիի ալգորիթմը:

**Ալգորիթմ**

Դիցուք տրված է $G$ էյլերյան գրաֆը:

**Քայլ 1:** Ընտրենք որևէ $u$ գագաթ և նրան հարևան $v$ գագաթ: Վերագրենք $uv$ կողին $1$ համար, այնուհետև հեռացնենք այդ կողը գրաֆից և անցնենք $v$ գագաթին:

**Քայլ 2:** Դիցուք $w$-ն այն գագաթն է, որում մենք գտնվում ենք նախորդ քայլը կատարելուց հետո և այդ քայլում որոշ կողի վերագրվել է $k$ համար: Ընտրենք $w$-ին կից ցանկացած կող, ընդ որում այն կողը, որը հանդիսանում է կամուրջ ստացված գրաֆում, ընտրում ենք այն դեպքում, եթե այլ ընտրության հնարավորություն չկա: Վերագրում ենք ընտրված կողին $k + 1$ համար, այնուհետև հեռացնում ենք այդ կողը գրաֆից և անցնում ենք այդ կողի մյուս գագաթին:

**Քայլ 3:** Կատարել քայլ 2-ը այնքան անգամ մինչև գրաֆում կող չմնա:

Համոզվենք, որ Ֆլյորիի ալգորիթմը իրոք կառուցում է էյլերյան ցիկլ: Նախ նկատենք, քանի որ $G$ գրաֆի յուրաքանչյուր գագաթի աստիճանը զույգ թիվ է, ուստի ալգորիթմը կարող է վերջացնի իր աշխատանքը միայն այն գագաթում, որից սկսել էր իր աշխատանքը: Այստեղից հետևում է, որ ալգորիթմը կառուցում է ինչ-որ $C$ ցիկլ: Մնում է ցույց տալ, որ այդ ցիկլը պարունակում է $G$ գրաֆի բոլոր կողերը, իսկ այն, որ այդ ցիկլը անցնում է $G$ գրաֆի յուրաքանչյուր կողով ճիշտ մեկ անգամ, հետևում է նրանից, որ յուրաքանչյուր կողով անցնելուց հետո ալգորիթմի համաձայն այդ կողը հեռացվում է գրաֆից: Ենթադրենք, որ $C$ ցիկլը պարունակում է $G$ գրաֆի ոչ բոլոր կողերը: Դիցուք $G'$-ը $G \backslash E(C)$ գրաֆի այն կապակցված բաղադրիչն է, որի համար $E(G') \neq \emptyset$: Դիտարկենք $C$ ցիկլի կողերի այն $A$ բազմությունը, որոնք կից են $G'$ գրաֆի գագաթներին: Նկատենք, որ $A \neq \emptyset$: Դիցուք $a$-ն $A$-ի այն կողն է, որը ալգորիթմի աշխատանքի ընթացքում ստացել է ամենամեծ համարը: Հեշտ է տեսնել, որ հեռացման պահին այդ $a$ կողը եղել է կամուրջ ստացված գրաֆում, իսկ դա հակասում է հերթական կողի ընտրությանը:

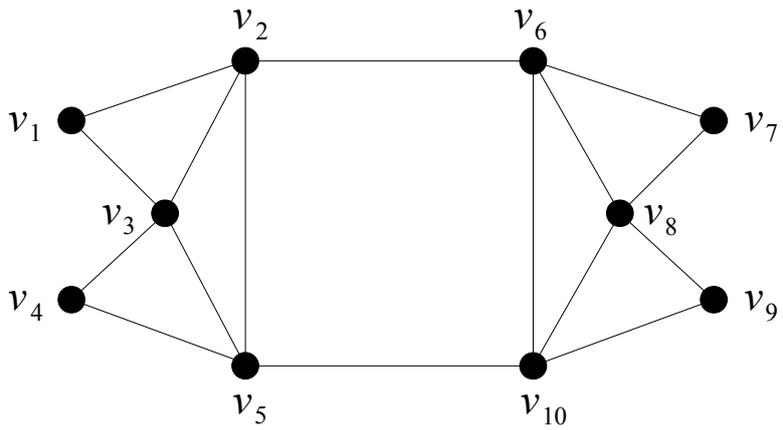

Նկ. 4.1.3

Բերենք ալգորիթմի աշխատանքը պարզաբանող մի օրինակ: Դիտարկենք նկ. 4.1.3-ում պատկերված $G$ էյլերյան գրաֆը: Նկ. 4.1.4-ում պատկերված է $G$ էյլերյան գրաֆի $C = v_1, v_2, v_6, v_7, v_8, v_6, v_{10}, v_8, v_9, v_{10}, v_5, v_2, v_3, v_5, v_4, v_3, v_1$ էյլերյան ցիկլը, որը ստացվել է Ֆլյորիի ալգորիթմի աշխատանքի արդյունքում: Նկատենք, որ $C$ ցիկլը կառուցելու ժամանակ անցնելով $v_1, v_2, v_6, v_7, v_8, v_6, v_{10}$ ճանապարհը մենք չենք կարող ընտրել $v_{10}v_5$ կողը, քանի որ այդ կողը հանդիսանում է կամուրջ ստացված գրաֆում, հետևաբար կարող ենք ընտրել $v_{10}v_8$ կամ $v_{10}v_9$ կողը ($C$-ում ընտրվել է $v_{10}v_8$-ը):

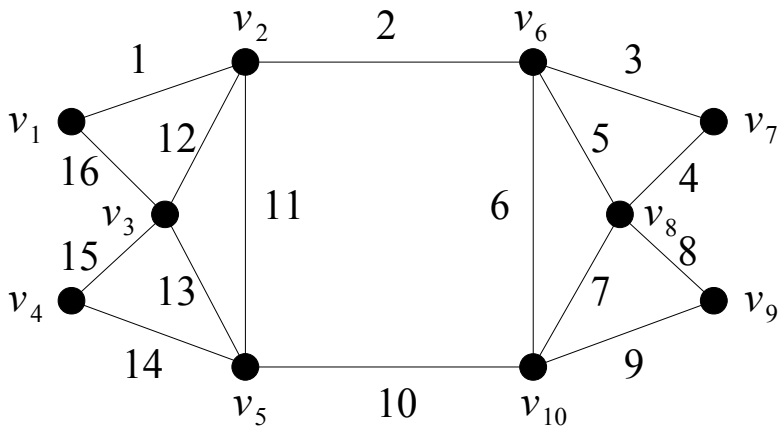

Նկ. 4.1.4

Պարզվում է, որ համարյա բոլոր գրաֆները էյլերյան չեն: Նախ տանք որոշ անհրաժեշտ սահմանումներ:

Դիցուք $\mathfrak{P}$-ն գրաֆների բազմության վրա որոշված հատկություն է: Օրինակ, որպես հատկություն կարող է հանդես գալ «կապակցված է» հատկությունը: Նշանակենք $G(n)$-ով $n$ գագաթ պարունակող բոլոր գրաֆների բազմությունը, իսկ $G_{\mathfrak{P}}(n)$-ով՝ $n$ գագաթ պարունակող բոլոր այն գրաֆների բազմությունը, որոնք օժտված են $\mathfrak{P}$ հատկությամբ:

**Սահմանում 4.1.4:** Կասենք, որ համարյա բոլոր գրաֆները օժտված են $\mathfrak{P}$ հատկությամբ, եթե

$$\lim_{n \to \infty} \frac{|G_{\mathfrak{P}}(n)|}{|G(n)|} = 1:$$

**Սահմանում 4.1.5:** Կասենք, որ համարյա բոլոր գրաֆները օժտված չեն $\mathfrak{P}$ հատկությամբ, եթե

$$\lim_{n \to \infty} \frac{|G_{\mathfrak{P}}(n)|}{|G(n)|} = 0:$$

**Թեորեմ 4.1.5 (Ռեյդ):** Համարյա բոլոր գրաֆները էյլերյան չեն:

**Ապացույց:** Դիտարկենք բոլոր գրաֆների բազմության վրա որոշված $\mathfrak{P}$ և $\mathfrak{P}'$ հատկությունները, որտեղ $\mathfrak{P}$-ն «էյլերյան է», իսկ $\mathfrak{P}'$-ը՝ «զույգ է» հատկություններն են: Համաձայն հետևանք 4.1.1-ի կստանանք՝ $G_{\mathfrak{P}}(n) \subseteq G_{\mathfrak{P}'}(n)$: Մյուս կողմից, ըստ թեորեմ 1.2.2-ի և 1.2.3-ի կստանանք՝ $|G(n)| = 2^{\binom{n}{2}}$ և $|G_{\mathfrak{P}'}(n)| = 2^{\binom{n-1}{2}}$: Այստեղից հետևում է, որ

$$|G_\mathfrak{P}(n)| \leq 2^{\binom{n-1}{2}} = 2^{\binom{n}{2}-n+1} = |G(n)| \cdot 2^{-n+1} \text{ և } \frac{|G_\mathfrak{P}(n)|}{|G(n)|} \leq \frac{1}{2^{n-1}},$$

հետևաբար

$$\lim_{n\to\infty} \frac{|G_\mathfrak{P}(n)|}{|G(n)|} = 0. \quad \blacksquare$$

Այժմ ներկայացնենք գրաֆում Էյլերյան ճանապարհի գոյության անհրաժեշտ և բավարար պայմանը։

**Թեորեմ 4.1.6:** Կապակցված $G$ գրաֆում գոյություն ունի Էյլերյան ճանապարհ այն և միայն այն դեպքում, երբ $G$ գրաֆի կենտ աստիճան ունեցող գագաթների քանակը երկուսից ավելի չէ։

**Ապացույց:** Նախ նկատենք, որ եթե $G$ կապակցված գրաֆում գոյություն ունի Էյլերյան ճանապարհի, ապա, ինչպես նշել ենք թեորեմ 4.1.2-ի ապացուցում, այդ ճանապարհի բոլոր գագաթները ունեն զույգ աստիճան, բացի, գուցե, ծայրակետերից։ Հետևաբար, $G$ գրաֆի կենտ աստիճան ունեցող գագաթների քանակը երկուսից ավելի չէ։

Ենթադրենք $G$ կապակցված գրաֆի կենտ աստիճան ունեցող գագաթների քանակը երկուսից ավել չէ։ Ցույց տանք, որ $G$ գրաֆը ունի Էյլերյան ճանապարհի։

Եթե $G$ կապակցված գրաֆը չունի կենտ աստիճան ունեցող գագաթներ, ապա, ըստ թեորեմ 4.1.2-ի, $G$ գրաֆը ունի Էյլերյան ցիկլ, որը նաև Էյլերյան ճանապարհ է։ Ըստ հետևանք 1.2.1-ի, $G$ գրաֆը չի կարող ունենալ կենտ աստիճան ունեցող միակ գագաթ։ Հետևաբար, թեորեմը ապացուցելու համար բավական է դիտարկել այն դեպքը, երբ $G$ գրաֆը պարունակում է կենտ աստիճան ունեցող երկու գագաթ։ Դիցուք այդ գագաթները $u$-ն և $v$-ն են։ Սահմանենք $G'$ գրաֆը հետևյալ կերպ.

$$V(G') = V(G) \cup \{w\}, \text{ որտեղ } w \notin V(G),$$

$$E(G') = E(G) \cup \{uw, wv\}.$$

Նկատենք, որ $G'$-ը կապակցված գրաֆ է, որի յուրաքանչյուր գագաթի աստիճանը զույգ թիվ է։ Հետևաբար, ըստ թեորեմ 4.1.2-ի $G'$ գրաֆը կպարունակի Էյլերյան ցիկլ։ Այդ ցիկլից դեն նետելով $w$ գագաթը կստանանք Էյլերյան $(u,v)$-ճանապարհի $G$ գրաֆում։ $\blacksquare$

Էյլերյան գրաֆներին ավելի մանրամասն կարելի է ծանոթանալ Ֆլեշների գրքում [14]:

## § 4.2. Համիլտոնյան ճանապարհներ և ցիկլեր

Դիցուք $G = (V, E)$-ն գրաֆ է։

**Սահմանում 4.2.1։** Կասենք, որ $G$ գրաֆում գոյություն ունի *համիլտոնյան ճանապարհ*, եթե $G$-ում գոյություն ունի պարզ ճանապարհ, որն անցնում է $G$-ի բոլոր գագաթներով։ Եթե այդ պարզ ճանապարհը պարզ ցիկլ է, ապա այն կանվանենք *համիլտոնյան ցիկլ*։

**Սահմանում 4.2.2։** $G$ գրաֆը կանվանենք *համիլտոնյան* գրաֆ, եթե այն պարունակում է համիլտոնյան ցիկլ։

Այստեղ նույնպես բնական է դիտարկել հետևյալ խնդիրը. ինչպիսի պայմանների պետք է բավարարի գրաֆը, որպեսզի այն ունենա համիլտոնյան ճանապարհ կամ ցիկլ (լինի համիլտոնյան)։ Պարզ է, որ եթե $G$ գրաֆը պարունակում է համիլտոնյան ճանապարհ կամ ցիկլ, ապա $G$-ն պետք է լինի կապակցված։ Սակայն ոչ բոլոր կապակցված գրաֆներն ունեն համիլտոնյան ճանապարհ կամ համիլտոնյան ցիկլ։ Այսպես, օրինակ, նկ. 4.2.1-ի աջ մասում պատկերված կապակցված գրաֆը ունի համիլտոնյան ճանապարհ (համոզվելու համար բավական է դիտարկել $w_1, u_1, w_2, u_2, w_3, u_3, w_4$ ճանապարհը), բայց չունի համիլտոնյան ցիկլ, իսկ ձախ մասում պատկերվածը՝ չունի համիլտոնյան ճանապարհի և, հետևաբար, նաև չունի համիլտոնյան ցիկլ։

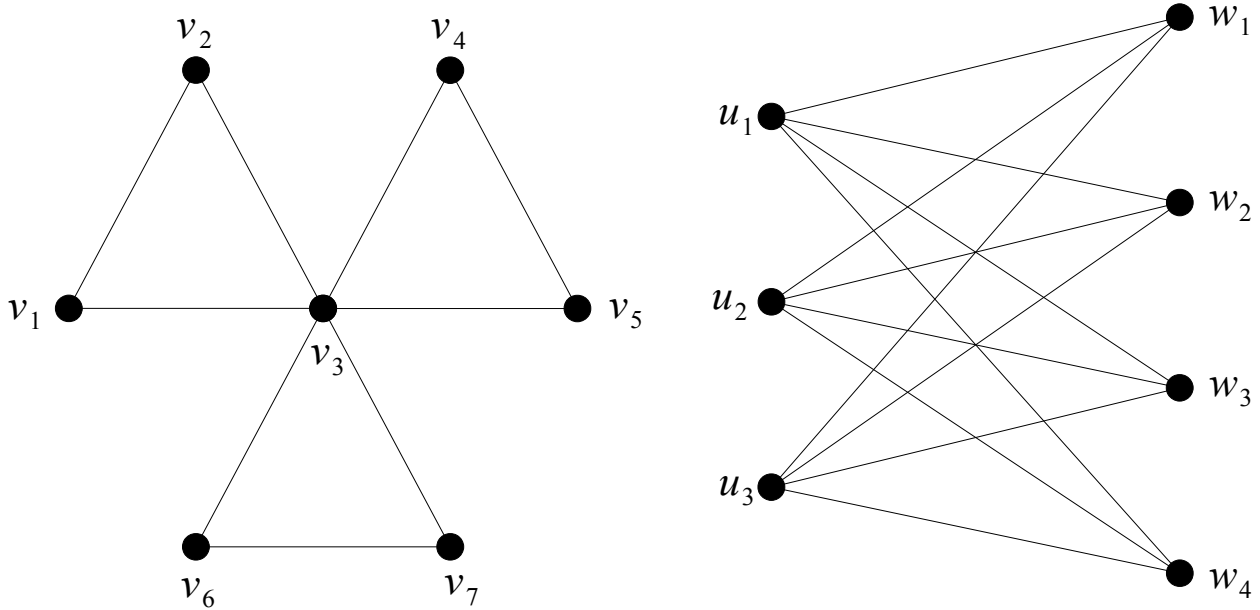

Նկ. 4.2.1

Նկ. 4.2.1-ի ձախ մասում պատկերված կապակցված գրաֆը հուշում է, որ կարելի է նույնիսկ պնդել, որ եթե $G$ գրաֆը պարունակում է համիլտոնյան ցիկլ (համիլտոնյան է), ապա $G$-ն 2-կապակցված գրաֆ է: Իրոք, եթե $G$ գրաֆը համիլտոնյան է, ապա ցանկացած $v \in V(G)$-ի համար $G - v$ գրաֆը պարունակում է համիլտոնյան ճանապարհի, ուստի $v$ գագաթը չի կարող լինել միակցման կետ: Սակայն, պետք է նշել, որ ոչ բոլոր 2-կապակցված գրաֆները համիլտոնյան են (օրինակ է Պետերսենի գրաֆը): Այժմ բերենք համիլտոնյան գրաֆների օրինակներ: Նկ. 4.2.2-ում պատկերված գրաֆները համիլտոնյան են և կետագծերով նշված են այդ գրաֆների համիլտոնյան ցիկլերը:

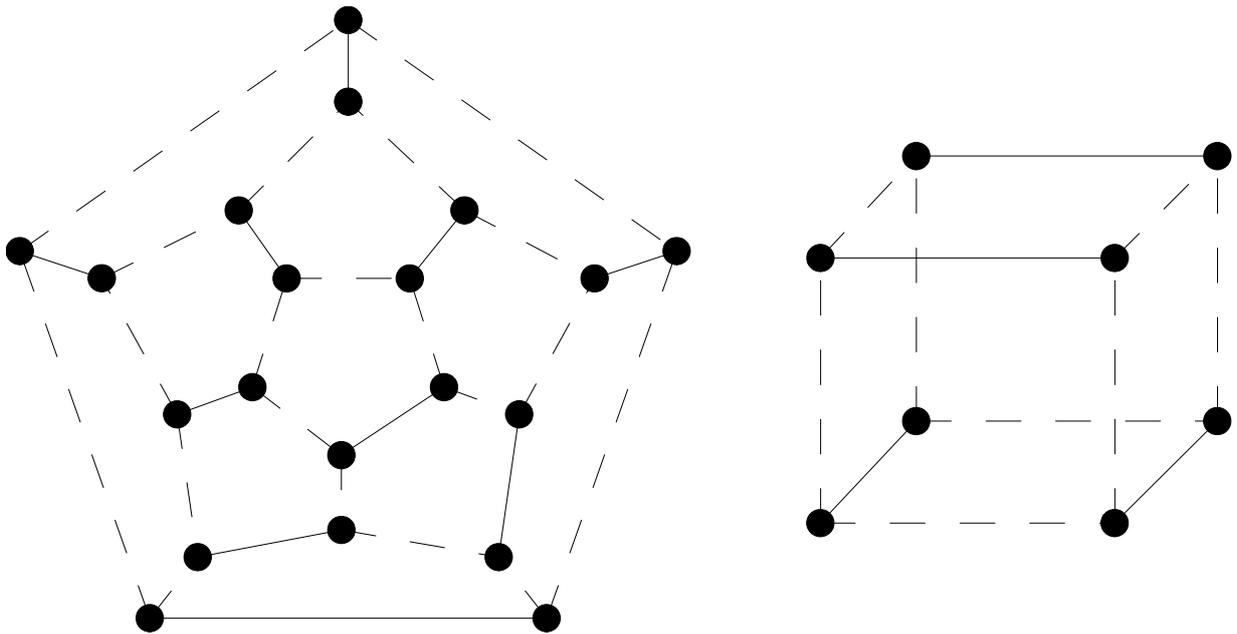

Նկ. 4.2.2

Գրաֆում համիլտոնյան ցիկլի գոյության անհրաժեշտ և բավարար պայմաններ հայտնի չեն, այդ պատճառով ստորև մենք կդիտարկենք համիլտոնյան ցիկլի գոյության որոշ անհրաժեշտ և որոշ բավարար պայմաններ: Նախ ապացուցենք Դիրակի լեմման, որը հետագայում կօգտագործենք համիլտոնյան ցիկլի գոյության որոշ բավարար պայմաններ ապացուցելու ժամանակ:

**Լեմմա 4.2.1:** Եթե $G$ գրաֆի համար տեղի ունի $\delta(G) \geq 2$ պայմանը, ապա $G$-ն պարունակում է առնվազն $\delta(G) + 1$ երկարություն ունեցող պարզ ցիկլ:

**Ապացույց:** Դիտարկենք $G$ գրաֆի առավելագույն երկարություն ունեցող որևէ պարզ ճանապարհի՝ $P = v_1, \ldots, v_s$: Քանի որ $P$-ի երկարությունը առավելագույնն է, ուստի այդ պարզ ճանապարհի $v_1$ և $v_s$ գագաթների բոլոր հարևանները պատկանում են $P$-ին (հակառակ դեպքում հնարավոր կլիներ նշել $P$-ից ավելի երկար պարզ ճանապարհի):

Մասնավորապես, դա նշանակում է, որ գագաթների $\{v_{\delta(G)+1}, v_{\delta(G)+2}, \ldots, v_s\}$ բազմության մեջ գոյություն ունի $v_{k_0}$ գագաթ, որը հարևան է $v_1$-ին։ Այստեղից հետևում է, որ $G$ գրաֆը պարունակում է առնվազն $\delta(G)+1$ երկարություն ունեցող $C = v_1, \ldots, v_{k_0}, v_1$ պարզ ցիկլ։ ∎

Ստորև կապացուցենք Դիրակի թեորեմը, որը համիլտոնյան ցիկլի գոյության պատմականորեն առաջին հայտնի բավարար պայմանն է։

**Թեորեմ 4.2.1 (Գ. Դիրակ):** Եթե $n$ գագաթ ($n \geq 3$) ունեցող $G$ գրաֆի համար տեղի ունի $\delta(G) \geq \frac{n}{2}$ պայմանը, ապա այն համիլտոնյան է։

**Ապացույց:** Նախ նկատենք, որ $n \geq 3$ պայմանը էական է այս թեորեմում։ Իրոք, լրիվ $K_2$ գրաֆը բավարարում է թեորեմի պայմանին, բայց այն համիլտոնյան գրաֆ չէ։

Ենթադրենք հակառակը, գոյություն ունի $n$ գագաթ ($n \geq 3$) ունեցող $G$ գրաֆ, որը բավարարում է $\delta(G) \geq \frac{n}{2}$ պայմանին, բայց $G$-ն համիլտոնյան գրաֆ չէ։ Նկատենք, որ եթե $u$-ն և $v$-ն $G$ գրաֆում կամայական երկու ոչ հարևան գագաթներ են, ապա $G + uv$ գրաֆը ևս կբավարարի թեորեմի պայմանին, այսինքն՝ $G + uv$ գրաֆը $n$ գագաթ ($n \geq 3$) ունեցող գրաֆ է, որը բավարարում է $\delta(G + uv) \geq \frac{n}{2}$ պայմանին։ Այստեղից հետևում է, որ $G$ հակաօրինակին ավելացնելով կողեր մենք կարող ենք ստանալ ապացուցվող թեորեմի մաքսիմալ հակաօրինակ, այսինքն՝ այնպիսի հակաօրինակ, որը բավարարում է թեորեմի պայմանին, որը համիլտոնյան գրաֆ չէ, բայց ցանկացած նոր կող ավելացնելուց ստացվող գրաֆը կլինի համիլտոնյան գրաֆ։

Ստորև կապացուցենք, որ գոյություն չունեն թեորեմի մաքսիմալ հակաօրինակներ։ Քանի որ ցանկացած հակաօրինակից կարելի է ստանալ մաքսիմալ հակաօրինակ, ապա այս պնդումից կստացվի, որ գոյություն չունեն նաև հակաօրինակներ, և, հետևաբար, թեորեմի պնդումը ճիշտ է։

Դիցուք $G$-ն մաքսիմալ հակաօրինակ է և $u$-ն և $v$-ն այդ գրաֆի կամայական երկու ոչ հարևան գագաթներ են։ Քանի որ $G$-ն մաքսիմալ հակաօրինակ է, ուստի $G + uv$ գրաֆը կլինի համիլտոնյան գրաֆ։ Այստեղից հետևում է, որ $G$-ն պարունակում է համիլտոնյան $(u, v)$-ճանապարհի։ Դիցուք այդ ճանապարհը $P = v_1, \ldots, v_n$-ն է, որտեղ $n \geq 3$ և $v_1 = u$, $v_n = v$։

Ցույց տանք, որ գոյություն ունի այնպիսի $i$ ինդեքս, որ $uv_{i+1} \in E(G)$ և $vv_i \in E(G)$։ Սահմանենք $S$ և $T$ բազմությունները հետևյալ կերպ.

$$S = \{i : uv_{i+1} \in E(G)\} \text{ և } T = \{i : vv_i \in E(G)\}։$$

Համոզվենք, որ $|S \cap T| \geq 1$: Քանի որ $|S \cup T| = |S| + |T| - |S \cap T|$ և $\delta(G) \geq \frac{n}{2}$, ուստի

$$|S \cup T| + |S \cap T| = |S| + |T| = d_G(u) + d_G(v) \geq n:$$

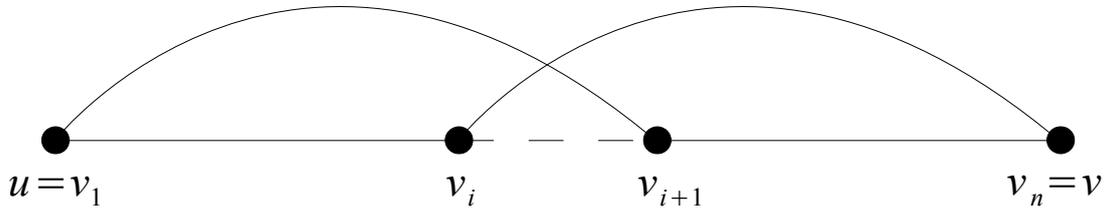

Նկ. 4.2.3

Մյուս կողմից, պարզ է որ $n \notin S \cup T$ և, հետևաբար, $|S \cup T| < n$: Այստեղից հետևում է, որ $|S \cap T| \geq 1$: Այժմ դիտարկենք $G$ գրաֆի $C = u, v_2, \ldots, v_{i-1}, v_i, v, v_{n-1}, \ldots, v_{i+2}, v_{i+1}, u$ պարզ ցիկլը (նկ. 4.2.3): Հեշտ է տեսնել, որ $C$-ն հանդիսանում է $G$ գրաֆի համիլտոնյան ցիկլ, ուստի $G$-ն համիլտոնյան գրաֆ է: ∎

Նկատենք, որ Դիրակի թեորեմը ապացուցելու ժամանակ $\delta(G) \geq \frac{n}{2}$ պայմանից մենք օգտվեցինք միայն $d_G(u) + d_G(v) \geq n$ ցույց տալու համար, երբ $uv \notin E(G)$: Այդ փաստը առաջին անգամ նշվել էր Օրէի կողմից. թեորեմ 4.2.1-ի ապացույցի դատողությունները կրկնելով, կարելի է ապացուցել հետևյալ թեորեմը.

**Թեորեմ 4.2.2 (О. Օրէ):** Եթե $n$ գագաթ ($n \geq 3$) ունեցող $G$ գրաֆում ցանկացած $u$ և $v$ ($u \neq v$) ոչ հարևան գագաթների համար տեղի ունի $d_G(u) + d_G(v) \geq n$ պայմանը, ապա $G$-ն համիլտոնյան գրաֆ է:

Նկատենք, որ թեորեմ 4.2.1-ում նշված՝ $\delta(G) \geq \frac{n}{2}$ և թեորեմ 4.2.2-ում նշված՝ ցանկացած $u$ և $v$ ($u \neq v$) ոչ հարևան գագաթների համար տեղի ունի $d_G(u) + d_G(v) \geq n$ պայմանները հնարավոր չէ լավացնել: Իրոք, դիտարկենք այն $G$ գրաֆը, որը ստացվում է $K_{\left\lfloor \frac{n+1}{2} \right\rfloor}$ և $K_{\left\lceil \frac{n+1}{2} \right\rceil}$ լրիվ գրաֆների մեկ գագաթը նույնացնելով: Հեշտ է տեսնել, որ $|V(G)| = n$, $\delta(G) = \left\lfloor \frac{n-1}{2} \right\rfloor$ և ցանկացած $u$ և $v$ ($u \neq v$) ոչ հարևան գագաթների համար տեղի ունի $d_G(u) + d_G(v) \geq n - 1$ պայմանը: Սակայն, քանի որ $G$-ն 2-կապակցված գրաֆ չէ, այն համիլտոնյան չէ:

Դիրակի և Օրէի թեորեմները հանդիսանում են համիլտոնյան գրաֆների տեսության դասական արդյունքներ, որոնք տալիս են գրաֆներում համիլտոնյան ցիկլի գոյության համար աստիճանային սահմանափակումներով բավարար պայմաններ: Համիլտոնյան գրաֆների տեսության մեկ այլ ուղղություն են ներկայացնում գրաֆներում համիլտոնյան ցիկլի գոյության այն բավարար պայմանները, որոնց դեպքում սահմանափակումները դրվում են գրաֆի ենթագրաֆների վրա։ Այս ուղղությունը համիլտոնյան գրաֆների տեսության մեջ հայտնի է որպես արգելված ենթագրաֆների ուղղություն։ Ստորև կձևակերպենք և կապացուցենք այս ուղղության դասական արդյունքներից մեկը՝ Գուդմանի և Հեդեթնիեմի թեորեմը։

**Թեորեմ 4.2.3 (Գուդման, Հեդեթնիեմի):** Եթե **2**-կապակցված $G$ գրաֆը չի պարունակում $K_{1,3}$ և $K_{1,3} + e$ ծնված ենթագրաֆներ, ապա $G$-ն համիլտոնյան գրաֆ է։

**Ապացույց:** Ենթադրենք հակառակը, գոյություն ունի **2**-կապակցված $G$ գրաֆ, որը չի պարունակում $K_{1,3}$ և $K_{1,3} + e$ ծնված ենթագրաֆներ, բայց $G$-ն համիլտոնյան գրաֆ չէ։ Քանի որ $G$-ն **2**-կապակցված գրաֆ է, ուստի համաձայն թեորեմ 3.2.1-ի $\delta(G) \geq 2$ և, հետևաբար, ըստ լեմմա 4.1.1-ի, $G$-ն կպարունակի պարզ ցիկլ։ Դիտարկենք $G$ գրաֆի ամենաերկար $C = v_1, \ldots, v_s, v_1$ պարզ ցիկլը։ Քանի որ $s < |V(G)|$ և $G$-ն **2**-կապակցված գրաֆ է, ուստի $G$ գրաֆում գոյություն ունի այնպիսի $v$ գագաթ, որ $v \in V(G) \setminus V(C)$ և $vv_i \in E(G)$ ($v_i \in V(C)$) (նկ. 4.2.4):

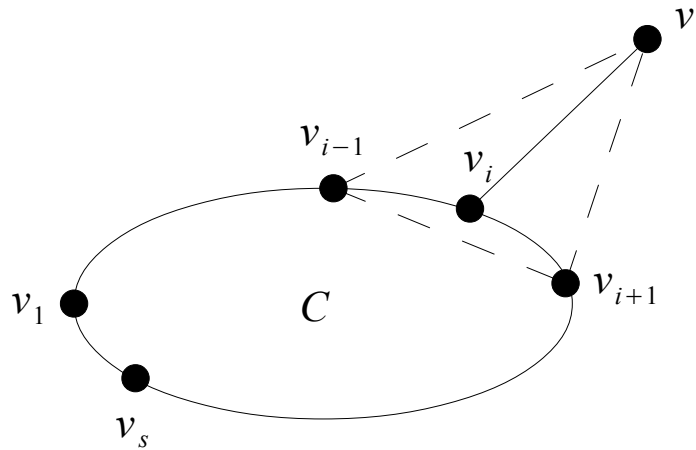

Նկ. 4.2.4

Դիցուք $S = \{v, v_{i-1}, v_i, v_{i+1}\}$։ Դիտարկենք $G$ գրաֆի $G[S]$ ծնված ենթագրաֆը։ Պարզ է, որ $G$ գրաֆում առկա են $vv_{i-1}, vv_{i+1}, v_{i-1}v_{i+1}$ կողերից առնվազն երկուսը (նկ. 4.2.4)։ Հետևաբար, $G$ գրաֆում առկա է $vv_{i-1}, vv_{i+1}$ կողերից մեկը, ինչը հակասում է $C$-ի ընտրությանը։ Այս հակասությունը ապացուցում է թեորեմը։ ∎

Համիլտոնյան գրաֆների տեսության մեջ առկա են նաև գրաֆներում համիլտոնյան ցիկլի գոյության այնպիսի բավարար պայմաններ, որոնց դեպքում սահմանափակումները դրվում են ոչ թե գագաթների աստիճանների կամ գրաֆի ենթագրաֆների վրա, այլ գրաֆի ուրիշ բնութագրիչների վրա։ Այդպիսին է, մասնավորապես, Խվատալի և Էրդյոշի թեորեմը: Մինչ այդ թեորեմի ձևակերպմանն անցնելը տանք մեկ անհրաժեշտ սահմանում: Դիցուք $G = (V, E)$-ն գրաֆ է և $I \subseteq V(G)$: Կասենք, որ $I$-ն հանդիսանում է գագաթների *անկախ բազմություն* $G$ գրաֆում, եթե $I$-ն չի պարունակում հարևան գագաթներ: Ամենաշատ գագաթներ պարունակող անկախ բազմություն հզորությունը նշանակենք $\alpha(G)$-ով:

**Թեորեմ 4.2.4 (Խվատալ, Էրդյոշ):** Եթե առնվազն երեք գագաթ պարունակող $G$ գրաֆում տեղի ունի $\alpha(G) \leq \varkappa(G)$ պայմանը, ապա $G$-ն համիլտոնյան գրաֆ է:

**Ապացույց:** Ենթադրենք հակառակը, գոյություն ունի առնվազն երեք գագաթ պարունակող $G$ գրաֆ, որի համար $\alpha(G) \leq \varkappa(G)$, բայց $G$-ն համիլտոնյան գրաֆ չէ:

Նախ նկատենք, որ եթե առնվազն երեք գագաթ պարունակող $G$ գրաֆում $\alpha(G) = 1$, ապա $G$-ն լրիվ գրաֆ է, որը համիլտոնյան գրաֆ է: Հետևաբար, մենք կարող ենք ենթադրել, որ $\alpha(G) \geq 2$: Քանի որ $\varkappa(G) \geq \alpha(G) \geq 2$, ուստի, համաձայն թեորեմ 3.2.1-ի, $\delta(G) \geq 2$: Ըստ լեմմա 4.2.1-ի, $G$-ն պարունակում է առնվազն $\delta(G) + 1$ երկարություն ունեցող պարզ ցիկլ: Դիտարկենք $G$ գրաֆի ամենաերկար $C$ պարզ ցիկլը: Պարզ է, որ $|C| \geq \delta(G) + 1$: Մյուս կողմից, համաձայն թեորեմ 3.2.1-ի, $\delta(G) \geq \varkappa(G)$ և, հետևաբար, $|C| \geq \varkappa(G) + 1$: Քանի որ $|C| < |V(G)|$, ուստի $G$ գրաֆում գոյություն ունի այնպիսի $v$ գագաթ, որ $v \in V(G) \setminus V(C)$: Դիցուք $\varkappa(G) = k$: Քանի որ $G$-ն $k$-կապակցված գրաֆ է, ուստի, համաձայն թեորեմ 3.4.5-ի, $G$-ն ունի $(v, V(C))$-հովհար: Այստեղից հետևում է, որ գոյություն ունեն առնվազն $k$ հատ $P_1, \dots, P_k$ պարզ ճանապարհներ, որոնք միացնում են $v$ գագաթը $C$ պարզ ցիկլին և այդ ճանապարհները հատվում են միայն $v$ գագաթում: Ենթադրենք $P_1, \dots, P_k$ ճանապարհները $C$ ցիկլի հետ հատվում են $v_1, \dots, v_k$ գագաթներում: Առանց ընդհանրությունը խախտելու կարող ենք ենթադրել, որ $C$ ցիկլի վրայով շարժվելիս $v_1, \dots, v_k$ գագաթները հերթականորեն հանդիպում են հենց այս հաջորդականությամբ: Ընտրենք $C$ ցիկլի վրա $u_1, \dots, u_k$ գագաթներ այնպես, որ յուրաքանչյուր $i$-ի համար, որտեղ $1 \leq i \leq k$, $u_i$ գագաթը $C$ ցիկլի վրա անմիջապես հաջորդում է $v_i$ գագաթին: Նկատենք, որ ցանկացած $i$-ի համար, որտեղ $1 \leq i \leq k$, $u_i v \notin E(G)$: Իրոք, հակառակ դեպքում գոյություն ունի $i_0$, որ $u_{i_0} v \in E(G)$: Այդ դեպքում $C$ ցիկլի $v_{i_0} u_{i_0}$ կողը փոխարինելով $P = P_{i_0}, u_{i_0}$ ճանապարհով մենք կստանանք $C$-ից ավելի երկար պարզ ցիկլ, ինչը հակասում է $C$-ի ընտրությանը: Ցույց տանք, որ կամայական $i, j$ զույգի համար, որտեղ $1 \leq i < j \leq k$, $u_i u_j \notin E(G)$: Ենթադրենք հակառակը, գոյություն ունի $i_0, j_0$ ($i_0 < j_0$), որ $u_{i_0}u_{j_0} \in E(G)$: Այդ դեպքում $C$ ցիկլից հեռացնելով $v_{i_0}u_{i_0}$ և $v_{j_0}u_{j_0}$ կողերը և ավելացնելով $P_{i_0}, P_{j_0}$ ճանապարհները և $u_{i_0}u_{j_0}$ կողը մենք կստանանք $C$-ից ավելի երկար պարզ ցիկլ, ինչը հակասում է $C$-ի ընտրությանը (նկ. 4.2.5):

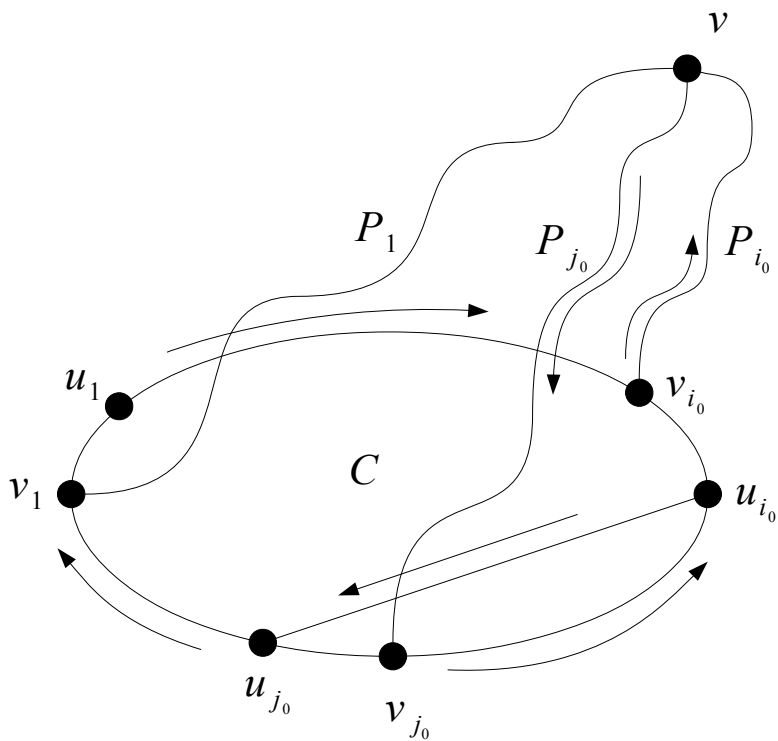

Նկ. 4.2.5

Այժմ դիտարկենք $G$ գրաֆի գագաթների $I = \{v, u_1, \ldots, u_k\}$ բազմությունը: Պարզ է, որ $I$-ն անկախ բազմություն է, ուստի $\alpha(G) \geq |I| = k+1$, ինչը հակասում է թեորեմի պայմանին: ∎

Նկատենք, որ թեորեմ 4.2.4-ում նշված՝ $\alpha(G) \leq \varkappa(G)$ պայմանը հնարավոր չէ լավացնել: Իրոք, դիտարկենք լրիվ երկկողմանի $K_{n,n+1}$ գրաֆը: Հեշտ է տեսնել, որ $\alpha(K_{n,n+1}) = n+1$ և $\varkappa(K_{n,n+1}) = n$, բայց $K_{n,n+1}$-ը համիլտոնյան գրաֆ չէ:

**Սահմանում 4.2.3**: $G$ գրաֆի $cl(G)$ *(համիլտոնյան) փակումը* գրաֆ է, որը ստացվում է $G$-ից հերթականորեն կողեր ավելացնելով, որոնք միացնում են երկու ոչ հարևան գագաթներ և որոնց աստիճանների գումարը մեծ կամ հավասար է այդ գրաֆի գագաթների քանակից, ընդ որում դա կատարվում է այնքան անգամ, ինչքան դա հնարավոր է:

Նկ. 4.2.6-ում պատկերված է $G$ գրաֆից այդ գրաֆի $cl(G)$ փակումը ստանալու ամբողջ ընթացքը:

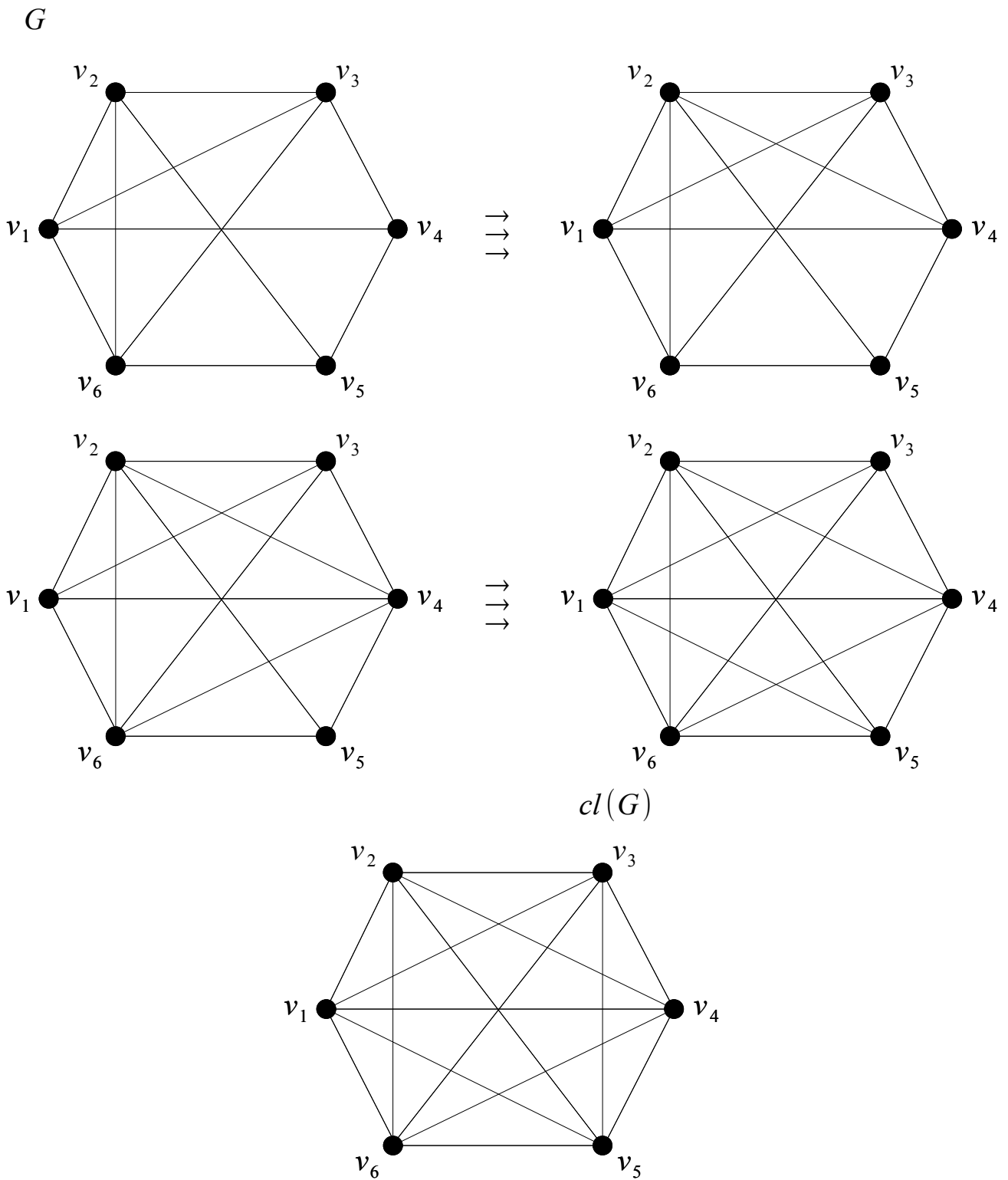

Նկ. 4.2.6

Բնական հարց է առաջանում. արդյո՞ք տրված $G$ գրաֆից այդ գրաֆի $cl(G)$ փակումը ստանալու ժամանակ կողերի ավելացման հերթականությունը կարևոր չէ: Այլ կերպ ասած, տրված $G$ գրաֆից միարժե՞ք է որոշվում նրա $cl(G)$ փակումը, թե ոչ: Համոզվենք, որ տրված $G$ գրաֆից այդ գրաֆի $cl(G)$ փակումը ստանալու ժամանակ կողերի ավելացման հերթականությունը կարևոր չէ: Ենթադրենք հակառակը, դիցուք $G_1$ գրաֆը ստացվել է $G$ գրաֆից $e_1, \ldots, e_r$ կողերը ավելացնելով, իսկ $G_2$ գրաֆը՝ $f_1, \ldots, f_s$ կողերը ավելացնելով։ Նախ նկատենք, որ եթե որևէ պահին ընթացիկ գրաֆի համար երկու ոչ հարևան $u$ և $v$ գագաթները, որոնց աստիճանների գումարը մեծ կամ հավասար է այդ գրաֆի գագաթների քանակից, հնարավոր է միացնել կողով, ապա դա պետք է կատարվի մինչև գրաֆի փակման կառուցման ավարտը։ Քանի որ սկզբնական $f_1$ կողը հնարավոր է ավելացնել $G$ գրաֆին, ուստի, ըստ վերը նշվածի, $f_1 \in E(G_1)$։ Համանման ձևով, եթե $f_1, \ldots, f_{i-1} \in E(G_1)$, ապա $f_i$-ն արդեն հնարավոր է ավելացնել $G_1$ գրաֆին, ուստի ըստ վերը նշվածի $f_i \in E(G_1)$։ Այստեղից հետևում է, որ $G_2 \subseteq G_1$։ Համանման ձևով կարելի է ցույց տալ, որ $G_1 \subseteq G_2$, ուստի $G_1 = G_2$։

Նշենք, որ գրաֆի փակման գաղափարը ներմուծվել է Բոնդիի և Խվատալի կողմից և դա պայմանավորված է հետևյալ արդյունքով։

**Թեորեմ 4.2.5 (Բոնդի, Խվատալ):** $G$ գրաֆը համիլտոնյան է այն և միայն այն դեպքում, երբ համիլտոնյան է $cl(G)$ գրաֆը։

**Ապացույց:** Թեորեմն ապացուցելու համար բավական է ցույց տալ, որ եթե $G$-ն $n$ գագաթ ունեցող գրաֆ է և $u, v$ ($u \neq v$) ոչ հարևան գագաթներ են, որոնց համար տեղի ունի $d_G(u) + d_G(v) \geq n$ պայմանը, ապա $G$-ն համիլտոնյան է այն և միայն այն դեպքում, երբ համիլտոնյան է $G + uv$ գրաֆը։ Նկատենք, որ եթե $G$-ն համիլտոնյան է, ապա ակնհայտ է, որ $G + uv$ գրաֆը նս կլինի համիլտոնյան։ Մյուս կողմի ապացույցը համընկնում է թեորեմ 4.2.1-ի ապացույցի հետ։ ∎

**Հետևանք 4.2.1:** Եթե $G$-ն $n$ գագաթ ($n \geq 3$) ունեցող գրաֆ է և $cl(G) = K_n$, ապա $G$-ն համիլտոնյան գրաֆ է։

Նկատենք, որ հետևանք 4.2.1-ից ստանում ենք, որ նկ. 4.2.6-ում պատկերված $G$ գրաֆը համիլտոնյան է, քանի որ $cl(G) = K_6$։ Օգտագործելով նշված հետևանքը, Խվատալին հաջողվեց ապացուցել հետևյալ թեորեմը։

**Թեորեմ 4.2.6 (Խվատալ):** Դիցուք $G$-ն $n$ գագաթ ($n \geq 3$) ունեցող գրաֆ է և $V(G) = \{v_1, \ldots, v_n\}$, որտեղ $d_G(v_1) \leq \cdots \leq d_G(v_n)$։ Այդ դեպքում եթե $i < \frac{n}{2}$ պայմանից հետևում է, որ $d_G(v_i) > i$ կամ $d_G(v_{n-i}) \geq n - i$, ապա $G$-ն համիլտոնյան գրաֆ է։

Ինչպես նշել ենք նախորդ պարագրաֆում, հայտնի է, որ համարյա բոլոր գրաֆները էլլերյան չեն։ Պարզվում է, որ գրաֆների համիլտոնյան լինելու հարցի վերաբերյալ իրավիճակը էապես տարբերվում է։

**Թեորեմ 4.2.7 (Պերեպելիցա):** Համարյա բոլոր գրաֆները համիլտոնյան են:

Այժմ սահմանենք գրաֆի կոշտության գաղափարը: Դիցուք $G = (V, E)$-ն գրաֆ է: $c(G)$-ով նշանակենք $G$ գրաֆի կապակցված բաղադրիչների քանակը:

**Սահմանում 4.2.4:** $G$ գրաֆը կոչվում է $t$-*կոշտ*, եթե $c(G - S) \geq 2$ պայմանին բավարարող ցանկացած $S \subseteq V(G)$-ի համար տեղի ունի $|S| \geq t \cdot c(G - S)$ անհավասարությունը: $G$ գրաֆի *կոշտություն* $\tau(G)$-ն կանվանենք այն ամենամեծ $t$-ն, որի դեպքում $G$ գրաֆը $t$-կոշտ է (համարենք, որ $\tau(K_n) = +\infty$, երբ $n \in \mathbb{N}$):

Հեշտ է տեսնել, որ եթե $m \leq n$, ապա $\tau(K_{m,n}) = \frac{m}{n}$: Նաև, դժվար չէ ցույց տալ, որ եթե $P$-ն Պետերսենի գրաֆ է, ապա $\tau(P) = \frac{4}{3}$:

Համարյա բոլոր վերը նշված թեորեմները հանդիսանում էին գրաֆներում համիլտոնյան ցիկլի գոյության բավարար պայմաններ: Ստորև մենք կձևակերպենք և կապացուցենք գրաֆներում համիլտոնյան ցիկլի գոյության մեկ անհրաժեշտ պայման:

**Թեորեմ 4.2.8:** Եթե $G$-ն համիլտոնյան գրաֆ է, ապա $\tau(G) \geq 1$:

**Ապացույց:** Դիտարկենք $c(G - S) \geq 2$ պայմանին բավարարող ցանկացած $S \subseteq V(G)$: Դիցուք $G_1, \ldots, G_{c(G-S)}$-ը $G - S$ գրաֆի կապակցված բաղադրիչներն են: Դիցուք $C$-ն $G$ գրաֆի որևէ համիլտոնյան ցիկլ է: $C$ ցիկլի վրայով շարժվելիս և դուրս գալով $G - S$ գրաֆի $G_1, \ldots, G_{c(G-S)}$ կապակցված բաղադրիչներից որևէ մեկից $C$ ցիկլը կարող է վերադառնալ միայն $S$-ի մեջ, ընդ որում տարբեր կապակցված բաղադրիչներից դուրս գալով $C$ ցիկլը այցելում է $S$-ի տարբեր գագաթներ: Այստեղից հետևում է, որ $|S| \geq c(G - S)$ և, հետևաբար, $\tau(G) \geq 1$: ∎

Գրաֆների կոշտության հետ է կապված համիլտոնյան գրաֆների տեսության հայտնի և բարդ հիպոթեզներից մեկը, որը ձևակերպել է Խվատալը:

**Հիպոթեզ 4.2.1:** Գոյություն ունի այնպիսի $t_0$ թիվ, որ ցանկացած $t_0$-կոշտ գրաֆ համիլտոնյան է:

Այս պարագրաֆի վերջում անդրադառնանք գրաֆում համիլտոնյան ճանապարհի գոյության խնդրին: Այստեղ նույնպես, ինչպես համիլտոնյան ցիկլի գոյության դեպքում, հայտնի չեն համիլտոնյան ճանապարհի գոյության անհրաժեշտ և բավարար պայմաններ: Ստորև մենք կդիտարկենք համիլտոնյան ճանապարհի գոյության որոշ բավարար պայմաններ:

**Թեորեմ 4.2.9 (Օ. Օրէ):** Եթե $n$ գագաթ ունեցող $G$ գրաֆում ցանկացած $u$ և $v$ ($u \neq v$) ոչ հարևան գագաթների համար տեղի ունի $d_G(u) + d_G(v) \geq n - 1$ պայմանը, ապա $G$ գրաֆում գոյություն ունի համիլտոնյան ճանապարհ:

**Ապացույց:** Նախ նկատենք, որ թեորեմի պայմաններին բավարարող $G$ գրաֆը կապակցված է: Իրոք, եթե մենք դիտարկենք $G$ գրաֆի ցանկացած $x$ և $y$ գագաթները, ապա հեշտ է տեսնել, որ $xy \in E(G)$ կամ $N_G(x) \cap N_G(y) \neq \emptyset$, ուստի $G$ գրաֆը կապակցված է:

Թեորեմն ապացուցելու համար կատարենք հակասող ենթադրություն՝ $G$-ն $n$ գագաթ ունեցող գրաֆ է, որում ցանկացած $u$ և $v$ ($u \neq v$) ոչ հարևան գագաթների համար տեղի ունի $d_G(u) + d_G(v) \geq n - 1$ պայմանը, բայց $G$-ն չի պարունակում համիլտոնյան ճանապարհ: Դիտարկենք $G$ գրաֆի առավելագույն երկարություն ունեցող որևէ պարզ ճանապարհ՝ $P = v_1, \ldots, v_s$: Քանի որ $P$-ի երկարությունն առավելագույնն է, ուստի այդ պարզ ճանապարհի $v_1$ և $v_s$ գագաթների բոլոր հարևանները պատկանում են $P$-ին (հակառակ դեպքում հնարավոր կլիներ նշել $P$-ից ավելի երկար պարզ ճանապարհ): Նկատենք, որ $v_1 v_s \notin E(G)$: Իրոք, եթե $v_1 v_s \in E(G)$, ապա $C = v_1, \ldots, v_s, v_1$-ն պարզ ցիկլ է և քանի որ $s < n$ և $G$-ն կապակցված գրաֆ է, ուստի $G$ գրաֆում գոյություն ունի $v$ գագաթ, որ $v \in V(G) \setminus V(C)$ և $vv_l \in E(G)$ ($v_l \in V(C)$): Այստեղից հետևում է, որ $P' = v, v_l, \ldots, v_s, v_1, \ldots, v_{l-1}$ պարզ ճանապարհը $P$-ից ավելի երկար պարզ ճանապարհ է $G$ գրաֆում, ինչը հակասում է $P$-ի ընտրությանը:

Ցույց տանք, որ գոյություն ունի այնպիսի $i$ ինդեքս, որ $v_1 v_i \in E(G)$ և $v_{i-1} v_s \in E(G)$: Իրոք, հակառակ դեպքում $P$ պարզ ճանապարհը կպարունակեր իրարից տարբեր հետևյալ գագաթները. $v_1$ գագաթը, $d_G(v_1)$ հատ գագաթ, որոնք հարևան են $v_1$ գագաթին և $d_G(v_s)$ հատ գագաթ, որոնք $v_s$ գագաթի հարևան գագաթների հաջորդներն են: Այսպիսով, քանի որ $d_G(v_1) + d_G(v_s) \geq n - 1$, ուստի $s \geq 1 + d_G(v_1) + d_G(v_s) \geq n$, որը հակասում է $s < n$ պայմանին: Այստեղից հետևում է, որ գոյություն ունի $i$ ինդեքս նշված հատկություններով: Դիտարկենք $G$ գրաֆի $C' = v_1, v_2, \ldots, v_{i-1}, v_s, v_{s-1}, \ldots, v_{i+1}, v_i, v_1$ պարզ ցիկլը: Քանի որ $s < n$ և $G$-ն կապակցված գրաֆ է, ուստի $G$ գրաֆում գոյություն ունի այնպիսի $u$ գագաթ, որ $u \in V(G) \setminus V(C')$ և $uv_t \in E(G)$ ($v_t \in V(C')$): Այստեղից հետևում է, որ $P'' = u, v_t, \ldots, v_s, v_1, \ldots, v_{t-1}$ պարզ ճանապարհը $P$-ից ավելի երկար պարզ ճանապարհ է $G$ գրաֆում, ինչը հակասում է $P$-ի ընտրությանը: Ստացված հակասությունը ապացուցում է թեորեմը: ∎

Նշենք, որ թեորեմ 4.2.1-ում նշված՝ ցանկացած $u$ և $v$ ($u \neq v$) ոչ հարևան գագաթների համար տեղի ունի $d_G(u) + d_G(v) \geq n - 1$ պայմանը հնարավոր չէ լավացնել: Իրոք, դիտարկենք $G = K_n \cup K_n$ գրաֆը: Հեշտ է տեսնել, որ $|V(G)| = 2n$ և ցանկացած $u$ և $v$ ($u \neq v$) ոչ հարևան գագաթների համար տեղի ունի $d_G(u) + d_G(v) \geq 2n - 2$ պայմանը: Սակայն, քանի որ $G$-ն կապակցված գրաֆ չէ, ուստի այն համիլտոնյան ճանապարհի պարունակել չի կարող:

Ստորև մենք կձևակերպենք և կապացուցենք համիլտոնյան ճանապարհի գոյության մեկ այլ բավարար պայման:

**Թեորեմ 4.2.10 (Խվատալ, Էրդյոշ):** Եթե $G$-ն $k$-կապակցված գրաֆ է և այդ գրաֆում տեղի ունի $\alpha(G) \leq k + 1$ պայմանը, ապա $G$ գրաֆում գոյություն ունի համիլտոնյան ճանապարհ:

**Ապացույց:** Ապացույցի համար սահմանենք $G'$ գրաֆը հետևյալ կերպ. $V(G') = V(G) \cup \{v\}, v \notin V(G)$ և $E(G') = E(G) \cup \{uv : u \in V(G)\}$: Հեշտ է տեսնել, որ $G'$-ը $(k+1)$-կապակցված գրաֆ է, որը բավարարում է թեորեմ 4.2.4-ի պայմանին, ուստի $G'$ գրաֆը կպարունակի համիլտոնյան ցիկլ: Դեն նետելով այդ ցիկլից $v$ գագաթը, կստանանք $G$ գրաֆի համիլտոնյան ճանապարհի: ∎

Հարկ է նշել, որ թեորեմ 4.2.10-ում նշված պայմանը ևս հնարավոր չէ լավացնել: Իրոք, դիտարկենք լրիվ երկկողմանի $K_{n,n+2}$ գրաֆը: Հեշտ է տեսնել, որ $K_{n,n+2}$-ը $n$-կապակցված գրաֆ է և $\alpha(K_{n,n+2}) = n + 2$, բայց $K_{n,n+2}$-ը չի պարունակում համիլտոնյան ճանապարհի:

Համիլտոնյան գրաֆների կարևորագույն կիրառություններից մեկը կապված է *շրջիկ գործակալի խնդրի* հետ: Այդ խնդիրը կայանում է հետևյալում. տրված են $1, \ldots, n$ բնակավայրերը, հայտնի են նրանց միջև $d_{ij}$ հեռավորությունները, $1 \leq i < j \leq n$: Այդ $n$ բնակավայրերից մեկում գտնվում է գործակալը, որը պետք է շրջագայի բոլոր բնակավայրերը՝ յուրաքանչյուրում լինելով մեկ անգամ և վերադառնա մեկնակետ: Ի՞նչ հերթականությամբ գործակալը պետք է այցելի այդ բնակավայրերը, որպեսզի անցած ճանապարհի երկարությունը լինի նվազագույնը: Այժմ վերաձևակերպենք այս խնդիրը: Դիցուք տրված է $K_n$ լրիվ գրաֆը, որի գագաթների բազմությունը $\{1, \ldots, n\}$-ն է, իսկ կողերի բազմությունը՝ $\{ij : 1 \leq i < j \leq n\}$: Ավելին, այդ լրիվ գրաֆը նաև կշռված գրաֆ է, այսինքն՝ այդ գրաֆի յուրաքանչյուր կող ունի $d_{ij}$ երկարություն (կշիռ), $1 \leq i < j \leq n$: Այդ դեպքում հեշտ է տեսնել, որ շրջիկ գործակալի խնդիրը կարելի է ձևակերպել այսպես.

գտնել $K_n$ կշռված լրիվ գրաֆի այն համիլտոնյան ցիկլը, որի երկարությունը ամենակարճն է։ Իհարկե, այս խնդրի լուծման համար կարելի է առաջարկել հատարկման եղանակը։ Օրինակ, կարելի է ֆիքսել $1, \dots, n$ գագաթներից մեկը, որից միշտ սկսել շրջանցումը և դիտարկել մնացած գագաթների բոլոր հնարավոր տեղափոխությունները։ Հետևաբար, դիտարկելով $(n-1)!$ հատ տարբերակներ և յուրաքանչյուրի համար պարզելով ստացված տեղափոխությանը համապատասխանում է համիլտոնյան ցիկլ թե ոչ, կարելի է առանձնացնել $K_n$ կշռված լրիվ գրաֆի բոլոր համիլտոնյան ցիկլերը։ Այնուհետև, $K_n$ կշռված լրիվ գրաֆի բոլոր համիլտոնյան ցիկլերից ընտրել ամենակարճը։ Սակայն պետք է նշել, որ այդքան գործողություն կատարելը, նույնիսկ $n$-ի փոքր արժեքների դեպքում գործնականում անհնար է։ Վերջում նշենք, որ շրջիկ գործակալի խնդրի լուծման արդյունավետ եղանակների գոյության հարցը հանդիսանում է դիսկրետ մաթեմատիկայի դժվար և դեռևս չլուծված խնդիրներից մեկը։

# Գլուխ 5

## Անկախ բազմություններ, զուգակցումներ, ֆակտորներ և ծածկույթներ

### § 5.1. Անկախ բազմություններ և ծածկույթներ

Դիցուք $G = (V, E)$-ն գրաֆ է և $S \subseteq V$:

**Սահմանում 5.1.1:** Կասենք, որ $S$-ը գագաթների *անկախ բազմություն է* $G$ գրաֆում, եթե $S$-ը չի պարունակում հարևան գագաթներ:

**Սահմանում 5.1.2:** Կասենք, որ $S$-ը *գագաթային ծածկույթ է* $G$ գրաֆում, եթե $G - S$ գրաֆը չի պարունակում կող:

**Սահմանում 5.1.3:** $G$ գրաֆում ամենաշատ գագաթներ պարունակող անկախ բազմությունները կանվանենք *առավելագույն անկախ բազմություններ*:

**Սահմանում 5.1.4:** $G$ գրաֆում ամենաքիչ գագաթներ պարունակող գագաթային ծածկույթները կանվանենք *նվազագույն գագաթային ծածկույթներ*:

$G$ գրաֆում առավելագույն անկախ բազմության հզորությունը նշանակենք $\alpha(G)$-ով, իսկ նվազագույն գագաթային ծածկույթի հզորությունը՝ $\beta(G)$-ով:

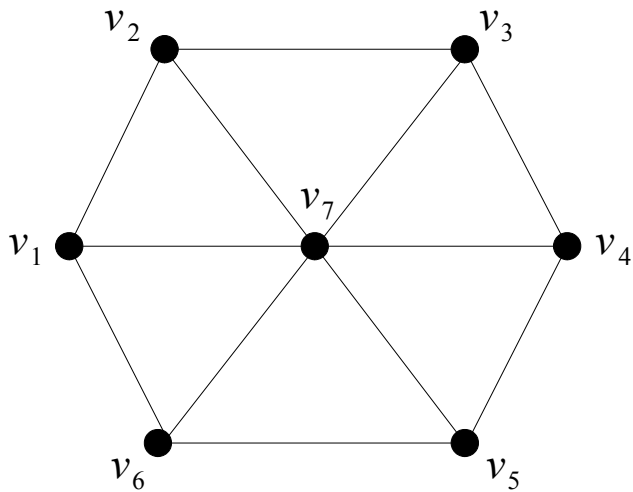

Նկ. 5.1.1

Նկ. 5.1.1-ում պատկերված $G$ գրաֆում $S_0 = \{v_1, v_3, v_5\}$ բազմությունը հանդիսանում է երեք հզորությամբ անկախ բազմություն, իսկ $S_1 = \{v_2, v_4, v_6, v_7\}$ բազմությունը՝ չորս հզորությամբ գագաթային ծածկույթ: Ավելին, կարելի է համոզվել, որ $\alpha(G) = 3$, իսկ $\beta(G) = 4$: Պարզվում է, որ ցանկացած $G$ գրաֆում $\alpha(G)$ և $\beta(G)$ պարամետրերը միմյանց հետ կապված են հետևյալ առնչությամբ:

**Թեորեմ 5.1.1 (Գալլաի):** Կամայական $G = (V, E)$ գրաֆում տեղի ունի $\alpha(G) + \beta(G) = |V|$ հավասարությունը:

**Ապացույց:** Դիցուք $S_0$-ն $G$ գրաֆում $\alpha(G)$ հզորությամբ անկախ բազմություն է: Նկատենք, որ այդ դեպքում $G$ գրաֆի ցանկացած կող կից է $V \backslash S_0$ բազմության գագաթի, և հետևաբար $V \backslash S_0$-ն հանդիսանում է $G$ գրաֆի գագաթային ծածկույթ, որտեղից հետևում է, որ

$$\beta(G) \leq |V \backslash S_0| = |V| - |S_0| = |V| - \alpha(G)$$

կամ

$$\alpha(G) + \beta(G) \leq |V|:$$

Մյուս կողմից, դիցուք $S_1$-ը $G$ գրաֆի $\beta(G)$ հզորությամբ գագաթային ծածկույթ է: Դա նշանակում է, որ $G - S_1$ գրաֆը չի պարունակում կող, որտեղից հետևում է, որ $V \backslash S_1$ բազմությունն իրենից ներկայացնում է գագաթների անկախ բազմություն: Սա նշանակում է, որ

$$\alpha(G) \geq |V \backslash S_1| = |V| - |S_1| = |V| - \beta(G)$$

կամ

$$\alpha(G) + \beta(G) \geq |V|:$$

$\alpha(G) + \beta(G) \leq |V|$ և $\alpha(G) + \beta(G) \geq |V|$ անհավասարություններից հետևում է, որ

$$\alpha(G) + \beta(G) = |V|: \blacksquare$$

Քանի որ ընդհանուր դեպքում տրված $G$ գրաֆի համար $\alpha(G)$ և $\beta(G)$ պարամետրերի գտնելու խնդիրները բարդ խնդիրներ են, բնական է դիտարկել տարբեր գնահատականներ այդ պարամետրերի համար: Այժմ մենք կներկայացնենք $\alpha(G)$-ի համար հայտնի գնահատականներից մեկը:

**Թեորեմ 5.1.2 (Կառո-Վեյի):** Կամայական $G$ գրաֆի համար տեղի ունի

$$\alpha(G) \geq \sum_{v \in V(G)} \frac{1}{1 + d_G(v)}$$

անհավասարությունը:

**Ապացույց:** Նախ նկատենք, որ եթե $G$ գրաֆը լրիվ գրաֆ է, ապա թեորեմը ճիշտ է: Ենթադրենք, որ $G$ գրաֆը լրիվ գրաֆ չէ:

Դիցուք $|V(G)| = n$։ Ապացույցը կատարենք մակածման եղանակով ըստ $n$-ի։ Հեշտ է տեսնել, որ թեորեմը ճիշտ է $n \leq 2$-ի դեպքում։ Ենթադրենք, թեորեմը ճիշտ է ցանկացած $G'$ գրաֆի համար, որի գագաթների քանակը $n$-ից ավելի չէ։ Դիտարկենք $n$ գագաթ ունեցող $G$ գրաֆը։ Դիցուք $x \in V(G)$ և $d_G(x) = \delta(G)$։ Քանի որ $G \neq K_n$, ուստի $\{x\} \cup N_G(x) \neq V(G)$։ Դիտարկենք $G' = G - x - N_G(x)$ գրաֆը։ Ըստ մակածման ենթադրության, $G'$ գրաֆի համար տեղի ունի $\alpha(G') \geq \sum_{v \in V(G')} \frac{1}{1+d_{G'}(v)}$ անհավասարությունը։ Դիցուք $S'$-ը առավելագույն անկախ բազմությունն է $G'$ գրաֆում։ Պարզ է, որ $S = \{x\} \cup S'$-ը անկախ բազմություն է $G$ գրաֆում։ Թեորեմի ապացույցը ավարտելու համար բավական է ցույց տալ, որ

$$\sum_{v \in V(G)} \frac{1}{1+d_G(v)} \leq \sum_{v \in V(G')} \frac{1}{1+d_{G'}(v)} + 1։$$

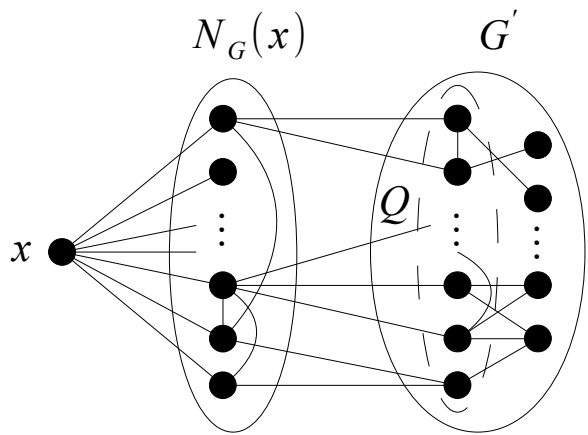

Նկ. 5.1.2

$Q$-ով նշանակենք $G'$ գրաֆի այն գագաթների բազմությունը, որոնք հարևան են $N_G(x)$-ի որևէ գագաթին $G$ գրաֆում (նկ. 5.1.2)։ Այդ դեպքում $G'$ գրաֆի համար տեղի ունի հետևյալը.

$$\sum_{v \in Q} \frac{1}{1+d_G(v)} \leq \sum_{v \in Q} \frac{1}{1+d_{G'}(v)} \text{ և } \sum_{v \in V(G') \setminus Q} \frac{1}{1+d_G(v)} = \sum_{v \in V(G') \setminus Q} \frac{1}{1+d_{G'}(v)}։$$

Այստեղից հետևում է, որ թեորեմն ապացուցելու համար մնում է ցույց տալ, որ

$$\frac{1}{1+d_G(x)} + \sum_{v \in N_G(x)} \frac{1}{1+d_G(v)} \leq 1։$$

Քանի որ $x$ գագաթը ընտրել ենք այնպես, որ $d_G(x) = \delta(G)$, ուստի ցանկացած $v \in N_G(x)$-ի համար տեղի ունի $d_G(x) \leq d_G(v)$ անհավասարությունը։ Հետևաբար,

$$\frac{1}{1+d_G(x)} + \sum_{v \in N_G(x)} \frac{1}{1+d_G(v)} \leq \frac{1}{1+d_G(x)} + \frac{d_G(x)}{1+d_G(x)} = 1: \quad \blacksquare$$

## § 5.2. Զուգակցումներ երկկողմանի գրաֆներում և min-max թեորեմներ

Դիցուք $G = (V, E)$-ն գրաֆ է և $M \subseteq E$:

**Սահմանում 5.2.1:** Կասենք, որ $M$-ը *կողերի անկախ բազմություն է* $G$ գրաֆում, եթե $M$-ը չի պարունակում հարևան կողեր:

Կողերի անկախ բազմությանն ընդունված է անվանել *զուգակցում*: Դիցուք $M$-ը $G$ գրաֆի զուգակցում է:

**Սահմանում 5.2.2:** Կասենք, որ $M$ զուգակցումը *փակուղային է*, եթե $G$ գրաֆում գոյություն չունի այնպիսի $e \notin M$ կող, որ $M \cup \{e\}$-ն լինի զուգակցում:

**Սահմանում 5.2.3:** Կասենք, որ $M$ զուգակցումը *առավելագույնն է*, եթե հզորությամբ նրանից մեծ զուգակցում $G$ գրաֆում չկա:

$G$ գրաֆում առավելագույն զուգակցման հզորությունը կնշանակենք $\alpha'(G)$-ով: Նկատենք, որ ցանկացած $G$ գրաֆում $\alpha'(G) \leq \frac{|V|}{2}$:

**Սահմանում 5.2.4:** Կասենք, որ $M$ զուգակցումը *կատարյալ է*, եթե այն պարունակում է $\frac{|V|}{2}$ կող:

Նկատենք, որ գրաֆի կատարյալ զուգակցումը հանդիսանում է առավելագույն զուգակցում: Նշենք, որ հակառակը ճիշտ չէ, քանի որ գրաֆը կարող է չպարունակել կատարյալ զուգակցում (օրինակ, եռանկյունը), մինչդեռ առավելագույն զուգակցում գոյություն ունի միշտ:

**Սահմանում 5.2.5:** Կասենք, որ $G$ գրաֆի $M$ զուգակցումը *հագեցնում է* $v$ գագաթը, եթե $M$ զուգակցումը պարունակում է $v$ գագաթին կից կող:

Նկատենք, որ ցանկացած $M$ զուգակցում հագեցնում է $2|M|$ գագաթ, և, հետևաբար, այն չի հագեցնում $|V| - 2|M|$ գագաթ: Ավելին, նկատենք, որ զուգակցումը հանդիսանում է կատարյալ զուգակցում այն և միայն այն դեպքում, երբ այն հագեցնում է գրաֆի բոլոր գագաթները:

Դիտարկենք բերված սահմանումները պարզաբանող օրինակ:

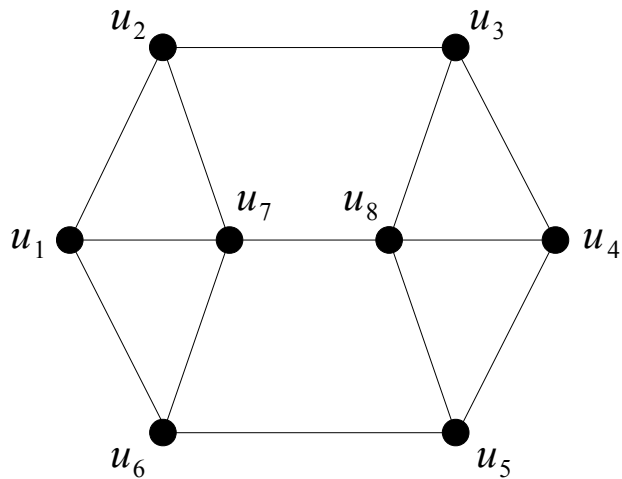

Նկ. 5.2.1

Նկ. 5.2.1-ում պատկերված $G$ գրաֆում կողերի $M = \{u_2u_3, u_7u_8\}$ բազմությունը հանդիսանում է զուգակցում։ $M$-ը հագեցնում է $u_2, u_3, u_7, u_8$ գագաթները, և չի հագեցնում՝ $u_1, u_4, u_5, u_6$ գագաթները։ Նկատենք, որ այն չի հանդիսանում փակուղային զուգակցում, քանի որ $M$-ին կարելի է ավելացնել $u_5u_6$ կողը և ստանալ ավելի մեծ զուգակցում։ Ավելին, նկատենք, որ $M' = \{u_2u_3, u_7u_8, u_5u_6\}$ զուգակցումն արդեն հանդիսանում է $G$ գրաֆի փակուղային զուգակցում, բայց այն չի հանդիսանում առավելագույն զուգակցում, քանի որ $G$ գրաֆում կողերի $\{u_1u_2, u_3u_4, u_5u_6, u_7u_8\}$ բազմությունը հանդիսանում է ավելի մեծ զուգակցում։ Նկատենք, որ վերջինս արդեն հանդիսանում է առավելագույն զուգակցում, քանի որ այն նաև կատարյալ զուգակցում է։

Դիցուք $M$-ը զուգակցում է, իսկ $P$-ն որևէ ճանապարհ է $G$ գրաֆում։

**Սահմանում 5.2.6:** Կասենք, որ $P$ ճանապարհը հանդիսանում է $M$-*հերթափոխ ճանապարհ*, եթե $P$ ճանապարհի կողերը հերթականորեն պատկանում են և չեն պատկանում $M$-ին։

**Սահմանում 5.2.7:** Կասենք, որ $M$-հերթափոխ $P$ ճանապարհը $M$-*ավելացնող է*, եթե այն միացնում է երկու գագաթ, որոնք հագեցած չեն $M$-ով։

Կրկին դիտարկենք նկ. 5.2.1-ում պատկերված $G$ գրաֆի $M = \{u_2u_3, u_7u_8\}$ զուգակցումը։ Նշենք, որ $u_1, u_2, u_3, u_8, u_7$ ճանապարհը հանդիսանում է $M$-հերթափոխ ճանապարհ, որը չի հանդիսանում $M$-ավելացնող ճանապարհ։ Մյուս կողմից, նշենք, որ $u_1, u_2, u_3, u_8, u_7, u_6$ ճանապարհը հանդիսանում է $M$-ավելացնող ճանապարհ։ Ավելին, $M$-ավելացնող ճանապարհի օրինակ է նաև $u_5, u_6$ մեկ երկարությամբ ճանապարհը։

$M$-ավելացնող ճանապարհների նշանակությունը և դերը երևում են հետևյալ

թեորեմի վրա:

**Թեորեմ 5.2.1 (Բերժ):** Որպեսզի $M$ զուգակցումը լինի առավելագույն, անհրաժեշտ է և բավարար, որ $G$ գրաֆը չպարունակի $M$-ավելացնող ճանապարհի:

**Ապացույց:** Նախ ենթադրենք, որ $M$ զուգակցումը առավելագույն է: Ցույց տանք, որ $G$ գրաֆում գոյություն չունեն $M$-ավելացնող ճանապարհներ:

Ենթադրենք հակառակը. դիցուք $G$ գրաֆը պարունակում է $M$-ավելացնող $P$ ճանապարհը: Դիտարկենք $G$ գրաֆի կողերի $M'$ ենթաբազմությունը, որը ստացվում է $M$-ից հետևյալ կերպ.

$$M' = (M \setminus E(P)) \cup (E(P) \setminus M):$$

Նկատենք, քանի որ $M$-ավելացնող $P$ ճանապարհի ծայրակետերը հագեցած չեն $M$-ով, $M'$-ը կհանդիսանա զուգակցում: Ավելին, նկատենք, որ $|M'| = |M| + 1 > |M|$, ինչը հակասում է $M$ զուգակցման առավելագույնը լինելուն: Հետևաբար, $G$ գրաֆում չկան $M$-ավելացնող ճանապարհներ:

Հիմա ենթադրենք, որ $G$ գրաֆում չկան $M$-ավելացնող ճանապարհներ: Ցույց տանք, որ $M$ զուգակցումը առավելագույն է:

Ենթադրենք հակառակը. դիցուք $M$ զուգակցումը առավելագույնը չէ: Սա նշանակում է, որ $G$ գրաֆում գոյություն ունի այնպիսի $M'$ զուգակցում, որ $|M'| > |M|$: Դիտարկենք $G$ գրաֆի $H$ ենթագրաֆը, որտեղ $V(H) = V(G)$ և $E(H) = (M \setminus M') \cup (M' \setminus M)$: Նկատենք, որ $H$ ենթագրաֆում $\Delta(H) \leq 2$, և, հետևաբար, նրա կապակցվածության յուրաքանչյուր բաղադրիչ իրենից ներկայացնում է ճանապարհ կամ ցիկլ: Քանի որ ցիկլի դեպքում նրա կողերը մեկումեջ պատկանում են $M$-ին և $M'$-ին, ապա ցիկլի երկարությունը կլինի զույգ:

Քանի որ $|M'| > |M|$, գոյություն կունենա $H$ ենթագրաֆի կապակցվածության բաղադրիչ, որում $M'$-ի կողերը ավելի շատ են, քան $M$-ի կողերը: Նկատենք, որ սա հնարավոր է միայն այն դեպքում, երբ կապակցվածության այդ բաղադրիչն իրենից ներկայացնում է կենտ երկարությամբ ճանապարհ, որի առաջին և վերջին կողերը $M'$-ից են: Նկատենք, որ կապակցվածության այդպիսի բաղադրիչը կլինի $M$-ավելացնող ճանապարհի, որն էլ կհակասի մեր ենթադրությանը: Հետևաբար, $M$ զուգակցումը առավելագույն է: ∎

Հիշենք, որ § 1.2-ում $G = (V, E)$ գրաֆի գագաթների $S \subseteq V$ ենթաբազմության համար սահմանեցինք $N_G(S)$ բազմությունը որպես

$$N_G(S) = \{u \in V \setminus S: \text{գոյություն ունի } v \in S, \text{որ } uv \in E\}:$$

Նկատենք, որ եթե $G$ գրաֆն երկկողմանի է, $V = V_1 \cup V_2$ նրա գագաթների բազմության համապատասխան տրոհումն է և $S \subseteq V_1$, ապա $N_G(S) \subseteq V_2$։ Ստորև կձևակերպենք և կապացուցենք Հոլլի թեորեմը։

**Թեորեմ 5.2.2 (Հոլլ):** Դիցուք $G = (V, E)$ գրաֆն երկկողմանի է և $V = V_1 \cup V_2$ գագաթների բազմության համապատասխան տրոհումն է։ Որպեսզի $G$ գրաֆը պարունակի $M$ զուգակցում, որը հագեցնում է $V_1$ բազմությունը, անհրաժեշտ է և բավարար, որ ցանկացած $S \subseteq V_1$ համար $|N_G(S)| \geq |S|$։

**Ապացույց:** Նախ նկատենք, որ եթե $G$ գրաֆում գոյություն ունի $V_1$ բազմությունը հագեցնող $M$ զուգակցում, ապա այդ զուգակցումը ցանկացած $S \subseteq V_1$ ենթաբազմության գագաթները կարտապատկերի $V_2$ ենթաբազմության իրարից տարբեր գագաթների, հետևաբար՝ $|N_G(S)| \geq |S|$։

Հակառակն ապացուցելու համար նկատենք, քանի որ $G$ երկկողմանի գրաֆում ցանկացած զուգակցում պարունակում է ոչ ավելի, քան $|V_1|$ կող, ապա պնդումն ապացուցելու համար բավական է ցույց տալ, որ $G$ գրաֆի կամայական առավելագույն զուգակցումը հագեցնում է $V_1$ բազմության բոլոր գագաթները։

Ենթադրենք հակառակը. այսինքն ենթադրենք, որ $G$ երկկողմանի գրաֆում ցանկացած $S \subseteq V_1$ համար $|N_G(S)| \geq |S|$, բայց նրա $M$ առավելագույն զուգակցումը չի հագեցնում $V_1$ բազմության բոլոր գագաթները։ Այդ դեպքում գոյություն կունենա $u \in V_1$ գագաթ, որը հագեցած չէ $M$ զուգակցումով։

Նշանակենք $S$-ով և $T$-ով $G$ գրաֆի գագաթների այն ենթաբազմությունները, համապատասխանաբար, $V_1$-ից և $V_2$-ից, որոնք հասանելի են $u$ գագաթից $M$-հերթափոխ ճանապարհով։ Նկատենք, որ $u \in S$։

Ցույց տանք, որ $M$ զուգակցումը հաստատում է փոխմիարժեք արտապատկերում $S \setminus \{u\}$ և $T$ բազմությունների գագաթների միջև։ Իրոք, $u$ գագաթից սկսվող $M$-հերթափոխ ճանապարհները հասնում են $V_2$ $M$ զուգակցմանը չպատկանող կողով, իսկ $V_1$՝ $M$ զուգակցմանը պատկանող կողով։ Հետևաբար, $S \setminus \{u\}$ բազմության ցանկացած գագաթ հասանելի է $T$ բազմության գագաթից $M$ զուգակցման կողով։ Քանի որ $M$ զուգակցումը առավելագույն էր, ապա համաձայն թեորեմ 5.2.1-ի, $G$ գրաֆում գոյություն չունեն $M$-ավելացնող ճանապարհներ, և, հետևաբար, $T$ բազմության ցանկացած գագաթ հագեցած է $M$ զուգակցումով, որտեղից հետևում է, որ եթե $M$-հերթափոխ ճանապարհը $u$ գագաթից հասել է $y \in T$ գագաթ, ապա $y$-ին կից $M$ զուգակցման կողը նրան կհամապատասխանեցնի $S \setminus \{u\}$ բազմության զագաթի: Ասվածից հետևում է, որ $M$ զուգակցումը հաստատում է փոխմիարժեք արտապատկերում $S \setminus \{u\}$ և $T$ բազմությունների զագաթների միջև, որը, մասնավորապես, նշանակում է, որ $|T| = |S \setminus \{u\}|$:

Մյուս կողմից, քանի որ $M$ զուգակցումը արտապատկերում է $T$ բազմության զագաթները $S \setminus \{u\}$ բազմության զագաթներին, ապա $T \subseteq N_G(S)$: Ցույց տանք, որ $T = N_G(S)$: Իրոք, եթե գոյություն ունենար $y \in N_G(S) \setminus T$ զագաթ, ապա այն հագեցած չէր լինի $M$ զուգակցման կողմով, և այն կառաջացներ $u$ զագաթից սկիզբ առնող $M$-հերթափոխ ճանապարհի դեպի $y$ զագաթ, ինչը կհակասեր $y \notin T$ պայմանին: Հետևաբար, $T = N_G(S)$: Արդյունքում՝

$$|N_G(S)| = |T| = |S \setminus \{u\}| = |S| - 1 < |S|$$

ինչը հակասում է թեորեմի պայմաններին: ∎

Դիցուք տրված է $S = \{s_1, s_2, \ldots, s_n\}$ բազմությունը և այդ բազմության ենթաբազմությունների $\mathfrak{F} = \{F_1, F_2, \ldots, F_m\}$ ընտանիքը:

**Սահմանում 5.2.8:** $S$ բազմության տարրերի $(s_{i_1}, s_{i_2}, \ldots, s_{i_m})$ $m$-յակը կանվանենք *տարբեր ներկայացուցիչների համակարգ* $\mathfrak{F}$ ընտանիքի համար, եթե $s_{i_1} \in F_1$, $s_{i_2} \in F_2$, …, $s_{i_m} \in F_m$ և $s_{i_p} \neq s_{i_q}$, երբ $p \neq q$:

Ստորև կձևակերպենք և կապացուցենք տարբեր ներկայացուցիչների համակարգի գոյության անհրաժեշտ և բավարար պայմանը:

**Թեորեմ 5.2.3 (Հոլլ):** Որպեսզի $S = \{s_1, s_2, \ldots, s_n\}$ բազմության ենթաբազմությունների $\mathfrak{F} = \{F_1, F_2, \ldots, F_m\}$ ընտանիքն ունենա տարբեր ներկայացուցիչների համակարգ, անհրաժեշտ է և բավարար, որ $\mathfrak{F}$ ընտանիքին պատկանող ցանկացած $k$ $F_{j_1}, F_{j_2}, \ldots, F_{j_k}$ բազմությունների համար տեղի ունենա

$$|F_{j_1} \cup F_{j_2} \cup \ldots \cup F_{j_k}| \geq k$$

պայմանը, որտեղ $1 \leq k \leq m$:

**Ապացույց:** Նախ նկատենք, որ թեորեմի անհրաժեշտությունն ակնհայտ է: Բավարարությունն ապացուցելու համար դիտարկենք $\mathfrak{F} = \{F_1, F_2, \ldots, F_m\}$ ընտանիքի § 1.5-ում սահմանված կցության $G(S, \mathfrak{F})$ գրաֆը: Նկատենք, որ այն երկկողմանի է, ավելին, այն բավարարում է թեորեմ 5.2.2-ի պայմաններին: Համաձայն այդ թեորեմի, $G(S, \mathfrak{F})$ գրաֆում գոյություն կունենա $\mathfrak{F}$-ը հագեցնող զուգակցում, որին կհամապատասխանի $\mathfrak{F}$ ընտանիքի տարբեր ներկայացուցիչների համակարգ: ∎

Դիտարկենք մեկ օրինակ։ Դիցուք $S = \{s_1, s_2, s_3, s_4\}$ և $\mathfrak{F} = \{F_1, F_2, F_3\}$, որտեղ $F_1 = \{s_2, s_3\}$, $F_2 = \{s_1, s_3, s_4\}$ և $F_3 = \{s_3, s_4\}$։

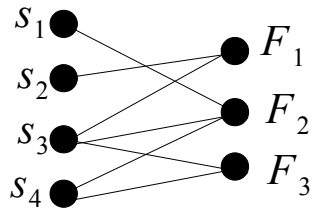

Նկ. 5.2.2

Նկ. 5.2.2-ում պատկերված է $\mathfrak{F}$ ընտանիքի կցության $G(S, \mathfrak{F})$ գրաֆը։ Նկատենք, որ այն պարունակում է $\mathfrak{F}$-ը հագեցնող երկու զուգակցում։ Մեկը $M = \{s_2F_1, s_1F_2, s_3F_3\}$, իսկ մյուսը՝ $M' = \{s_2F_1, s_1F_2, s_4F_3\}$։ $M$-ին համապատասխանում է տարբեր ներկայացուցիչների ($s_2, s_1, s_3$) համակարգը, իսկ $M'$-ին՝ ($s_2, s_1, s_4$) համակարգը։

**Թեորեմ 5.2.4 (Ֆրոբենիուս):** Ցանկացած երկկողմանի $r$-համասեռ գրաֆ ($r \in \mathbb{N}$) պարունակում է կատարյալ զուգակցում։

**Ապացույց:** Դիցուք $G$-ն երկկողմանի $r$-համասեռ գրաֆ է, $r \in \mathbb{N}$, և դիցուք $V = V_1 \cup V_2$ նրա գագաթների բազմության համապատասխան տրոհումն է։ Նկատենք, որ $G$ գրաֆի կողերի քանակը հավասար է $r|V_1|$։ Մյուս կողմից, այն հավասար է $r|V_2|$, և քանի որ $r \geq 1$, ապա $|V_1| = |V_2|$։ Այստեղից հետևում է, որ թեորեմն ապացուցելու համար բավական է ցույց տալ, որ $G$-ն պարունակում է $V_1$ բազմությունը հագեցնող զուգակցում։

Դիտարկենք ցանկացած $S \subseteq V_1$։ Նկատենք, որ $S$-ից դուրս եկող կողերի քանակը հավասար է $r|S|$։ Մյուս կողմից, $S$-ից դուրս եկող ցանկացած կող հարևան է $N_G(S)$ բազմության գագաթի, հետևաբար այդպիսի կողերի քանակը չի գերազանցում $r|N_G(S)|$-ը։ Այստեղից հետևում է, որ $|N_G(S)| \geq |S|$, քանի որ $r \geq 1$։ Հաշվի առնելով, որ $S$-ը ընտրված էր կամայապես, թեորեմ 5.2.2-ից կստանանք, որ $G$-ն պարունակում է $V_1$ բազմությունը հագեցնող զուգակցում։ ∎

Ստորև կձևակերպենք և կապացուցենք երեք թեորեմ, որոնք պատկանում են այսպես կոչված *մինմաքս թեորեմների* շարքին։ Նշենք, որ թեորեմն անվանում են *մինմաքս*, եթե այն պնդում է, որ ինչ-որ մի պարամետրի մաքսիմում հնարավոր արժեքը հավասար է մեկ այլ պարամետրի մինիմում հնարավոր արժեքին։

Նկատենք, որ ցանկացած $G$ գրաֆում $\beta(G) \geq \alpha'(G)$։ Սա հետևում է այն բանից, որ $\alpha'(G)$ հատ անկախ կող պարունակող զուգակցման յուրաքանչյուր կող ծածկելու համար անհրաժեշտ է առնվազն մեկ գագաթ։ Ստորև կապացուցենք Քյոնիգ-Էգերվարի թեորեմը, որն առաջարկում է վերը նշված պարամետրերի հավասարության բավարար պայման։

**Թեորեմ 5.2.5 (Քյոնիգ-Էգերվարի):** Ցանկացած երկկողմանի $G$ գրաֆում $\beta(G) = \alpha'(G)$։

**Ապացույց:** Քանի որ ցանկացած $G$ գրաֆում $\beta(G) \geq \alpha'(G)$, ապա պնդումն ապացուցելու համար բավական է ցույց տալ, որ եթե $U$-ն $\beta(G)$ հզորությամբ զազաքային ծածկույթ է $G$ երկկողմանի գրաֆում, ապա նրանում գոյություն ունի $\beta(G)$ հզորությամբ զուգակցում։

Դիցուք $V = V_1 \cup V_2$ $G$ երկկողմանի գրաֆի զազաքների բազմության համապատասխան տրոհումն է, և $U$-ն $\beta(G)$ հզորությամբ զազաքային ծածկույթ է։ Նշանակենք՝

$$R = U \cap V_1 \text{ և } T = U \cap V_2,$$

և դիցուք $H = G[R \cup (V_2 \setminus T)]$, $H' = G[T \cup (V_1 \setminus R)]$։ Նկատենք, որ $H$-ը և $H'$-ը երկկողմանի գրաֆներ են։

Ցույց տանք, որ $H$-ը պարունակում է $R$-ը հագեցնող զուգակցում, և $H'$-ը պարունակում է $T$-ն հագեցնող զուգակցում։ Նկատենք, որ եթե այս երկու պնդումն ապացուցենք, ապա, քանի որ $H$ և $H'$ գրաֆների կողերի բազմությունները չեն հատվում, ապա այդ զուգակցումների միավորումը կհանդիսանա $G$ երկկողմանի գրաֆի $|R| + |T| = |U| = \beta(G)$ հզորությամբ զուգակցում, ինչը կապացուցի թեորեմը։

Քանի որ $R \cup T$-ն զազաքային ծածկույթ է, ապա $G$ գրաֆում չկա կող, որը միացնում է $V_2 \setminus T$ բազմության զազաքը $V_1 \setminus R$ բազմության զազաքին։ Դիտարկենք ցանկացած $S \subseteq R$ և $N_H(S) \subseteq V_2 \setminus T$ բազմությունները։ Եթե $|N_H(S)| < |S|$, ապա $U$ զազաքային ծածկույթի մեջ $S$-ը փոխարինելով $N_H(S)$-ով, մենք կստանանք $|U| = \beta(G)$-ից ավելի փոքր հզորությամբ զազաքային ծածկույթ $G$ գրաֆում ($N_H(S)$-ը ծածկում է $S$-ից դուրս եկող բոլոր այն կողերը, որոնք ծածկված չեն $T$-ով), ինչը կիակասի $U$-ի ընտրությանը։ Հետևաբար, ցանկացած $S \subseteq R$ համար $|N_H(S)| \geq |S|$, որտեղից, համաձայն թեորեմ 5.2.2-ի, կստանանք, որ $H$-ը պարունակում է $R$-ը հագեցնող զուգակցում։ Համանման դատողություններով կարելի է ապացուցել, որ $H'$-ը պարունակում է $T$-ն հագեցնող զուգակցում։ ∎

Ապացուցված թեորեմն ունի մեկ հետաքրքիր մեկնաբանություն։ Դիցուք $A$-ն $m \times n$ կարգի մատրից է, որի տարրերը զրո կամ մեկ են։ Այդպիսի մատրիցում *շարք* ասելով կհասկանանք ցանկացած տող կամ սյուն, և $A$ մատրիցում երկու մեկ կհամարենք *անկախ*, եթե նրանք գտնվում են տարբեր տողերում և տարբեր սյուներում (չեն պատկանում միևնույն շարքին): $\psi(A)$-ով նշանակենք նվազագույն թվով շարքերի քանակը $A$ մատրիցում, որոնք ընդգրկում են $A$-ի բոլոր մեկերը, և $\Psi(A)$-ով նշանակենք առավելագույն թվով անկախ մեկերի քանակը: Պարզ է, որ $\Psi(A) \leq \psi(A)$:

**Թեորեմ 5.2.6:** Զրոներից և մեկերից կազմված ցանկացած $A$ մատրիցում տեղի ունի $\Psi(A) = \psi(A)$ հավասարությունը:

**Ապացույց:** Դիցուք $A$-ն զրոներից և մեկերից կազմված ցանկացած $m \times n$ կարգի մատրից է: Դիտարկենք $G$ գրաֆը, որի գագաթների բազմությունն է $V = \{b_1, \ldots, b_m, c_1, \ldots, c_n\}$-ը, իսկ կողերը ստացվում են հետևյալ կերպ. $b_i c_j \in E(G)$ այն և միայն այն դեպքում, երբ $A$ մատրիցում $i$-րդ տողի և $j$-րդ սյան հատման վանդակում գրված է մեկ:

Նկատենք, որ $G$ գրաֆը երկկողմանի է: Ավելին, $A$ մատրիցի ցանկացած թվով անկախ մեկերի համապատասխանում է նույն հզորությամբ զուգակցում $G$ գրաֆում և հակառակը, հետևաբար՝ $\Psi(A) = \alpha'(G)$: Եվ վերջապես, $A$ մատրիցի ցանկացած թվով շարքերի, որոնք ընդգրկում են մատրիցի բոլոր մեկերը, համապատասխանում է $G$ գրաֆի նույն հզորությամբ գագաթային ծածկույթ և հակառակը, հետևաբար՝ $\psi(A) = \beta(G)$: Հաշվի առնելով թեորեմ 5.2.5-ը, կստանանք՝

$$\Psi(A) = \alpha'(G) = \beta(G) = \psi(A): \blacksquare$$

Դիտարկենք վերջին երկու թեորեմները պարզաբանող օրինակ: Դիցուք՝

$$A = \begin{pmatrix} 1 & 0 & 1 & 0 & 0 \\ 0 & 0 & 0 & 1 & 1 \\ 0 & 0 & 1 & 0 & 1 \\ 0 & 1 & 0 & 1 & 0 \end{pmatrix}:$$

Թեորեմ 5.2.6-ի ապացույցում $A$ մատրիցին համապատասխանեցվող երկկողմանի $G$ գրաֆը պատկերված է ստորև.

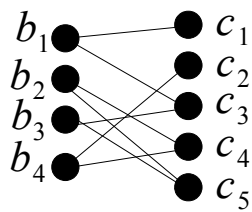

Նկ. 5.2.3

Նկատենք, որ $A$ մատրիցում կան չորս անկախ մեկ. $a_{11}, a_{24}, a_{33}, a_{42}$, որոնց համապատասխանում է $G$ գրաֆի $\{b_1 c_1, b_2 c_4, b_3 c_3, b_4 c_2\}$ զուգակցումը: Հեռացնելով $A$ մատրիցի բոլոր տողերը, մենք կոչնչացնենք նրա բոլոր մեկերը: Այս շարքերին կիսամապատասխանի $G$ գրաֆի $\{b_1, b_2, b_3, b_4\}$ զագաթային ծածկույթը։ Նկատենք, որ բերված օրինակում

$$\Psi(A) = \alpha'(G) = \beta(G) = \psi(A) = 4:$$

Դիցուք $P$-ն ցանկացած վերջավոր բազմություն է, և $\alpha$-ն նրա վրա տրված բինար հարաբերություն է։ Հիշենք, որ $\alpha$-ն անվանում են *կարգի հարաբերություն* տրված $P$ բազմության վրա, եթե այն բավարարում է հետևյալ երեք պայմաններին.

1. ցանկացած $x \in P$ համար $x\alpha x$;
2. ցանկացած $x, y \in P$ համար եթե $x\alpha y$ և $y\alpha x$, ապա $x = y$;
3. ցանկացած $x, y, z \in P$ համար եթե $x\alpha y$ և $y\alpha z$, ապա $x\alpha z$:

Կարգի հարաբերությունը սովորաբար նշանակում են ≤ սիմվոլով, իսկ այն բազմությունը, որի վրա տրված է կարգի հարաբերությունը, անվանում են *կարգավորված բազմություն*։ Եթե $x \leq y$ և $x \neq y$, ապա սովորաբար գրում են $x < y$։ Հիշենք նաև, որ $x \in P$ տարրն անվանում են *մինիմալ* (*մաքսիմալ*), եթե $P$-ում գոյություն չունի $y$ տարր այնպես, որ $y < x$ ($x < y$)։ Դժվար չէ համոզվել, որ ցանկացած վերջավոր բազմության վրա տրված կարգի հարաբերության նկատմամբ միշտ գոյություն ունեն մինիմալ և մաքսիմալ տարրեր։

$(P, \leq)$ կարգավորված բազմության $x$ և $y$ տարրերն անվանում են *համեմատելի*, եթե $x \leq y$ կամ $y \leq x$ (հակառակ դեպքում այդ տարրերը կանվանենք *անհամեմատելի*)։ Վերջավոր $(P, \leq)$ կարգավորված բազմության համար $\Phi(P)$-ով նշանակենք $P$-ում առավելագույն թվով զույգ առ զույգ անհամեմատելի տարրերի քանակը։

Ընդունված է ասել, որ $P$ բազմության իրարից տարբեր $x_1, x_2, \ldots, x_k$ տարրերը կազմում են *շղթա*, եթե $x_1 < x_2, x_2 < x_3, \ldots, x_{k-1} < x_k$։ $\varphi(P)$-ով նշանակենք վերջավոր $(P, \leq)$ կարգավորված բազմության նվազագույն թվով չհատվող շղթաների քանակը, որոնք պարունակում են $P$ բազմության բոլոր տարրերը։ Ստորև կձևակերպենք և կապացուցենք Դիլվորթի թեորեմը, որը պնդում է, որ վերը սահմանված երկու պարամետրերն իրականում հավասար են ցանկացած վերջավոր $(P, \leq)$ կարգավորված բազմության համար։

**Թեորեմ 5.2.7 (Դիլվորթ):** Ցանկացած վերջավոր $(P, \leq)$ կարգավորված բազմության համար տեղի ունի $\Phi(P) = \varphi(P)$ հավասարությունը։

**Ապացույց:** Նախ նկատենք, քանի որ $P$ բազմության ցանկացած շղթայի վրա ցանկացած երկու տարր համեմատելի են, ապա $\Phi(P) \leq \varphi(P)$։ Հակառակ անհավասարությունը ցույց տալու համար ենթադրենք, որ $\Phi(P) = \mathbf{n}$, և մակածման եղանակով ըստ $P$ բազմության տարրերի քանակի ցույց տանք, որ $P$-ն հնարավոր է տրոհել $\mathbf{n}$ չհատվող շղթաների:

Նկատենք, որ պնդումն ակնհայտ է, երբ $|P| = 1$: Ենթադրենք, որ նշված պնդումը ճիշտ է բոլոր այն կարգավորված $Q$ բազմությունների համար, որոնք բավարարում են $|Q| < |P|$ պայմանին, և դիտարկենք $P$ կարգավորված բազմությունը:

Քննարկենք երկու դեպք:

Դեպք 1: $P$ բազմության մեջ գոյություն ունի $\mathbf{n}$ հատ զույգ առ զույգ անհամեմատելի տարրեր պարունակող այնպիսի $U$ բազմություն, որը չի ներառում $P$ բազմության բոլոր մաքսիմալ տարրերը, և ոչ էլ $P$ բազմության բոլոր մինիմալ տարրերը:

Դիտարկենք $P^+$ և $P^-$ բազմություններն, որոնք սահմանվում են հետևյալ կերպ.

$$P^+ = \{p \in P : u \leqslant p \text{ ինչ որ մի } u \in U - \text{ից}\},$$

$$P^- = \{p \in P : p \leqslant u \text{ ինչ որ մի } u \in U - \text{ից}\}:$$

Հաշվի առնելով $U$ բազմության ընտրությունը, կստանանք

$$P^+ \neq P, \quad P^- \neq P \text{ և } P = P^+ \cup P^-, U = P^+ \cap P^-:$$

Համաձայն մակածման ենթադրության, $P^+$ և $P^-$ բազմություններից յուրաքանչյուրը հնարավոր է տրոհել $\mathbf{n}$ չհատվող շղթաների: Սոսնձելով այդ շղթաները $U$ բազմությանը պատկանող տարրերում, մենք կստանանք $P$ բազմության տրոհում $\mathbf{n}$ չհատվող շղթաների:

Դեպք 2: $P$ բազմության մեջ ցանկացած $\mathbf{n}$ հատ զույգ առ զույգ անհամեմատելի տարրեր պարունակող $U$ բազմություն ներառում է $P$ բազմության բոլոր մաքսիմալ տարրերը կամ $P$ բազմության բոլոր մինիմալ տարրերը:

Հետևաբար, $P$ բազմության մեջ գոյություն ունեն ոչ ավելի, քան երկու $\mathbf{n}$ հատ զույգ առ զույգ անհամեմատելի տարրեր պարունակող բազմություն, որոնցից մեկը կներառի $P$ բազմության բոլոր մաքսիմալ տարրերը, իսկ մյուսը $P$ բազմության բոլոր մինիմալ տարրերը: Դիտարկենք $P$ բազմության ցանկացած $a$ մինիմալ և $b$ մաքսիմալ տարրեր, որոնք բավարարում են $a \leqslant b$: Համաձայն մակածման ենթադրության, $P \setminus \{a, b\}$ բազմությունը հնարավոր է տրոհել $\mathbf{n} - \mathbf{1}$ չհատվող շղթաների: Այդ շղթաներին ավելացնելով $a \leqslant b$ շղթան, կստանանք $P$ բազմության որոնելի տրոհումը $\mathbf{n}$ շղթաների: ∎

# § 5.3. Զուգակցումներ գրաֆներում և Տատտի թեորեմը

$G$ գրաֆի համար $o(G)$-ով նշանակենք $G$-ի կապակցվածության այն բաղադրիչների քանակը, որոնք պարունակում են կենտ թվով գագաթներ: Ստորև կձևակերպենք և կապացուցենք Տատտի թեորեմը, որը նկարագրում է կատարյալ զուգակցում պարունակող գրաֆները: Հիշեցնենք, որ զուգակցումը կոչվում է կատարյալ, եթե այն հագեցնում է գրաֆի բոլոր գագաթները:

**Թեորեմ 5.3.1 (Տատտ):** Որպեսզի $G$ գրաֆը պարունակի կատարյալ զուգակցում, անհրաժեշտ է և բավարար, որ ցանկացած $S \subseteq V(G)$ համար տեղի ունենա $o(G - S) \leq |S|$ պայմանը:

**Ապացույց:** Նախ ենթադրենք, որ $G$ գրաֆը պարունակում է $M$ կատարյալ զուգակցումը: Դիտարկենք ցանկացած $S \subseteq V(G)$: Դիցուք $o(G - S) = t$, և ենթադրենք, որ $G_1, \ldots, G_t$-ն $G - S$ գրաֆի կապակցվածության այն բաղադրիչներն են, որոնք պարունակում են կենտ թվով գագաթներ: Նկատենք, որ քանի որ $i = 1, \ldots, t$ համար $G_i$ գրաֆի գագաթների քանակը կենտ է, ապա $M \cap E(G_i)$ զուգակցումը չի կարող հագեցնել $G_i$ գրաֆի բոլոր գագաթները, հետևաբար $G_i$-ում գոյություն ունի գոնե մեկ գագաթ, որը հագեցնող և $M$ կատարյալ զուգակցմանը պատկանող կողը կից է $S$-ին պատկանող գագաթի: Քանի որ $M$-ը զուգակցում է, ապա $S$-ին պատկանող վերոհիշյալ գագաթները տարբեր են $G_1, \ldots, G_t$-ի համար: Հետևաբար, $S$-ը պարունակում է առնվազն $t$ գագաթ, որը նշանակում է, որ $o(G - S) = t \leq |S|$:

Բավարարությունն ապացուցելու համար ցույց տանք, որ բոլոր $G$ գրաֆներում, որոնցում ցանկացած $S \subseteq V(G)$ համար տեղի ունի $o(G - S) \leq |S|$ անհավասարությունը, միշտ գոյություն ունի կատարյալ զուգակցում:

Ենթադրենք հակառակը, դիցուք գոյություն ունեն հակաօրինակներ, և դիցուք $G$-ն նրանցից մեկն է: Վերցնենք $S = \emptyset$: Նկատենք, որ
$$o(G) = o(G - S) \leq |S| = |\emptyset| = 0,$$
որը նշանակում է, որ $G$ գրաֆի կապակցվածության բոլոր բաղադրիչների գագաթների քանակները զույգ թվեր են, որտեղից, մասնավորապես, հետևում է, որ $G$ գրաֆի գագաթների քանակը ևս զույգ է:

Նկատենք, որ եթե $u$-ն և $v$-ն $G$ գրաֆում կամայական երկու ոչ հարևան գագաթներ են, ապա ցանկացած $S \subseteq V(G)$ համար $o(G + uv - S) \leq o(G - S) \leq |S|$: Այստեղից հետևում է, որ $G$ հակաօրինակին ավելացնելով կողեր մենք կարող ենք ստանալ ապացուցվող թեորեմի մաքսիմալ հակաօրինակ, այսինքն՝ այնպիսի հակաօրինակ, որը բավարարում է $o(G - S) \leq |S|$ պայմանին ցանկացած $S \subseteq V(G)$ համար, որում գոյություն չունի կատարյալ զուգակցում, բայց ցանկացած նոր կողի ավելացումից ստացվող գրաֆում արդեն գոյություն ունի կատարյալ զուգակցում:

Ստորև կապացուցենք, որ գոյություն չունեն թեորեմի մաքսիմալ հակաօրինակներ: Քանի որ ցանկացած հակաօրինակից կարելի է ստանալ մաքսիմալ հակաօրինակ, ապա այս պնդումից կստացվի, որ գոյություն չունեն նաև հակաօրինակներ, և, հետևաբար, թեորեմի պնդումը ճիշտ է:

Դիցուք $G$-ն մաքսիմալ հակաօրինակ է: $U$-ով նշանակենք $G$ գրաֆի այն գագաթների բազմությունը, որոնք հարևան են մնացած բոլոր գագաթներին: Հատուկ նշենք, որ $U$-ն կարող է լինել դատարկ: Ցույց տանք, որ հնարավոր չէ, որ $G - U$ գրաֆի կապակցվածության բոլոր բաղադրիչները լինեն լրիվ գրաֆներ:

Իրոք, եթե $G - U$ գրաֆի բոլոր կապակցվածության բաղադրիչները լրիվ գրաֆներ են, ապա, քանի որ $o(G - U) \leq |U|$, ապա մենք $G$ գրաֆում կկառուցենք կատարյալ զուգակցում հետևյալ եղանակով (նկ. 5.3.1).

1. $G - U$ գրաֆի բոլոր զույգ թվով գագաթներ պարունակող կապակցվածության բաղադրիչներում վերցնենք կատարյալ զուգակցում,
2. $G - U$ գրաֆի բոլոր կենտ թվով գագաթներ պարունակող կապակցվածության բաղադրիչներում վերցնենք մաքսիմում զուգակցում, նկատենք, որ այն չի հագեցնում ճիշտ մեկ գագաթ, այդ մեկ գագաթը վերցնենք $U$ բազմության գագաթներից մեկի հետ որպես կառուցվելիք կատարյալ զուգակցման կող, քանի որ $o(G - U) \leq |U|$, ապա $U$ բազմության գագաթները մենք միշտ կարող ենք ընտրել իրարից տարբեր,
3. և վերջապես $U$ բազմության մնացած $|U| - o(G - U)$ գագաթները վերցնենք զույգերով կամայապես նկատենք, քանի որ $G$ գրաֆի գագաթների քանակը զույգ է, ապա $|U| - o(G - U)$ թիվը ևս զույգ է, և, հետևաբար, մենք միշտ դա կարող ենք անել:

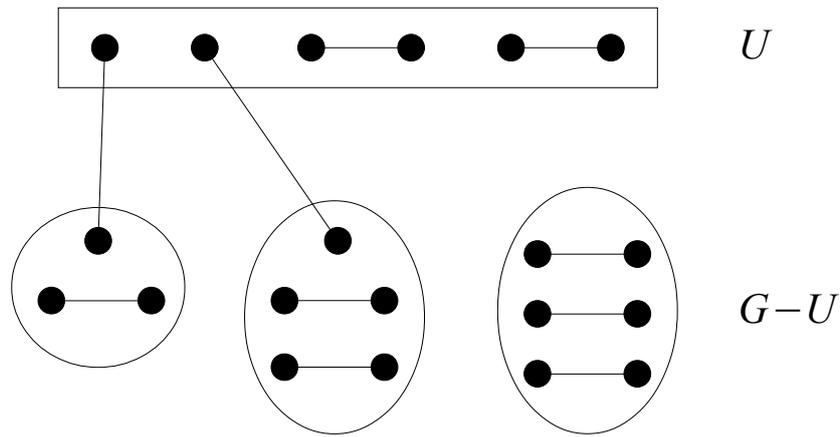

Նկ. 5.3.1

Ասվածից հետևում է, որ $G - U$ գրաֆի կապակցվածության ոչ բոլոր բաղադրիչներն են լրիվ գրաֆներ։ Սա նշանակում է, որ $G - U$ գրաֆի կապակցվածության ինչ-որ բաղադրիչում կան երկու ոչ հարևան $x$ և $y$ գագաթներ։

Ցույց տանք, որ մենք այդ ոչ հարևան գագաթները կարող ենք ընտրել այնպես, որ նրանք ունենան ընդհանուր հարևան։ Իրոք, դիտարկենք $G - U$ գրաֆի կապակցվածության բաղադրիչում $x$ և $y$ գագաթները միացնող ճանապարհը։ Դիցուք այդ ճանապարհը $x = z_0, z_1, \ldots, z_k = y$-ն է։ Քանի որ $x$ և $y$ գագաթները հարևան չեն, ապա $k \geq 2$։

Պնդումն ապացուցենք մակածման եղանակով ըստ $k$-ի։ Եթե $k = 2$, ապա $z_1$ գագաթը ոչ հարևան $x$ և $y$ գագաթների ընդհանուր հարևան է։ Ենթադրենք, որ պնդումն արդեն ապացուցել ենք $k - 1$ համար, և դիտարկենք $z_1$ գագաթը։ Եթե $z_1$ գագաթը հարևան է $y$ գագաթին, ապա կրկին $z_1$ գագաթը կլինի $x$ և $y$ գագաթների ընդհանուր հարևանը, իսկ եթե $z_1$ գագաթը հարևան չէ $y$ գագաթին, ապա նկատենք, որ ոչ հարևան $z_1$ և $y$ գագաթները միացված են $k - 1$ երկարությամբ ճանապարհով, որտեղից, համաձայն մակածման ենթադրության, կարող ենք ընտրել երկու ոչ հարևան գագաթ, որոնք ունեն ընդհանուր հարևան։

Ասվածից հետևում է, որ առանց ընդհանրությունը խախտելու մենք կարող ենք ենթադրել, որ $G - U$ գրաֆի ոչ հարևան $x$ և $y$ գագաթներն ունեն ընդհանուր հարևան $z$ գագաթ։ Քանի որ $z \notin U$, ապա, համաձայն $U$ բազմության սահմանման, $G$ գրաֆում գոյություն ունի այնպիսի $w$ գագաթ, որ $z$ և $w$ գագաթները հարևան չեն։

Հիշենք, որ $G$ հակաօրինակը մաքսիմալ էր։ Սա, մասնավորապես, նշանակում է, որ ցանկացած նոր կող ավելացնելուց ստացված գրաֆն արդեն ունի կատարյալ զույգակցում։

Դիտարկենք $G + xy$ և $G + zw$ գրաֆները, և դիցուք $M_1$-ը և $M_2$-ը, համապատասխանաբար, այդ գրաֆներում կատարյալ զուգակցումներ են: Թեորեմի ապացույցն ավարտելու համար, մենք, օգտվելով $M_1$ և $M_2$ կատարյալ զուգակցումներից, ցույց կտանք, որ $G + xy + zw$ գրաֆում գոյություն ունի կատարյալ զուգակցում, որը չի պարունակում $xy$ և $zw$ կողերը: Նկատենք, որ սա հակասություն է, քանի որ $G + xy + zw$ գրաֆի $xy$ և $zw$ կողերը չպարունակող կատարյալ զուգակցումը կլինի $G$ գրաֆի կատարյալ զուգակցում, ինչը կհակասի այն բանին, որ $G$-ն հակաօրինակ է և, հետևաբար, չի պարունակում կատարյալ զուգակցում:

Նախ նկատենք, որ $G + xy + zw$ գրաֆի ցանկացած գագաթ կամ կից է $M_1 \cap M_2$-ին պատկանող կողի, կամ մեկ կողի $M_1 \setminus M_2$-ից և մեկ կողի $M_2 \setminus M_1$-ից: Ասվածից հետևում է, որ $M_1 \cup M_2$-ին պատկանող կողերը կամ պատկանում են $M_1 \cap M_2$-ին կամ կազմում են ցիկլեր, ընդ որում վերջիններս ունեն զույգ երկարություն, քանի որ այդ ցիկլերի կողերը մեկընդմեջ պատկանում են $M_1 \setminus M_2$-ին և $M_2 \setminus M_1$-ին: Ավելին, նկատենք, որ $xy \in M_1 \setminus M_2$ և $zw \in M_2 \setminus M_1$, և, հետևաբար, $xy$ և $zw$ կողերը պատկանում են $M_1 \cup M_2$-ի զույգ ցիկլերին:

Թեորեմի ապացույցն ավարտելու համար քննարկենք երկու դեպք:

Դեպք 1: $xy$ և $zw$ կողերը պատկանում են $M_1 \cup M_2$-ի տարբեր ցիկլերին:

Այս դեպքում $G + xy + zw$ գրաֆի $xy$ և $zw$ կողերը չպարունակող կատարյալ զուգակցումը կարելի է ստանալ հետևյալ կերպ. վերցնենք $M_1 \cap M_2$-ին պատկանող կողերը և դրանց ավելացնենք $M_1 \cup M_2$-ի զույգ ցիկլերի կատարյալ զուգակցումները: Նկատենք, քանի որ զույգ ցիկլերը ունեն երկու կատարյալ զուգակցում, մենք միշտ կարող ենք խուսափել $xy$ և $zw$ կողերը վերցնելուց:

Դեպք 2: $xy$ և $zw$ կողերը պատկանում են $M_1 \cup M_2$-ի միևնույն $C$ ցիկլին:

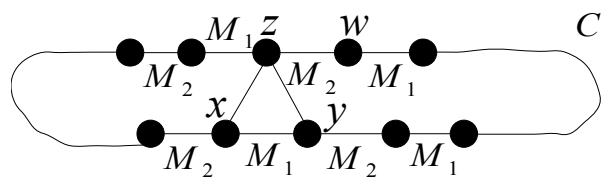

Նկ. 5.3.2

Մենք կենթադրենք, որ $z$ գագաթից $M_2$-ին պատկանող կողից սկսելով և $C$ ցիկլով շարժվելով, մենք առաջինը հանդիպում ենք $y$ գագաթին, որից հետո նոր $x$ գագաթին (նկ. 5.3.2): Նշենք, որ մյուս դեպքը քննարկվում է համանման ձևով:

Այս դեպքում $G + xy + zw$ գրաֆի $xy$ և $zw$ կողերը չպարունակող կատարյալ զուգակցումը կարելի ստանալ հետևյալ կերպ. վերցնենք $M_1 \cap M_2$-ին պատկանող կողերը և դրանց ավելացնենք $M_1 \cup M_2$-ի բոլոր զույգ ցիկլերի բացի $C$-ից կատարյալ զուգակցումները, $C$ ցիկլի ձախ մասում վերցնենք $M_1$-ին պատկանող կողերը, աջ մասում՝ $M_2$-ին պատկանող կողերը և ավելացնենք $yz$ կողը (նկ. 5.3.2)։ ∎

Օգտվելով Տատտի թեորեմից, ձևակերպենք և ապացուցենք Պետերսենի թեորեմը, որը տալիս է խորանարդ գրաֆում կատարյալ զուգակցման գոյության բավարար պայման։

**Թեորեմ 5.3.2 (Պետերսեն)։** Դիցուք $G$-ն խորանարդ գրաֆ է, որը պարունակում է ոչ ավելի, քան երկու կամուրջ։ Այդ դեպքում նրանում գոյություն ունի կատարյալ զուգակցում։

**Ապացույց։** Դիտարկենք ցանկացած $S \subseteq V(G)$, և դիցուք $o(G - S) = t$, և $G_1, \ldots, G_t$-ն $G - S$ գրաֆի կապակցվածության այն բաղադրիչներն են, որոնք պարունակում են կենտ թվով գագաթներ։ Նկատենք, որ համաձայն թեորեմ 5.3.1-ի, պնդումն ապացուցելու համար բավական է ցույց տալ, որ $t \leq |S|$։

Նկատենք, որ $i = 1, \ldots, t$-ի համար տեղի ունի
$$3|V(G_i)| = 2|E(G_i)| + |\partial_G(V(G_i))|$$
հավասարությունը, որտեղից հաշվի առնելով, որ $|V(G_i)|$-ն կենտ է, կստանանք, որ կենտ է նաև $|\partial_G(V(G_i))|$-ն։ Քանի որ ըստ ենթադրության $G$ գրաֆում գոյություն ունի ոչ ավելի, քան երկու կամուրջ, մենք կստանանք, որ $i = 1, \ldots, t$ համար $|\partial_G(V(G_i))|$ թվերից ամենաշատը երկուսն են հավասար մեկի, իսկ մնացածը առնվազն երեք են։

Դիտարկենք $G$ գրաֆի այն կողերը, որոնք միացնում են $S$ բազմության գագաթները $G - S$ գրաֆի կենտ թվով գագաթներ պարունակող $G_1, \ldots, G_t$ կապակցվածության բաղադրիչներին։ Նկատենք, որ համաձայն վերը ասվածի, այդ կողերի քանակը առնվազն $3(t - 2) + 2$-է։ Մյուս կողմից, քանի որ $G$-ն խորանարդ գրաֆ է, ապա $S$ բազմության գագաթների աստիճանն երեք է, և, հետևաբար, $S$ բազմության գագաթները կից են ոչ ավելի, քան $3|S|$ կողի։ Արդյունքում՝
$$3t - 4 = 3(t - 2) + 2 \leq 3|S|$$
կամ
$$t \leq |S| + 1 :$$

Հաշվի առնելով, որ խորանարդ գրաֆները պարունակում են զույգ թվով գագաթներ, մենք կստանանք, որ $t$-ն և $|S|$-ը ունեն նույն զույգությունը, հետևաբար, $t = |S| + 1$ հավասարությունը հնարավոր չէ։ Այստեղից հետևում է, որ $t \leq |S|$։ ∎

Հատուկ նշենք, որ Պետերսենի վերը նշված թեորեմում խորանարդ գրաֆի ոչ ավելի, քան երկու կամուրջ չպարունակելու պայմանը հնարավոր չէ թուլացնել։

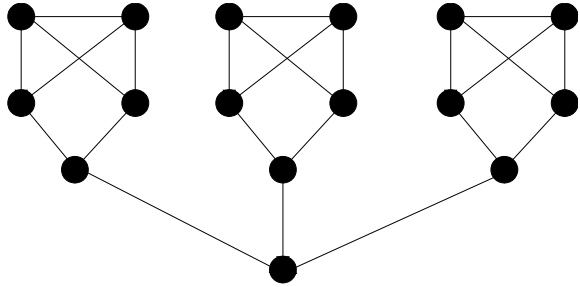

Նկ. 5.3.3

Նկ. 5.3.3-ում բերված է երեք կամուրջ պարունակող *Սիլվեստրի գրաֆի* օրինակը, որը չունի կատարյալ զուգակցում։

Ստորև կձևակերպենք և կապացուցենք Տատտ-Բերժի բանաձևն, որը կարելի է մեկնաբանել նաև որպես մինմաքս թեորեմ։ Ինչպես արդեն նշել ենք § 5.2-ում, $G$ գրաֆի ցանկացած $M$ զուգակցում հագեցնում է $2|M|$ գագաթ, և, հետևաբար, այն չի հագեցնում $|V| - 2|M|$ գագաթ։ Ասվածից հետևում է, որ ցանկացած $G$ գրաֆում գագաթների նվազագույն թիվը, որոնք չեն ծածկվում զուգակցումով, հավասար է $|V| - 2\alpha'(G)$-ի։ Տատտ-Բերժի բանաձևն պնդում է, որ այդ մինիմումը հավասար է ստորև բերված մաքսիմումին։

**Թեորեմ 5.3.3 (Տատտ-Բերժ):** Կամայական $G$ գրաֆում տեղի ունի հետևյալ հավասարությունը.

$$|V| - 2\alpha'(G) = \max_{S \subseteq V(G)}(o(G - S) - |S|)։$$

**Ապացույց:** Նախ նկատենք, որ եթե $M$-ը $G$ գրաֆի ցանկացած զուգակցում է, և $S \subseteq V(G)$, ապա $M$-ին պատկանող ամենաշատը $|S|$ կող կարող է միացնել $S$-ին պատկանող գագաթը $G - S$ գրաֆի կենտ թվով գագաթներ պարունակող կապակցվածության բաղադրիչին պատկանող գագաթին, հետևաբար $G - S$ գրաֆի առնվազն $o(G - S) - |S|$ հատ կենտ թվով գագաթներ պարունակող կապակցվածության բաղադրիչում գոյություն կունենան $M$-ով չհագեցած գագաթներ, որը, նշանակում է, որ $G$ գրաֆում $M$-ով չհագեցած գագաթների քանակն առնվազն $o(G - S) - |S|$-է։ Նկատենք, քանի որ $M$-ը և $S$-ն ընտրված էին կամայապես, ասվածից հետևում է, որ

$$|V| - 2\alpha'(G) \geq \max_{S \subseteq V(G)}(o(G - S) - |S|)։$$

Հակառակ անհավասարությունն ապացուցելու համար նշանակենք

$$d = \max_{S \subseteq V(G)} (o(G - S) - |S|)$$

և ցույց տանք, որ $G$ գրաֆում գոյություն ունի զուգակցում, որը չի հագեցնում $G$ գրաֆի ոչ ավելի, քան $d$ գագաթ։

Դիտարկենք $G'$ գրաֆը, որն իրենից ներկայացնում է $G$ գրաֆի և $K_d$ գրաֆի ($d$ գագաթ պարունակող լրիվ գրաֆի) գումարը (§ 1.3)։ Ցույց տանք, որ $G'$ գրաֆը պարունակում է $M'$ կատարյալ զուգակցում։ Նկատենք, որ սա կապացուցի թեորեմը, քանի որ եթե մենք $G'$ գրաֆից հեռացնենք $V(K_d)$-ն, ապա $M'$-ի մնացած կողերը կկազմեն $G$ գրաֆի զուգակցում, որը չի հագեցնում $G$ գրաֆի ոչ ավելի, քան $d$ գագաթ։

Դիտարկենք ցանկացած $S' \subseteq V(G')$, և քննարկենք հետևյալ երեք դեպքերը։

Դեպք 1: $S' = \emptyset$:

Նկատենք, որ ցանկացած $S \subseteq V(G)$ համար $|V(G)|$-ն և $o(G-S) - |S|$ տարբերությունն ունեն միևնույն զույգությունը, հետևաբար, $|V(G)|$-ի և $d$-ի զույգությունը նույնպես համընկնում է, որտեղից հետևում է, որ $|V(G')|$-ը զույգ թիվ է։ Հաշվի առնելով, որ $G'$-ը կապակցված գրաֆ է, կստանանք՝

$$o(G' - S') = o(G') = 0 = |S'|:$$

Դեպք 2: $S' \neq \emptyset$ և $V(K_d) \not\subseteq S'$:

Նկատենք, որ այս դեպքում $G' - S'$ գրաֆը կապակցված է, հետևաբար

$$o(G' - S') \leq 1 \leq |S'|:$$

Դեպք 3: $V(K_d) \subseteq S'$:

Նշանակենք $S = S' \setminus V(K_d)$, և նկատենք, որ $S \subseteq V(G)$։ $d$-ի սահմանումից հետևում է, որ

$$o(G - S) - |S| \leq d,$$

որտեղից հաշվի առնելով, որ $o(G' - S') = o(G - S)$, կստանանք՝

$$o(G' - S') = o(G - S) \leq d + |S| = |S'|:$$

Քննարկված երեք դեպքերի արդյունքում տեսանք, որ ցանկացած $S' \subseteq V(G')$, համար $o(G' - S') \leq |S'|$։ Համաձայն թեորեմ 5.3.1-ի, $G'$ գրաֆը պարունակում է կատարյալ զուգակցում։ ∎

Ստորև կձևակերպենք և կապացուցենք § 5.1-ում ապացուցված Գալլայի հավասարության կողային տարբերակը։ Այդ նպատակով տանք մեկ սահմանում։

Դիցուք $G = (V, E)$-ն գրաֆ է և $L \subseteq E$:

**Սահմանում 5.3.1:** Կասենք, որ $L$-ը հանդիսանում է *կողային ծածկույթ* $G$ գրաֆում, եթե $G$ գրաֆի ցանկացած զագաթ կից է $L$-ին պատկանող գոնե մեկ կողի:

**Սահմանում 5.3.2:** $G$ գրաֆում ամենաքիչ կողեր պարունակող կողային ծածկույթներին կանվանենք *նվազագույն կողային ծածկույթներ*:

$G$ գրաֆում նվազագույն կողային ծածկույթի հզորությունը նշանակենք $\beta'(G)$-ով: Նկատենք, որ գրաֆում գոյություն ունի կողային ծածկույթ այն և միայն այն դեպքում, երբ գրաֆում չկան մեկուսացված զագաթներ: Սա, մասնավորապես, նշանակում է, որ $\beta'(G)$ պարամետրը սահմանված է միայն այն $G$ գրաֆների համար, որոնցում ցանկացած զագաթի աստիճանն առնվազն մեկ է:

Կրկին դիտարկենք նկ. 5.1.1-ում պատկերված $G$ գրաֆը, և նկատենք, որ նրանում կողերի $\{v_1v_2, v_3v_4, v_5v_6, v_5v_7\}$ բազմությունը հանդիսանում է կողային ծածկույթ, ավելին, դժվար չէ համոզվել, որ այդ գրաֆում $\beta'(G) = 4$:

**Թեորեմ 5.3.4 (Գալլաի):** Մեկուսացված զագաթներ չպարունակող ցանկացած $G$ գրաֆում տեղի ունի $\alpha'(G) + \beta'(G) = |V|$ հավասարությունը:

**Ապացույց:** Դիցուք $M$-ը $G$ գրաֆի առավելագույն զուգակցում է: Նկատենք, որ այն ծածկում է $2|M| = 2\alpha'(G)$ զագաթ: Վերցնենք $G$ գրաֆի մնացած $|V| - 2\alpha'(G)$ զագաթներին կից կողեր (յուրաքանչյուր զագաթին մեկ կից կող), և դիտարկենք $G$ գրաֆի կողերի $N$ բազմությունը, որը ստացվում է այդ կողերին միավորելով $M$-ը: Նկատենք, որ $N$-ը կողային ծածկույթ է, և

$$\beta'(G) \leq |N| = \alpha'(G) + |V| - 2\alpha'(G) = |V| - \alpha'(G),$$

կամ

$$\beta'(G) + \alpha'(G) \leq |V|:$$

Մյուս կողմից, դիցուք $L$-ը $G$ գրաֆի $\beta'(G)$ հատ կող պարունակող կողային ծածկույթ է: Նկատենք, որ եթե $uv \in L$, ապա հնարավոր չէ, որ $u$-ն և $v$-ն միաժամանակ կից լինեն $L \setminus \{uv\}$ բազմությանը պատկանող կողերի: Իրոք, հակառակ դեպքում $L \setminus \{uv\}$-ն կհանդիսանար $G$ գրաֆի $\beta'(G) - 1$ հատ կող պարունակող կողային ծածկույթ, ինչը կհակասեր $\beta'(G)$-ի սահմանմանը: Այստեղից հետևում է, որ եթե դիտարկենք $G$ գրաֆի $H = (V, L)$ ենթագրաֆը, ապա $H$-ի կապակցվածության բաղադրիչներն իրենցից ներկայացնում են աստղեր (§ 1.2): Եթե $H$-ի կապակցվածության բաղադրիչների քանակը նշանակենք $k$-ով, ապա, համաձայն թեորեմ 2.3.1-ի,

$$\beta'(G) = |L| = |V| - k:$$

Դիտարկենք $G$ գրաֆի $M'$ զուգակցումը, որը ստացվում է, եթե $H$-ի կապակցվածության բաղադրիչներից յուրաքանչյուրից վերցնենք մեկական կող։ Նկատենք, որ

$$\alpha'(G) \geq |M'| = k = |V| - \beta'(G),$$

կամ

$$\beta'(G) + \alpha'(G) \geq |V|:$$

$\beta'(G) + \alpha'(G) \leq |V|$ և $\beta'(G) + \alpha'(G) \geq |V|$ անհավասարություններից հետևում է, որ

$$\beta'(G) + \alpha'(G) = |V|: \blacksquare$$

## § 5.4. Ֆակտորներ և տարբեր ֆակտորիզացիաներ

Դիցուք $G = (V, E)$-ն գրաֆ է, և $H$-ը նրա ցանկացած ենթագրաֆ է։

**Սահմանում 5.4.1:** Կասենք, որ $H$-ը հանդիսանում է $G$ գրաֆի *ֆակտոր*, եթե այն $G$ գրաֆի կմախքային ենթագրաֆ է։

Դիցուք $f$-ը արտապատկերում է, որը բավարարում է $f: V(G) \to \mathbb{Z}_+$ պայմանին։

**Սահմանում 5.4.2:** Կասենք, որ $G$ գրաֆի $H$ ֆակտորը հանդիսանում է $G$ գրաֆի $f$-*ֆակտոր*, եթե ցանկացած $v \in V(G)$ համար $d_H(v) = f(v)$։

Դիցուք $k \in \mathbb{Z}_+$։

**Սահմանում 5.4.3:** Եթե ցանկացած $v \in V(G)$ համար $f(v) = k$, ապա $G$ գրաֆի $f$-ֆակտորին կանվանենք $k$-*ֆակտոր*։

Փաստորեն, գրաֆի $k$-ֆակտորները նրա կմախքային $k$-համասեռ ենթագրաֆներն են։ Դիտարկենք բերված սահմանումները պարզաբանող օրինակներ։

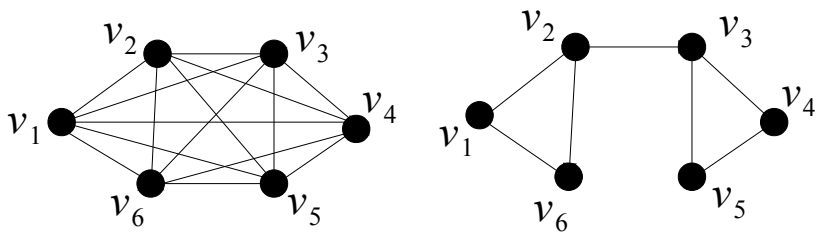

Նկ. 5.4.1

Նկ. 5.4.1-ի ձախ մասում պատկերված է $K_6$ վեց գագաթ պարունակող լրիվ գրաֆը, իսկ աջ մասում՝ նրա $H$ ենթագրաֆը, որը հանդիսանում է $f$-ֆակտոր, որտեղ $f(v_1) =$ $f(v_4) = f(v_5) = f(v_6) = 2$, $f(v_2) = f(v_3) = 3$։ Ավելին, ստորև բերված նկ. 5.4.2-ում պատկերված են $K_6$ գրաֆի $F_1$ և $F_2$, համապատասխանաբար, 1- և 2-ֆակտորներ:

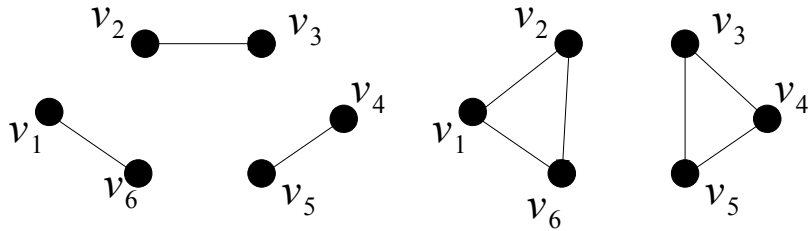

Նկ. 5.4.2

Վերջապես, նկատենք, որ եթե դիտարկենք $F_1$-ի և $F_2$-ի լրացում գրաֆները $K_6$-ում, ապա նրանք կհանդիսանան $K_6$ գրաֆի, համապատասխանաբար, 4- և 3-ֆակտորներ:

Գրաֆի 1-ֆակտորն իրենից ներկայացնում է կմախքային 1-համասեռ ենթագրաֆ: Դժվար չէ համոզվել, որ այդ ենթագրաֆի կողերի բազմությունն իրենից ներկայացնում է կատարյալ զուգակցում: Ավելին, ճիշտ է նաև հակառակը. եթե ունենք գրաֆի կատարյալ զուգակցում, ապա եթե դիտարկենք գրաֆի կմախքային ենթագրաֆը, որի կողերի բազմությունը այդ տրված կատարյալ զուգակցումն է, ապա մենք կստանանք գրաֆի 1-ֆակտոր: Ասվածից հետևում է, որ թեորեմ 5.3.1-ը կարելի է նաև դիտել որպես 1-ֆակտոր պարունակող գրաֆների նկարագիր:

Այստեղից առաջանում է բնական հարց. հնարավո՞ր է առաջարկել եղանակ, որը կպարզի, թե տրված գրաֆը պարունակում է արդյոք $f$-ֆակտոր: Ստորև կնկարագրենք Տատտի ալգորիթմը, որը թույլ է տալիս տրված գրաֆի $f$-ֆակտոր պարունակելու հարցը հանգեցնել մեկ այլ գրաֆի 1-ֆակտոր պարունակելու հարցին:

Դիցուք տրված է $G$ գրաֆը և $f:V(G) \to \mathbb{Z}_+$ արտապատկերումը: Նախ նկատենք, որ եթե $G$ գրաֆում գոյություն ունի $v$ գագաթ, որի համար $f(v) > d_G(v)$, ապա $G$-ն չի կարող ունենալ $f$-ֆակտոր, հետևաբար, առանց ընդհանրությունը խախտելու, կարող ենք ենթադրել, որ $G$ գրաֆի ցանկացած $v$ գագաթի համար $f(v) \le d_G(v)$: Նշանակենք $e(v) = d_G(v) - f(v)$ և նկատենք, որ $e(v) \ge 0$: Դիտարկենք $H$ գրաֆը, որը ստացվում է $G$-ից նրա ցանկացած $v$ գագաթ փոխարինելով $K_{d_G(v),e(v)}$ -լրիվ երկկողմանի գրաֆով (§ 1.2), որի մի կողմը՝ $A(v)$-ն, պարունակում է $d_G(v)$ գագաթ, իսկ մյուս կողմը՝ $B(v)$-ն, $e(v)$ գագաթ, և $G$ գրաֆի ցանկացած $vw$ կողի համար $A(v)$-ին պատկանող մեկ գագաթ միացնենք կողով $A(w)$-ին պատկանող մեկ գագաթի հետ այնպես, որ $A$-գագաթները մասնակցեն ճիշտ մեկ այդպիսի կողի մեջ:

Ստորև բերված նկ. 5.4.3-ի ձախ մասում պատկերված է $G$ գրաֆը և նրա գագաթների բազմության վրա որոշված $f: V(G) \to \mathbb{Z}_+$ արտապատկերումը, իսկ աջ մասում պատկերված է $G$ գրաֆին համապատասխանող $H$ գրաֆը և կետագծերով ցույց է տրված $H$-ի այն **1**-ֆակտորը, որին համապատասխանում է $G$ գրաֆի $f$-ֆակտորը:

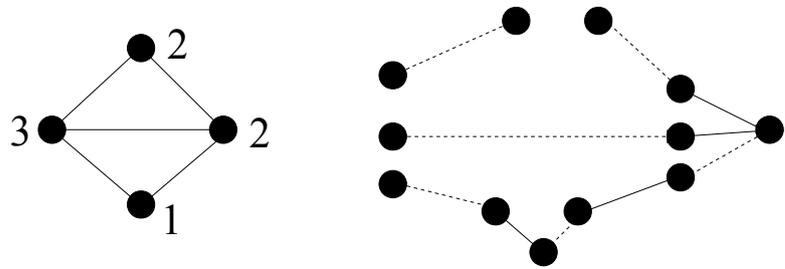

Նկ. 5.4.3

**Թեորեմ 5.4.1 (Տատտ):** $G$ գրաֆը պարունակում է $f$-ֆակտոր այն և միայն այն դեպքում, երբ վերը նշված եղանակով կառուցված $H$ գրաֆը պարունակում է **1**-ֆակտոր:

**Ապացույց:** Եթե $G$ գրաֆը պարունակում է $f$-ֆակտոր, ապա $H$ գրաֆի այդ $f$-ֆակտորին համապատասխանող կողերը կազմում են զուգակցում, որը $G$ գրաֆի ցանկացած $v$ գագաթի համար $A(v)$-ից չի հագեցնում ճիշտ $e(v)$ գագաթ: Յուրաքանչյուր $v$ գագաթի համար այդ զուգակցմանն ավելացնենք $K_{d_G(v), e(v)}$-լրիվ երկկողմանի գրաֆի զուգակցումը, որը հագեցնում է այդ չհագեցած գագաթները և $B(v)$-ն: Արդյունքում կստանանք $H$ գրաֆի **1**-ֆակտոր:

Հակառակը, եթե $H$ գրաֆն ունի **1**-ֆակտոր, ապա $G$ գրաֆի յուրաքանչյուր $v$ գագաթի համար նրանից հեռացնելով $B(v)$-ին պատկանող գագաթներին կից կողերը, մենք կստանանք $H$ գրաֆի կողերի բազմություն, որոնց, դժվար չէ տեսնել, համապատասխանում է $G$ գրաֆի $f$-ֆակտոր: ∎

Դիցուք $G = (V, E)$-ն գրաֆ է, և $f$-ը արտապատկերում է, որը բավարարում է $f: V(G) \to \mathbb{Z}_+$ պայմանին:

**Սահմանում 5.4.4:** $G$ գրաֆը կանվանենք $f$-*ֆակտորիզացվող*, եթե այն հնարավոր է տրոհել զույգ առ զույգ չհատվող $f$-ֆակտորների:

Նշենք, որ գրաֆի $f$-ֆակտորների տրոհման բուն պրոցեսը հաճախ անվանում են $f$-*ֆակտորիզացիա* կամ, կրճատ, *ֆակտորիզացիա*:

Դիցուք $k \in \mathbb{Z}_+$: Այն դեպքում, երբ $G$ գրաֆի ցանկացած $v \in V(G)$ համար $f(v) = k$, $f$-ֆակտորիզացիան մենք կանվանենք $k$-*ֆակտորիզացիա*:

Ստորև պատկերված է չորս գագաթանի $K_4$ լրիվ գրաֆը և բերված է նրա **1**-ֆակտորիզացիայի օրինակ:

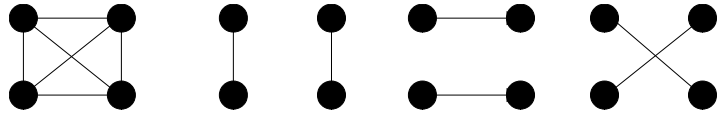

Նկ. 5.4.4

Ստորև կնշենք **1**-ֆակտորիզացվող գրաֆների դասի օրինակ:

**Թեորեմ 5.4.2:** Կամայական երկկողմանի $r$-համասեռ ($r \in \mathbb{N}$) գրաֆ **1**-ֆակտորիզացվող է:

**Ապացույց:** Դիցուք $G$-ն երկկողմանի $r$-համասեռ գրաֆ է: Համաձայն թեորեմ 5.2.4-ի, այն պարունակում է $F_1$ կատարյալ զուգակցում: Դիտարկենք $G - F_1$ գրաֆը: Նկատենք, որ այն երկկողմանի $(r-1)$-համասեռ գրաֆ է: Համաձայն թեորեմ 5.2.4-ի, այն պարունակում է $F_2$ կատարյալ զուգակցում: Դիտարկենք $G - F_1 - F_2$ գրաֆը: Նկատենք, որ այն երկկողմանի $(r-2)$-համասեռ գրաֆ է: Նշված քայլերը կիրառելով $r$ անգամ, մենք կստանանք $G$ երկկողմանի $r$-համասեռ գրաֆի կողերի բազմության տրոհում զույգ առ զույգ չհատվող կատարյալ զուգակցումների, որոնցից դժվար չէ ստանալ $G$ գրաֆի **1**-ֆակտորիզացիա: ∎

Նշենք, որ գրաֆների **1**-ֆակտորիզացիաներին առնչվող ոչ բոլոր խնդիրներն են լուծված: Մասնավորապես, չի լուծված հանրահայտ **1**-ֆակտորիզացիայի հիպոթեզը, որը ձևակերպված է ստորև:

**Հիպոթեզ 5.4.1:** Դիցուք $G$-ն $2n$ գագաթ պարունակող $r$-համասեռ գրաֆ է: Այդ դեպքում. եթե $n$-ը կենտ է և $r \geq n$ կամ $n$-ը զույգ է և $r \geq n-1$, ապա $G$-ն **1**-ֆակտորիզացվող է:

Հիմա կձևակերպենք և կապացուցենք Պետերսենի թեորեմը, որը նկարագրում է **2**-ֆակտորիզացվող գրաֆները:

**Թեորեմ 5.4.3:** Որպեսզի $G$ գրաֆը լինի **2**-ֆակտորիզացվող, անհրաժեշտ է և բավարար, որ այն լինի զույգ համասեռ:

**Ապացույց:** Դիցուք $G$ գրաֆը **2**-ֆակտորիզացվող է: Այդ դեպքում այն հնարավոր է տրոհել զույգ առ զույգ չհատվող $H_1, \ldots, H_k$ **2**-ֆակտորների: Վերցնենք $G$ գրաֆի ցանկացած $v$ գագաթ: Նկատենք, որ $i = 1, 2, \ldots, k$ համար $d_{H_i}(v) = 2$, և հետևաբար $d_G(v) = 2k$: Սա նշանակում է, որ $G$ գրաֆը $2k$-համասեռ է:

Հակառակն ապացուցելու համար ենթադրենք, որ $G$ գրաֆը $2k$-համասեռ է, $k \in \mathbb{N}$:

Ցույց տանք, որ այն հնարավոր է տրոհել $k$ հատ 2-ֆակտորների:

Նկատենք, որ պնդումն ապացուցելու համար բավական է ցույց տալ, որ $G$ գրաֆն ունի 2-ֆակտոր: Իրոք, եթե $G$ գրաֆից հեռացնենք այդ 2-ֆակտորի կողերը, ապա կստանանք $(2k-2)$-համասեռ գրաֆ: Ստացված գրաֆը կրկին կպարունակի 2-ֆակտոր: Այն կրկին կհեռացնենք: Նկարագրված քայլերը $k$ անգամ կրկնելուց հետո, մենք կստանանք $G$ գրաֆի 2-ֆակտորիզացիա:

Հետևաբար, բավական է ցույց տալ, որ $2k$-համասեռ $G$ գրաֆն ունի 2-ֆակտոր: Այս պնդումն ապացուցելու համար բավական է ապացուցել այն կապակցված գրաֆների համար: Իրոք, վերցնելով $G$ գրաֆի կապակցվածության յուրաքանչյուր բաղադրիչում մեկական 2-ֆակտոր, մենք կստանանք $G$ գրաֆի 2-ֆակտոր:

Հետևաբար, կարող ենք նաև ենթադրել, որ $G$ գրաֆը կապակցված է: Դիցուք, $G$ գրաֆի գագաթներն են $v_1, v_2, \ldots, v_n$-ը: Քանի որ $G$ գրաֆը կապակցված է և ցանկացած գագաթի աստիճանը զույգ թիվ է, ապա, համաձայն թեորեմ 4.1.2-ի, $G$ գրաֆում գոյություն ունի $C$ Էյլերյան ցիկլ: Քանի որ $G$-ում ցանկացած գագաթի աստիճանը $2k$ է, մենք, $C$ ցիկլով շարժվելուց, յուրաքանչյուր գագաթ կհանդիպենք $k$ անգամ, և ամեն անգամ տրված գագաթը հանդիպելուց մենք մեկ անգամ մուտք ենք գործում գագաթ, և մեկ անգամ՝ նրանից դուրս գալիս:

Դիտարկենք օժանդակ $H$ գրաֆը, որտեղ $V(H) = \{u_1, u_2, \ldots, u_n, w_1, w_2, \ldots, w_n\}$ և $u_i w_j \in E(H)$ այն և միայն այն դեպքում, երբ $C$ ցիկլով շարժվելուց մենք $v_i$ գագաթից անմիջապես հետո անցնում ենք $v_j$ գագաթ:

Նկատենք, որ ըստ սահմանման, $H$ գրաֆն երկկողմանի է, ավելին, քանի որ $G$ գրաֆի յուրաքանչյուր գագաթ մենք մուտք ենք գործում $k$ անգամ, և նրանից դուրս գալիս $k$ անգամ, ապա $H$ գրաֆը $k$-համասեռ է: Համաձայն թեորեմ 5.2.4-ի, $H$ գրաֆը պարունակում է $F$ կատարյալ զուգակցում:

Ցույց տանք, որ $H$ գրաֆի $F$ կատարյալ զուգակցումից մենք կարող ենք կառուցել $G$ գրաֆի 2-ֆակտոր: Դրա համար բավական է ցույց տալ, որ $F$-ի միջոցով մենք $G$ գրաֆի ցանկացած $v_i$ գագաթին կարող ենք համապատասխանեցնել նրան կից երկու կող:

Դիտարկենք $v_i$ գագաթին համապատասխան $H$ գրաֆի $u_i$ և $w_i$ գագաթները: Քանի որ $F$-ը կատարյալ զուգակցում է, ապա այն կպարունակի $u_i$ գագաթին կից $u_i w_j$ կող, և $w_i$ գագաթին կից $u_k w_i$ կող: Համաձայն $H$ գրաֆի սահմանման, սա նշանակում է, որ $C$ ցիկլով շարժվելուց մենք $v_k$ գագաթից անմիջապես հետո անցնում ենք $v_i$ գագաթ, որից հետո անցնում ենք $v_j$ գագաթ: $G$ գրաֆի $F$ կատարյալ զուգակցմանը համապատասխանող 2-ֆակտորը կառուցելիս վերցնենք $v_i$ գագաթին կից $v_kv_i$ և $v_iv_j$ կողերը: ∎

Դիտարկենք վերջին թեորեմի ապացույցը պարզաբանող օրինակ:

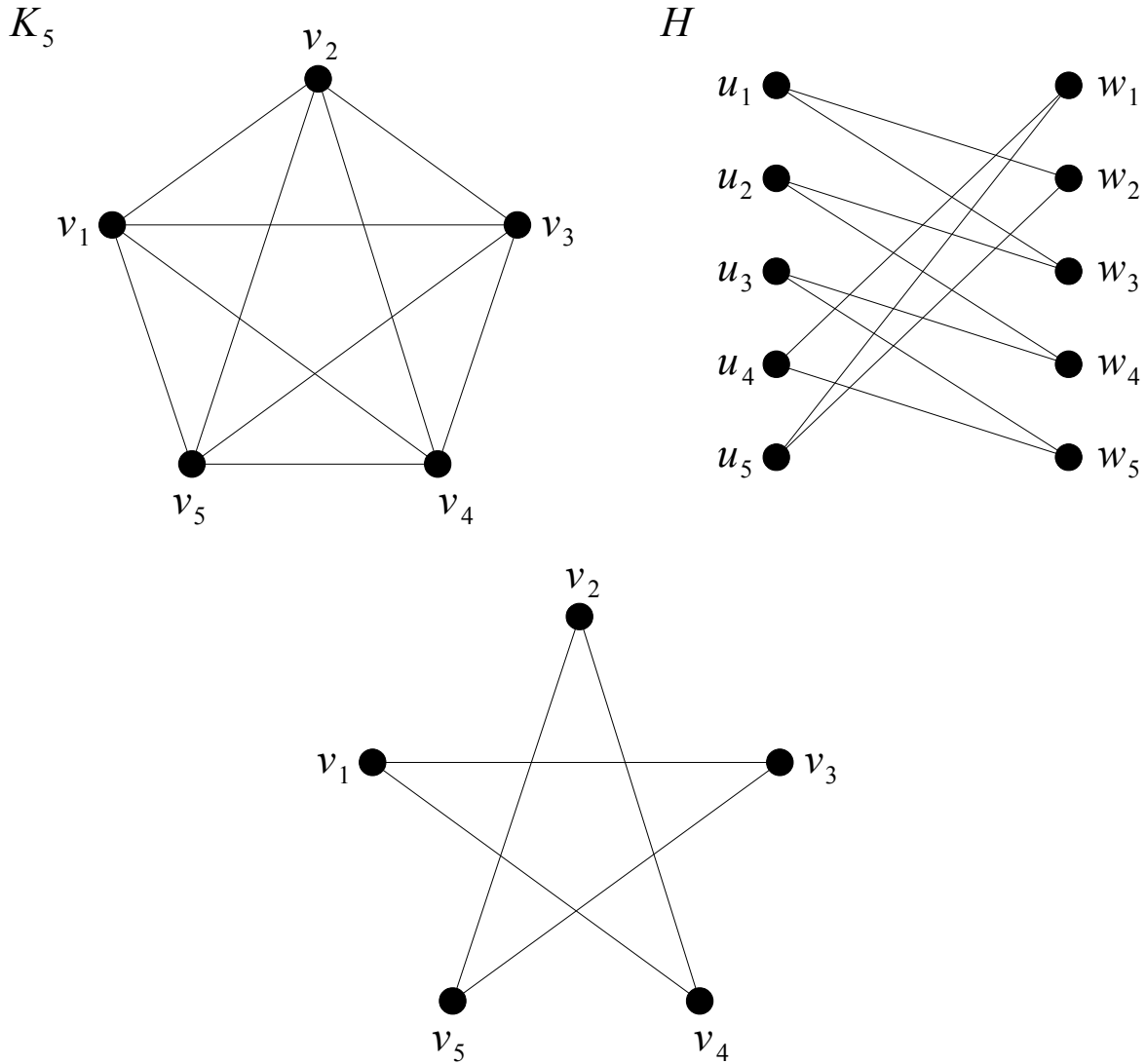

Նկ. 5.4.5

Նկ. 5.4.5-ում պատկերված է $K_5$ լրիվ գրաֆը, որը 4-համասեռ գրաֆ է, և $K_5$ գրաֆի $C = v_1, v_2, v_3, v_4, v_5, v_1, v_3, v_5, v_2, v_4, v_1$ էյլերյան ցիկլին համապատասխան օժանդակ $H$ երկկողմանի գրաֆը: Նաև նկ. 5.4.5-ում պատկերված է $G$ գրաֆի այն 2-ֆակտորը, որը համապատասխանում է $H$ գրաֆի $F = \{u_4w_1, u_5w_2, u_1w_3, u_2w_4, u_3w_5\}$ կատարյալ զուգակցմանը:

Այժմ դիտարկենք գրաֆների մի այլ տիպի ֆակտորիզացիա: Դիցուք $G$-ն $(n,m)$-գրաֆ է: Այդ դեպքում պարզ է, որ $G$-ն կարելի է ներկայացնել կողերով չհատվող կմախքային անտառների միավորման տեսքով, օրինակ, վերցնելով $G$ գրաֆի յուրաքանչյուր կող և մնացած մեկուսացված գագաթները որպես կմախքային անտառ։ Պարզ է նաև, որ այս դեպքում մենք կստանանք $G$ գրաֆի ներկայացում $m$ հատ կողերով չհատվող կմախքային անտառների միավորման տեսքով։ Առաջանում է բնական խնդիր. տրված $G$ գրաֆի համար գտնել նվազագույն $a(G)$ թիվը, որի դեպքում $G$ գրաֆը կարելի է ներկայացնել $a(G)$ հատ կողերով չհատվող կմախքային անտառների միավորման տեսքով։ Այդ $a(G)$ թիվը կանվանենք $G$ գրաֆի *անտառների տրոհման թիվ*։ Օրինակ, հեշտ է տեսնել, որ $a(K_4) = 2$, իսկ $a(K_5) = 3$ (նկ. 5.4.6):

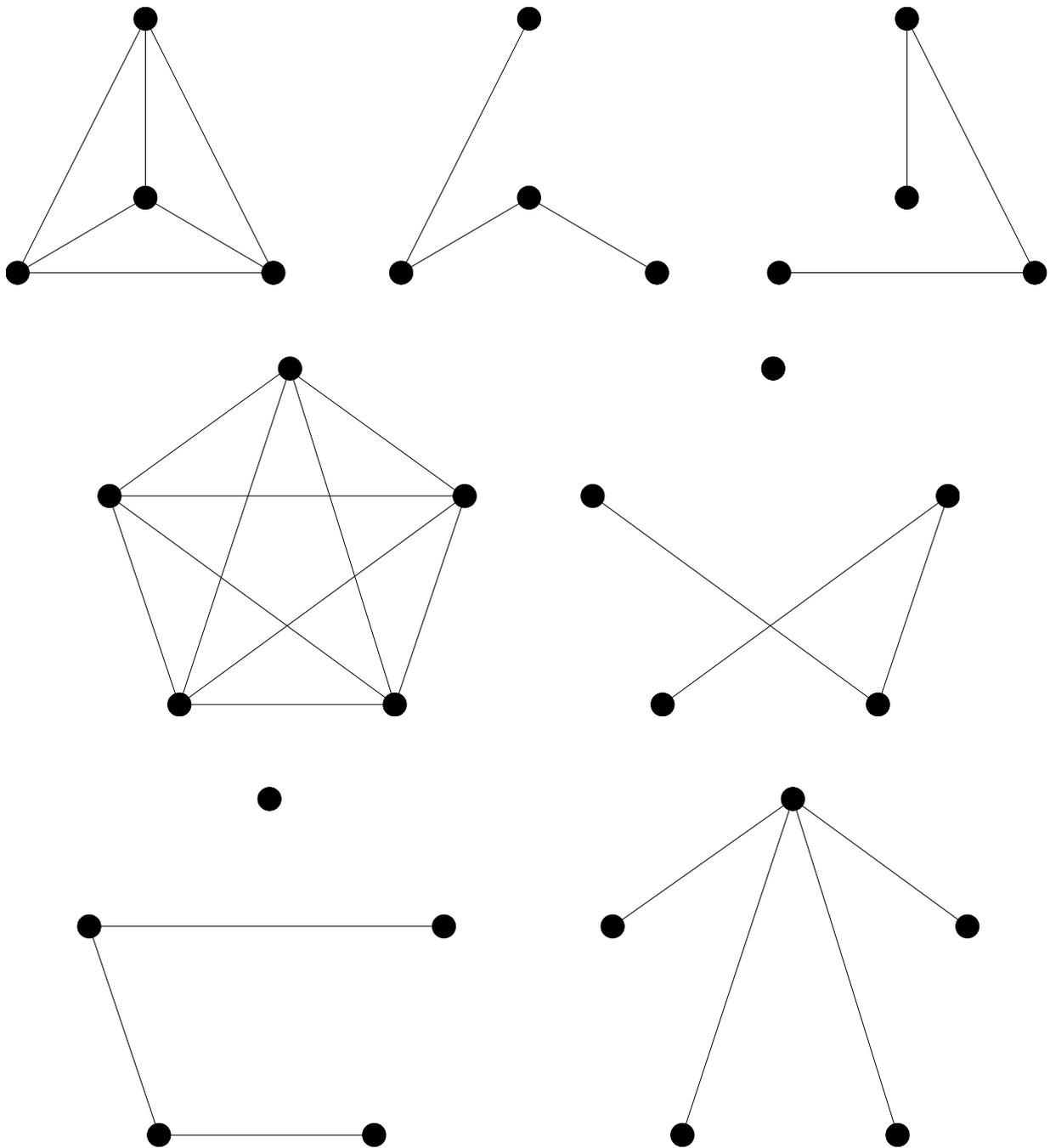

Նկ. 5.4.6

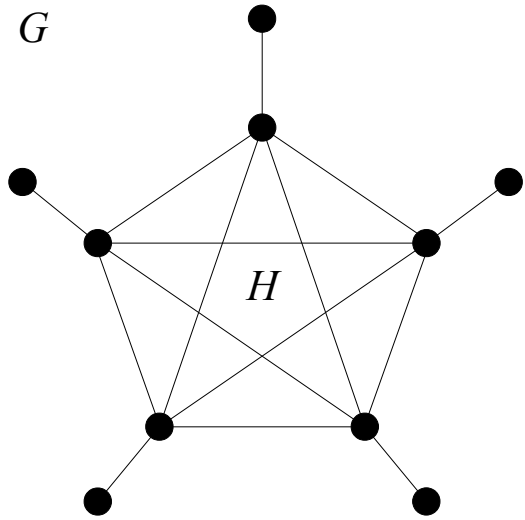

Նկ. 5.4.7

Նախ նկատենք, որ ցանկացած $G$ ($|E(G)| \geq 1$) գրաֆի համար տեղի ունի $a(G) \geq \max\limits_{\substack{H \subseteq G, \\ |V(H)| > 1}} \left\lceil \frac{|E(H)|}{|V(H)|-1} \right\rceil$ անհավասարությունը։ Իրոք, եթե $H \subseteq G$ և $|V(H)| > 1$, ապա $H$ գրաֆի կամայական կմախքային անտառի կողերի քանակը մեծ չէ $|V(H)| - 1$-ից, հետևաբար, ամենաքիչը $\left\lceil \frac{|E(H)|}{|V(H)|-1} \right\rceil$ հատ կմախքային անտառ է անհրաժեշտ $H$ գրաֆի կողերով չհատվող կմախքային անտառներով ներկայացման համար, ուստի $a(H) \geq \left\lceil \frac{|E(H)|}{|V(H)|-1} \right\rceil$։ Մյուս կողմից պարզ է, որ ցանկացած $H \subseteq G$-ի համար $a(G) \geq a(H)$, ինչից հետևում է, որ $a(G) \geq \max\limits_{\substack{H \subseteq G, \\ |V(H)| > 1}} \left\lceil \frac{|E(H)|}{|V(H)|-1} \right\rceil$։ Նշենք, որ այս բանաձևում աջ մասի առավելագույն արժեքը կարող է հասանելի լինի ոչ թե $G$ գրաֆի դեպքում, այլ նրա կողերով հարուստ որևէ $H$ ենթագրաֆի դեպքում։ Այսպես, օրինակ, նկ. 5.4.7-ում պատկերված $G$ գրաֆը $(10, 15)$-գրաֆ է, իսկ նրա $H$ ($H = K_5$) ենթագրաֆը՝ $(5, 10)$-գրաֆ է և $a(H) \geq \left\lceil \frac{|E(H)|}{|V(H)|-1} \right\rceil = \left\lceil \frac{10}{4} \right\rceil = 3 > 2 = \left\lceil \frac{15}{10} \right\rceil = \left\lceil \frac{|E(G)|}{|V(G)|-1} \right\rceil$։

Ինչպես նշել ենք, ցանկացած $G$ գրաֆի համար տեղի ունի $a(G) \geq \max\limits_{\substack{H \subseteq G, \\ |V(H)| > 1}} \left\lceil \frac{|E(H)|}{|V(H)|-1} \right\rceil$ անհավասարությունը։ Պարզվում է, այս ստորին գնահատականը միշտ հասանելի է։

**Թեորեմ 5.4.4 (Նեշ-Վիլյամս):** Կամայական $G$ ($|E(G)| \geq 1$) գրաֆի համար տեղի ունի

$$a(G) = \max\limits_{\substack{H \subseteq G, \\ |V(H)| > 1}} \left\lceil \frac{|E(H)|}{|V(H)|-1} \right\rceil$$

հավասարությունը։

**Ապացույց:** Ենթադրենք հակառակը, դիցուք գոյություն ունի $G'$ ($|E(G')| \geq 1$) գրաֆ,

որի համար $a(G') > \max_{\substack{H \subseteq G', \\ |V(H)|>1}} \left\lceil \frac{|E(H)|}{|V(H)|-1} \right\rceil$: Բոլոր հակաօրինակներից ընտրենք այն $G$ ($|E(G)| \geq 1$) գրաֆը, որի համար $a(G) > \max_{\substack{H \subseteq G, \\ |V(H)|>1}} \left\lceil \frac{|E(H)|}{|V(H)|-1} \right\rceil$, և $|V(G)| + |E(G)|$-ն ընդունում է նվազագույն արժեքը: Պարզ է, որ այդ դեպքում $G$-ն (նվազագույն հակաօրինակը) կլինի կապակցված գրաֆ և $a(G-e) < a(G)$ անհավասարությունը տեղի կունենա ցանկացած $e \in E(G)$-ի համար: Նախ ապացուցենք, որ եթե $G$-ն նվազագույն հակաօրինակ է, ապա ցանկացած $e \in E(G)$-ի համար, $G-e$ գրաֆի ցանկացած ներկայացում $a(G)-1$ հատ կողերով չհատվող կմախքային անտառների միավորման տեսքով իրենից ներկայացնում է $G$ գրաֆի կողերով չհատվող $a(G)-1$ հատ կմախքային ծառերի միավորում: Ենթադրենք հակառակը. դիցուք գոյություն ունի $e = uv \in E(G)$-ի, որ $G-e$ գրաֆը ունի ներկայացում կողերով չհատվող $a(G)-1$ հատ կմախքային անտառների միավորման տեսքով, որտեղ ոչ բոլոր կմախքային անտառները կմախքային ծառեր են: Դիցուք $E(G-e) = E_1 \cup \cdots \cup E_{a(G)-1}$ և $G_i = (V_i, E_i)$, որտեղ $V_i = \{x : x \in V(G)$ և կից է որևէ կողի $E_i$–ից$\}$ և $G_i$-ն անտառ է ($1 \leq i \leq a(G)-1$): Առանց ընդհանրությունը խախտելու կարող ենք ենթադրել, որ $G_1$-ը $G$ գրաֆի կմախքային ծառ չէ:

Քանի որ $G$-ն նվազագույն հակաօրինակ է, ուստի $G_1 + e$-ն պարունակում է ցիկլ: Դիցուք $T$-ն $G_1$ գրաֆի այն կապակցվածության բաղադրիչն է, որը պարունակում է $e$ կողի $u$ և $v$ ցագաթները: Դիցուք $K = G[V(T)]$: Պարզ է, որ $e \in E(K)$ և քանի որ $G_1$-ը $G$ գրաֆի կմախքային ծառ չէ, ուստի $V(T) \neq V(G)$: Մյուս կողմից, քանի որ $G$-ն կապակցված գրաֆ է և $E(G) \setminus E(K) \neq \emptyset$, ուստի $E(K) = A_1 \cup \cdots \cup A_{a(G)-1}$ և $H_i = (U_i, A_i)$, որտեղ $U_i = \{y : y \in V(T)$ և կից է որևէ կողի $A_i$–ից$\}$ և $H_i$-ն անտառ է ($1 \leq i \leq a(G)-1$): Այժմ դիտարկենք հավաքածուների հետևյալ $S$ բազմությունը. բոլոր $(E_1', \ldots, E_{a(G)-1}', \{e'\})$ տեսքի հավաքածուներն, որոնց համար $E(G) = E_1' \cup \cdots \cup E_{a(G)-1}' \cup \{e'\}$-ն $G$ գրաֆի կողերի տրոհում է $a(G)$ հատ անտառների, որի դեպքում $G_1'$ գրաֆի (այստեղ $G_1'$ գրաֆը սահմանվում է նույն ձևով, ինչպես $G_1$-ը) կապակցվածության բաղադրիչը հանդիսանում է $K$-ի կմախքային ծառ և $e' \in E(K)$: Հեշտ է տեսնել, որ $(E_1, \ldots, E_{a(G)-1}, \{e\}) \in S$, ուստի $S \neq \emptyset$: Ընտրենք $S$-ից այն $(\overline{E}_1, \ldots, \overline{E}_{a(G)-1}, \{\overline{e}\})$ հավաքածուն, որի դեպքում

$$\sum_{i=1}^{a(G)-1} |\overline{E}_i \cap A_i|$$

արտահայտությունը ստանում է իր առավելագույն արժեքը: Քանի որ $\overline{e} \in E(K)$, ուստի

$\bar{e} \in A_t$, որտեղ $1 \leq t \leq a(G) - 1$։ Պարզ է, որ $\bar{G}_t + \bar{e}$ գրաֆը (այստեղ $\bar{G}_t$ գրաֆը սահմանվում է նույնպես, ինչպես $G_t$-ն) պարունակում է $C$ ցիկլ և $\bar{e} \in E(C)$։ Եթե $t = 1$, ապա $E(C) \subseteq E(K)$։ Եթե $t \neq 1$ և $E(C) \nsubseteq E(K)$, ապա գոյություն ունի $f \in E(C)$, որի մի զազաթը $V(K)$-ից է, իսկ մյուսը՝ $V(G) \setminus V(K)$-ից է (քանի որ $K$-ն $G$-ի ծնված ենթագրաֆ է)։ Քանի որ $\bar{G}_1$ գրաֆը հանդիսանում է $K$-ի կմախքային ծառ և $\bar{G}_1 + f$-ը անտառ է, ուստի $(\bar{E}_1 + f, \ldots, \bar{E}_t + \bar{e} - f, \ldots, \bar{E}_{a(G)-1})$ հավաքածուին համապատասխանում է $G$ գրաֆի տրոհում կողերով չհատվող $a(G) - 1$ հատ կմախքային անտառների, որը հակասություն է։ Այսպիսով, մենք կարող ենք ենթադրել, որ $E(C) \subseteq E(K)$։ Քանի որ $H_t$-ն անտառ է, ուստի գոյություն ունի $f \in E(C) \setminus A_t \subseteq E(K)$։ Այժմ դիտարկենք $(\bar{E}_1, \ldots, \bar{E}_t + \bar{e} - f, \ldots, \bar{E}_{a(G)-1}, \{f\}) \in S$ հավաքածուն։ Պարզ է, որ այս դեպքում $|\bar{E}_t \cap A_t| < |(\bar{E}_t + \bar{e} - f) \cap A_t|$, որը հակասում է $\sum_{i=1}^{a(G)-1} |\bar{E}_i \cap A_i|$ արտահայտության արժեքի առավելագույն լինելուն։ Այստեղից հետևում է, որ եթե $G$-ն նվազագույն հակաօրինակ է, ապա ցանկացած $e \in E(G)$-ի համար, $G - e$ գրաֆի ցանկացած ներկայացում կողերով չհատվող $a(G) - 1$ հատ կմախքային անտառների միավորման տեսքով իրենից ներկայացնում է $G$ գրաֆի կողերով չհատվող $a(G) - 1$ հատ կմախքային ծառերի միավորում։ Այժմ դիտարկենք հետևյալ հավասարությունը.

$$|E(G)| - 1 = |E(G - e)| = (|V(G)| - 1)(a(G) - 1),$$

որտեղից ստացվում է հետևյալ հակասությունը.

$$a(G) > \left\lceil \frac{|E(G)|}{|V(G)| - 1} \right\rceil = \left\lceil a(G) - 1 + \frac{1}{|V(G)| - 1} \right\rceil = a(G)։$$

∎

**Հետևանք 5.4.1:** Ցանկացած $m$ և $n$ բնական թվերի համար տեղի ունեն

$$a(K_n) = \left\lceil \frac{n}{2} \right\rceil \text{ և } a(K_{m,n}) = \left\lceil \frac{mn}{m+n-1} \right\rceil$$

հավասարությունները։

# Գլուխ 6

## Աստիճանային հավաքածուներ

### § 6.1. Պսևդոգրաֆների և մուլտիգրաֆների աստիճանային հավաքածուներ

Դիցուք $G$-ն պսևդոգրաֆ է և $V(G) = \{v_1, \ldots, v_n\}$-ը այդ պսևդոգրաֆի գագաթների բազմությունն է: Պարզ է, որ յուրաքանչյուր $G$ պսևդոգրաֆին կհամապատասխանի $d = (d_1, \ldots, d_n)$ ամբողջ ոչ բացասական թվերի հավաքածուն, որտեղ $d_i = d_G(v_i)$ ($1 \leq i \leq n$): Այդ $d = (d_1, \ldots, d_n)$ հավաքածուն կանվանենք $G$ պսևդոգրաֆի *աստիճանային հավաքածու*: Հիշենք, որ պսևդոգրաֆի յուրաքանչյուր օղակը այդ պսևդոգրաֆի գագաթի աստիճանն ավելացնում է երկուսով: Նկատենք, որ տարբեր պսևդոգրաֆներ կարող են ունենալ նույն աստիճանային հավաքածուն: Այսպես, օրինակ, նկ. 6.1.1-ում պատկերված պսևդոգրաֆները տարբեր են, սակայն նրանցից յուրաքանչյուրին համապատասխանում է $(4, 4, 2, 2, 1, 1)$ հավաքածուն:

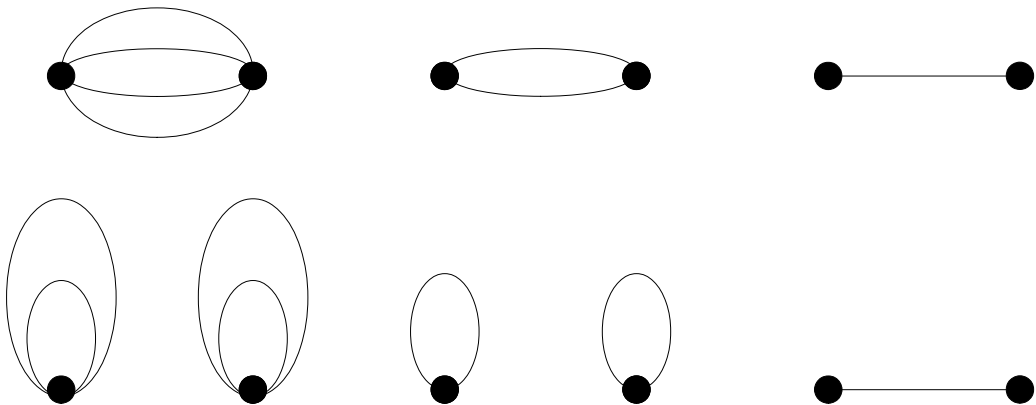

Նկ. 6.1.1

Այժմ ենթադրենք, որ տրված է $d = (d_1, \ldots, d_n)$ ամբողջ ոչ բացասական թվերի հավաքածուն: Կասենք, որ այդ $d = (d_1, \ldots, d_n)$ հավաքածուն *իրացվում է պսևդոգրաֆում*, եթե գոյություն ունի $G$ պսևդոգրաֆ, որի աստիճանային հավաքածուն $d = (d_1, \ldots, d_n)$-ն է: Նկատենք, որ ըստ դիտողություն 1.2.1-ի, եթե $d = (d_1, \ldots, d_n)$ հավաքածուն իրացվում է պսևդոգրաֆում, ապա $\sum_{i=1}^{n} d_i$-ն զույգ թիվ է։ Հակիմին ցույց է տվել, որ այս պարզ անհրաժեշտ պայմանը նաև հանդիսանում է բավարար պայման։

**Թեորեմ 6.1.1:** $d = (d_1, \ldots, d_n)$ ամբողջ ոչ բացասական թվերի հաջաբածուն իրացվում է պսևդոգրաֆում այն և միայն այն դեպքում, երբ $\sum_{i=1}^{n} d_i$-ն զույգ թիվ է։

**Ապացույց:** Ինչպես նշել ենք, ըստ դիտողություն 1.2.1-ի, $\sum_{i=1}^{n} d_i$-ի զույգ թիվ լինելը անհրաժեշտ պայման է պսևդոգրաֆում $d = (d_1, \ldots, d_n)$-ի իրացվելիության համար։

Ցույց տանք, որ եթե $\sum_{i=1}^{n} d_i$-ն զույգ թիվ է, ապա գոյություն ունի $G$ պսևդոգրաֆ, որի աստիճանային հաջաբածուն $d = (d_1, \ldots, d_n)$-ն է։ Մենք կկառուցենք $G$ պսևդոգրաֆ, որի համար $V(G) = \{v_1, \ldots, v_n\}$ և $d_G(v_i) = d_i$ ($1 \leq i \leq n$)։ Քանի որ $\sum_{i=1}^{n} d_i$-ն զույգ թիվ է, ապա $|\{i: d_i - \text{ն կենտ թիվ է}\}|$ թիվը ևս զույգ է։ Տրոհենք $\{v_i: d_i - \text{ն կենտ թիվ է}\}$ բազմության գագաթները զույգերի և յուրաքանչյուր զույգի գագաթները միացնենք կողով։ Այնուհետև յուրաքանչյուր $v_i$ գագաթին ավելացնենք $\left\lfloor \frac{d_i}{2} \right\rfloor$ հատ օղակներ ($1 \leq i \leq n$)։ Հեշտ է տեսնել, որ ստացված $G$ պսևդոգրաֆի աստիճանային հաջաբածուն $d = (d_1, \ldots, d_n)$-ն է։ ∎

Օրինակ, դիտարկենք $(5, 5, 3, 2, 2, 2, 1)$ հաջաբածուն։ Ըստ թեորեմ 6.1.1-ի, այն իրացվում է պսևդոգրաֆում։ Կառուցենք $G$ պսևդոգրաֆ թեորեմ 6.1.1-ի ապացույցում բերված եղանակով։ Դիցուք $V(G) = \{v_1, v_2, v_3, v_4, v_5, v_6, v_7\}$։ Ստորև պատկերված է այդ $G$ պսևդոգրաֆը։

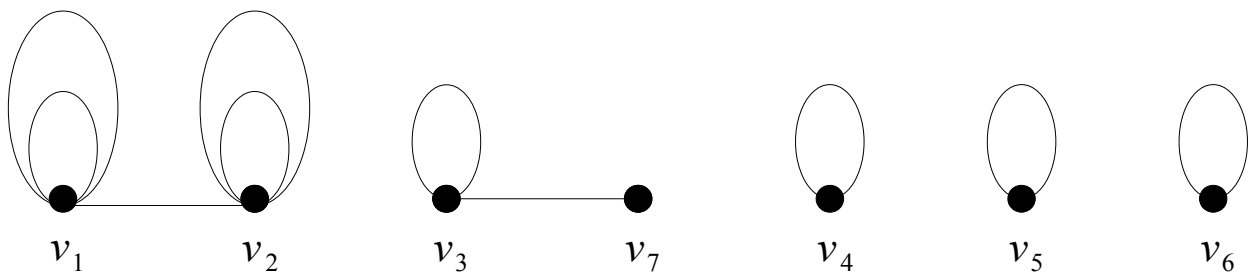

Նկ. 6.1.2

Դիցուք $G$-ն մուլտիգրաֆ է և $V(G) = \{v_1, \ldots, v_n\}$-ը այդ մուլտիգրաֆի գագաթների բազմությունն է։ Պարզ է, որ յուրաքանչյուր $G$ մուլտիգրաֆին կհամապատասխանի $d = (d_1, \ldots, d_n)$ ամբողջ ոչ բացասական թվերի հաջաբածուն, որտեղ $d_i = d_G(v_i)$ ($1 \leq i \leq n$)։ Այդ $d = (d_1, \ldots, d_n)$ հաջաբածուն կանվանենք $G$ մուլտիգրաֆի աստիճանային հաջաբածու։ Նկատենք, որ տարբեր մուլտիգրաֆներ կարող են ունենալ նույն աստիճանային հաջաբածուն։ Այսպես, օրինակ, նկ. 6.1.3-ում պատկերված մուլտիգրաֆները տարբեր են, սակայն նրանցից յուրաքանչյուրին համապատասխանում է $(5, 3, 2, 2, 2)$ հաջաբածուն:

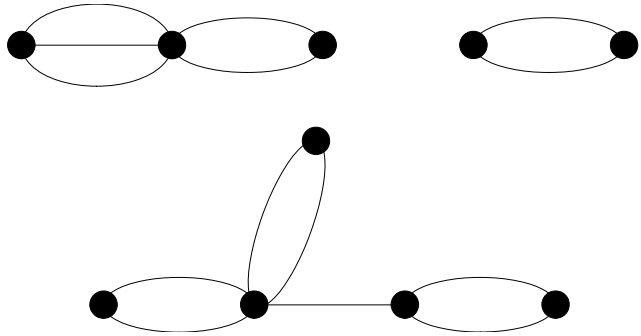

Նկ. 6.1.3

Այժմ ենթադրենք, որ տրված է $d = (d_1, …, d_n)$ ամբողջ ոչ բացասական թվերի հաջաբածուն: Կասենք, որ այդ $d = (d_1, …, d_n)$ հաջաբածուն *իրացվում է մուլտիգրաֆում*, եթե գոյություն ունի $G$ մուլտիգրաֆ, որի աստիճանային հաջաբածուն $d = (d_1, …, d_n)$-ն է: Նկատենք, որ ըստ դիտողություն 1.2.1-ի, եթե $d = (d_1, …, d_n)$ հաջաբածուն իրացվում է մուլտիգրաֆում, ապա $\sum_{i=1}^{n} d_i$-ն զույգ թիվ է: Նշենք, որ այս պայմանը մուլտիգրաֆների դեպքում չի հանդիսանում նաև բավարար պայման: Իրոք, եթե դիտարկենք $(5, 2, 1)$ հաջաբածուն, ապա հեշտ է տեսնել, որ այս հաջաբածուն չի իրացվում մուլտիգրաֆներում:

**Թեորեմ 6.1.2 (Ս. Հակիմի):** $d = (d_1, …, d_n)$ ամբողջ ոչ բացասական թվերի հաջաբածուն, որտեղ $d_1 \geq d_2 \geq \cdots \geq d_n$, իրացվում է մուլտիգրաֆում այն և միայն այն դեպքում, երբ $\sum_{i=1}^{n} d_i$-ն զույգ թիվ է և $d_1 \leq \sum_{i=2}^{n} d_i$:

**Ապացույց:** Ինչպես նշել ենք, ըստ դիտողություն 1.2.1-ի, $\sum_{i=1}^{n} d_i$-ի զույգ թիվ լինելը անհրաժեշտ պայման է մուլտիգրաֆում $d = (d_1, …, d_n)$-ի իրացվելիության համար: Ցույց տանք, որ անհրաժեշտ պայման է հանդիսանում նաև $d_1 \leq \sum_{i=2}^{n} d_i$ անհավասարությունը: Իրոք, եթե ենթադրենք հակառակը՝ $d_1 > \sum_{i=2}^{n} d_i$, ապա հեշտ է տեսնել, որ այդ դեպքում $v_1$ գագաթը, որի աստիճանը $d_1$ է, կից է առնվազն մի կողի, որի մյուս գագաթը տարբեր է $v_2, …, v_n$ գագաթներից, ինչը հակասություն է:

Այժմ ցույց տանք, որ եթե $\sum_{i=1}^{n} d_i$-ն զույգ թիվ է և $d_1 \leq \sum_{i=2}^{n} d_i$, ապա գոյություն ունի $G$ մուլտիգրաֆ, որի աստիճանային հաջաբածուն $d = (d_1, …, d_n)$-ն է: Նախ կառուցենք օժանդակ $H$ մուլտիգրաֆը, որի համար $V(H) = \{v_1, …, v_n\}$ և $d_H(v_1) = \sum_{i=2}^{n} d_i$, $d_H(v_i) = d_i$ ($2 \leq i \leq n$) (նկ. 6.1.4):

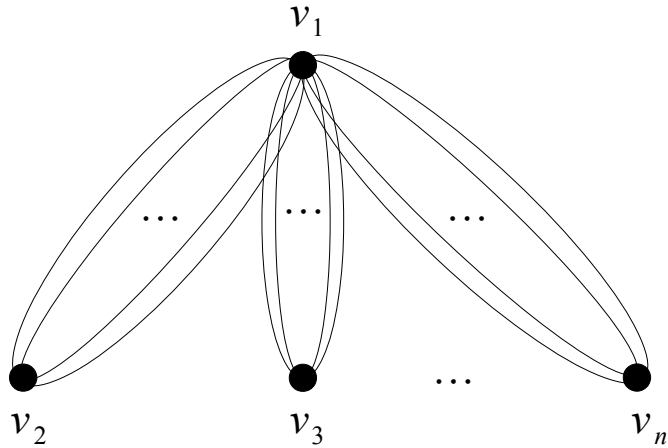

Նկ. 6.1.4

Եթե $d_1 = \sum_{i=2}^{n} d_i$, ապա մենք կվերցնենք $G = H$։ Ենթադրենք $d_1 < \sum_{i=2}^{n} d_i$։ Ցույց տանք, որ $\sum_{i=2}^{n} d_i \equiv d_1 (mod\ 2)$։ Իրոք, քանի որ $\sum_{i=1}^{n} d_i$-ն զույգ թիվ է, ուստի $\sum_{i=1}^{n} d_i - 2d_1 = \sum_{i=2}^{n} d_i - d_1 \equiv 0 (mod\ 2)$։ Այստեղից հետևում է, որ $\sum_{i=2}^{n} d_i \equiv d_1 (mod\ 2)$։ Եթե $d_1 < \sum_{i=2}^{n} d_i$ և $v_1$ գագաթը հարևան է երկու տարբեր $v_{i_0}$ և $v_{i_1}$ գագաթներին $H$ մուլտիգրաֆում, ապա կառուցենք $H_1$ մուլտիգրաֆը հետևյալ կերպ. $H_1 = (H - v_1 v_{i_0} - v_1 v_{i_1}) + v_{i_0} v_{i_1}$։ Հեշտ է տեսնել, որ $H_1$ մուլտիգրաֆում $v_1$ գագաթի աստիճանը $\sum_{i=2}^{n} d_i - 2$ է, իսկ մնացած գագաթների աստիճանները չեն փոխվել։ Եթե $d_1 = \sum_{i=2}^{n} d_i - 2$, ապա մենք կվերցնենք $G = H_1$։ Եթե $d_1 < \sum_{i=2}^{n} d_i - 2$ և $v_1$ գագաթը հարևան է երկու տարբեր $v_{j_0}$ և $v_{j_1}$ գագաթներին $H_1$ մուլտիգրաֆում, ապա կառուցենք $H_2$ մուլտիգրաֆը հետևյալ կերպ. $H_2 = (H_1 - v_1 v_{j_0} - v_1 v_{j_1}) + v_{j_0} v_{j_1}$։ Համանման ձևով մենք կկառուցենք $H_1, \dots, H_l$ մուլտիգրաֆները։ Պարզ է, որ $H_l$ մուլտիգրաֆում $v_1$ գագաթի աստիճանը $\sum_{i=2}^{n} d_i - 2l$ է, իսկ մնացած գագաթների աստիճանները չեն փոխվել։ Եթե $d_1 = \sum_{i=2}^{n} d_i - 2l$, ապա մենք կվերցնենք $G = H_l$։ Ցույց տանք, որ եթե $H_l$ մուլտիգրաֆում $v_1$ գագաթը չունի երկու տարբեր հարևան գագաթներ, ապա $d_1 = \sum_{i=2}^{n} d_i - 2l$։ Իրոք, եթե $H_l$ մուլտիգրաֆում $v_1$ գագաթը չունի երկու տարբեր հարևան գագաթներ և $d_1 < \sum_{i=2}^{n} d_i - 2l$, ապա այդ մուլտիգրաֆում գոյություն կունենա $v_t$ գագաթ, որին կից են $v_1$ գագաթից դուրս եկող բոլոր կողերը, և, հետևաբար, $d_t > d_1$, որը հակասում է $d_1 \geq d_i\ (1 \leq i \leq n)$ պայմանին։ Այստեղից հետևում է, որ, վերցնելով $G = H_l$-ի, մենք կստանանք մուլտիգրաֆ, որի աստիճանային հավաքածուն $d = (d_1, \dots, d_n)$-ն է։ ∎

Օրինակ, դիտարկենք $(5, 4, 3, 3, 2, 2, 1)$ հավաքածուն։ Ըստ թեորեմ 6.1.2-ի այն իրացվում է մուլտիգրաֆում։ Կառուցենք $G$ մուլտիգրաֆ թեորեմ 6.1.2-ի ապացուցում բերված եղանակով։ Դիցուք $V(G) = \{v_1, v_2, v_3, v_4, v_5, v_6, v_7\}$։ Ստորև պատկերված է այդ $G$ մուլտիգրաֆը ստանալու ամբողջ ընթացքը։

**H**

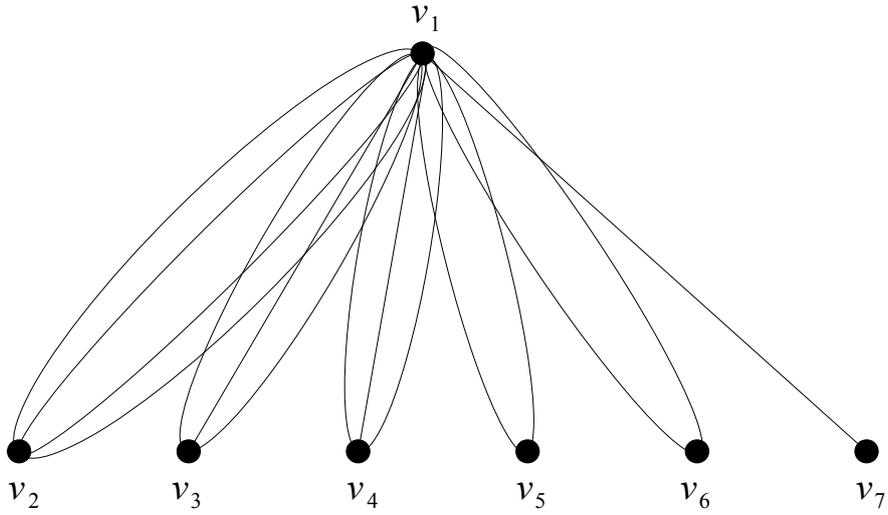

**H₁**

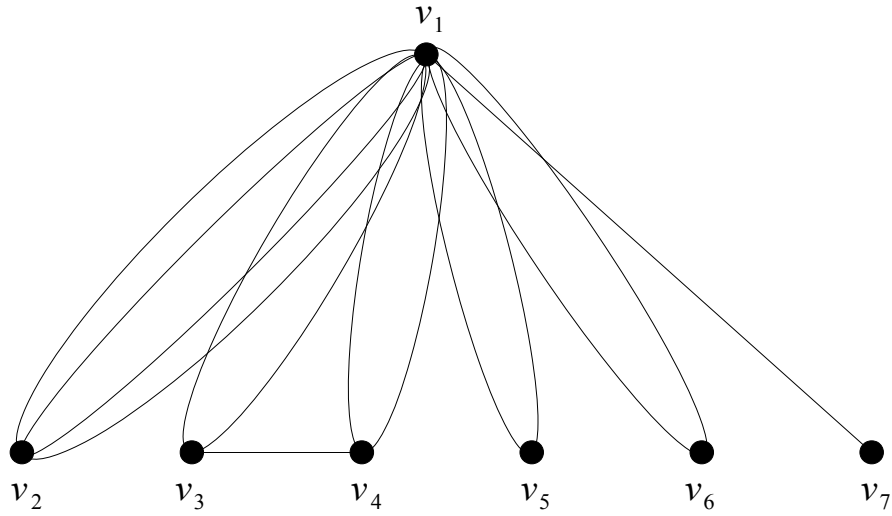

**H₂**

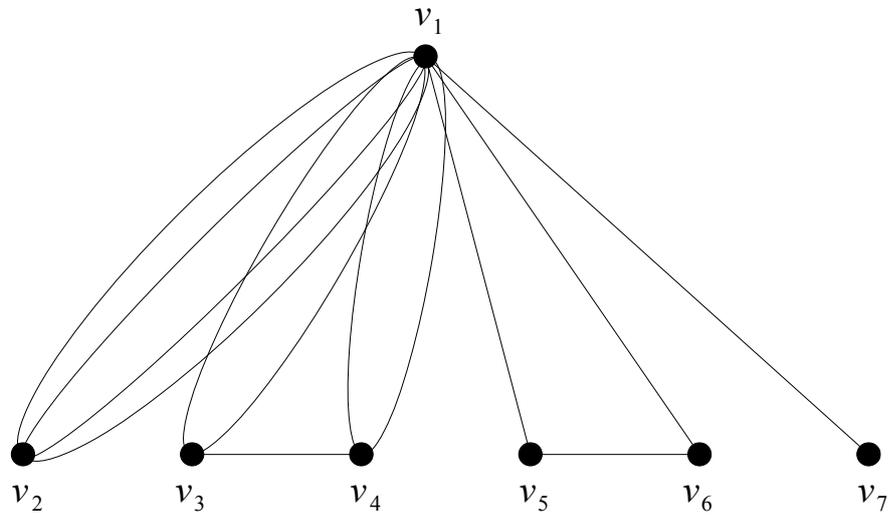

**$H_3$**

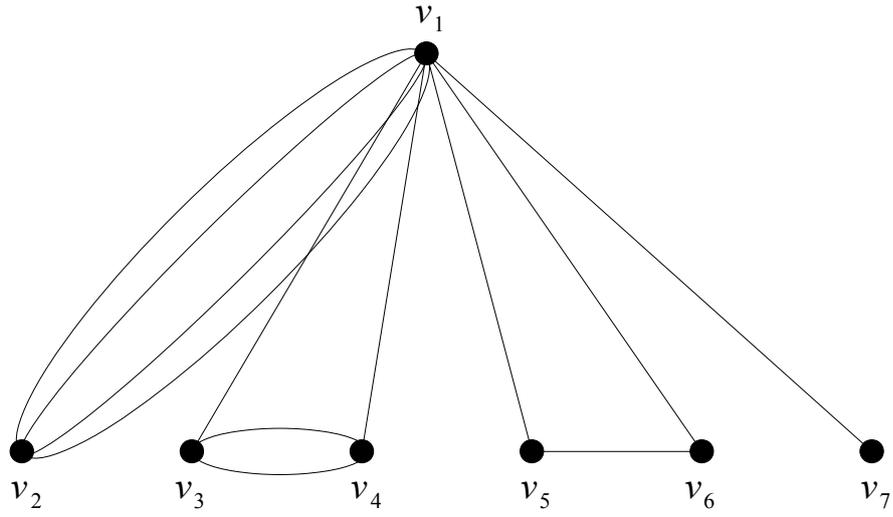

**$H_4$**

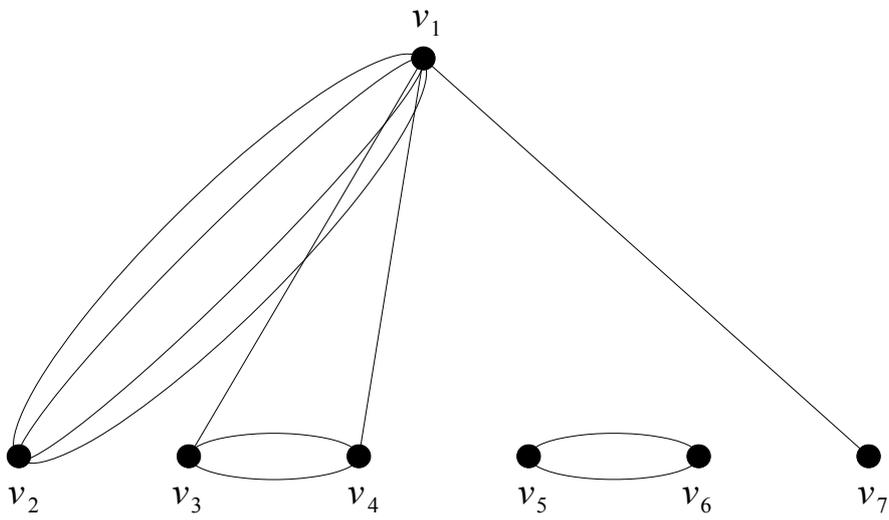

**$G = H_5$**

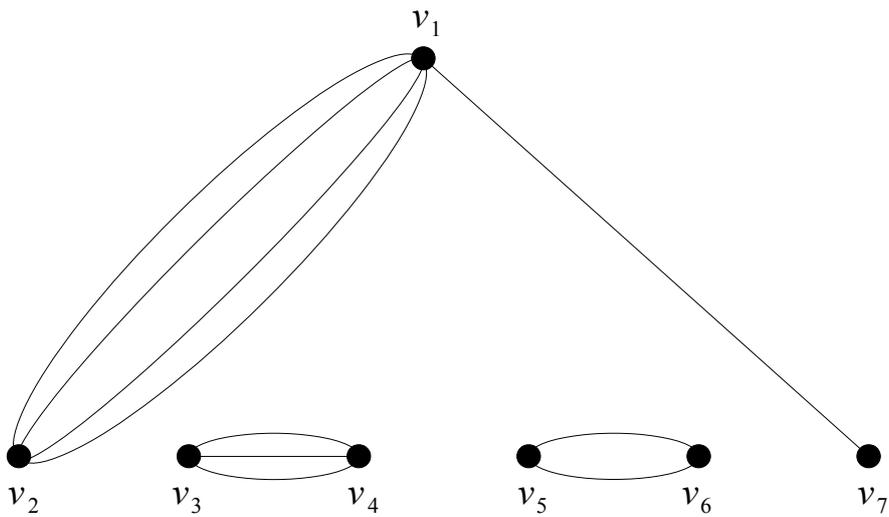

Նկ. 6.1.5

## § 6.2. Գրաֆային հավաքածուներ

Դիցուք $G$-ն գրաֆ է և $V(G) = \{v_1, \ldots, v_n\}$-ը այդ գրաֆի գագաթների բազմությունն է։ Պարզ է, որ յուրաքանչյուր $G$ գրաֆին կհամապատասխանի $d = (d_1, \ldots, d_n)$ ամբողջ ոչ բացասական թվերի հավաքածուն, որտեղ $d_i = d_G(v_i)$ $(1 \leq i \leq n)$։ Այդ $d = (d_1, \ldots, d_n)$ հավաքածուն կանվանենք $G$ գրաֆի աստիճանային հավաքածու։ Նկատենք, որ տարբեր գրաֆներ կարող են ունենալ նույն աստիճանային հավաքածուն։ Այսպես, օրինակ, նկ. 6.2.1-ում պատկերված գրաֆները տարբեր են, սակայն նրանցից յուրաքանչյուրին համապատասխանում է $(2, 2, 2, 1, 1)$ հավաքածուն։

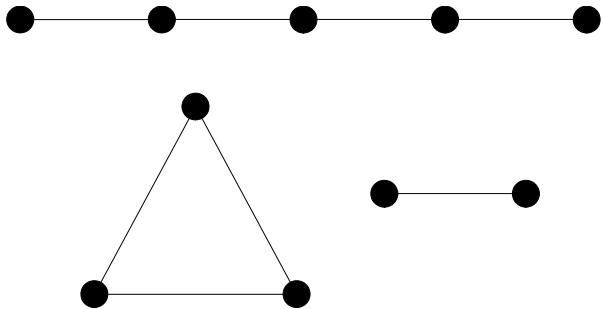

Նկ. 6.2.1

Այժմ ենթադրենք, որ տրված է $d = (d_1, \ldots, d_n)$ ամբողջ ոչ բացասական թվերի հավաքածուն։

**Սահմանում 6.2.1:** $d = (d_1, \ldots, d_n)$ հավաքածուն կանվանենք *գրաֆային*, եթե գոյություն ունի $G$ գրաֆ, որի աստիճանային հավաքածուն $d = (d_1, \ldots, d_n)$-ն է։

Նկատենք, որ ըստ դիտողության 1.2.1-ի, եթե $d = (d_1, \ldots, d_n)$ հավաքածուն գրաֆային է, ապա $\sum_{i=1}^{n} d_i$-ն զույգ թիվ է։ Հեշտ է տեսնել նաև, որ եթե $d = (d_1, \ldots, d_n)$ հավաքածուն գրաֆային է, ապա $0 \leq d_i \leq n-1$ $(1 \leq i \leq n)$։ Նշենք, որ այս պայմանը գրաֆների դեպքում չի հանդիսանում նաև բավարար պայման։ Իրոք, եթե դիտարկենք $(3, 3, 3, 1)$ հավաքածուն, ապա հեշտ է տեսնել, որ այս հավաքածուն գրաֆային չէ։ Առաջանում է բնական հարց, թե ո՞ր հավաքածուներն են գրաֆային և որոնք՝ ոչ, և եթե հավաքածուն գրաֆային է, ապա ինչպե՞ս կառուցել այդ աստիճանային հավաքածուն ունեցող գրաֆը։ Առաջին հարցի պատասխանը տրվել է Էրդյոշի և Գալլայի կողմից, իսկ երկրորդինը՝ Հավելի և Հակիմի կողմից։

Թեորեմ 6.2.1 (Պ. Էրդյոշի, Տ. Գալլայի): $d = (d_1, \dots, d_n)$ ամբողջ ոչ բացասական թվերի հաջորդություն, որտեղ $d_1 \geq d_2 \geq \dots \geq d_n$, գրաֆային է այն և միայն այն դեպքում, երբ $\sum_{i=1}^{n} d_i$-ն զույգ թիվ է և ցանկացած $k$-ի համար ($1 \leq k \leq n-1$) տեղի ունի

$$\sum_{i=1}^{k} d_i \leq k(k-1) + \sum_{i=k+1}^{n} \min\{k, d_i\}$$

պայմանը:

Այժմ ձևակերպենք և ապացուցենք Հավելի և Հակիմի թեորեմը:

**Թեորեմ 6.2.2:** Եթե $n = 1$, ապա միակ գրաֆային հաջորդություն $d = (0)$-ն է: Եթե $n \geq 2$, ապա $d = (d_1, \dots, d_n)$ ամբողջ ոչ բացասական թվերի հաջորդություն, որտեղ $d_1 \geq d_2 \geq \dots \geq d_n$, գրաֆային է այն և միայն այն դեպքում, երբ գրաֆային է

$$d' = (d_2 - 1, d_3 - 1, \dots, d_{d_1+1} - 1, d_{d_1+2}, \dots, d_n)$$

հաջորդությունը:

**Ապացույց:** Նախ ցույց տանք, որ եթե $d' = (d_2 - 1, d_3 - 1, \dots, d_{d_1+1} - 1, d_{d_1+2}, \dots, d_n)$-ը գրաֆային հաջորդություն է, ապա $d = (d_1, \dots, d_n)$-ն ևս գրաֆային է: Իրոք, եթե $d'$ հաջորդությունը գրաֆային է, ապա գոյություն ունի $G'$ գրաֆ, որի աստիճանային հաջորդություն $d'$-ն է, և եթե մենք $G'$ գրաֆին ավելացնենք մեկ գագաթ և միացնենք այն $d_2 - 1, d_3 - 1, \dots, d_{d_1+1} - 1$ աստիճաններ ունեցող $d_1$ հատ գագաթների հետ, ապա կստացվի $G$ գրաֆ, որի աստիճանային հաջորդություն $d$-ն է:

Ցույց տանք, որ եթե $d = (d_1, \dots, d_n)$-ն գրաֆային հաջորդություն է, ապա $d' = (d_2 - 1, d_3 - 1, \dots, d_{d_1+1} - 1, d_{d_1+2}, \dots, d_n)$ հաջորդություն ևս գրաֆային է: Քանի որ $d$-ն գրաֆային է, ապա գոյություն ունի $G$ գրաֆ, որի աստիճանային հաջորդություն $d$-ն է: Դիցուք $w \in V(G)$ և $d_G(w) = d_1$: $S$-ով նշանակենք $G$ գրաֆի $d_2, d_3, \dots, d_{d_1+1}$ աստիճաններ ունեցող և $d_1$ հատ գագաթ պարունակող բազմությունը: Եթե $N_G(w) = S$, ապա հեշտ է տեսնել, որ $G - w$ գրաֆի աստիճանային հաջորդություն կլինի $d'$-ը: Ենթադրենք $N_G(w) \neq S$: Այս դեպքում մենք $G$ գրաֆից կառուցենք $G'$ գրաֆի, որի աստիճանային հաջորդություն $d$-ն է և $|N_G(w) \cap S| < |N_{G'}(w) \cap S|$: Եթե $N_{G'}(w) = S$, ապա պարզ է որ $G' - w$ գրաֆի աստիճանային հաջորդություն կլինի $d'$-ը: Եթե $N_{G'}(w) \neq S$, ապա մենք $G'$ գրաֆից կառուցենք $G''$ գրաֆի, որի աստիճանային հաջորդություն $d$-ն է և $|N_{G'}(w) \cap S| < |N_{G''}(w) \cap S|$: Պարզ է, որ համանման անցումներ կատարելով, մենք ամենաշատը $d_1$ քայլից կկառուցենք $G^*$ գրաֆ, որի աստիճանային հաջորդություն ևս կլինի $d$-ն, և $N_{G^*}(w) =$ $S$, իսկ $G^* - w$ գրաֆի աստիճանային հավաքածուն կլինի $d'$-ը: Այժմ ցույց տանք, որ այդպիսի անցումները հնարավոր են: Եթե $N_G(w) \neq S$, ապա քանի որ $d_G(w) = d_1 = |S|$, ուստի $G$ գրաֆում գոյություն ունեն այնպիսի $x \in S$ և $z \notin S$, որ $wx \notin E(G)$ և $wz \in E(G)$: Ըստ $S$-ի սահմանման՝ $d_G(x) \geq d_G(z)$: Դիցուք $T = \{w, x, z\}$: Ցույց տանք, որ $G$ գրաֆում գոյություն ունի այնպիսի $y$ գագաթ, որ $y \notin T$ և $xy \in E(G)$, $yz \notin E(G)$ (նկ. 6.2.2): Սահմանենք $\varepsilon$ թիվը հետևյալ կերպ.

$$\varepsilon = \begin{cases} 0, & \text{եթե } xz \notin E(G), \\ 1, & \text{հակառակ դեպքում:} \end{cases}$$

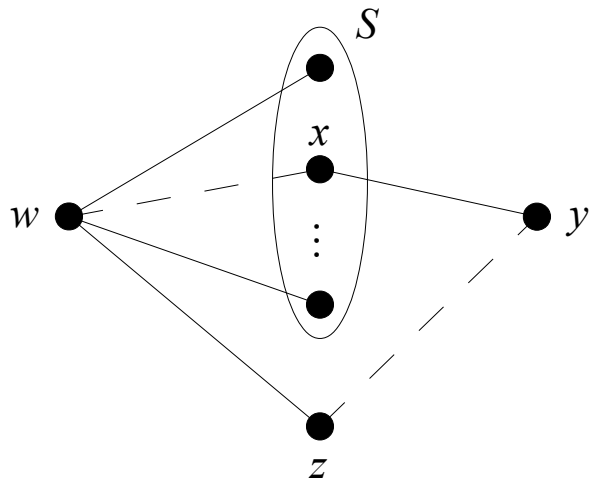

Նկ. 6.2.2

Այժմ $y$ գագաթի գոյությունը $G$ գրաֆում հետևում է նրանից, որ $x$ գագաթը $T$-ից դուրս գտնվող մասում ունի $d_G(x) - \varepsilon$ հատ հարևան գագաթներ, $z$ գագաթը $T$-ից դուրս գտնվող մասում ունի $d_G(z) - 1 - \varepsilon$ հատ հարևան գագաթներ և $d_G(x) \geq d_G(z)$: $G$ գրաֆից անցնենք $G'$ գրաֆի հետևյալ կերպ. $G' = G - wz - xy + wx + yz$: Հեշտ է տեսնել, որ $G'$ գրաֆի աստիճանային հավաքածուն ևս $d$-ն է և $|N_G(w) \cap S| < |N_{G'}(w) \cap S|$: ∎

Թեորեմ 6.2.2-ը հնարավորություն է տալիս առաջարկել ալգորիթմ, որը թույլ է տալիս կառուցել տրված հավաքածուի համար նրան համապատասխան գրաֆը կամ տալիս է բացասական պատասխան, եթե այդ հավաքածուն գրաֆային չէ:

### Ալգորիթմ

Դիցուք տրված է $d = (d_1, \ldots, d_n)$ ամբողջ ոչ բացասական թվերի հավաքածուն, որտեղ $d_1 \geq d_2 \geq \cdots \geq d_n$:

**Քայլ 1:** $d$ հավաքածուի համար կառուցել ըստ թեորեմ 6.2.2-ի $d'$ հավաքածուն:

**Քայլ 2:** Դասավորել $d'$ հավաքածուի տարրերը չաճման կարգով և ստացված հավաքածուն անվանել $d^{(1)}$:

**Քայլ 3:** $d^{(1)}$ հավաքածուի համար կառուցել ըստ թեորեմ 6.2.2-ի $d''$ հավաքածուն:

**Քայլ 4:** Դասավորել $d''$ հավաքածուի տարրերը չամման կարգով և ստացված հավաքածուն անվանել $d^{(2)}$:

**Քայլ 5:** Կատարել նմանատիպ քայլեր մինչև չգտնվի բացասական տարր պարունակող հավաքածու (հավաքածուն գրաֆային չէ) կամ հավաքածուի բոլոր տարրերը լինեն զրոներ (հավաքածուն գրաֆային է):

Օրինակ, դիտարկենք $d = (5, 5, 3, 3, 2, 2, 2)$ հավաքածուն: Կիրառենք նկարագրված ալգորիթմը պարզելու համար, գրաֆային է արդյոք $(5, 5, 3, 3, 2, 2, 2)$ հավաքածուն, թե ոչ: Դրա համար ալգորիթմի աշխատանքի ընթացքում ստացվող յուրաքանչյուր հավաքածուի համար պահենք հետևյալ տիպի աղյուսակը.

| $V$ | $v_1$ | $v_2$ | $v_3$ | $v_4$ | $v_5$ | $v_6$ | $v_7$ |
|---|---|---|---|---|---|---|---|
| $d$ | 5 | 5 | 3 | 3 | 2 | 2 | 2 |

Քայլ 1-ից հետո կստանանք հետևյալ աղյուսակը.

| $V'$ | $v_2$ | $v_3$ | $v_4$ | $v_5$ | $v_6$ | $v_7$ |
|---|---|---|---|---|---|---|
| $d'$ | 4 | 2 | 2 | 1 | 1 | 2 |

Քայլ 2-ից հետո կստանանք հետևյալ աղյուսակը.

| $V^{(1)}$ | $v_2$ | $v_3$ | $v_4$ | $v_7$ | $v_5$ | $v_6$ |
|---|---|---|---|---|---|---|
| $d^{(1)}$ | 4 | 2 | 2 | 2 | 1 | 1 |

Քայլ 3-ից հետո կստանանք հետևյալ աղյուսակը.

| $V''$ | $v_3$ | $v_4$ | $v_7$ | $v_5$ | $v_6$ |
|---|---|---|---|---|---|
| $d''$ | 1 | 1 | 1 | 0 | 1 |

Քայլ 4-ից հետո կստանանք հետևյալ աղյուսակը.

| $V^{(2)}$ | $v_3$ | $v_4$ | $v_7$ | $v_6$ | $v_5$ |
|---|---|---|---|---|---|
| $d^{(2)}$ | 1 | 1 | 1 | 1 | 0 |

Քայլ 5-ից հետո կստանանք հետևյալ աղյուսակը.

| $V'''$ | $v_4$ | $v_7$ | $v_6$ | $v_5$ |
|---|---|---|---|---|
| $d'''$ | 0 | 1 | 1 | 0 |

Քայլ 6-ից հետո կստանանք հետևյալ աղյուսակը.

| $V^{(3)}$ | $v_7$ | $v_6$ | $v_4$ | $v_5$ |
|---|---|---|---|---|
| $d^{(3)}$ | 1 | 1 | 0 | 0 |

Քայլ 7-ից հետո կստանանք հետևյալ աղյուսակը.

| $V''''$ | $v_6$ | $v_4$ | $v_5$ |
|---|---|---|---|
| $d''''$ | 0 | 0 | 0 |

Քանի որ $d'''' = (0, 0, 0)$ հավաքածուի բոլոր տարրերը զրոներ են, ուստի $d = (5, 5, 3, 3, 2, 2, 2)$ հավաքածուն գրաֆային է: Այժմ կառուցենք $d = (5, 5, 3, 3, 2, 2, 2)$ հավաքածուին համապատասխան գրաֆը: Կառուցումը կատարվում է վերջից (նկ. 6.2.3):

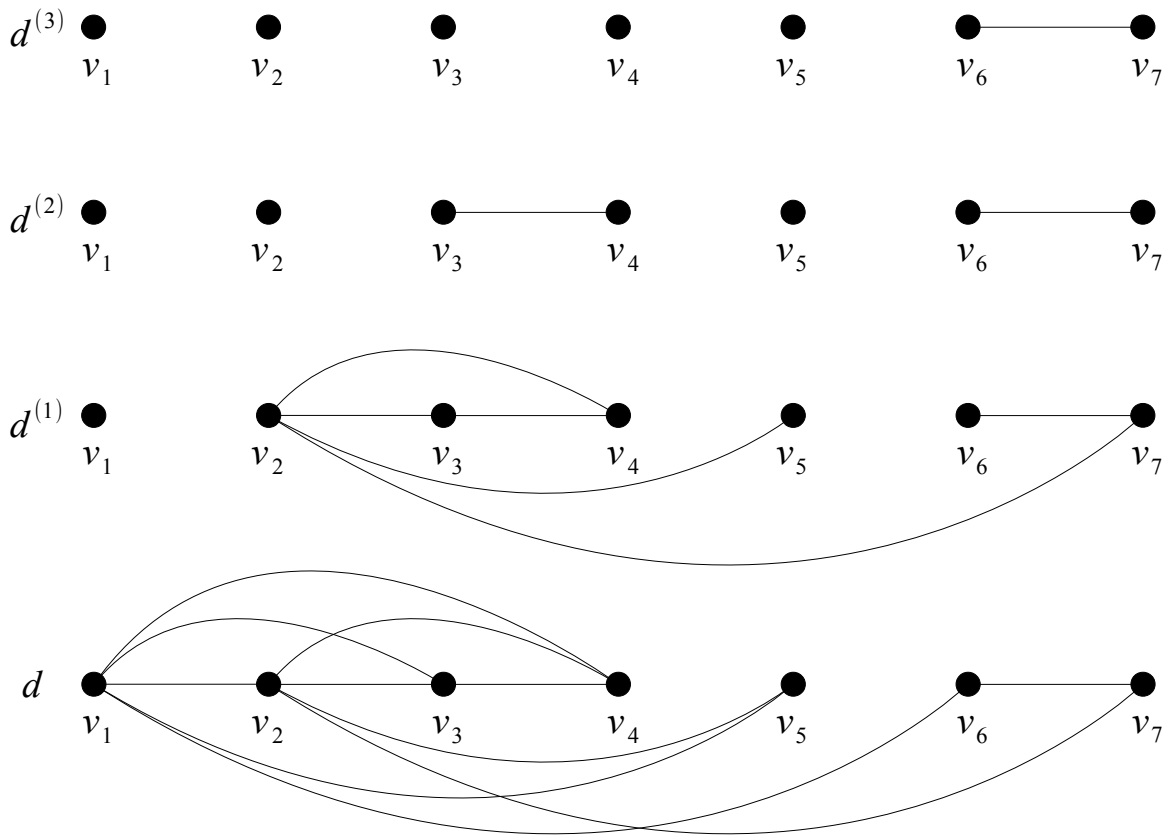

Նկ. 6.2.3

Այժմ դիտարկենք այն հարցը, թե երբ $d = (d_1, ..., d_n)$ գրաֆային հավաքածուին համապատասխանում է կապակցված գրաֆ և, մասնավորապես, ծառ:

Կասենք, որ $d = (d_1, ..., d_n)$ հավաքածուն *իրացվում է կապակցված գրաֆում*, եթե գոյություն ունի $G$ կապակցված գրաֆ, որի աստիճանային հավաքածուն $d = (d_1, ..., d_n)$-ն է:

**Թեորեմ 6.2.3:** Եթե $n \geq 2$, ապա $d = (d_1, ..., d_n)$ գրաֆային հավաքածուն, որտեղ $d_1 \geq d_2 \geq \cdots \geq d_n$, իրացվում է կապակցված գրաֆում այն և միայն այն դեպքում, երբ $d_n > 0$ և տեղի ունի $\sum_{i=1}^{n} d_i \geq 2(n-1)$ պայմանը:

**Ապացույց:** Նախ ցույց տանք, որ եթե $d = (d_1, ..., d_n)$-ն գրաֆային հավաքածու է, որտեղ $d_1 \geq d_2 \geq \cdots \geq d_n$, իրացվում է կապակցված գրաֆում, ապա $d_n > 0$ և տեղի ունի $\sum_{i=1}^{n} d_i \geq 2(n-1)$ պայմանը: Իրոք, քանի որ կապակցված գրաֆը չի պարունակում մեկուսացված գագաթներ, կստանանք, որ $d_n > 0$: Մյուս կողմից քանի որ ցանկացած կապակցված գրաֆ պարունակում է կմախքային ծառ ըստ թեորեմ 2.3.4-ի, ուստի ըստ թեորեմ 1.2.1-ի կստանանք, որ $\sum_{i=1}^{n} d_i \geq 2(n-1)$:

Ցույց տանք, որ եթե $d = (d_1, ..., d_n)$-ն գրաֆային հավաքածու է, որտեղ $d_1 \geq d_2 \geq \cdots \geq d_n > 0$ և տեղի ունի $\sum_{i=1}^{n} d_i \geq 2(n-1)$ պայմանը, ապա այդ հավաքածուն իրացվում է կապակցված գրաֆում: Ապացույցը կատարենք մակածման եղանակով ըստ $n$-ի: Եթե $n = 2$, ապա պարզ է, որ այդ պայմաններին բավարարող հավաքածուն միակն է և այն $d = (1, 1)$-ն է, որին համապատասխանում է երկու իրար հարևան գագաթներ պարունակող կապակցված գրաֆը: Ենթադրենք, $n \geq 3$ և պնդումը ճիշտ է ցանկացած $d' = (d'_1, ..., d'_{n'})$ գրաֆային հավաքածուի համար, երբ $n' < n$: Դիտարկենք $d = (d_1, ..., d_n)$ գրաֆային հավաքածուն: Տրոհենք ապացույցը երկու դեպքի:

**Դեպք 1:** $d_n = 1$:

Քանի որ $n \geq 3$, ուստի $d_1 \geq 2$: Դիտարկենք $d' = (d_1 - 1, d_2, ..., d_{n-1})$ հավաքածուն: Դժվար չէ ցույց տալ, որ եթե $d = (d_1, ..., d_n)$-ն գրաֆային է, ապա $d' = (d_1 - 1, d_2, ..., d_{n-1})$ հավաքածուն ևս կլինի գրաֆային: Նկատենք, որ $d_{n-1} > 0$ և $\sum_{i=1}^{n-1} d'_i = \sum_{i=1}^{n} d_i - 2 \geq 2(n-2)$: Այստեղից, ըստ մակածման ենթադրության, հետևում է, որ $d'$-ը իրացվում է կապակցված գրաֆում: Դիցուք $G'$-ը կապակցված գրաֆ է, որի աստիճանային հավաքածուն $d' = (d_1 - 1, d_2, ..., d_{n-1})$-ն է: Կառուցենք կապակցված $G$ գրաֆ, որի աստիճանային հավաքածուն $d = (d_1, ..., d_n)$-ն է, հետևյալ կերպ. $G'$ գրաֆին ավելացնենք նոր գագաթ և միացնենք այն $d_1 - 1$ աստիճան ունեցող գագաթի հետ:

Դեպք 2: $d_n = l \geq 2$:

Դիտարկենք $d' = (d_1 - 1, d_2 - 1, \ldots, d_l - 1, \ldots, d_{n-1})$ հավաքածուն: Դժվար չէ ցույց տալ, որ եթե $d = (d_1, \ldots, d_n)$-ն գրաֆային է, ապա $d' = (d_1 - 1, d_2 - 1, \ldots, d_l - 1, \ldots, d_{n-1})$ հավաքածուն ևս կլինի գրաֆային: Նկատենք, որ $d_{n-1} > 0$ և

$$\sum_{i=1}^{n-1} d'_i = \sum_{i=1}^{n} d_i - 2l \geq nl - 2l = l(n-2) \geq 2(n-2):$$

Այստեղից, ըստ մակածման ենթադրության, հետևում է, որ $d'$-ը իրացվում է կապակցված գրաֆում: Դիցուք $G'$-ը կապակցված գրաֆ է, որի աստիճանային հավաքածուն $d' = (d_1 - 1, d_2 - 1, \ldots, d_l - 1, \ldots, d_{n-1})$-ն է: Կառուցենք կապակցված $G$ գրաֆ, որի աստիճանային հավաքածուն $d = (d_1, \ldots, d_n)$-ն է, հետևյալ կերպ. $G'$ գրաֆին ավելացնենք նոր գագաթ և միացնենք այն $d_1 - 1, d_2 - 1, \ldots, d_l - 1$ աստիճան ունեցող գագաթների հետ: ∎

Կասենք, որ $d = (d_1, \ldots, d_n)$ հավաքածուն *իրացվում է ծառում*, եթե գոյություն ունի $T$ ծառ, որի աստիճանային հավաքածուն $d = (d_1, \ldots, d_n)$-ն է: Թեորեմ 5.2.3-ի ապացույցի դեպք 1-ի դատողությունները կրկնելով, կարելի է ապացուցել հետևյալ թեորեմը:

**Թեորեմ 6.2.4:** Եթե $n \geq 2$, ապա $d = (d_1, \ldots, d_n)$ հավաքածուն, որտեղ $d_1 \geq d_2 \geq \cdots \geq d_n$, իրացվում է ծառում այն և միայն այն դեպքում, երբ $d_n > 0$ և տեղի ունի $\sum_{i=1}^{n} d_i = 2(n-1)$ պայմանը:

## § 6.3. Տրոհվող գրաֆների և կատարյալ զույգակցում պարունակող գրաֆների աստիճանային հավաքածուներ

**Սահմանում 6.3.1:** $G$ գրաֆը կոչվում է *տրոհվող*, եթե այդ գրաֆի գագաթների $V(G)$ բազմությունը կարելի է տրոհել $I$ և $C$ ենթաբազմությունների այնպես, որ $G[I]$-ն պարունակում է միայն մեկուսացված գագաթներ, իսկ $G[C]$-ն՝ լրիվ գրաֆ է:

**Թեորեմ 6.3.1 (Պ. Համմեր, Բ. Սիմեոնե):** Եթե $d = (d_1, \ldots, d_n)$-ն գրաֆային հավաքածու է, որտեղ $d_1 \geq d_2 \geq \cdots \geq d_n$, և $G$-ն ցանկացած գրաֆ է, որի աստիճանային հավաքածուն $d = (d_1, \ldots, d_n)$-ն է, ապա $G$-ն տրոհվող գրաֆ է այն և միայն այն դեպքում, երբ տեղի ունի

$$\sum_{i=1}^{k} d_i = k(k-1) + \sum_{i=k+1}^{n} d_i$$

պայմանը, որտեղ $k = k(d) = \max\{i: d_i \geq i - 1\}$:

**Ապացույց:** Դիցուք $G$-ն տրոհվող գրաֆ է: Տրոհենք $G$ գրաֆի գագաթների բազմությունը $I$ և $C$ ենթաբազմությունների այնպես, որ $G[C]$-ն պարունակի ամենաշատ քանակությամբ զույգ առ զույգ հարևան գագաթներ: Պարզ է, որ այդ դեպքում, եթե $v \in I$, $u \in C$ և $|C| = l$, ապա $d_G(u) \geq l - 1$ և $d_G(v) < l$: Այստեղից հետևում է, որ $l = k$: Քանի որ $G[C]$-ն լրիվ գրաֆ է, իսկ $G[I]$-ն անկախ բազմություն է, ուստի ճիշտ է $\sum_{i=1}^{k} d_i = k(k-1) + \sum_{i=k+1}^{n} d_i$ պայմանը, որտեղ $k = k(d) = \max\{i: d_i \geq i - 1\}$:

Այժմ ենթադրենք $G$-ն գրաֆ է, որի աստիճանային հավաքածուն $d = (d_1, \ldots, d_n)$-ն է: Դիցուք $V(G) = \{v_1, \ldots, v_n\}$, $C = \{v_1, \ldots, v_k\}$ և $I = \{v_{k+1}, \ldots, v_n\}$: Տրոհենք $\sum_{i=1}^{k} d_i$ գումարը երկու մասի. $\sum_{i=1}^{k} d_i = A + B$, որտեղ $A$-ն այդ գումարի մեջ $v_i v_j$ ($v_i, v_j \in C$) տեսքի կողերի ներդրումն է, իսկ $B$-ն` այդ գումարի մեջ $v_i v_j$ ($v_i \in C, v_j \in I$) տեսքի կողերի ներդրումն է: Հեշտ է տեսնել, որ $A \leq k(k-1)$ և $B \leq \sum_{i=k+1}^{n} d_i$: Մյուս կողմից պարզ է, որ $\sum_{i=1}^{k} d_i = k(k-1) + \sum_{i=k+1}^{n} d_i$ պայմանը տեղի ունի այն և միայն այն դեպքում, երբ $A = k(k-1)$ և $B = \sum_{i=k+1}^{n} d_i$, իսկ դա նշանակում է, որ $G[C]$-ն լրիվ գրաֆ է, $G[I]$-ն անկախ բազմություն է և, հետևաբար, $G$-ն տրոհվող գրաֆ է: ∎

Թեորեմ 6.3.1-ից հետևում է, որ եթե $d = (d_1, \ldots, d_n)$ գրաֆային հավաքածուին համապատասխան գրաֆը տրոհվող գրաֆ է, ապա այդ հավաքածուին համապատասխանող բոլոր գրաֆները տրոհվող են:

Հետևյալ թեորեմի պնդումն առաջարկվել է որպես հիպոթեզ 1970-ին Գյունբաումի կողմից և ապացուցվել Կունդուի կողմից 1973-ին և Լովասի կողմից՝ 1974-ին:

**Թեորեմ 6.3.2 (Կունդու, Լովաս):** Դիցուք $d_1, \ldots, d_n$-ը ամբողջ թվեր են: Որպեսզի գոյություն ունենա կատարյալ զուգակցում պարունակող $G$ գրաֆ, որի աստիճանային հավաքածուն $d = (d_1, \ldots, d_n)$-ն է, անհրաժեշտ է և բավարար, որ $n$-ը լինի զույգ, և $d = (d_1, \ldots, d_n)$, $d' = (d_1 - 1, \ldots, d_n - 1)$ հավաքածուները լինեն գրաֆային:

**Ապացույց:** Ենթադրենք, որ $G$ կատարյալ զուգակցում պարունակող գրաֆի աստիճանային հավաքածուն $d = (d_1, \ldots, d_n)$-ն է, և դիցուք $M$-ը $G$ գրաֆի որևէ կատարյալ զուգակցում է: Նկատենք, որ այդ դեպքում $n$-ը զույգ է, $d = (d_1, \ldots, d_n)$-ը գրաֆային է և $G - M$ գրաֆի աստիճանային հավաքածուն $d' = (d_1 - 1, \ldots, d_n - 1)$-ն է, որտեղից հետևում է, որ գրաֆային է նաև $d'$ հավաքածուն։

Հիմա ենթադրենք, որ $n$-ը զույգ է, $d = (d_1, \ldots, d_n)$ և $d' = (d_1 - 1, \ldots, d_n - 1)$ հավաքածուները գրաֆային են, և ցույց տանք, որ այդ դեպքում գոյություն ունի կատարյալ զուգակցում պարունակող $G$ գրաֆ, որի աստիճանային հավաքածուն $d = (d_1, \ldots, d_n)$-ն է։

Քանի որ $d = (d_1, \ldots, d_n)$ հավաքածուն գրաֆային է, ապա գոյություն ունի $G$ գրաֆ, որի աստիճանային հավաքածուն $d$-ն է։ Ենթադրենք, որ $V(G) = \{v_1, \ldots, v_n\}$։ Քանի որ գրաֆային է նաև $d'$ հավաքածուն, ապա գոյություն ունի $G'$ գրաֆ, որի աստիճանային հավաքածուն $d'$-ն է։ Առանց ընդհանրությունը խախտելու, կարող ենք ենթադրել, որ $V(G') = V(G) = \{v_1, \ldots, v_n\}$։

Դիտարկենք վերոհիշյալ պայմաններին բավարարող բոլոր $G$ և $G'$ գրաֆները, և նրանցից ընտրենք $G$ և $G'$ գրաֆներն այնպես, որ $|(E(G) \setminus E(G')) \cup (E(G') \setminus E(G))|$ ամենափոքրն է։ Ցույց տանք, որ այսպիսի ընտրության դեպքում $G' \subseteq G$, այսինքն $G'$-ը $G$ գրաֆի ենթագրաֆ է։ Նկատենք, որ վերջինս կապացուցի թեորեմը, քանի որ եթե $G' \subseteq G$, ապա $E(G) \setminus E(G')$-ը կլինի $G$ գրաֆի կատարյալ զուգակցում։

Ենթադրենք, որ $G' \not\subseteq G$ և դիտարկենք $V(G)$ բազմությանը պատկանող $v$ գագաթը, որը կից է առավելագույն թվով կողերի $E(G') \setminus E(G)$-ից։ Դիցուք այդ թիվը $r$-է։ Այդ դեպքում, ակնհայտ է, որ $v$-ն կից է $r + 1$ կողի $E(G) \setminus E(G')$-ից։ Դիտարկենք $V(G)$ բազմությանը պատկանող ցանկացած $z$ գագաթ, որտեղ $zv \in E(G') \setminus E(G)$, և ընտրենք որևէ $w \in V(G)$ գագաթ այնպես, որ $zw \in E(G) \setminus E(G')$։ Նկատենք, որ այդպիսի $w$ գագաթ գոյություն ունի, քանի որ $d_G(w) = d_{G'}(w) + 1$։

Ցույց տանք, որ ցանկացած $y \in V(G) \setminus \{v, w\}$ համար, եթե $vy \in E(G) \setminus E(G')$, ապա $wy \in E(G) \setminus E(G')$։ Նախ ցույց տանք, որ $wy \in E(G)$։ Ենթադրենք, որ $wy \notin E(G)$։ Դիտարկենք $H$ գրաֆը, որն ստացվում է $G$-ից հետևյալ կերպ. հեռացնենք $vy$ և $zw$ կողերը $G$-ից, և ավելացնենք $vz$ և $wy$ կողերը։ Նկատենք, որ $H$ գրաֆի աստիճանները համընկնում են $G$ գրաֆի աստիճանների հետ, բայց

$$|(E(H) \setminus E(G')) \cup (E(G') \setminus E(H))| < |(E(G) \setminus E(G')) \cup (E(G') \setminus E(G))|,$$

ինչը հակասում է $G$ գրաֆի ընտրությանը։ Հետևաբար, $wy \in E(G)$։ Հիմա ցույց տանք, որ $wy \notin E(G')$։ Ենթադրենք, որ $wy \in E(G')$։ Դիտարկենք $G''$ գրաֆը, որն ստացվում է $G'$-ից հետևյալ կերպ. հեռացնենք $vz$ և $wy$ կողերը $G'$-ից, և ավելացնենք $vy$ և $zw$ կողերը։ Նկատենք, որ $G''$ գրաֆի աստիճանները համընկնում են $G'$ գրաֆի աստիճանների հետ,

բայց

$$|(E(G)\backslash E(G''))\cup(E(G'')\backslash E(G))|<|(E(G)\backslash E(G'))\cup(E(G')\backslash E(G))|,$$

ինչը հակասում է $G'$ գրաֆի ընտրությանը: Հետևաբար, $wy\notin E(G')$ և $wy\in E(G)\backslash E(G')$:

Նկատենք, որ $E(G)\backslash E(G')$ բազմությանը պատկանող այն կողերի քանակը, որոնք կից են $w$ գագաթին, ավելին է, քան նույն բազմությանը պատկանող այն կողերի քանակը, որոնք կից են $v$ գագաթին: Իրոք, ցանկացած $vy\in E(G)\backslash E(G')$ կողի համար $wy\in E(G)\backslash E(G')$ ($y\neq v,w$): Ավելին, մենք ունենք նաև $zw\in E(G)\backslash E(G')$ կողը ($vw$ կողը կարող է պատկանել կամ չպատկանել $E(G)\backslash E(G')$-ին, բայց այն միևնույն ներդրումն է ունենում $d_G(v)$-ում և $d_G(w)$-ում երկու դեպքում էլ): Սա հակասում է $v$ գագաթի ընտրությանն, ինչն ապացուցում է թեորեմը: ∎

## Գլուխ 7

## Հարթ գրաֆներ

### § 7.1. Հարթ գրաֆների պարզագույն հատկությունները

Դիցուք $G = (V, E)$-ն գրաֆ է:

**Սահմանում 7.1.1:** $G$ գրաֆը կոչվում է *հարթ*, եթե այն կարելի է այնպես պատկերել հարթության վրա, որ ցանկացած կող չունենա ինքնահատում և ցանկացած երկու կողեր չունենան ընդհանուր կետեր, բացի զագաթներից:

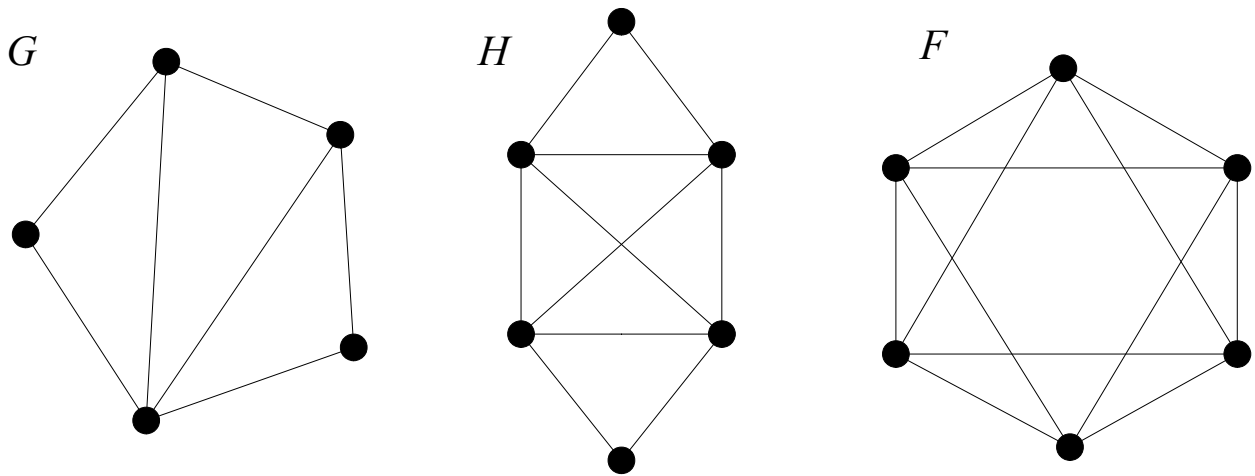

Նկ. 7.1.1

Դիտարկենք նկ. 7.1.1-ում պատկերված $G, H$ և $F$ գրաֆները: Հեշտ է տեսնել, որ $G$-ն հարթ գրաֆ է: Համոզվենք, որ $H$ և $F$ գրաֆները ևս հարթ են: Իրոք, դիտարկենք նկ. 7.1.2-ում նշված $H$ և $F$ գրաֆների պատկերումները: Հեշտ է տեսնել, որ այդ պատկերումներն այնպիսին են, որ ցանկացած կող չունի ինքնահատում և ցանկացած երկու կողեր չունենան ընդհանուր կետեր, բացի զագաթներից: Թեն դիտարկված գրաֆները հարթ են, սակայն գոյություն ունեն նաև գրաֆներ, որոնց հնարավոր չէ պատկերել հարթության վրա այնպես, որ կողերը չհատվեն: Պարզվում է, որ լրիվ երկկողմանի $K_{3,3}$ գրաֆը հնարավոր չէ պատկերել հարթության վրա այնպես, որ կողերը չհատվեն: Այն, որ $K_{3,3}$ գրաֆը հարթ չէ, մենք կապացուցենք քիչ անց, իսկ այժմ ցույց տանք, որ ցանկացած գրաֆ միշտ հնարավոր է այնպես պատկերել $\mathbb{R}^3$-ում, որ ցանկացած կող չունենա ինքնահատում և ցանկացած երկու կողեր չունենան ընդհանուր կետեր, բացի գագաթներից:

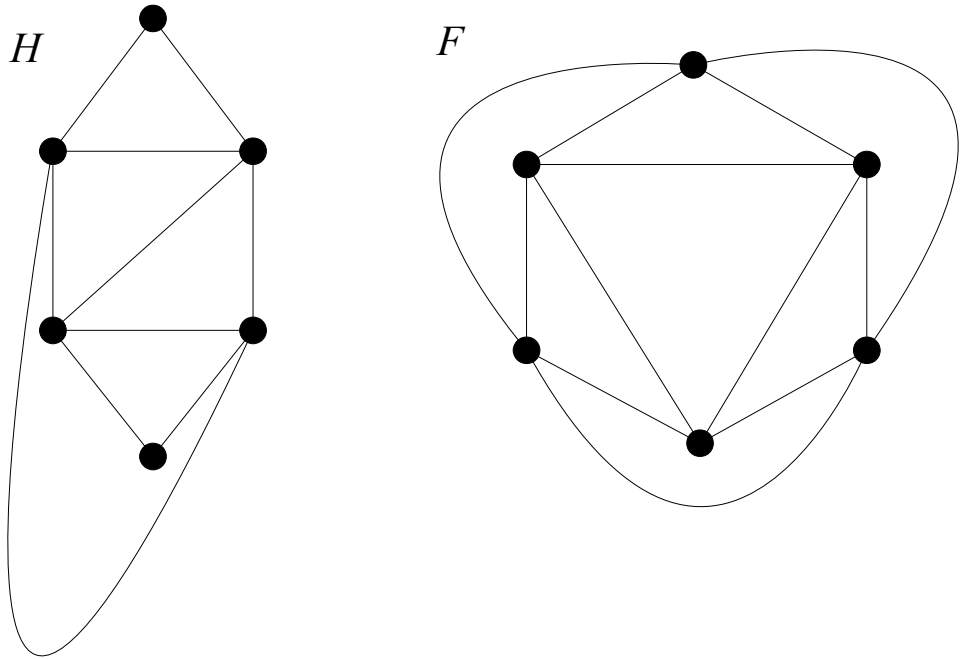

Նկ. 7.1.2

**Թեորեմ 7.1.1:** Կամայական $G$ գրաֆ կարելի է այնպես պատկերել տարածության մեջ, որ ցանկացած կող չունենա ինքնահատում և ցանկացած երկու կողեր չունենան ընդհանուր կետեր, բացի գագաթներից:

**Ապացույց:** Ապացույցի համար դիտարկենք $Ox$ առանցքը և նշենք $Ox$-ի վրա իրարից տարբեր $|V(G)|$ հատ կետեր: Այժմ դիտարկենք $Ox$-ով անցնող հարթությունները և նրանց մեջ ֆիքսենք $|E(G)|$ հատ հարթություն: Այդ հարթությունները ֆիքսելուց հետո յուրաքանչյուր հարթության վրա տանենք մեկական կող կիսաշրջանի տեսքով: Հեշտ է տեսնել, որ այդ պատկերման դեպքում թեորեմի պնդումը դառնում է ակներև: ∎

Դիցուք $G$ հարթ գրաֆը պատկերված է հարթության վրա այնպես, որ ցանկացած կող չունի ինքնահատում և ցանկացած երկու կողեր չունեն ընդհանուր կետեր, բացի գագաթներից: $G$ գրաֆի *նիստ* կոչվում է հարթության այն մաքսիմալ տիրույթը, որտեղ ցանկացած երկու կետ կարող են միացվել գրաֆի կողերը չհատող անընդհատ գծով: Նիստերից մեկը անսահմանափակ է, որը նաև անվանում են *անվերջ* նիստ, մյուսները *սահմանափակ* են: Նիստի *եզր* կանվանենք այդ նիստին պատկանող գագաթների և կողերի բազմությունը: Երկու նիստեր կանվանենք *հարևան*, եթե նրանք ունեն ընդհանուր կող: Դիտարկենք նկ. 7.1.3-ում բերված $G$ գրաֆը: Հեշտ է տեսնել, որ այդ գրաֆի $f_1$-ից $f_8$ նիստերը սահմանափակ են, իսկ $f_9$՝ անսահմանափակ է: Նշենք նաև, որ այդ գրաֆի $f_2$ և $f_3$ նիստերը հարևան են, իսկ $f_5$ և $f_8$-ը՝ հարևան չեն:

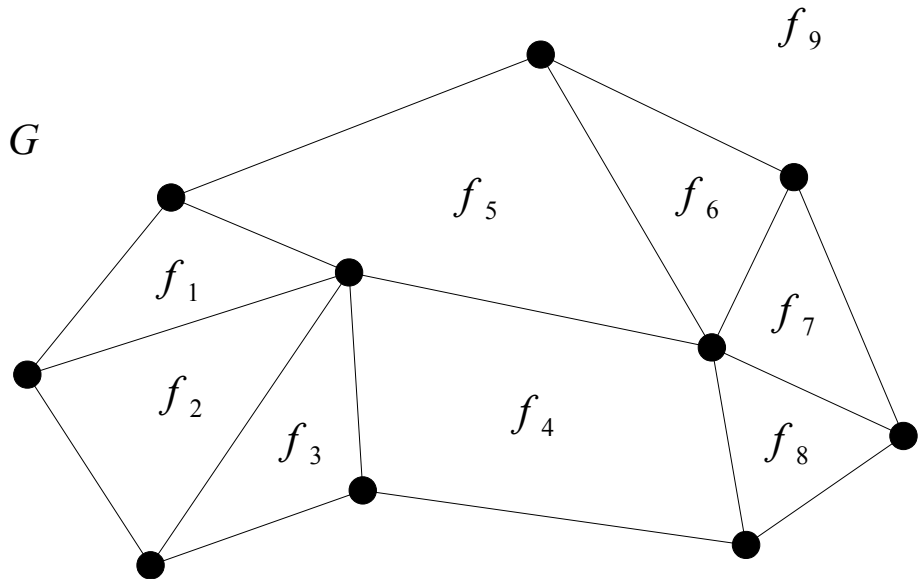

Նկ. 7.1.3

Պարզվում է, որ ցանկացած կապակցված հարթ գրաֆում այդ գրաֆի գագաթների, կողերի և նիստերի քանակների միջև կապ գոյություն ունի, որը հայտնաբերվել է Էյլերի կողմից և հայտնի է որպես Էյլերի բանաձև:

**Թեորեմ 7.1.2 (Լ. Էյլեր):** Եթե $G$-ն կապակցված հարթ $(n, m)$-գրաֆ է, որն ունի $r$ հատ նիստ, ապա տեղի ունի

$$n - m + r = 2$$

հավասարությունը:

**Ապացույց:** Ենթադրենք, որ տրված է $G$ կապակցված հարթ գրաֆը, որը պատկերված է հարթության վրա այնպես, որ ցանկացած կող չունի ինքնահատում և ցանկացած երկու կողեր չունեն ընդհանուր կետեր, բացի գագաթներից:

Նկատենք, որ ցանկացած ծառ հարթ գրաֆ է, որն ունի միայն անսահմանափակ նիստ: Եթե $G$-ն չի պարունակում ցիկլ, ապա այն ծառ է ($r = 1$) և ըստ թեորեմ 2.3.1-ի $n = m + 1$, ուստի $n - m + r = 2$: Այստեղից հետևում է, որ թեորեմը ճիշտ է ծառերի համար:

Այժմ ենթադրենք, որ $G$-ն ծառ չէ: Ըստ թեորեմ 2.3.4-ի $G$-ն ունի $T$ կմախքային ծառ: Նշանակենք այդ $T$ ծառի կողերի քանակը $m_T$-ով, իսկ նիստերի քանակը՝ $r_T$-ով: Քանի որ $T$-ն կմախքային ծառ է, ուստի $m_T = n - 1$ և $r_T = 1$: Այդ կմախքային ծառը ստանալու ժամանակ մենք $G$ գրաֆի ցիկլերից հեռացնում ենք կողեր, պահպանելով կապակցվածությունը: Ցանկացած այդպիսի կող հեռացնելուց տեղի ունի հետևյալը.

1. գագաթների քանակը չի փոխվում,
2. կողերի քանակը մեկով պակասում է,
3. նիստերի քանակը մեկով պակասում է, քանի որ կողը հեռացնելուց հետո այդ կողին հարևան երկու նիստերը միաձուլվում են:

Այստեղից հետևում է, որ $m - r = m_T - r_T$: Մյուս կողմից, ինչպես նշել ենք, $m_T = n - 1$ և $r_T = 1$, ուստի $m - r = m_T - r_T = n - 2$ և, հետևաբար, $n - m + r = 2$: ∎

**Հետևանք 7.1.1**: Լրիվ երկկողմանի $K_{3,3}$ գրաֆը հարթ չէ:

**Ապացույց**: Իրոք, ենթադրենք հակառակը. $K_{3,3}$ գրաֆը հարթ է: Նկատենք, որ $K_{3,3}$-ը կապակցված $(6, 9)$-գրաֆ է: Այդ դեպքում, համաձայն թեորեմ 7.1.2-ի, այն կունենա $r = 2 - 6 + 9 = 5$ նիստ: Քանի որ $K_{3,3}$-ը երկկողմանի գրաֆ է, ուստի նրա յուրաքանչյուր նիստի եզրը պարունակում է առնվազն չորս կող: Այստեղից, հաշվի առնելով, որ յուրաքանչյուր կող մասնակցում է ամենաշատը երկու նիստում, կստանանք, որ $9 = |E(K_{3,3})| \geq \frac{4 \cdot 5}{2} = 10$, ինչը հնարավոր չէ: ∎

**Հետևանք 7.1.2**: Լրիվ $K_5$ գրաֆը հարթ չէ:

**Ապացույց**: Իրոք, ենթադրենք հակառակը. $K_5$ գրաֆը հարթ է: Նկատենք, որ $K_5$-ը կապակցված $(5, 10)$-գրաֆ է: Այդ դեպքում, համաձայն թեորեմ 7.1.2-ի, այն կունենա $r = 2 - 5 + 10 = 7$ նիստ: Քանի որ $K_5$ գրաֆի յուրաքանչյուր նիստի եզրը պարունակում է առնվազն երեք կող, և հաշվի առնելով, որ յուրաքանչյուր կող մասնակցում է ամենաշատը երկու նիստում, կստանանք, որ $10 = |E(K_5)| \geq \frac{3 \cdot 7}{2} > 10$, ինչը հնարավոր չէ: ∎

**Դիտողություն 7.1.1**: Թեորեմ 7.1.2-ից նաև հետևում է, որ ցանկացած կապակցված հարթ $(n, m)$-գրաֆի համար նիստերի քանակը կախված չի լինի հարթության վրա այդ գրաֆի պատկերման ձևից, եթե ցանկացած կող չունենա ինքնահատում և ցանկացած երկու կողեր չունենան ընդհանուր կետեր, բացի գագաթներից: Կամայական այդպիսի պատկերման մեջ նիստերի քանակը կլինի $m - n + 2$:

Հարթ գրաֆների կարևոր հատկություններից է այդ գրաֆներում գագաթների և կողերի քանակների միջև առկա գծային կապը:

**Թեորեմ 7.1.3**: Եթե $G$-ն կապակցված հարթ $(n, m)$-գրաֆ է ($n \geq 3$), ապա տեղի ունի

$m \leq 3n - 6$ անհավասարությունը։

**Ապացույց:** Իրոք, քանի որ յուրաքանչյուր նիստ սահմանափակված է առնվազն երեք կողերով (բացառությամբ այն դեպքի, երբ $G$-ն երեք գագաթ ունեցող ծառ է, որի դեպքում $m \leq 3n - 6$ անհավասարությունը տեղի ունի) և յուրաքանչյուր կող մասնակցում է ամենաշատը երկու նիստում, կստանանք, որ $3r \leq 2m$: Մյուս կողմից, համաձայն թեորեմ 7.1.2-ի, ստանում ենք հետևյալը.

$$2 = n - m + r \leq n - m + \frac{2m}{3} = n - \frac{m}{3},$$

որտեղից ստացվում է $m \leq 3n - 6$ անհավասարությունը։ ∎

**Հետևանք 7.1.3**: Ցանկացած հարթ գրաֆում գոյություն ունի գագաթ, որի աստիճանը հինգից ավել չէ։

**Ապացույց:** Իրոք, ենթադրենք հակառակը. գոյություն ունի այնպիսի հարթ $G$ գրաֆ, որի ցանկացած $v \in V(G)$-ի համար $d_G(v) \geq 6$: Այստեղից և թեորեմ 1.2.1-ից հետևում է, որ $|E(G)| = \frac{1}{2}\sum_{v \in V(G)} d_G(v) \geq 3|V(G)|$, ինչը հակասում է թեորեմ 7.1.3-ին։ ∎

**Սահմանում 7.1.3:** $G$ հարթ գրաֆը կոչվում է *մաքսիմալ հարթ* գրաֆ, եթե այդ գրաֆին ցանկացած նոր կող ավելացնելուց ստացվող գրաֆը հարթ չէ։

Հեշտ է տեսնել, որ $G$-ն մաքսիմալ հարթ գրաֆ է այն և միայն այն դեպքում, երբ այդ գրաֆի յուրաքանչյուր նիստ եռանկյուն է։

Ստորև պատկերված է մաքսիմալ հարթ գրաֆի օրինակ.

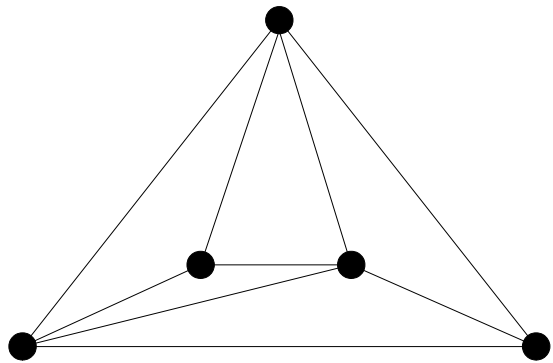

Նկ. 7.1.4

Նշենք նաև առանց ապացույցի Ուիտնիի թեորեմը մաքսիմալ հարթ գրաֆների մասին։

**Թեորեմ 7.1.4:** Առնվազն չորս գագաթ պարունակող ցանկացած մաքսիմալ հարթ գրաֆ **3-կապակցված է**։

Հայտնի է, որ հարթության կետերի և սֆերայի կետերի միջև գոյություն ունի փոխմիարժեք համապատասխանություն։ Պարզվում է, հարթ գրաֆների և սֆերայի վրա պատկերվող գրաֆների միջև ոս գոյություն ունի կապ։ Այստեղ հասկանում ենք, որ գրաֆը պատկերվող է սֆերայի վրա, եթե այդ գրաֆի գագաթներին կարող ենք համապատասխանեցնել սֆերայի կետեր (տարբեր գագաթներին տարբեր կետեր), և եթե երկու գագաթ կազմում են կող գրաֆում, ապա նրանց համապատասխան կետերը միացվում են անընդհատ կորով առանց ինքնահատումների, որը չի անցնում մեկ այլ գագաթին համապատասխան կետով, և տարբեր կողերին համապատասխանող անընդհատ կորերը չունեն ընդհանուր կետեր, բացի գագաթներից:

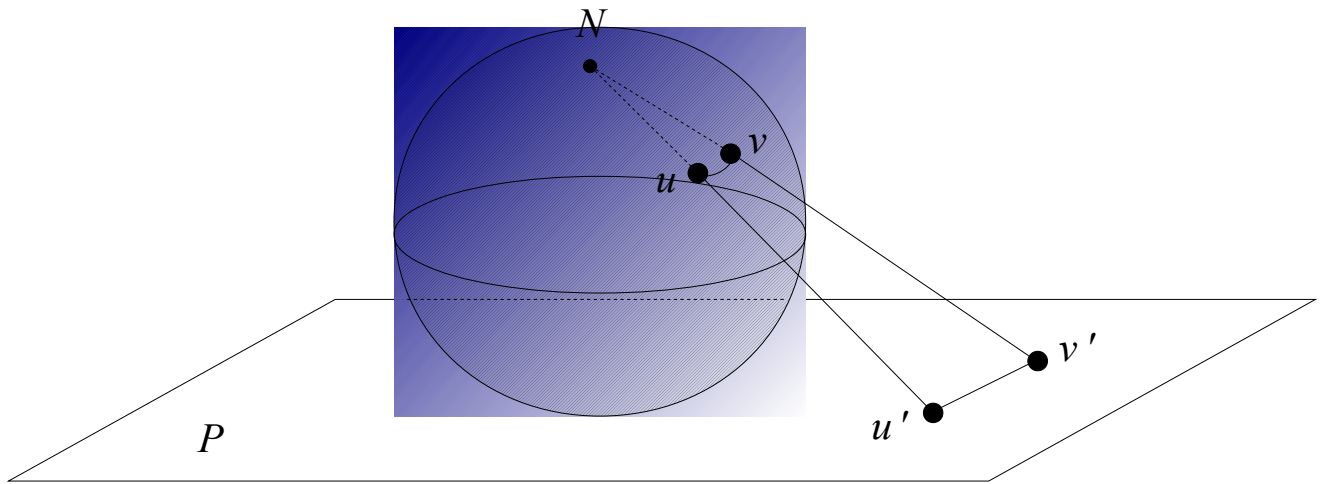

Նկ. 7.1.5

**Թեորեմ 7.1.5:** Գրաֆը սֆերայի վրա պատկերվող է այն և միայն այն դեպքում, երբ այն հարթ է:

**Ապացույց:** Այս թեորեմը ապացուցելու համար բավական է դիտարկել ստերեոգրաֆիկ պրոյեկցիան (նկ. 7.1.5)։ Դիցուք $G$ գրաֆը պատկերված է սֆերայի վրա: Տանենք այդ սֆերային շոշափող $P$ հարթություն այնպիսի կետում, որի տրամագծորեն հակառակ $N$ կետը («հյուսիսային բևեռը») չգտնվի $G$ գրաֆի կողի վրա կամ հանդիսանա այդ գրաֆի գագաթ: Այժմ դիտարկենք $G'$ գրաֆը, որը ստացվել է $N$ կետից $P$ հարթություն վրա $G$ գրաֆի ստերեոգրաֆիկ պրոյեկցիայի միջոցով: Քանի որ գոյություն ունի փոխմիարժեք համապատասխանություն սֆերայի $N$ կետից տարբեր կետերի և նրանց ստերեոգրաֆիկ պրոյեկցիաների միջև, ուստի $G'$ գրաֆը պատկերված է $P$ հարթության վրա այնպես, որ ցանկացած կող չունի ինքնահատում և ցանկացած երկու կողեր չունեն ընդհանուր կետեր, բացի գագաթներից: Այստեղից հետևում է, որ $G'$ գրաֆը հարթ է և իզոմորֆ է $G$-ին:

Համանման ձևով ապացուցվում է հակառակ պնդումը, հաշվի առնելով վերը նշված փոխմիարժեք համապատասխանությունը: ∎

**Թեորեմ 7.1.6:** Եթե $G$ գրաֆը պատկերված է հարթության վրա այնպես, որ ցանկացած կող չունի ինքնահատում և ցանկացած երկու կողեր չունեն ընդհանուր կետեր, բացի գագաթներից, և $E$-ն այդ գրաֆի որևէ նիստի եզրի կողերի բազմությունն է, ապա $G$ գրաֆը կարելի է այնպես պատկերել հարթության վրա, որ ցանկացած կող չունենա ինքնահատում և ցանկացած երկու կողեր չունենան ընդհանուր կետեր, բացի գագաթներից, և $E$-ն հանդիսանա անսահմանափակ նիստի եզրի կողերի բազմություն:

**Ապացույց:** Թեորեմն ապացուցելու համար բավական է դիտարկել ստերեոգրաֆիկ պրոյեկցիան:

Դիցուք $E$-ն $f$ նիստի եզրի կողերի բազմությունն է: Եթե $f$-ը անսահմանափակ նիստ է, ապա ակնհայտ է, որ $E$-ն անսահմանափակ նիստի եզրի կողերի բազմություն է: Ենթադրենք, $f$-ը սահմանափակ նիստ է: Ընտրենք այդ $f$ նիստի որևէ ներքին $N$ կետ և անվանենք այն «հյուսիսային բևեռ»: Այնուհետև պատկերենք $G$ գրաֆը սֆերայի վրա վերը նշված եղանակով: Տանենք այդ սֆերային շոշափող $Q$ հարթություն այն կետում, որը տրամագծորեն հակառակ է «հյուսիսային բևեռին»: Այժմ դիտարկենք $G'$ գրաֆը, որը ստացվել է $N$ կետից $Q$ հարթության վրա $G$ գրաֆի ստերեոգրաֆիկ պրոյեկցիայի միջոցով: Պարզ է, որ $G'$ գրաֆը իզոմորֆ է $G$-ին: Քանի որ գոյություն ունի փոխմիարժեք համապատասխանություն սֆերայի $N$ կետից տարբեր կետերի և նրանց ստերեոգրաֆիկ պրոյեկցիաների միջև, ուստի $G'$ գրաֆը պատկերված է $Q$ հարթության վրա այնպես, որ ցանկացած կող չունի ինքնահատում և ցանկացած երկու կողեր չունեն ընդհանուր կետեր, բացի գագաթներից: Այստեղից հետևում է, որ $G'$ գրաֆը հարթ է: Մյուս կողմից, հեշտ է տեսնել, որ այդ պատկերման դեպքում $f$-ը կհանդիսանա անսահմանափակ նիստ, իսկ $E$-ն` այդ նիստի եզրի կողերի բազմություն: ∎

Նշենք առանց ապացույցի ևս մի թեորեմ, որն ապացուցվել է Ֆարիի կողմից և հայտնի է նրա անունով:

**Թեորեմ 7.1.7:** Ցանկացած հարթ գրաֆ կարելի է այնպես պատկերել հարթության վրա, որ ցանկացած կող հանդիսանա հատված և ցանկացած երկու կողեր չունենան ընդհանուր կետեր, բացի գագաթներից:

Այս պարագրաֆի վերջում բերենք հայտնի **3-**, **4-** և **5-**համասեռ հարթ գրաֆների օրինակներ (նկ. 7.1.6):

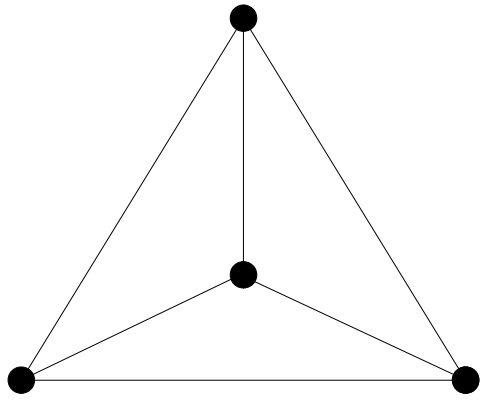
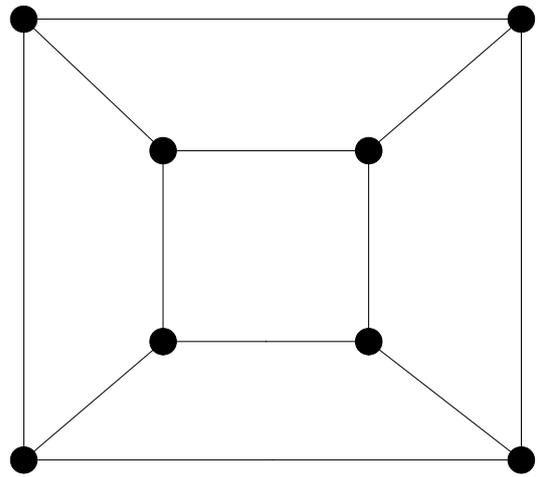
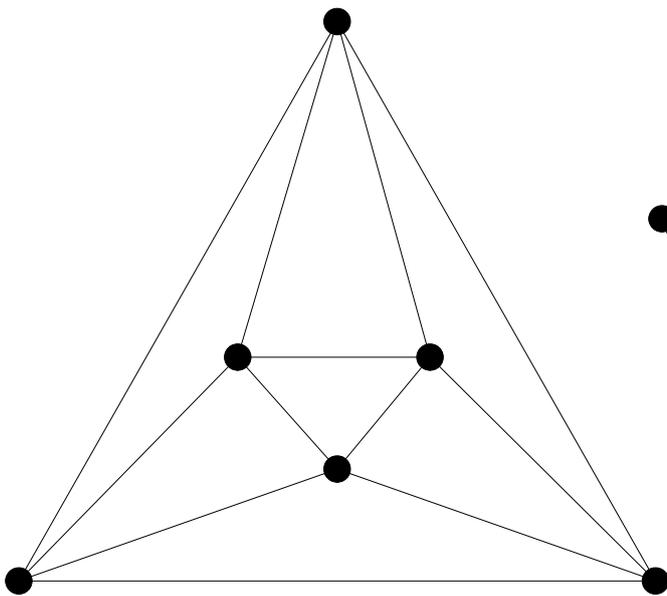
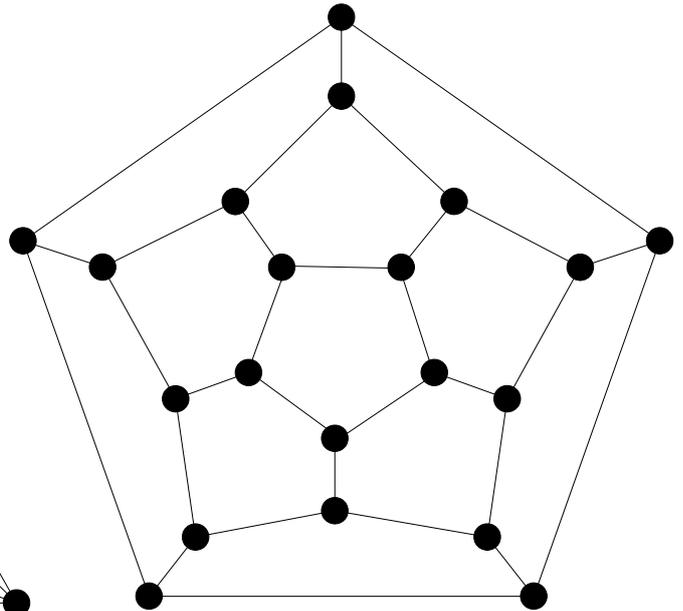
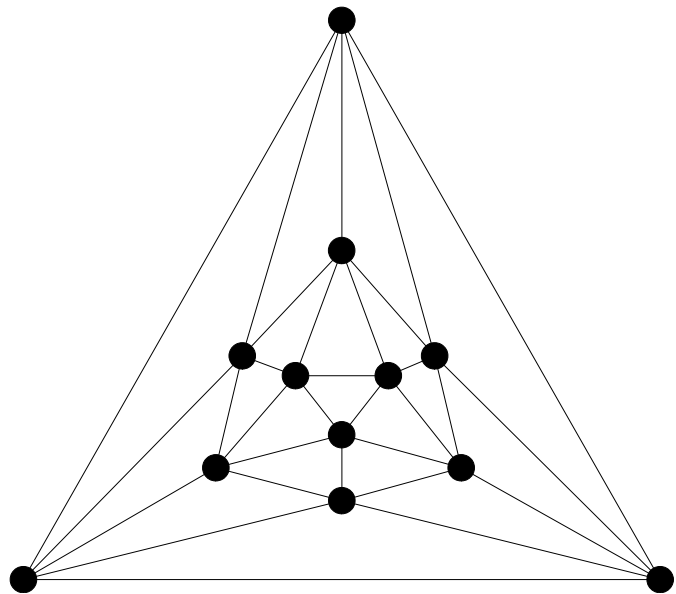

Նկ. 7.1.6

# § 7.2. Հարթ գրաֆների նկարագրությունը. Պոնտրյագին-Կուրատովսկու և Վագների թեորեմները

Դիցուք $G = (V, E)$-ն գրաֆ է։

**Սահմանում 7.2.1:** Կասենք, որ $H$ գրաֆը հանդիսանում է $G$ գրաֆի *ենթատրոհում*, եթե $H$ գրաֆը ստացվում է $G$ գրաֆից կողի տրոհում գործողության հաջորդական կիրառումների միջոցով։

Եթե $\delta(G) \geq 3$ և $H$ գրաֆը $G$ գրաֆի ենթատրոհում է, ապա $H$ գրաֆի այն $v$ գագաթները, որոնց համար $d_H(v) \geq 3$, կանվանենք *իրական գագաթներ*։ Նկատենք, որ իրական գագաթները հանդիսանում են $G$ գրաֆի գագաթների պատկերներ $H$ գրաֆում։

Նկ. 7.2.1-ում պատկերված են $G$ և $H$ գրաֆները։ Հեշտ է տեսնել, որ $H$ գրաֆը հանդիսանում է $G$ գրաֆի ենթատրոհում։ Նկատենք նաև, որ $H$ գրաֆի $u_1, u_2, u_3$ և $w_1, w_2, w_3$ գագաթները իրական են։

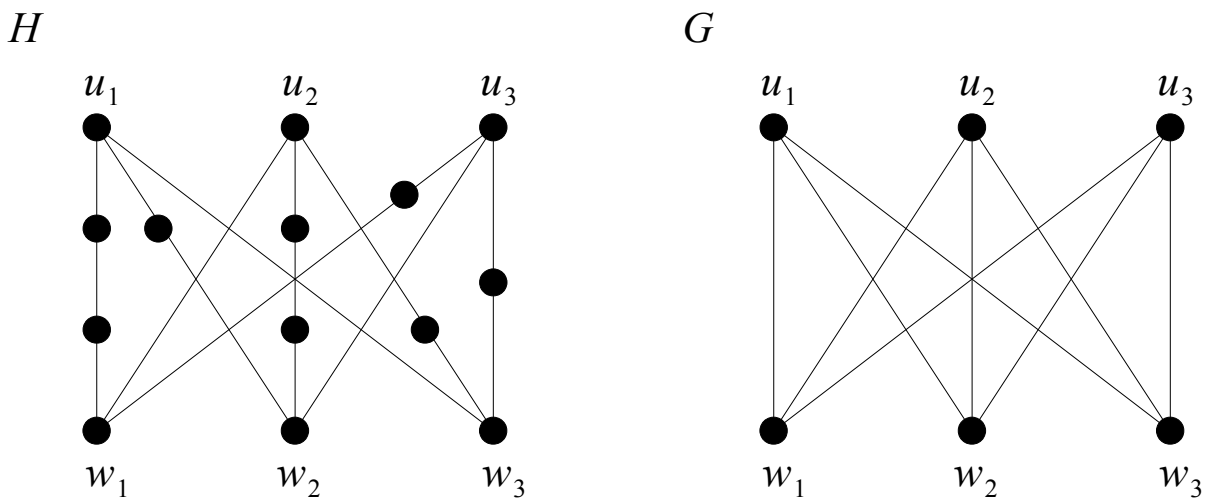

Նկ. 7.2.1

**Սահմանում 7.2.2:** Կասենք, որ $H$ գրաֆը հանդիսանում է $G$ գրաֆի *մինոր*, եթե $H$ գրաֆը ստացվում է $G$ գրաֆի ենթագրաֆից կողի կծկում գործողության հաջորդական կիրառումների միջոցով։

Նկ. 7.2.2-ում պատկերված են $G$ և $H$ գրաֆները։ Ցույց տանք, որ $H$ գրաֆը հանդիսանում է $G$ գրաֆի մինոր։ Իրոք, դիտարկենք $G$ գրաֆը և այդ գրաֆի $u_1v_1, u_2v_2, u_3v_3, u_4v_4$ և $u_5v_5$ կողերը։ Այդ կողերի հաջորդական կծկումների արդյունքում կառաջանան $H$ գրաֆի $w_1, w_2, w_3, w_4$ և $w_5$ գագաթները, իսկ ստացված գրաֆը կլինի $H$-ը, որը իզոմորֆ է $K_5$ լրիվ գրաֆին։

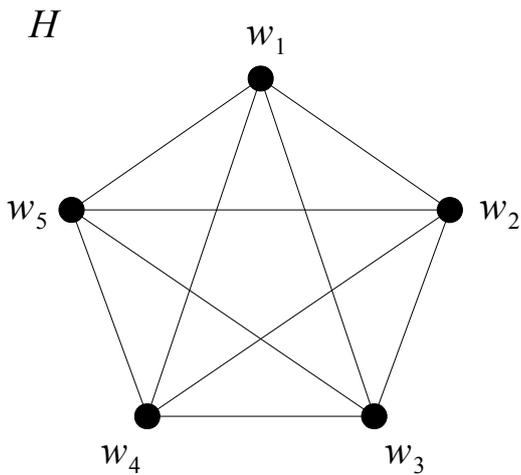
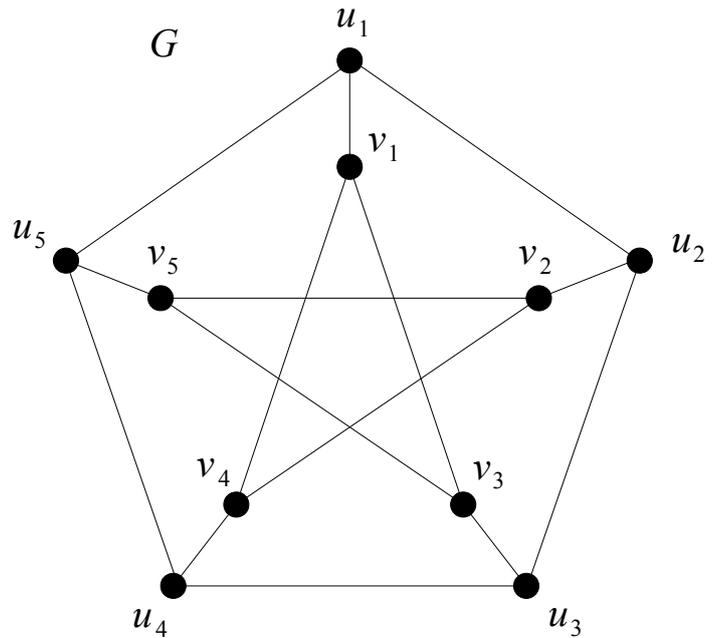

Նկ. 7.2.2

Այս պարագրաֆում մեր հիմնական նպատակը կլինի տալ հարթ գրաֆների նկարագրությունը։ Քանի որ կողի տրոհման գործողությունը չի խախտում գրաֆի հարթ լինելու հատկությունը, ուստի մենք կփորձենք գտնել այն տոպոլոգիական մինիմալ ոչ հարթ գրաֆները, որոնք չեն հանդիսանում այլ ոչ հարթ գրաֆների ենթատրոհումներ։ Մյուս կողմից պարզ է, որ ցանկացած գրաֆ, որը պարունակում է ենթագրաֆ, որը $K_{3,3}$-ի կամ $K_5$-ի ենթատրոհում է, հարթ լինել չի կարող։ Հետագայում, պարզության համար, $G$ գրաֆի այն ենթագրաֆը, որը $K_{3,3}$-ի կամ $K_5$-ի ենթատրոհում է, կանվանենք $G$ գրաֆի $PK$-*ենթագրաֆ*։ *Մինիմալ ոչ հարթ գրաֆ* կանվանենք այն $G$ գրաֆը, որը հարթ չէ, բայց ցանկացած $e \in E(G)$ կողի համար $G - e$ գրաֆը հարթ է։ Պարզ է, որ մինիմալ ոչ հարթ գրաֆը կապակցված գրաֆ է։

Այժմ անցնենք օժանդակ լեմմաների ապացուցմանը։

**Լեմմա 7.2.1:** Ցանկացած մինիմալ ոչ հարթ գրաֆ 2-կապակցված է։

**Ապացույց:** Դիցուք $G$-ն մինիմալ ոչ հարթ գրաֆ է։ Քանի որ $G$-ն կապակցված գրաֆ է, ուստի լեմման ապացուցելու համար բավական է ցույց տալ, որ $G$-ն չի պարունակում միակցման կետեր։ Դիցուք $v$ գագաթը $G$ գրաֆի միակցման կետ է։ Ենթադրենք $H_1, H_2, \ldots, H_k$-ն $G - v$ գրաֆի կապակցված բաղադրիչներն են։ Դիտարկենք $G$ գրաֆի $G_1, G_2, \ldots, G_k$ ենթագրաֆները, որտեղ $G_i = G[V(H_i) \cup \{v\}]$ ($1 \leq i \leq k$)։ Քանի որ $G$-ն մինիմալ ոչ հարթ գրաֆ է, ուստի $G_1, G_2, \ldots, G_k$ ենթագրաֆները հարթ են։ Համաձայն թեորեմ 7.1.6-ի, $G_1, G_2, \ldots, G_k$ ենթագրաֆները կարելի է այնպես պատկերել հարթության վրա, որ ցանկացած կող չունենա ինքնահատում և ցանկացած երկու կողեր չունենան ընդհանուր կետեր, բացի գագաթներից, և $v$ գագաթը հանդիսանա անսահմանափակ նիստի եզրի գագաթ այդ բոլոր ենթագրաֆներում: Պատկերելով այդ ենթագրաֆները $\frac{360°}{k}$-ից փոքր չափի և $v$-ն որպես միակ ընդհանուր կետ ունեցող $k$ հատ տարբեր անկյունների ներսում, մենք կստանանք հարթության վրա $G$ գրաֆի այնպիսի պատկերում, որի դեպքում ցանկացած կող չունի ինքնահատում և ցանկացած երկու կողեր չունեն ընդհանուր կետեր, բացի գագաթներից: Այստեղից հետևում է, որ $G$ գրաֆը հարթ է, ինչը հակասություն է: ∎

**Լեմմա 7.2.2:** Դիցուք $G$-ն գրաֆ է, $S = \{x, y\}$ և $c(G - S) \geq 2$: Ենթադրենք, $G_1$ և $G_2$-ը $G$ գրաֆի ենթագրաֆներ են, որոնց համար $V(G_1) \cap V(G_2) = S$ և $E(G_1) \cup E(G_2) = E(G)$: Այդ դեպքում, եթե $G$ գրաֆը հարթ չէ, ապա $H_1$ $\left(H_1 = \begin{cases} G_1 + xy, & \text{եթե } xy \notin E(G), \\ G_1, & \text{եթե } xy \in E(G) \end{cases}\right)$ և $H_2$ $\Big(H_2 =$ $\begin{cases} G_2 + xy, & \text{եթե } xy \notin E(G), \\ G_2, & \text{եթե } xy \in E(G) \end{cases}\Big)$ գրաֆներից առնվազն մեկը ևս հարթ չէ:

**Ապացույց:** Ենթադրենք, $H_1$ և $H_2$-ը հարթ են: Այդ դեպքում, համաձայն թեորեմ 7.1.6-ի, $H_1$ և $H_2$ գրաֆները կարելի է այնպես պատկերել հարթության վրա, որ ցանկացած կող չունենա ինքնահատում և ցանկացած երկու կողեր չունենան ընդհանուր կետեր, բացի գագաթներից, և $xy$ կողը հանդիսանա անսահմանափակ նիստի եզրի կող այդ գրաֆներում: Պատկերելով $H_1$ և $H_2$ գրաֆները, նույնացնենք այդ պատկերների $xy$ կողը, և, եթե $xy \notin E(G)$, ապա նետենք այդ կողը: Մենք կստանանք հարթության վրա $G$ գրաֆի այնպիսի պատկերում, որի դեպքում ցանկացած կող չունի ինքնահատում և ցանկացած երկու կողեր չունեն ընդհանուր կետեր, բացի գագաթներից: Այստեղից հետևում է, որ $G$ գրաֆը հարթ է, իսկ դա հակասություն է: ∎

**Լեմմա 7.2.3:** Եթե $G$-ն ամենաքիչ կողեր պարունակող այնպիսի գրաֆ է, որը հարթ չէ և չունի $PK$-ենթագրաֆ, ապա $G$-ն 3-կապակցված է:

**Ապացույց:** Նախ նկատենք, որ $G$ գրաֆից կող հեռացնելու արդյունքում գրաֆում չի առաջանում $PK$-ենթագրաֆ: Քանի որ ցանկացած $e \in E(G)$ կողի համար $G - e$ ենթագրաֆը հարթ է, ուստի $G$-ն մինիմալ ոչ հարթ գրաֆ է: Համաձայն լեմմա 7.2.1-ի, այդ $G$ գրաֆը 2-կապակցված է: Դիցուք $S = \{x, y\}$ և $c(G - S) \geq 2$: Քանի որ $G$-ն հարթ չէ, ուստի, ըստ լեմմա 7.2.2-ի, $H_1$ և $H_2$ գրաֆներից առնվազն մեկը հարթ չէ: Ենթադրենք, որ $H_1$-ը հարթ չէ: Քանի որ $|E(H_1)| < |E(G)|$, ուստի $H_1$-ը պարունակում է $PK$-ենթագրաֆ:

Պարզ է, որ այդ ենթագրաֆը նաև մասնակցում է $G$ գրաֆում, բացի, մի գուցե, $xy$ կողից։ Մյուս կողմից, քանի որ $G$ գրաֆը 2-կապակցված է, ուստի $H_2$-ում գոյություն կունենա $(x,y)$-ճանապարհի, որը $G$ գրաֆի $PK$-ենթագրաֆում կհամապատասխանի $xy$ կողին։ Այստեղից հետևում է, որ $G$ գրաֆը պարունակում է $PK$-ենթագրաֆ, իսկ դա հակասություն է։ ∎

**Լեմմա 7.2.4:** Եթե $G/e$-ն ($e \in E(G)$) պարունակում է $PK$-ենթագրաֆ, ապա $G$-ն ևս պարունակում է $PK$-ենթագրաֆ։

**Ապացույց:** Դիցուք $e = xy$ և $G' = G/e$։ Դիցուք նաև $H$-ը $G'$ գրաֆի $PK$-ենթագրաֆ է, ընդ որում $z$-ը $G$ գրաֆի $e = xy$ կողի կծկման արդյունքում ստացված գագաթն է։ Եթե $z$-ը $H$ գրաֆի իրական գագաթ չէ, ապա ակնհայտ է, որ $G$ գրաֆը ևս պարունակում է $PK$-ենթագրաֆ։ Եթե $z$-ը $H$ գրաֆի իրական գագաթ է և այդ գագաթին կից ամենաշատը մեկ կողն է կից $x$ գագաթին $G$ գրաֆում, ապա $z$ գագաթից մենք կարող ենք անցնել $xy$ կողի $G$ գրաֆում և այդ դեպքում $y$ գագաթը կհանդիսանա $G$ գրաֆի $PK$-ենթագրաֆի իրական գագաթ։ Այսպիսով, այս դեպքում ևս $G$ գրաֆը պարունակում է $PK$-ենթագրաֆ։

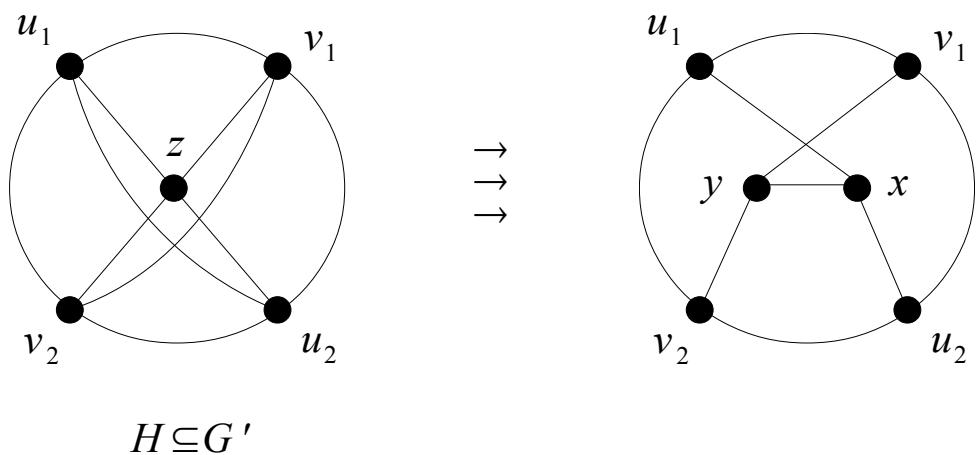

$H \subseteq G'$

Նկ. 7.2.3

Միակ դեպքը, որը մենք չենք դիտարկել, այն է, որ $H$-ը հանդիսանում է $K_5$ գրաֆի ենթատրոհում, $z$-ը $H$ գրաֆի իրական գագաթ է և $x$ և $y$ գագաթներից յուրաքանչյուրը $G$ գրաֆում կից է $z$ գագաթին կից չորս կողերից երկուսին (նկ. 7.2.3)։ Դիցուք $u_1$ և $u_2$-ը $H$ գրաֆի այն իրական գագաթներն են, որոնք $z$ գագաթով սկսվող և $G$ գրաֆում $x$ գագաթին կից կողով անցնող ճանապարհների մյուս ծայրակետերն են, իսկ $v_1$ և $v_2$-ը՝ $H$ գրաֆի այն իրական գագաթներն են, որոնք $z$ գագաթով սկսվող և $G$ գրաֆում $y$ գագաթին կից կողով անցնող ճանապարհների մյուս ծայրակետերն են։ Այդ դեպքում, դեն նետելով $H$ գրաֆից ($u_1, u_2$)-ճանապարհը և ($v_1, v_2$)-ճանապարհը և անցնելով $z$ գագաթից $xy$ կողի, մենք կստանանք $G$ գրաֆում $K_{3,3}$ գրաֆի ենթատրոհում, որտեղ $y, u_1, u_2$-ը մի կողմի իրական գագաթներ են, իսկ $x, v_1, v_2$-ը՝ մյուս կողմի (նկ. 7.2.3): ∎

Նախորդ պարագրաֆում մենք նշել ենք Ֆարիի թեորեմը, որը պնդում է, որ ցանկացած հարթ $G$ գրաֆ կարելի է այնպես պատկերել հարթության վրա, որ ցանկացած կող հանդիսանա հատված և ցանկացած երկու կողեր չունենան ընդհանուր կետեր, բացի գագաթներից: Եթե նաև պահանջենք, որ յուրաքանչյուր սահմանափակ նիստ այդ պատկերման դեպքում հանդիսանա ուռուցիկ բազմանկյուն, ապա կասենք, որ $G$ գրաֆն ունի *ուռուցիկ պատկերում*: Այսպես, օրինակ, նկ. 7.1.6-ում բերված են որոշ 3-համասեռ, 4-համասեռ և 5-համասեռ հարթ գրաֆների ուռուցիկ պատկերումներ: Նշենք, որ ոչ բոլոր հարթ գրաֆներն ունեն ուռուցիկ պատկերումներ. օրինակ, հեշտ է տեսնել, որ 2-կապակցված հարթ $K_{2,n}$ գրաֆը չունի այդպիսի պատկերում, երբ $n \geq 4$: Սակայն, պարզվում է, որ բոլոր 3-կապակցված հարթ գրաֆներն ունեն այդպիսի պատկերում: Մենք կբերենք այդ փաստի Տոմասենի ապացույցը՝ միաժամանակ ապացուցելով հարթ գրաֆների նկարագրությունը:

**Թեորեմ 7.2.1 (Տատտ):** Եթե $G$-ն 3-կապակցված գրաֆ է և այն չի պարունակում $PK$-ենթագրաֆ, ապա $G$ գրաֆը ունի ուռուցիկ պատկերում:

**Ապացույց:** Ապացույցը կատարենք մակածման եղանակով ըստ $|V(G)|$-ի: Նախ նկատենք, որ չորսից ոչ ավելի գագաթ ունեցող գրաֆների համար թեորեմի պնդումը ճիշտ է: Իրոք, այդ պայմանին բավարարող միակ 3-կապակցված գրաֆը $K_4$-ն է, իսկ այդ գրաֆի ուռուցիկ պատկերումը բերված է նկ. 7.1.7-ում: Ենթադրենք, $|V(G)| \geq 5$ և թեորեմի պնդումը ճիշտ է ցանկացած 3-կապակցված $G'$ գրաֆի համար, որը չի պարունակում $PK$-ենթագրաֆ, երբ $|V(G')| < |V(G)|$: Դիտարկենք 3-կապակցված $G$ գրաֆը: Համաձայն թեորեմ 3.3.6-ի, $G$ գրաֆում գոյություն ունի $e = xy$ կող, որ $G/e$ գրաֆը ևս 3-կապակցված գրաֆ է, ընդ որում $z$-ը $G$ գրաֆի $e = xy$ կողի կծկման արդյունքում ստացված գագաթն է: Համաձայն լեմմա 7.2.4-ի, $G/e$ գրաֆը ևս չի պարունակում $PK$-ենթագրաֆ և, հետևաբար, ըստ մակածման ենթադրության, $G/e$ գրաֆը ունի ուռուցիկ պատկերում: Դիտարկենք այդ պատկերումը: Պարզ է, որ $G/e - z$ գրաֆի պատկերում գոյություն ունի նիստ (այդ նիստը կարող է լինել նաև անսահմանափակ նիստը), որը պարունակում է $z$ գագաթը: Դիցուք այդ նիստը սահմանափակված է $C$ պարզ ցիկլով:

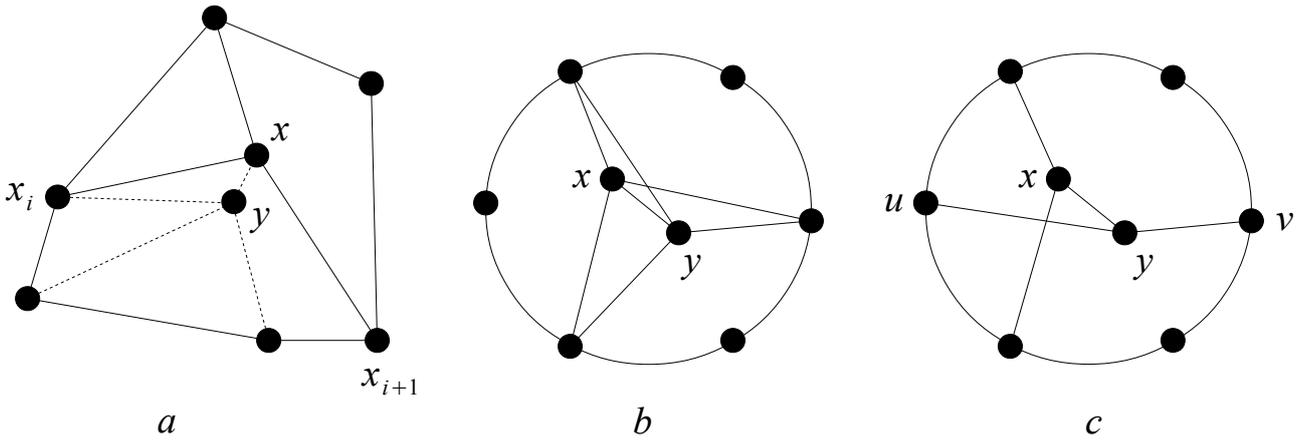

նկ. 7.2.4

Քանի որ $G/e$ գրաֆը ունի ուռուցիկ պատկերում, ուստի $z$ գագաթը միացված է իր հարևան գագաթների հետ հատվածներով։ Դիցուք $x$ գագաթի հարևան $x_1, \ldots, x_k$ ($k \geq 3$) գագաթները $C$ պարզ ցիկլի վրայով շարժվելիս հերթականորեն հանդիպում են հենց այդ հաջորդականությամբ։ Եթե $y$ գագաթի բոլոր հարևան գագաթները պատկանում են $C$ պարզ ցիկլի որևէ $(x_i, x_{i+1})$-ճանապարհին, ապա մենք կարող ենք ստանալ $G$ գրաֆի ուռուցիկ պատկերում, տեղադրելով $x$ գագաթը $z$ գագաթի կետում $G/e$ գրաֆի ուռուցիկ պատկերման դեպքում, իսկ $y$ գագաթը տեղադրելով $x$ գագաթին մոտ կետում՝ $xx_i$ և $xx_{i+1}$ կողերի արանքում (նկ. 7.2.4 a)։ Հակառակ դեպքում՝ $x$ և $y$ գագաթները ունեն երեք ընդհանուր հարևան գագաթներ (նկ. 7.2.4 b) կամ $y$ գագաթը ունի երկու $u$ և $v$ հարևան գագաթներ $C$ պարզ ցիկլի վրա, որոնք պատկանում են $C - x_i - x_{i+1}$ ենթագրաֆի կապակցվածության տարբեր բաղադրիչներին (նկ. 7.2.4 c)։ Առաջին դեպքում հեշտ է տեսնել, որ $G$ գրաֆը պարունակում է $K_5$ գրաֆի ենթատրոհում, իսկ երկրորդ դեպքում՝ $C$-ն, $uyv$, $x_i x x_{i+1}$ ճանապարհները և $xy$ կողը կազմում են $K_{3,3}$ գրաֆի ենթատրոհում։ ∎

Այժմ մենք կարող ենք ձևակերպել հարթ գրաֆ լինելու հայտանիշը։

**Թեորեմ 7.2.2 (Պոնտրյագին, Կուրատովսկի):** $G$ գրաֆը հարթ է այն և միայն այն դեպքում, երբ այն չի պարունակում ենթագրաֆ, որը $K_{3,3}$-ի կամ $K_5$-ի ենթատրոհում է։

Նկատենք, որ այս թեորեմի ապացույցը բխում է հետևանք 7.1.1, 7.1.2, լեմմա 7.2.3 և թեորեմ 7.2.1-ից։

Նշենք նաև առանց ապացույցի հարթ գրաֆների ևս մի նկարագրում։

**Թեորեմ 7.2.3 (Վագներ):** $G$ գրաֆը հարթ է այն և միայն այն դեպքում, երբ $G$-ն չի պարունակում $K_{3,3}$ կամ $K_5$ որպես մինոր։

Նկատենք, որ եթե $G$ գրաֆը պարունակում է ենթագրաֆ, որը $K_{3,3}$ կամ $K_5$-ի ենթատրոհում է, ապա ակնհայտ է, որ այն նաև կպարունակի $K_{3,3}$ կամ $K_5$ որպես մինոր։ Մյուս կողմից պարզվում է, հակառակ պնդումը ևս ճիշտ է. եթե $G$ գրաֆը պարունակում է $K_{3,3}$ կամ $K_5$ որպես մինոր, ապա այդ գրաֆը նաև կպարունակի ենթագրաֆ, որը $K_{3,3}$ կամ $K_5$-ի ենթատրոհում է։ Այսպիսով, թեորեմ 7.2.2-ը և 7.2.3-ը համարժեք են։

## § 7.3. Արտաքին հարթ գրաֆները և Հարարի-Չարտրանդի թեորեմը

Դիցուք $G = (V, E)$-ն հարթ գրաֆ է։

**Սահմանում 7.3.1**։ $G$ հարթ գրաֆը կոչվում է *արտաքին հարթ* գրաֆ, եթե այն կարելի է այնպես պատկերել հարթության վրա, որ ցանկացած կող չունենա ինքնահատում, ցանկացած երկու կողեր չունենան ընդհանուր կետեր (բացի գագաթներից), և բոլոր գագաթները պատկանեն միևնույն նիստին։

Նկատենք, որ թեորեմ 7.1.6-ից հետևում է, որ մենք միշտ կարող ենք ենթադրել, որ արտաքին հարթ գրաֆում բոլոր գագաթները պարունակող նիստը անսահմանափակ նիստն է։

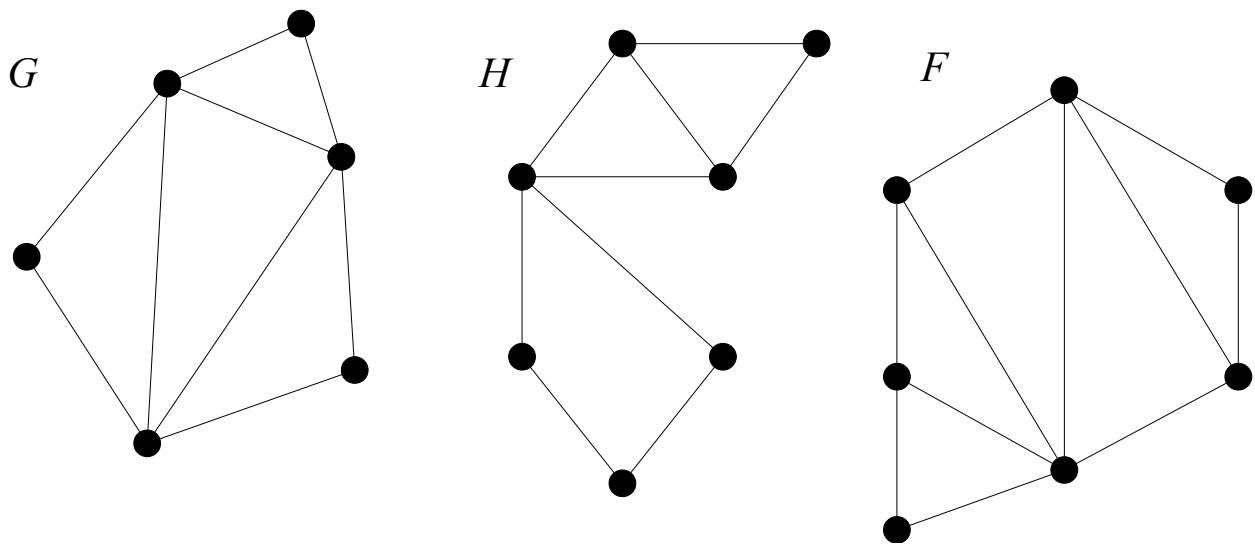

Նկ. 7.3.1

Դիտարկենք նկ. 7.3.1-ում պատկերված $G, H$ և $F$ գրաֆները։ Հեշտ է տեսնել, որ $G, H$ և $F$ գրաֆները արտաքին հարթ են, քանի որ այդ հարթ գրաֆների բոլոր գագաթները պատկանում են անսահմանափակ նիստին։ Չնայած նկ. 7.3.1-ում պատկերված գրաֆներին, գոյություն ունեն նաև գրաֆներ, որոնք արտաքին հարթ չեն։ Այսպես, օրինակ, հեշտ է տեսնել, որ նկ. 7.3.2-ում պատկերված երկու գրաֆները արտաքին հարթ չեն։

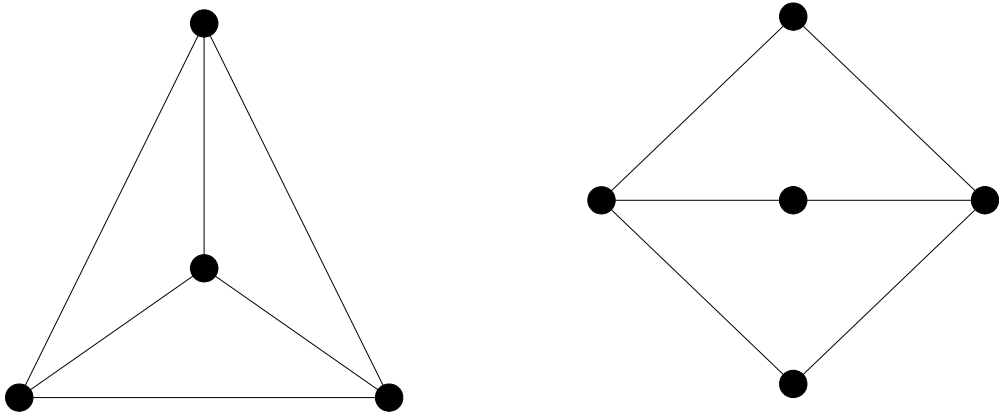

Նկ. 7.3.2

Ստորև կձևակերպենք և կապացուցենք արտաքին հարթ գրաֆների մի նկարագրում:

**Թեորեմ 7.3.1:** $G$ գրաֆը արտաքին հարթ է այն և միայն այն դեպքում, երբ $G + K_1$ գրաֆը հարթ է:

**Ապացույց:** Նախ ենթադրենք, $G$ գրաֆը արտաքին հարթ է և այն պատկերված է հարթության վրա այնպես, որ ցանկացած կող չունի ինքնահատում, ցանկացած երկու կողեր չունեն ընդհանուր կետեր (բացի գագաթներից), և բոլոր գագաթները պատկանում են անսահմանափակ նիստին: Այդ դեպքում, տեղադրելով $v$ գագաթը ($v \notin V(G)$) անսահմանափակ նիստի որևէ կետում և միացնելով այդ գագաթը $G$ գրաֆի բոլոր գագաթների հետ այնպես, որ միացման արդյունքում չառաջանան կողերի հատման կետեր, մենք կստանանք, որ $G + K_1$ գրաֆը հարթ է:

Այժմ ենթադրենք, որ $G + K_1$ գրաֆը հարթ է: Այդ դեպքում $G + K_1$ գրաֆը պարունակում է $v$ գագաթ, որը հարևան է $G$ գրաֆի բոլոր գագաթներին: Դիցուք $G + K_1$ հարթ գրաֆը պատկերված է հարթության վրա այնպես, որ ցանկացած կող չունի ինքնահատում և ցանկացած երկու կողեր չունեն ընդհանուր կետեր (բացի գագաթներից): Հեռացնելով այդ պատկերման մեջ մասնակցող $v$ գագաթը և նրան կից կողերը մենք կստանանք, որ $G$ գրաֆի բոլոր գագաթները պատկանում են միևնույն նիստին, ուստի $G$-ն արտաքին հարթ գրաֆ է: ∎

Ինչպես նշել ենք, $K_{2,3}$ և $K_4$ գրաֆները հարթ են, բայց արտաքին հարթ չեն (նկ. 7.3.2): Պարզվում է, այդ գրաֆները կարևոր դեր ունեն արտաքին հարթ գրաֆների նկարագրության մեջ:

**Թեորեմ 7.3.2 (Հարարի, Չարտրանդ):** $G$ գրաֆն արտաքին հարթ է այն և միայն այն դեպքում, երբ այն չի պարունակում ենթագրաֆ, որը $K_{2,3}$-ի կամ $K_4$-ի ենթատրոհում է:

**Ապացույց։** Նախ ենթադրենք. $G$ գրաֆը արտաքին հարթ է և պարունակում $H$ ենթագրաֆ, որը $K_{2,3}$-ի կամ $K_4$-ի ենթատրոհում է։ Համաձայն թեորեմ 7.3.1-ի, $G + K_1$ գրաֆը հարթ է։ Քանի որ $H + K_1$ գրաֆը հանդիսանում է $K_{3,3}$-ի կամ $K_5$-ի ենթատրոհում, ուստի $G + K_1$ գրաֆը կպարունակի ենթագրաֆ, որը $K_{3,3}$-ի կամ $K_5$-ի ենթատրոհում է, հետևաբար, ըստ թեորեմ 7.2.2-ի, $G + K_1$ գրաֆը հարթ չէ, իսկ դա հակասություն է։

Այժմ ենթադրենք, որ գոյություն ունի $G$ գրաֆ, որը չի պարունակում ենթագրաֆ, որը $K_{2,3}$-ի կամ $K_4$-ի ենթատրոհում է, բայց $G$-ն արտաքին հարթ չէ։ Համաձայն թեորեմ 7.3.1-ի, $G + K_1$ գրաֆը հարթ չէ և, հետևաբար, ըստ թեորեմ 7.2.2-ի, $G + K_1$ գրաֆը կպարունակի ենթագրաֆ, որը $K_{3,3}$-ի կամ $K_5$-ի ենթատրոհում է։ Հեշտ է տեսնել, որ այդ դեպքում $G$ գրաֆը կպարունակի ենթագրաֆ, որը $K_{2,3}$-ի կամ $K_4$-ի ենթատրոհում է, ինչը հակասություն է։ ∎

Նշենք նաև առանց ապացույցի արտաքին հարթ գրաֆների ևս մի նկարագրում:

**Թեորեմ 7.3.3 (Հարարի, Չարտրանդ)։** $G$ գրաֆը արտաքին հարթ է այն և միայն այն դեպքում, երբ $G$-ն չի պարունակում $K_{2,3}$ կամ $K_4$ որպես մինոր։

**Սահմանում 7.3.2։** $G$ արտաքին հարթ գրաֆը կոչվում է մաքսիմալ արտաքին հարթ գրաֆ, եթե այդ գրաֆին ցանկացած նոր կող ավելացնելուց ստացվող գրաֆը արտաքին հարթ չէ։

Դժվար չէ համոզվել, որ նկ. 7.3.1-ում պատկերված $G$ և $F$ գրաֆները մաքսիմալ արտաքին հարթ գրաֆներ են, իսկ $H$-ը՝ մաքսիմալ արտաքին հարթ գրաֆ չէ։ Հեշտ է տեսնել նաև, որ $G$-ն մաքսիմալ արտաքին հարթ գրաֆ է այն և միայն այն դեպքում, երբ այդ գրաֆի յուրաքանչյուր նիստ եռանկյուն է, բացի, մի գուցե, անսահմանափակ նիստից։

Այժմ ապացուցենք մի օժանդակ լեմմա։

**Լեմմա 7.3.1։** Կամայական արտաքին հարթ գրաֆում գոյություն ունի գագաթ, որի աստիճանը երկուսից ավելի չէ։

**Ապացույց։** Նախ նկատենք, որ լեմմայի պնդումն ակնհայտ է չորսից ոչ ավելի գագաթ ունեցող արտաքին հարթ գրաֆներում։ Դիցուք $G$-ն արտաքին հարթ գրաֆ է և $|V(G)| \geq 5$։ Այդ գրաֆին ավելացնենք կողեր այնպես, որ ստացված $G'$ գրաֆը լինի մաքսիմալ արտաքին հարթ գրաֆ։ Պարզ է, որ այդ դեպքում անսահմանափակ նիստի եզրը կհանդիսանա համիլտոնյան ցիկլ։ Դիցուք այդ համիլտոնյան ցիկլը $C$-ն է։ Դիտարկենք $C$ ցիկլի բոլոր $xy$ լարերը և ընտրենք այն մեկը, որի դեպքում $C$ ցիկլի $(x, y)$-ճանապարհով և $xy$ կողով սահմանափակված հարթության տիրույթը պարունակում է նվազագույն քանակությամբ սահմանային նիստեր։ Դիցուք այդ լարը $uv$-ն է։ Հեշտ է տեսնել, որ այդ դեպքում նիստերի քանակը հավասար է մեկի։ Այստեղից հետևում է, որ այդ նիստի եզրը կպարունակի $w \in V(G')$ գագաթ, որի համար $d_{G'}(w) \leq 2$ և, հետևաբար, $d_G(w) \leq 2$։ ∎

Ինչպես հարթ գրաֆների դեպքում, այդպես էլ արտաքին հարթ գրաֆների դեպքում, գոյություն ունի կապ գրաֆի գագաթների և կողերի քանակների միջև։

**Թեորեմ 7.3.4**: Եթե $G$-ն արտաքին հարթ $(n, m)$-գրաֆ է ($n \geq 2$), ապա տեղի ունի $m \leq 2n - 3$ անհավասարությունը։

**Ապացույց**: Թեորեմի ապացույցը կատարենք մակածման եղանակով ըստ $n$-ի։ Թեորեմի պնդումն ակնհայտ է $n \leq 3$ դեպքում։ Ենթադրենք, որ $n \geq 4$ և թեորեմը ճիշտ է ցանկացած $G'$ արտաքին հարթ գրաֆի համար, երբ $|V(G')| < n$։ Դիտարկենք $G$ արտաքին հարթ $(n, m)$-գրաֆը։ Համաձայն լեմմա 7.3.1-ի, $G$ գրաֆում գոյություն ունի $v$ գագաթ, որի համար $d_G(v) \leq 2$։ Դիտարկենք $G' = G - v$ գրաֆը։ Պարզ է, որ $G'$ գրաֆը արտաքին հարթ է և $|V(G')| = n - 1$։ Ըստ մակածման ենթադրության, $|E(G')| \leq 2(n-1) - 3 = 2n - 5$։ Այստեղից հետևում է, որ $m = |E(G)| \leq |E(G')| + 2 \leq 2n - 3$։ ∎

## § 7.4. Գրաֆները մակերևույթների վրա, գրաֆի սեռը և խաչումների թիվը

Նախորդ պարագրաֆներում մենք հիմնականում դիտարկում էինք հարթ գրաֆներ, ինչպես նաև նշեցինք, որ գոյություն ունեն գրաֆներ, որոնք հարթ չեն (օրինակ, $K_n$ լրիվ գրաֆը, երբ $n \geq 5$)։ Նկատեցինք նաև, որ հարթ գրաֆներն այն և միայն այն գրաֆներն են, որոնք հնարավոր է պատկերել սֆերայի վրա այնպես, որ տարբեր կողերին համապատասխանող անընդհատ կորերը չունենան ընդհանուր կետեր, բացի գագաթներից։ Ավելի ընդհանուր ձևով կարող ենք ասել, որ $G$ գրաֆը կարելի է պատկերել $S$ մակերևույթի վրա, եթե այդ գրաֆի գագաթներին կարող ենք համապատասխանեցնել $S$-ի կետեր (տարբեր գագաթներին տարբեր կետեր), այնպես, որ եթե երկու գագաթ կազմում են կող գրաֆում, ապա դրանց համապատասխան կետերը միացվում են ինքնահատում չունեցող անընդհատ կորով, որը չի անցնում մեկ այլ գագաթին համապատասխան կետով, և տարբեր կողերին համապատասխանող անընդհատ կորերը չունեն ընդհանուր կետեր, բացի գագաթներից։ Բնական հարց է առաջանում. հնարավո՞ր է արդյոք ոչ հարթ գրաֆները պատկերել այլ մակերևույթների վրա։ Այդ հարցին դրական պատասխան է տրվել Քյոնիգի կողմից։ Դրա համար ցանկացած գրաֆը պատկերենք սֆերայի վրա, թույլ տալով, որ կողերը հատվեն, բայց թույլ չտալով, որ երեք կող հատվեն միևնույն կետում։ Այնուհետև կողերի յուրաքանչյուր հատման համար այդ սֆերային ավելացնենք բռնակ և հատվող կողերից մեկը անցկացնենք բռնակի վրայով, իսկ մյուսը՝ բռնակի տակ։ Պարզ է, որ արդյունքում մենք կպատկերենք գրաֆը որոշակի քանակությամբ բռնակներ ունեցող սֆերայի վրա։ Իհարկե, այդպիսի պատկերման դեպքում մենք կօգտագործենք ավելի շատ բռնակներ, քան անհրաժեշտ է, և, հետևաբար, բնական է դիտարկել ցանկացած գրաֆի համար բռնակների նվազագույն քանակը, որը անհրաժեշտ է այդ գրաֆի պատկերման համար բռնակներ ունեցող սֆերայի վրա։

Դիցուք $k \in \mathbb{Z}_+$։ $\mathbb{S}_k$-ով նշանակենք *k բռնակ ունեցող սֆերան*։ Ստորև պատկերված են $\mathbb{S}_1$, $\mathbb{S}_2$ և $\mathbb{S}_3$-ը։

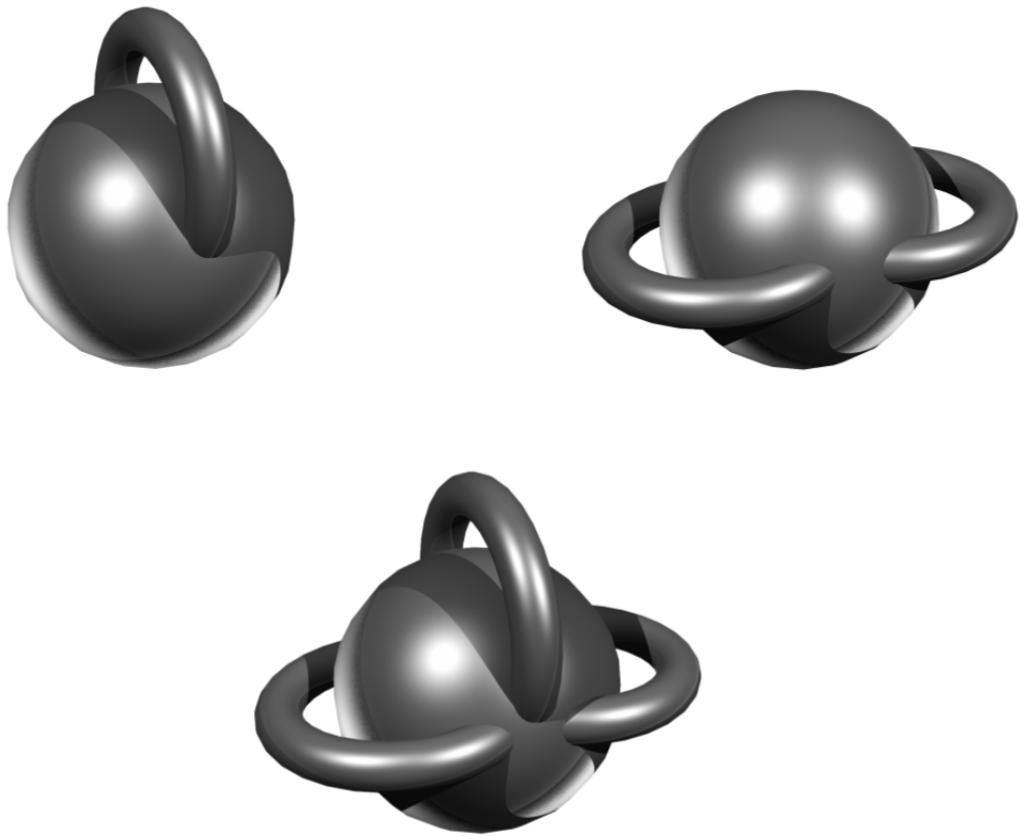

Նկ. 7.4.1 (Այս տարածական պատկերները նկարեց Մարտիրոսյան Հայկը)

**Սահմանում 7.4.1:** $G$ գրաֆի $\gamma(G)$ *սեռ* կանվանենք այն նվազագույն $k$-ն, որի դեպքում այդ գրաֆը կարելի է պատկերել $\mathbb{S}_k$-ի վրա։

Թեորեմ 7.1.5-ից հետևում է, որ $\gamma(G) = 0$ այն և միայն այն դեպքում, երբ $G$-ն հարթ գրաֆ է: Այստեղից և հետևանքներ 7.1.1, 7.1.2-ից բխում է, որ $\gamma(K_{3,3}) \geq 1$ և $\gamma(K_5) \geq 1$: Մյուս կողմից, նկ. 7.4.2-ում բերված են $K_{3,3}$-ի և $K_5$-ի պատկերումները $\mathbb{S}_1$-ի (տորի) վրա: Հետևաբար, $\gamma(K_{3,3}) = 1$ և $\gamma(K_5) = 1$: Այն գրաֆները, որոնք հնարավոր չէ պատկերել սֆերայի ($\mathbb{S}_0$-ի) վրա, սակայն հնարավոր է պատկերել տորի ($\mathbb{S}_1$-ի) վրա, կոչվում են տորոիդալ գրաֆներ: Օրինակ, այդպիսին են $K_{3,3}$ և $K_5$ գրաֆները:

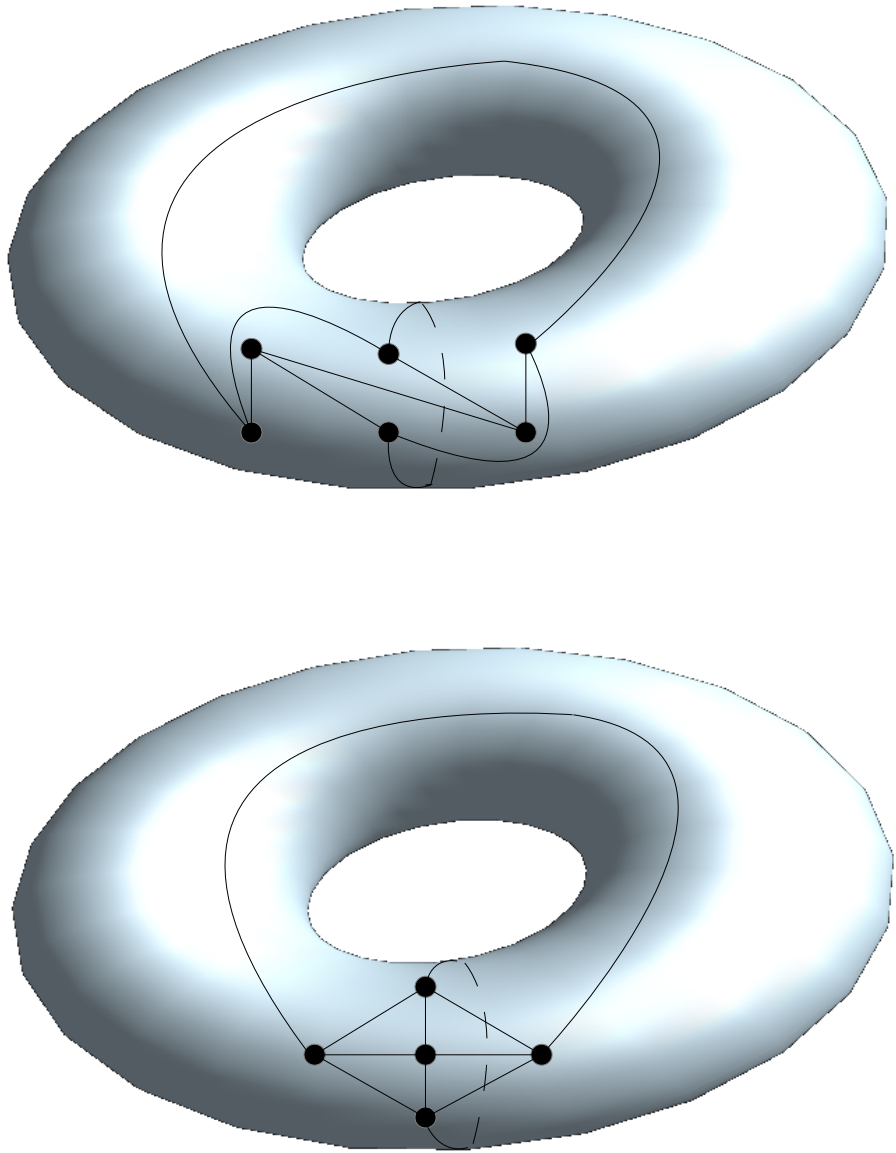

Նկ. 7.4.2

Հայտնի է, որ գոյություն ունի Էյլերի բանաձևի ընդհանրացում կամայական սեռ ունեցող $G$ կապակցված գրաֆի համար: Նկատենք, որ $G$ կապակցված գրաֆը $\mathbb{S}_{\gamma(G)}$-ի վրա պատկերման դեպքում այդ մակերևույթը կտրողի սահմանային մաքսիմալ տիրույթների (այստեղ բացակայում է անսահմանային տիրույթը, որը առկա էր հարթության վրա պատկերման դեպքում): Մյուս կողմից, Քյոնիգզը ցույց է տվել, որ այդ տիրույթները կլինեն կապակցված տիրույթներ (այդ տիրույթների ցանկացած երկու կետեր կարող են միացվել գրաֆի կողերը չհատող անընդհատ գծով), ուստի մենք կարող ենք անվանել այդ տիրույթները $G$ կապակցված գրաֆի *նիստեր* $\mathbb{S}_{\gamma(G)}$-ի վրա պատկերման դեպքում:

**Թեորեմ 7.4.1 (Լ. Էյլեր):** Եթե $G$-ն կապակցված $(n, m)$-գրաֆ է, որն ունի $r$ հատ նիստ $\mathbb{S}_{\gamma(G)}$-ի վրա պատկերման դեպքում, ապա տեղի ունի

$$n - m + r = 2 - 2 \cdot \gamma(G)$$

հավասարությունը:

Ստորև մենք առանց ապացույցի կնշենք գրաֆի սեռի հետ կապված որոշ արդյունքներ:

**Թեորեմ 7.4.2 (Ռինգել, Յանգս):** Ցանկացած $n \geq 3$ բնական թվի համար տեղի ունի

$$\gamma(K_n) = \left\lceil \frac{(n-3)(n-4)}{12} \right\rceil$$

հավասարությունը:

**Թեորեմ 7.4.3 (Ռինգել):** Ցանկացած $m \geq 2$ և $n \geq 2$ բնական թվերի համար տեղի ունի

$$\gamma(K_{m,n}) = \left\lceil \frac{(m-2)(n-2)}{4} \right\rceil$$

հավասարությունը:

**Թեորեմ 7.4.4 (Ռինգել, Բայնեկե, Հարարի):** Ցանկացած $n \geq 2$ բնական թվի համար տեղի ունի

$$\gamma(Q_n) = 1 + (n-4) \cdot 2^{n-3}$$

հավասարությունը:

Նշենք նաև, որ տրված $G$ գրաֆի համար $\gamma(G)$ պարամետրը գտնելու խնդիրը պատկանում է գրաֆների տեսության բարդ խնդիրների դասին [35]: Տարբեր մակերևույթների վրա տրված գրաֆներին ավելի մանրամասն կարելի է ծանոթանալ Մոհարի և Տոմասենի գրքում [28]:

**Սահմանում 7.4.2:** $G$ գրաֆի $cr(G)$ *խաչումների թիվ* կանվանենք այդ գրաֆի հարթության վրա պատկերելու դեպքում կողերի հատումների նվազագույն քանակը:

Հեշտ է տեսնել, որ $cr(G) = 0$ այն և միայն այն դեպքում, երբ $G$-ն հարթ գրաֆ է: Այստեղից և հետևանքներ 7.1.1, 7.1.2-ից բխում է, որ $cr(K_{3,3}) \geq 1$ և $cr(K_5) \geq 1$: Մյուս կողմից, նկ. 7.4.3-ում $K_{3,3}$ և $K_5$ գրաֆները պատկերված են այնպես, որ այդ գրաֆների կողերը հատվում են ճիշտ մեկ անգամ: Հետևաբար, $cr(K_{3,3}) = 1$ և $cr(K_5) = 1$:

Ստորև մենք կձևակերպենք և կապացուցենք գրաֆի խաչումների թվի մի պարզագույն ստորին գնահատական:

**Թեորեմ 7.4.5:** Եթե $G$-ն $(n, m)$-գրաֆ է ($n \geq 3$), ապա տեղի ունի $cr(G) \geq m - 3n + 6$ անհավասարությունը:

**Ապացույց:** Իրոք, դիտարկենք հարթության վրա $G$ գրաֆի կողերի հատումների նվազագույն քանակով պատկերումը: Եթե $G$ գրաֆը հարթ է, ապա համաձայն թեորեմ 7.1.3-ի, $m \leq 3n - 6$, ուստի $m - 3n + 6 \leq 0$: Եթե $G$ գրաֆը հարթ չէ, ապա կողերի յուրաքանչյուր հատման դեպքում այդ հատումը առաջացնող կողերից մեկը հեռացնենք գրաֆից: Պարզ է, որ արդյունքում կունենանք $G'$ հարթ գրաֆ, որի համար, ըստ թեորեմ 7.1.3-ի, ստույգ է $m - cr(G) \leq |E(G')| \leq 3n - 6$ անհավասարությունը: Այստեղից հետևում է, որ $cr(G) \geq m - 3n + 6$: ∎

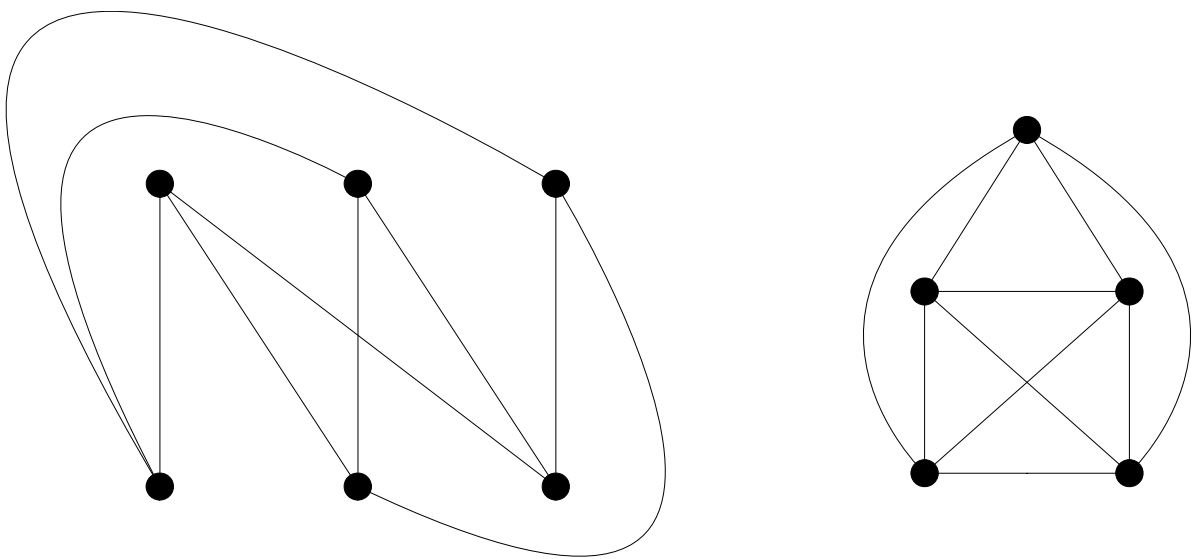

Նկ. 7.4.3

Ստորև մենք կնշենք լրիվ և լրիվ երկկողմանի գրաֆների խաչումների թվերի վերին գնահատականներ:

**Թեորեմ 7.4.6 (Գայ):** Ցանկացած $n$ բնական թվի համար տեղի ունի

$$cr(K_n) \leq \frac{1}{4} \left\lfloor \frac{n}{2} \right\rfloor \left\lfloor \frac{n-1}{2} \right\rfloor \left\lfloor \frac{n-2}{2} \right\rfloor \left\lfloor \frac{n-3}{2} \right\rfloor$$

անհավասարությունը:

Հայտնի է, որ այս վերին գնահատականը ճշգրիտ է, երբ $n \leq 12$:

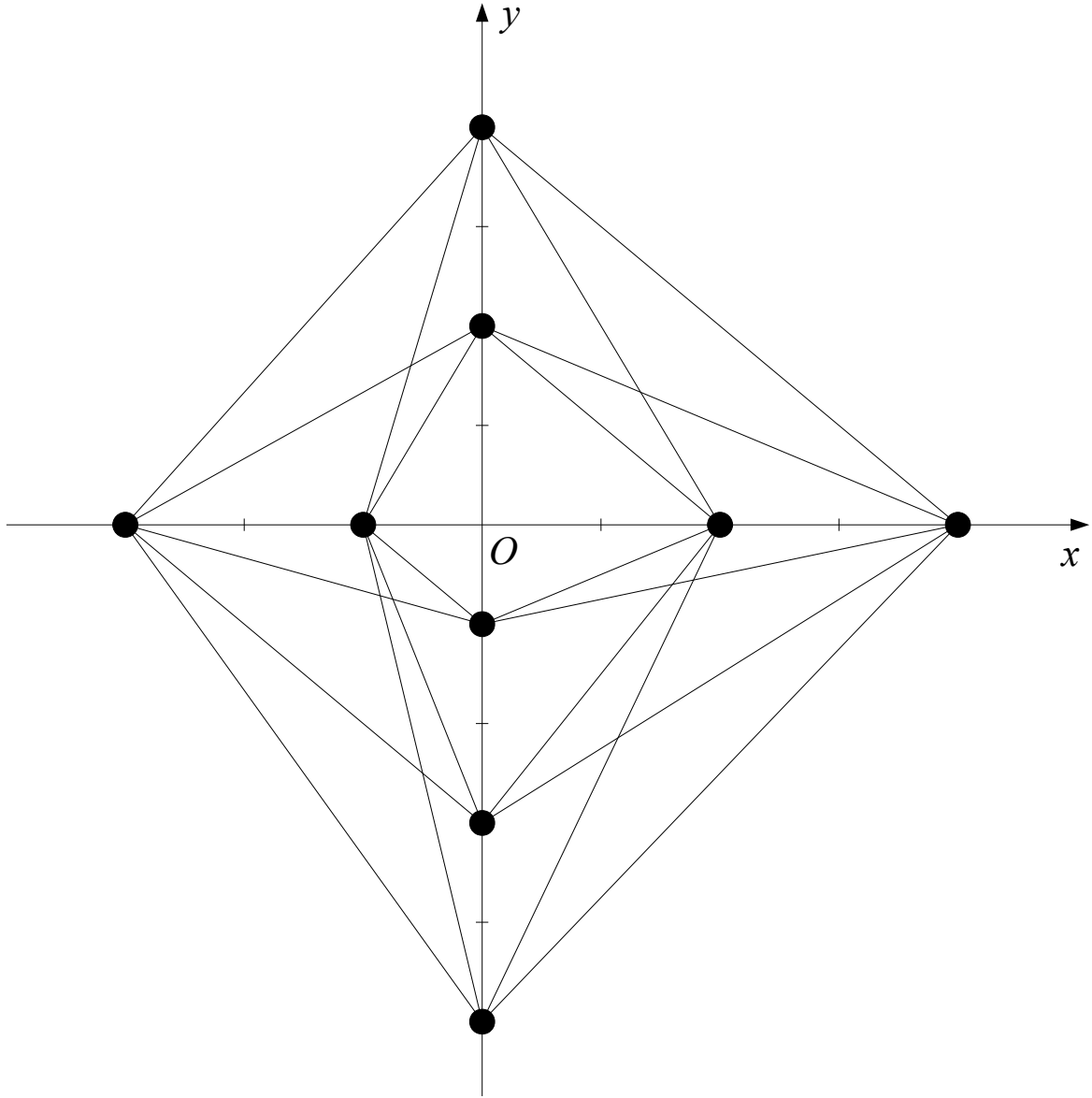

Նկ. 7.4.4

**Թեորեմ 7.4.7 (Զարանկևիչ):** Ցանկացած $m$ և $n$ բնական թվերի համար տեղի ունի

$$cr(K_{m,n}) \leq \left\lfloor \frac{m}{2} \right\rfloor \left\lfloor \frac{m-1}{2} \right\rfloor \left\lfloor \frac{n}{2} \right\rfloor \left\lfloor \frac{n-1}{2} \right\rfloor$$

անհավասարությունը:

**Ապացույց:** Դիցուք $V(K_{m,n}) = \{u_1, \ldots, u_m, v_1, \ldots, v_n\}$ և $E(K_{m,n}) = \{u_i v_j : 1 \leq i \leq m, 1 \leq j \leq n\}$:

Տեղադրենք $u_i$ գագաթը հարթության $(i \cdot (-1)^i, 0)$ կոորդինատներ ունեցող կետում, որտեղ $1 \leq i \leq m$, իսկ $v_j$ գագաթը՝ $(0, j \cdot (-1)^j)$ կոորդինատներ ունեցող կետում, որտեղ $1 \leq j \leq n$: Այնուհետև $u_i$ գագաթը միացնենք $v_j$ գագաթի հետ ուղիղ գծով, որտեղ $1 \leq i \leq m, 1 \leq j \leq n$: Դժվար չէ տեսնել, որ այս պատկերման դեպքում կողերի հատումների քանակը կլինի հավասար $\left\lfloor\frac{m}{2}\right\rfloor\left\lfloor\frac{m-1}{2}\right\rfloor\left\lfloor\frac{n}{2}\right\rfloor\left\lfloor\frac{n-1}{2}\right\rfloor$ թվին, ուստի $cr(K_{m,n}) \leq \left\lfloor\frac{m}{2}\right\rfloor\left\lfloor\frac{m-1}{2}\right\rfloor\left\lfloor\frac{n}{2}\right\rfloor\left\lfloor\frac{n-1}{2}\right\rfloor$: ∎

Նկ. 7.4.4-ում պատկերված է թեորեմ 7.4.7-ի ապացույցում կառուցված պատկերը $K_{4,5}$ գրաֆի դեպքում:

Հայտնի է, որ այս վերին գնահատականը ճշգրիտ է, երբ $min\{m,n\} \leq 6$ կամ, երբ $7 \leq m \leq 8$ և $7 \leq n \leq 10$:

Գրաֆների խաչումների թվերի հետ են կապված գրաֆների տեսության հայտնի և բարդ երկու հիպոթեզները, որոնք ձևակերպել է Զարանկիչը:

**Հիպոթեզ 7.4.1:** Ցանկացած $m$ և $n$ բնական թվերի համար տեղի ունի
$$cr(K_{m,n}) = \left\lfloor\frac{m}{2}\right\rfloor\left\lfloor\frac{m-1}{2}\right\rfloor\left\lfloor\frac{n}{2}\right\rfloor\left\lfloor\frac{n-1}{2}\right\rfloor$$
հավասարությունը:

**Հիպոթեզ 7.4.2:** Ցանկացած $n$ բնական թվի համար տեղի ունի
$$cr(K_n) = \frac{1}{4}\left\lfloor\frac{n}{2}\right\rfloor\left\lfloor\frac{n-1}{2}\right\rfloor\left\lfloor\frac{n-2}{2}\right\rfloor\left\lfloor\frac{n-3}{2}\right\rfloor$$
հավասարությունը:

Նշենք նաև, որ կամայական $G$ գրաֆի համար $cr(G)$ պարամետրը գտնելու խնդիրը պատկանում է գրաֆների տեսության բարդ խնդիրների դասին [15]:

Գլուխ 8

Գրաֆների ներկումներ

§ 8.1. Գրաֆների գագաթային ներկումներ

Դիցուք $G = (V, E)$-ն գրաֆ է:

**Սահմանում 8.1.1:** $G$ գրաֆի *գագաթային $k$-ներկում* կոչվում է $\alpha: V(G) \to \{1, \ldots, k\}$ արտապատկերումը, իսկ՝ $1, \ldots, k$ թվերը կոչվում են *գույներ*:

**Սահմանում 8.1.2:** $G$ գրաֆի $\alpha$ գագաթային $k$-ներկումը կոչվում է *ճիշտ գագաթային $k$-ներկում*, եթե ցանկացած $uv \in E(G)$-ի համար ստույգ է $\alpha(u) \neq \alpha(v)$ պայմանը: Այլ կերպ ասած, ճիշտ գագաթային ներկումն այնպիսի ներկում է, որի դեպքում հարևան գագաթները ներկվում են տարբեր գույներով:

**Սահմանում 8.1.3:** $G$ գրաֆը կոչվում է $k$-*ներկելի*, եթե գոյություն ունի $G$ գրաֆի ճիշտ գագաթային $k$-ներկում: Այն նվազագույն $k$-ն, որի դեպքում $G$-ն $k$-ներկելի է, կոչվում է $G$ գրաֆի *քրոմատիկ թիվ*: $\chi(G)$-ով նշանակենք $G$ գրաֆի քրոմատիկ թիվը:

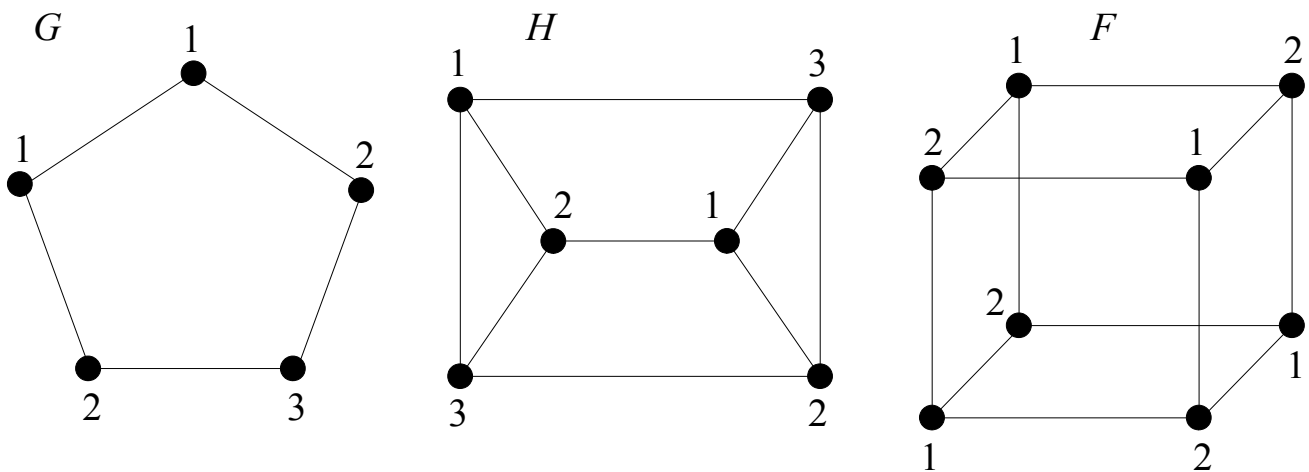

Նկ. 8.1.1

Դիտարկենք նկ. 8.1.1-ում պատկերված $G, H$ և $F$ գրաֆները: Հեշտ է տեսնել, որ այդ նկարում բերված է $G$ գրաֆի գագաթային 3-ներկումը, որը, սակայն, ճիշտ գագաթային 3-ներկում չէ: Մյուս կողմից, հեշտ է տեսնել նաև, որ նկ. 8.1.1-ում պատկերված են $H$ գրաֆի ճիշտ գագաթային 3-ներկում և $F$ գրաֆի ճիշտ գագաթային 2-ներկում: Ավելին, դժվար չէ համոզվել, որ $\chi(H) = 3$ և $\chi(F) = 2$:

Նկատենք, որ $\chi(G) = 1$ այն և միայն այն դեպքում, երբ $G$ գրաֆի բոլոր գագաթները մեկուսացված են և $\chi(G) = 2$ այն և միայն այն դեպքում, երբ $G$-ն առնվազն մեկ կող ունեցող երկկողմանի գրաֆ է: Քանի որ կենտ երկարություն ունեցող պարզ ցիկլը, ըստ թեորեմ 2.2.1-ի, չի հանդիսանում երկկողմանի գրաֆ, ուստի այն ճիշտ ներկելու համար անհրաժեշտ է առնվազն երեք գույն: Մյուս կողմից հեշտ է կառուցել կենտ երկարություն ունեցող պարզ ցիկլի ճիշտ գագաթային 3-ներկում: Այսպիսով, ստանում ենք հետևյալը. ցանկացած $n \geq 3$-ի համար տեղի ունի

$$\chi(C_n) = \begin{cases} 2, \text{ եթե } n-\text{ը զույգ է,} \\ 3, \text{ եթե } n-\text{ը կենտ է,} \end{cases}$$

հավասարությունը:

Այժմ դիտարկենք լրիվ գրաֆները: Քանի որ $K_n$ լրիվ գրաֆի ցանկացած երկու գագաթ հարևան են, ապա $\chi(K_n) = n$:

**Սահմանում 8.1.4:** $G$ գրաֆի $\omega(G)$ *խտությունը* այն ամենամեծ $n$ թիվն է, որի համար $G$ գրաֆն ունի $n$ գագաթ ունեցող լրիվ ենթագրաֆ:

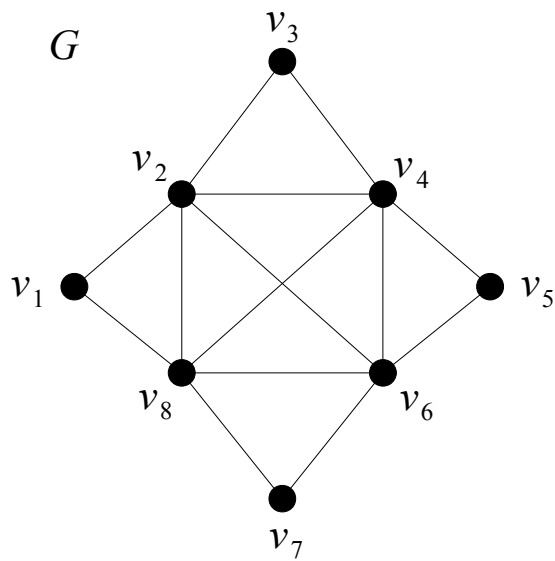

Նկ. 8.1.2

Դիտարկենք նկ. 8.1.2-ում պատկերված $G$ գրաֆը: Հեշտ է տեսնել, որ $G$ գրաֆի $G[\{v_2, v_4, v_6, v_8\}]$ ենթագրաֆը չորս գագաթ ունեցող լրիվ գրաֆ է, ուստի $\omega(G) \geq 4$: Մյուս կողմից, դժվար չէ համոզվել, որ այդ գրաֆը չի պարունակում հինգ գագաթ ունեցող լրիվ ենթագրաֆ, հետևաբար, $\omega(G) = 4$:

Հայտնի են գրաֆի քրոմատիկ թվի տարբեր գնահատականներ: Նշենք դրանցից մի քանիսը:

**Թեորեմ 8.1.1:** Ցանկացած $G$ գրաֆի համար տեղի ունի
$$\omega(G) \leq \chi(G) \leq \Delta(G) + 1$$
անհավասարությունը։

**Ապացույց:** Քանի որ $G$-ն պարունակում է $\omega(G)$ գագաթ ունեցող լրիվ ենթագրաֆ, ուստի այդ ենթագրաֆի գագաթները պետք է ներկվեն տարբեր գույներով և, հետևաբար, $\chi(G) \geq \omega(G)$։

Ցույց տանք, որ $\chi(G) \leq \Delta(G) + 1$։ Ապացույցը կատարենք մակածման եղանակով ըստ $|V(G)|$-ի։ Եթե $|V(G)| = 1$, ապա $\chi(G) = 1$ և $\Delta(G) = 0$։ Ենթադրենք, $\chi(G') \leq \Delta(G') + 1$ անվահասարությունը ճիշտ է ցանկացած $G'$ գրաֆի համար, երբ $|V(G')| < |V(G)|$։ Դիտարկենք $G$ գրաֆը։ Վերցնենք $G$ գրաֆի որևէ $v$ գագաթ և դիտարկենք $G' = G - v$ գրաֆը։ Ըստ մակածման ենթադրության $\chi(G') \leq \Delta(G') + 1 \leq \Delta(G) + 1$։ Մյուս կողմից, քանի որ $d_G(v) \leq \Delta(G)$, ուստի $1, \ldots, \Delta(G) + 1$ գույների մեջ կգտնվի մեկը, որը չի օգտագործվել որպես $v$ գագաթի հարևան գագաթի գույն և, հետևաբար, $v$ գագաթը կարելի է ներկել այդ գույնով։ Այսպիսով, $\chi(G) \leq \Delta(G) + 1$։ ■

Նշենք, որ թեորեմ 8.1.1-ում բերված ստորին և վերին գնահատականները հասանելի են։ Իրոք, դիտարկենք առնվազն մեկ կող ունեցող երկկողմանի $G$ գրաֆը. համաձայն թեորեմ 2.2.1-ի, այն չի պարունակում եռանկյուն, ուստի $\chi(G) = \omega(G) = 2$։ Այժմ դիտարկենք $K_n$ լրիվ գրաֆը և կենտ երկարություն ունեցող $C_{2n+1}$ պարզ ցիկլը. ցանկացած $n$ բնական թվի համար տեղի ունեն $\chi(K_n) = \Delta(K_n) + 1 = n$ և $\chi(C_{2n+1}) = \Delta(C_{2n+1}) + 1 = 3$ հավասարությունները։ Մյուս կողմից, գոյություն ունեն գրաֆներ, որոնց համար թեորեմ 8.1.1-ում բերված անհավասարությունները խիստ են։ Այսպես, օրինակ, հեշտ է տեսնել, որ $3 = \chi(C_5) > \omega(C_5) = 2$ և $2 = \chi(C_4) < \Delta(C_4) + 1 = 3$։

Նշենք նաև, որ թեորեմ 8.1.1-ի ստորին և վերին գնահատականների միջև թվաբանականի հետ կապված է գրաֆների տեսության հայտնի և բարդ հիպոթեզներից մեկը, որը ձևակերպել է Ռիդը։

**Հիպոթեզ 8.1.1:** Ցանկացած $G$ գրաֆի համար տեղի ունի
$$\chi(G) \leq \left\lceil \frac{\omega(G) + \Delta(G) + 1}{2} \right\rceil$$
անհավասարությունը։

**Թեորեմ 8.1.2:** Ցանկացած $G$ գրաֆի համար տեղի ունի

$$\frac{|V(G)|}{\alpha(G)} \leq \chi(G) \leq |V(G)| - \alpha(G) + 1$$

անհավասարությունը:

**Ապացույց:** Դիտարկենք $G$ գրաֆի որևէ ճիշտ գազաթային $\chi(G)$-ներկում: Պարզ է, որ

$$V(G) = V_1 \cup V_2 \cup \cdots \cup V_{\chi(G)} \text{ և } V_i \cap V_j = \emptyset, \text{ երբ } 1 \leq i \neq j \leq \chi(G),$$

որտեղ $V_i$-ն $i$-րդ գույնով ներկված գազաթների բազմությունն է $(1 \leq i \leq \chi(G))$: Քանի որ ներկումը ճիշտ գազաթային $\chi(G)$-ներկում է, ուստի յուրաքանչյուր $i$-ի համար $V_i$-ն անկախ բազմություն է $(1 \leq i \leq \chi(G))$: Հետևաբար, ցանկացած $i$-ի համար $(1 \leq i \leq \chi(G))$ ստույգ է $|V_i| \leq \alpha(G)$ անհավասարությունը: Այսպիսով,

$$|V(G)| = \sum_{i=1}^{\chi(G)} |V_i| \leq \chi(G) \cdot \alpha(G),$$

ուստի $\chi(G) \geq \frac{|V(G)|}{\alpha(G)}$:

Դիցուք $S$-ը $G$ գրաֆի որևէ առավելագույն անկախ բազմություն է: Պարզ է, որ $|S| = \alpha(G)$: Դիտարկենք $G - S$ գրաֆը և ներկենք այդ գրաֆի գազաթները $\chi(G - S)$ գույներով այնպես, որ հարևան գազաթները ներկված լինեն տարբեր գույներով: Այնուհետև $G$ գրաֆի $S$-ի գազաթները ներկենք նոր գույնով: Պարզ է, որ կունենանք հետևյալը.

$$\chi(G) \leq \chi(G - S) + 1 \leq |V(G)| - |S| + 1 = |V(G)| - \alpha(G) + 1: \blacksquare$$

Այստեղ նույնպես, ինչպես և թեորեմ 8.1.1-ում, գոյություն ունեն գրաֆներ, որոնց համար թեորեմ 8.1.2-ում բերված ստորին և վերին գնահատականները հասանելի են, և նաև գոյություն ունեն գրաֆներ, որոնց համար այդ թեորեմում բերված անհավասարությունները խիստ են:

Ստորև կնկարագրենք գրաֆների ճիշտ գազաթային ներկում կառուցելու մի պարզագույն ալգորիթմ:

### Ներկման պարզագույն ալգորիթմ

Դիցուք տրված են $G$ գրաֆը և նրա գազաթների $v_1, v_2, \ldots, v_n$ հաջորդականությունը: Հերթով ներկենք $G$ գրաֆի $v_1, v_2, \ldots, v_n$ գազաթները, վերագրելով $v_i$ գազաթին այն ամենափոքր գույնը, որը չի մասնակցում այդ գազաթին հարևան ներկված գազաթների վրա:

Ալգորիթմի նկարագրությունից հետևում է, որ նրա աշխատանքի արդյունքում ստացվում է $G$ գրաֆի ճիշտ գագաթային ներկում: Ավելին, նկատենք, որ կիրառելով այս ալգորիթմը ցանկացած $G$ գրաֆի և նրա գագաթների որևէ հաջորդականության համար, հեշտ է ստանալ թեորեմ 8.1.1-ում բերված $\chi(G) \leq \Delta(G) + 1$ գնահատականը: Իրոք, քանի որ $G$ գրաֆի գագաթների հաջորդականության մեջ ցանկացած գագաթ ունի ամենաշատը $\Delta(G)$ հատ ներկված հարևան գագաթներ, ապա ներկման պարզագույն ալգորիթմի աշխատանքի արդյունքում օգտագործվող գույների քանակը մեծ չէ $(\Delta(G) + 1)$-ից:

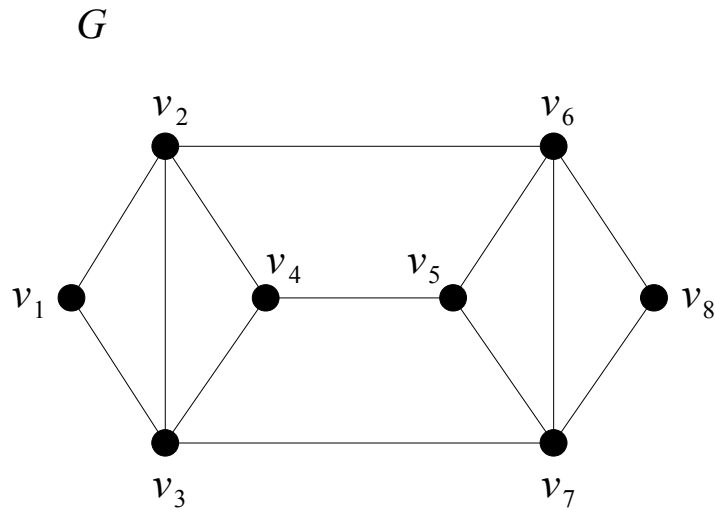

Նկ. 8.1.3

Բերենք ներկման պարզագույն ալգորիթմի աշխատանքը պարզաբանող մի քանի օրինակներ: Դիտարկենք նկ. 8.1.3-ում պատկերված $G$ գրաֆը:

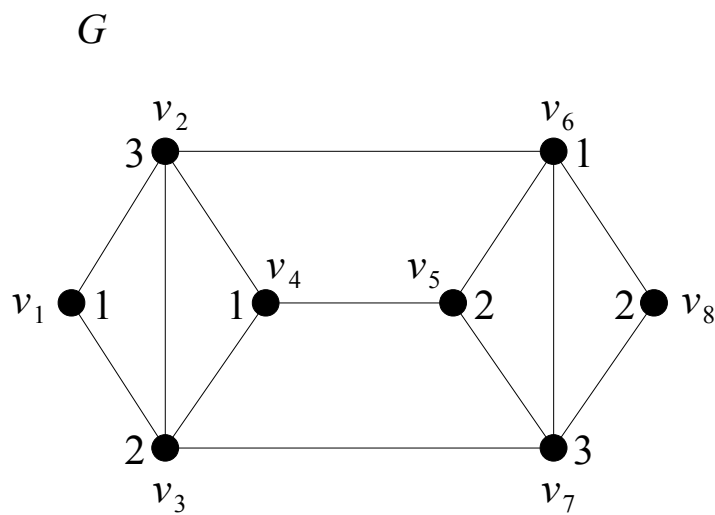

Նկ. 8.1.4

Նկ. 8.1.4-ում պատկերված է $G$ գրաֆի ճիշտ գագաթային 3-ներկումը, որը ստացվել է $G$ գրաֆի և նրա գագաթների $v_1, v_4, v_5, v_6, v_8, v_7, v_3, v_2$ հաջորդականության նկատմամբ ներկման պարզագույն ալգորիթմի կիրառման արդյունքում։ Քանի որ այդ $G$ գրաֆի համար $\omega(G) = 3$, ապա $\chi(G) = 3$։ Այսպիսով, $G$ գրաֆի գագաթների $v_1, v_4, v_5, v_6, v_8, v_7, v_3, v_2$ հաջորդականության դեպքում ներկման պարզագույն ալգորիթմը կառուցում է $G$ գրաֆի ճիշտ գագաթային $\chi(G)$-ներկում։

Այժմ դիտարկենք նկ. 8.1.5-ում պատկերված $H$ գրաֆը։

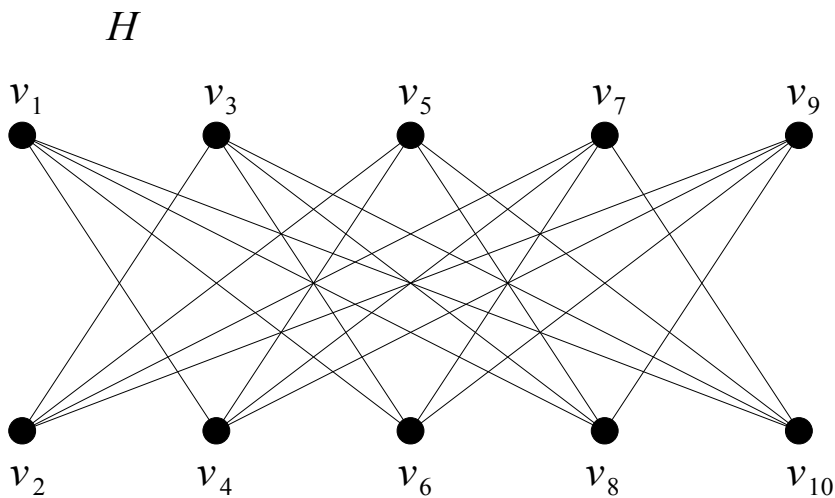

Նկ. 8.1.5

Նկ. 8.1.6-ում պատկերված է $H$ գրաֆի ճիշտ գագաթային $5$-ներկումը, որը ստացվել է $H$ գրաֆի և նրա գագաթների $v_1, v_2, v_3, v_4, v_5, v_6, v_7, v_8, v_9, v_{10}$ հաջորդականության նկատմամբ ներկման պարզագույն ալգորիթմի կիրառման արդյունքում։ Մյուս կողմից, քանի որ $H$-ը երկկողմանի գրաֆ է, ապա $\chi(H) = 2$։ Այսպիսով, $H$ գրաֆի գագաթների $v_1, v_2, v_3, v_4, v_5, v_6, v_7, v_8, v_9, v_{10}$ հաջորդականության դեպքում ներկման պարզագույն ալգորիթմը կառուցում է $H$ գրաֆի ճիշտ գագաթային $(\chi(H) + 3)$-ներկում։

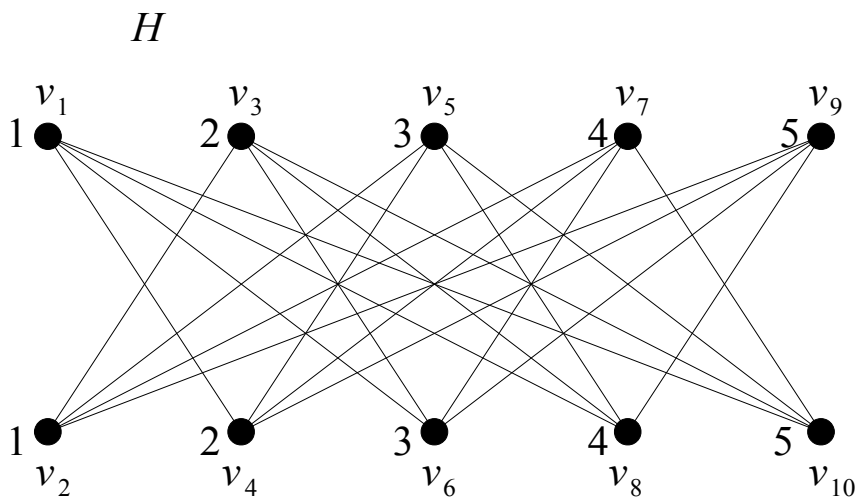

Նկ. 8.1.6

Հեշտ է տեսնել, որ ներկման պարզագույն ալգորիթմի աշխատանքի արդյունքը էապես կախված է գրաֆի գագաթների հաջորդականությունից: Իհարկե, գույների նվազագույն քանակով ներկում գտնելու համար կարելի է առաջարկել հատարկման եղանակը: Օրինակ, կարելի է դիտարկել գրաֆի գագաթների բոլոր հնարավոր տեղափոխությունները որպես ներկման պարզագույն ալգորիթմի գագաթների հաջորդականություններ: Հետնաբար, եթե գրաֆի գագաթների քանակը $n$ է, ապա, դիտարկելով $n!$ հատ տարբերակներ և յուրաքանչյուրի համար կիրառելով ներկման պարզագույն ալգորիթմը, կարելի է գտնել այդ գրաֆի ներկման պարզագույն ալգորիթմով ստացվող բոլոր ճիշտ գագաթային ներկումները: Այնուհետն, գրաֆի բոլոր ստացված ճիշտ գագաթային ներկումներից ընտրել գույների նվազագույն քանակությամբ որևէ ներկում: Սակայն պետք է նշել, որ այդքան գործողություն կատարելը, նույնիսկ $n$-ի փոքր արժեքների դեպքում գործնականում անհնար է: Մյուս կողմից, գոյություն ունեն գրաֆների դասեր, որոնց համար հայտնի են գագաթների հաջորդականությունները, որոնց դեպքում ներկման պարզագույն ալգորիթմը կառուցում է գույների նվազագույն քանակությամբ ներկում: Այդպիսին են, օրինակ, § 1.5-ում սահմանված միջակայքերի գրաֆները:

**Թեորեմ 8.1.3:** Եթե $G$-ն միջակայքերի գրաֆ է, ապա $\chi(G) = \omega(G)$:

**Ապացույց:** Քանի որ $G$-ն միջակայքերի գրաֆ է, ապա իրական թվերի $\mathbb{R}$ բազմության վրա գոյություն ունի $\mathfrak{F} = \{I_1, I_2, \ldots, I_n\}$ փակ հատվածների ընտանիք, որի հատումների $\Omega(\mathbb{R}, \mathfrak{F})$ գրաֆը $G$-ն է: Դիտարկենք $G$ գրաֆի գագաթների այն հաջորդականությունը, որը ստացվում է $I_1, I_2, \ldots, I_n$ փակ հատվածները ըստ ձախ ծայրակետի հանդիպման դասավորելով: Այնուհետն այդ գրաֆի համար կիրառենք ներկման պարզագույն ալգորիթմը: Դիցուք այդ ալգորիթմի ավարտի պահին $v$ գագաթը ստացել է առավելագույն $k$ գույնը: Քանի որ $v$ գագաթը ալգորիթմի աշխատանքի ընթացքում չի ստացել $1, \ldots, k-1$ գույներ, ուստի $v$ գագաթին համապատասխան միջակայքի ձախ $a$ ծայրակետը պատկանում է նան $1, \ldots, k-1$ գույներն արդեն ունեցող գագաթներին համապատասխան միջակայքերին: Բոլոր այդ միջակայքերը պարունակում են $a$ կետը, ուստի $v$ գագաթը և նրա $1, \ldots, k-1$ գույներ ունեցող հարևան գագաթները կազմում են լրիվ ենթագրաֆ $G$ գրաֆում: Այսպիսով, $\omega(G) \geq k \geq \chi(G)$: Մյուս կողմից, ըստ թեորեմ 8.1.1-ի, $\chi(G) \geq \omega(G)$ և, հետնաբար, $\chi(G) = \omega(G)$: ∎

Այժմ ձևակերպենք և ապացուցենք Վելշի և Պաուելլի թեորեմը:

**Թեորեմ 8.1.4:** Եթե $G$-ն գրաֆ է, որի աստիճանային հավաքածուն $d = (d_1, \ldots, d_n)$-ն է, որտեղ $d_1 \geq \cdots \geq d_n$, ապա $\chi(G) \leq 1 + \max_{1 \leq i \leq n} \min\{d_i, i-1\}$:

**Ապացույց:** Թեորեմը ապացուցելու համար մենք կկիրառենք ներկման պարզագույն ալգորիթմը $d_1, \ldots, d_n$ աստիճաններ ունեցող գագաթների հաջորդականության նկատմամբ: Հեշտ է տեսնել, որ $i$-րդ գագաթը ներկելու ժամանակ նրա հարևան գագաթներից ոչ ավելի, քան $\min\{d_i, i-1\}$-ը արդեն ստացել են գույներ: Հետևաբար, $i$-րդ գագաթի գույնը ավելի չէ $1 + \min\{d_i, i-1\}$-ից: Այսպիսով, ներկման պարզագույն ալգորիթմի աշխատանքի արդյունքում օգտագործվող գույների քանակը չի գերազանցում $1 + \max_{1 \leq i \leq n} \min\{d_i, i-1\}$-ը: ∎

**Սահմանում 8.1.5:** Եթե $G$ գրաֆի համար տեղի ունի $\chi(G) = k$ հավասարությունը, ապա $G$ գրաֆը կոչվում է $k$-*քրոմատիկ* գրաֆ:

**Սահմանում 8.1.6:** Եթե $G$ գրաֆի համար տեղի ունի $\chi(G) = k$ հավասարությունը, սակայն $G$ գրաֆի $G$-ից տարբեր ցանկացած $H$ ենթագրաֆի համար տեղի ունի $\chi(H) < k$ անհավասարությունը, ապա $G$ գրաֆը կոչվում է $k$-*կրիտիկական* գրաֆ:

Նկատենք, որ ցանկացած $k$-քրոմատիկ գրաֆ ունի $k$-կրիտիկական ենթագրաֆ: Նշենք գրաֆի քրոմատիկ թվի մի վերին գնահատական, որը ստացվում է այդ գրաֆի կրիտիկական ենթագրաֆները ուսումնասիրելու արդյունքում: Ստորև մենք կձևակերպենք և կապացուցենք Սեկերեշի և Վիլֆի թեորեմը:

**Թեորեմ 8.1.5:** Ցանկացած $G$ գրաֆի համար տեղի ունի
$$\chi(G) \leq 1 + \max_{H \subseteq G} \delta(H)$$
անհավասարությունը:

**Ապացույց:** Դիցուք $k = \chi(G)$ և $H'$-ը $G$ գրաֆի $k$-կրիտիկական ենթագրաֆն է: Նախ ցույց տանք, որ $\delta(H') \geq k - 1$:

Քանի որ $H'$-ը $G$ գրաֆի $k$-կրիտիկական ենթագրաֆ է, ուստի ցանկացած $v \in V(H')$ համար $H' - v$ գրաֆը ունի ճիշտ գագաթային $(k-1)$-ներկում: Եթե $H'$-ում գոյություն ունի $v' \in V(H')$, որ $d_{H'}(v') < k - 1$, ապա $H' - v'$ գրաֆն ունի ճիշտ գագաթային $(k-1)$-ներկում, որի դեպքում $v'$ գագաթի համար գոյություն կունենա $1, \ldots, k-1$ գույներից մեկը, որը չի օգտագործվում նրա հարևան գագաթների համար և, հետևաբար, $v'$ գագաթը կարելի է ներկել այդ գույնով: Այսպիսով, $\chi(H') \leq k - 1$, ինչը հակասում է $H'$-ի $k$-կրիտիկական գրաֆ լինելուն: Այստեղից հետևում է, որ $\delta(H') \geq k - 1$: Հետևաբար,

$$\chi(G) - 1 \leq \delta(H') \leq \max_{H \subseteq G} \delta(H). \quad \blacksquare$$

Նկատենք, որ թեորեմ 8.1.5-ից ևս հետևում է թեորեմ 8.1.1-ում բերված $\chi(G) \leq \Delta(G) + 1$ գնահատականը, քանի որ ցանկացած $G$ գրաֆի համար տեղի ունի $\max_{H \subseteq G} \delta(H) \leq \Delta(G)$ անհավասարությունը։ Մյուս կողմից, արդեն նշել ենք, որ այդ վերին գնահատականները հասանելի են $K_n$ լրիվ գրաֆի և կենտ երկարություն ունեցող $C_{2n+1}$ պարզ ցիկլի համար, քանի որ $\chi(K_n) = \Delta(K_n) + 1 = n$ և $\chi(C_{2n+1}) = \Delta(C_{2n+1}) + 1 = 3$։ Պարզվում է այդ գրաֆները միակն են, որոնց համար այդ վերին գնահատականը հասանելի է։ Սույն փաստն առաջին անգամ նշվել է Բրուքսի կողմից:

**Թեորեմ 8.1.6 (Բրուքս):** Եթե $G$-ն կապակցված գրաֆ է, որը լրիվ գրաֆ կամ կենտ երկարություն ունեցող պարզ ցիկլ չէ, ապա տեղի ունի $\chi(G) \leq \Delta(G)$ անհավասարությունը։

**Ապացույց:** Նախ նկատենք, որ եթե $G$-ն կապակցված գրաֆ է և $\Delta(G) \leq 2$, ապա թեորեմն ակնհայտ է։

Ենթադրենք, $G$-ն կապակցված գրաֆ է և $\Delta(G) \geq 3$։ Ցույց տանք, որ եթե $G$-ն կապակցված գրաֆ է և լրիվ գրաֆ չէ, ապա տեղի ունի $\chi(G) \leq \Delta(G)$ անհավասարությունը։ Ենթադրենք հակառակը․ գոյություն ունեն $H$ կապակցված գրաֆներ, որոնք լրիվ չեն, $\Delta(H) \geq 3$ և $\chi(H) = \Delta(H) + 1$։ Ընտրենք այդ գրաֆներից նվազագույն քանակությամբ գագաթներ ունեցող $G$ գրաֆը։

Վերցնենք $G$ գրաֆի որևէ $v$ գագաթ և դիտարկենք $G' = G - v$ գրաֆը։ $G$ գրաֆի ընտրությունից հետևում է, որ $G'$ գրաֆն ունի ճիշտ գագաթային $\Delta(G)$-ներկում։ Հետևաբար, $d_G(v) = \Delta(G)$ (եթե $d_G(v) < \Delta(G)$, ապա $G'$ գրաֆի ճիշտ գագաթային $\Delta(G)$-ներկման դեպքում $v$ գագաթի համար գոյություն կունենա $1, \ldots, \Delta(G)$ գույներից մեկը, որը չի օգտագործվում նրա հարևան գագաթների համար և, հետևաբար, ներկելով $v$ գագաթը այդ գույնով, մենք կստանանք $G$ գրաֆի ճիշտ գագաթային $\Delta(G)$-ներկում, ինչը հակասում է $\chi(G) = \Delta(G) + 1$ պայմանին)։ Ավելին, կարելի է պնդել հետևյալը.

**Հատկություն 1:** $G'$ գրաֆի ցանկացած ճիշտ գագաթային $\Delta(G)$-ներկման դեպքում $v$ գագաթի հարևան գագաթների գույները զույգ առ զույգ տարբեր են։

Իրոք, եթե գոյություն ունենար $G'$ գրաֆի ճիշտ գագաթային $\Delta(G)$-ներկում, որի դեպքում $v$ գագաթի հարևան գագաթների գույները կրկնվում են, ապա գոյություն կունենար $1, \ldots, \Delta(G)$ գույներից մեկը, որը չի օգտագործվում նրա հարևան գագաթների համար և, հետևաբար, ներկելով $v$ գագաթը այդ գույնով, մենք նորից կստանանք $G$ գրաֆի ճիշտ գագաթային $\Delta(G)$-ներկում:

Դիցուք $\Delta = \Delta(G)$ և $N_G(v) = \{v_1, v_2, \ldots, v_\Delta\}$: Առանց ընդհանրությունը խախտելու կարող ենք ենթադրել, որ $G'$ գրաֆի $\alpha$ ճիշտ գագաթային $\Delta$-ներկման դեպքում $v_1, v_2, \ldots, v_\Delta$ գագաթների գույներն են $1, 2, \ldots, \Delta$-ն: Դիցուք $G_{ij} = G'[S]$, որտեղ $S = \{w : w \in V(G') \text{ և } (\alpha(w) = i \text{ կամ } \alpha(w) = j)\}$ $(1 \leq i \neq j \leq \Delta)$:

Այժմ ապացուցենք հետևյալը.

**Հատկություն 2:** $v_i$ և $v_j$ գագաթները պատկանում են $G_{ij}$ գրաֆի միևնույն կապակցված բաղադրիչին $(1 \leq i \neq j \leq \Delta)$:

Իրոք, եթե $v_i$ և $v_j$ գագաթները պատկանեն $G_{ij}$ գրաֆի տարբեր կապակցված բաղադրիչներին, ապա վերաներկելով $v_i$ գագաթը պարունակող կապակցված բաղադրիչում $i$ գույնով ներկված գագաթները $j$ գույնով, իսկ $j$-ով ներկված գագաթները՝ $i$ գույնով, կստանանք $G'$ գրաֆի $\alpha'$ ճիշտ գագաթային $\Delta$-ներկում, որի դեպքում $v_i$ և $v_j$ գագաթները ներկված են նույն գույնով, ինչը հակասում է հատկություն 1-ին:

Դիցուք $C_{ij}$-ն $G_{ij}$ գրաֆի $v_i$ և $v_j$ գագաթները պարունակող կապակցված բաղադրիչն է:

**Հատկություն 3:** $C_{ij}$-ն $v_i$ և $v_j$ գագաթները միացնող պարզ ճանապարհ է $(1 \leq i \neq j \leq \Delta)$:

Ենթադրենք $d_{C_{ij}}(v_i) > 1$: Պարզ է, որ այդ դեպքում $v_i$ գագաթը հարևան է առնվազն երկու $j$-ով ներկված գագաթների հետ: Քանի որ $d_{G'}(v_i) \leq \Delta - 1$, ուստի $v_i$ գագաթը կարելի է վերաներկել $k$ գույնով, որտեղ $k \neq i, j$, այնպես որ ստանանք $G'$ գրաֆի $\alpha''$ ճիշտ գագաթային $\Delta$-ներկում, որի դեպքում $v_i$ և $v_k$ գագաթները ներկված են նույն գույնով, ինչը հակասում է հատկություն 1-ին: Համանման ձևով կարելի է ցույց տալ, որ $d_{C_{ij}}(v_j) = 1$:

Այժմ համոզվենք, որ $C_{ij}$-ի մնացած բոլոր գագաթների աստիճանները 2 են: Ենթադրենք հակառակը, և դիցուք $u$-ն $C_{ij}$-ում $v_i$-ից $v_j$ գագաթները միացնող ճանապարհի առաջին գագաթն է, որի համար $d_{C_{ij}}(u) > 2$: Եթե $\alpha(u) = i$, ապա $u$ գագաթը հարևան է առնվազն երեք $j$-ով ներկված գագաթներին, իսկ եթե $\alpha(u) = j$, ապա $u$ գագաթը հարևան է առնվազն երեք $i$-ով ներկված գագաթներին: Քանի որ $d_{G'}(u) \leq \Delta$, ուստի $u$ գագաթը կարելի է վերաներկել $k$ գույնով, որտեղ $k \neq i, j$, այնպես որ ստանանք $G'$ գրաֆի $\alpha'''$ ճիշտ գագաթային $\Delta$-ներկում, որի դեպքում $v_i$ և $v_j$ գագաթները կպատկանեն $G_{ij}$ գրաֆի տարբեր կապակցված բաղադրիչներին, ինչը հակասում է

հատկություն 2-ին:

**Հատկություն 4:** $C_{ij}$ և $C_{ik}$ պարզ ճանապարհներն ունեն միայն $v_i$ ընդհանուր գագաթ $(i, j, k = 1, 2, \ldots, \Delta, i \neq j \neq k \neq i)$:

Ենթադրենք հակառակը, և դիցուք $u$-ն $v_i$-ից տարբեր $C_{ij}$ և $C_{ik}$ պարզ ճանապարհների ընդհանուր գագաթն է: Պարզ է, որ այդ դեպքում $\alpha(u) = i$ և $u$ գագաթը հարևան է առնվազն երկու $j$-ով և երկու $k$-ով ներկված գագաթներին: Քանի որ $d_{G'}(u) \leq \Delta$, ուստի $u$ գագաթը կարելի է վերաներկել $l$ գույնով, որտեղ $l \neq i, j, k$, այնպես որ ստանանք $G'$ գրաֆի $\alpha''''$ ճիշտ գագաթային $\Delta$-ներկում, որի դեպքում $v_i$ և $v_j$ գագաթները կպատկանեն $G_{ij}$ գրաֆի տարբեր կապակցված բաղադրիչներին, ինչը հակասում է հատկություն 2-ին:

Քանի որ $G$-ն կապակցված գրաֆ է և լրիվ գրաֆ չէ, ուստի գոյություն ունեն երկու գագաթներ, որոնք հարևան չեն, բայց ունեն ընդհանուր հարևան գագաթ: Առանց ընդհանրությունը խախտելու կարող ենք ենթադրել, որ այդ գագաթները $v_1$ և $v_2$-ն են: Դիտարկենք $C_{12}$ պարզ ճանապարհի $u$ ($u \neq v_2$) գագաթը, որը հարևան է $v_1$-ին: Քանի որ $\Delta \geq 3$, ուստի $G'$ գրաֆում գոյություն ունի $C_{13}$ պարզ ճանապարհ: Այդ դեպքում վերաներկելով $C_{13}$ պարզ ճանապարհի $1$ գույնով ներկված գագաթները $3$ գույնով, իսկ $3$-ով ներկված գագաթները՝ $1$ գույնով, կստանանք $G'$ գրաֆի $\alpha'''''$ ճիշտ գագաթային $\Delta$-ներկում, որի դեպքում $v_1$ գագաթը ներկված է $3$ գույնով, իսկ $v_3$-ը՝ $1$ գույնով: Այստեղից հետևում է, որ նոր առաջացած $C'_{12}$ և $C'_{23}$ կապակցված բաղադրիչները ունեն ընդհանուր $u$ ($u \neq v_2$) գագաթ, ինչը հակասում է հատկություն 4-ին: ∎

Մենք արդեն նշել ենք, որ գոյություն ունեն $G$ գրաֆներ, որոնց համար $\chi(G) = \omega(G)$: Սյուս կողմից պարզվում է, որ $\chi(G)$ և $\omega(G)$ պարամետրերի միջև տարբերությունը կարող է լինել մեծ ցանկացած նախապես տրված թվից: Այդ խնդիրը առաջին անգամ դրվել է Դիրակի կողմից. գոյություն ունի արդյոք եռանկյուն չպարունակող $G$ գրաֆ, որի քրոմատիկ թիվն ինչքան ասեք մեծ է: Դրական պատասխանը խնդրին տրվել է Զիկովի, Տատտի և Միցելկու կողմից: Այստեղ մենք կանգ կառնենք Միցելկու կողմից կառուցված այդպիսի գրաֆների վրա:

Դիցուք $G$-ն գրաֆ է և $V(G) = \{v_1, v_2, \ldots, v_n\}$: Սահմանենք $G$ գրաֆի *Միցելկու $\mu(G)$ գրաֆը* հետևյալ կերպ.

$$V(\mu(G)) = \{v_1, v_2, \ldots, v_n, u_1, u_2, \ldots, u_n\} \cup \{w\},$$

$$E(\mu(G)) = E(G) \cup \{u_i v_j : v_j \in N_G(v_i), 1 \leq i \leq n, 1 \leq j \leq n\} \cup \{w u_i : 1 \leq i \leq n\}:$$

Օրինակ, նկ. 8.1.7-ում պատկերված է $C_5$ գրաֆի Միցելսկու $\mu(C_5)$ գրաֆը, որը նաև հայտնի է *Գրոթշի գրաֆ* անունով։

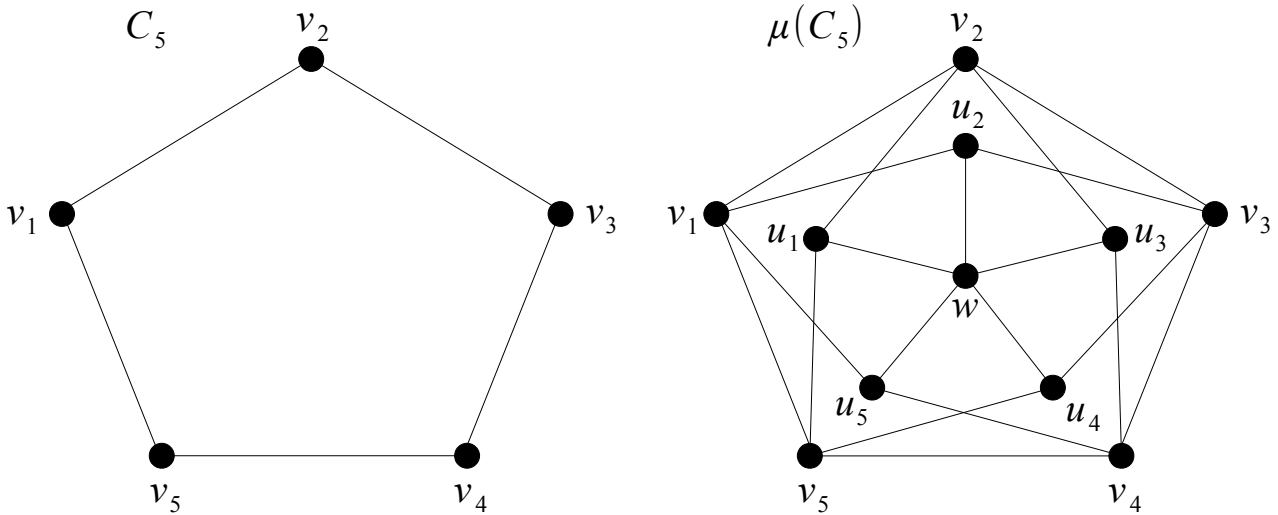

Նկ. 8.1.7

Ստորև մենք կձևակերպենք և կապացուցենք Միցելսկու թեորեմը։

**Թեորեմ 8.1.7:** Եթե $G$-ն եռանկյուն չպարունակող $k$-քրոմատիկ գրաֆ է, ապա $\mu(G)$-ն եռանկյուն չպարունակող $(k+1)$-քրոմատիկ գրաֆ է։

**Ապացույց:** Դիցուք $G$-ն եռանկյուն չպարունակող գրաֆ է, $\chi(G) = k$ և $V(G) = \{v_1, v_2, \ldots, v_n\}$։ Դիտարկենք $\mu(G)$ գրաֆը։ Քանի որ $\mu(G)$ գրաֆի գագաթների $\{u_1, u_2, \ldots, u_n\}$ բազմությունն անկախ է, ուստի $u_i$ գագաթը պարունակող եռանկյունը պետք է երկու հարևան գագաթ ունենա $V(G)$-ից։ Ըստ $\mu(G)$ գրաֆի կառուցման, այդ երկու գագաթները կարող են լինել միայն $v_i$ գագաթին հարևան գագաթներից։ Այստեղից հետևում է, որ եթե $G$-ն եռանկյուն չպարունակող գրաֆ է, ապա $\mu(G)$ գրաֆը ևս չի պարունակում եռանկյուն։

Այժմ ցույց տանք, որ $\chi(\mu(G)) = k + 1$։ Դիցուք $\alpha$-ն $G$ գրաֆի ճիշտ գագաթային $k$-ներկում է։ Սահմանենք $\mu(G)$ գրաֆի $\beta$ գագաթային $(k+1)$-ներկումը հետևյալ կերպ. $\beta(v_i) = \beta(u_i) = \alpha(v_i)$, երբ $1 \leq i \leq n$, և $\beta(w) = k + 1$։ Հեշտ է տեսնել, որ $\beta$-ն $G$ գրաֆի ճիշտ գագաթային $(k+1)$-ներկում է, ուստի $\chi(\mu(G)) \leq k + 1$։ Համոզվենք, որ $\chi(\mu(G)) \geq k + 1$։ Ենթադրենք հակառակը. $\mu(G)$ գրաֆն ունի ճիշտ գագաթային $k$-ներկում։ Դիցուք այդ ներկումը $\gamma$-ն է։ Առանց ընդհանրությունը խախտելու կարող ենք ենթադրել, որ $\gamma(w) = k$։ Այստեղից հետևում է, որ $\gamma(u_i) \in \{1, \ldots, k-1\}$, երբ $1 \leq i \leq n$։ Դիցուք $A = \{v_i : \gamma(v_i) = k, 1 \leq i \leq n\}$։ Կառուցենք $G$ գրաֆի ճիշտ գագաթային $(k-1)$-ներկումը հետևյալ կերպ. յուրաքանչյուր $v_i \in A$-ի համար վերաներկենք այդ գագաթը $\gamma(u_i)$ գույնով, որտեղ $1 \leq i \leq n$: Քանի որ $\gamma$-ն ճիշտ գագաթային ներկում է, ուստի $A$-ն անկախ բազմություն է: Այսպիսով, թեորեմն ապացուցելու համար բավական է ցույց տալ, որ $v_i v'$ կողի կից $v_i$ և $v'$ գագաթները ներկված են տարբեր գույներով, որտեղ $v' \in V(G) \setminus A$: Եթե $v_i v' \in E(G)$, ապա, ըստ $\mu(G)$ գրաֆի կառուցման, $u_i v' \in E(\mu(G))$ և, հետևաբար, $\gamma(u_i) \neq \gamma(v')$: Այժմ, հեռացնելով $\mu(G)$ գրաֆից $u_1, u_2, \ldots, u_n$ գագաթները և $w$ գագաթը, կստանանք $G$ գրաֆի ճիշտ գագաթային $(k-1)$-ներկում, որի գոյությունը հակասում է $\chi(G) = k$ պայմանին: Այսպիսով, $\chi(\mu(G)) = k+1$ և, հետևաբար, $\mu(G)$-ն եռանկյուն չպարունակող $(k+1)$-քրոմատիկ գրաֆ է: ∎

Տարբեր գնահատականներ գրաֆների և նրանց լրացումների քրոմատիկ թվերի գումարի և արտադրյալի համար առաջին անգամ տրվել են Նորդհաուս և Գադումայի կողմից:

**Թեորեմ 8.1.8:** $n$ գագաթ ունեցող ցանկացած $G$ գրաֆի և նրա $\overline{G}$ լրացման համար տեղի ունեն

$$2\sqrt{n} \leq \chi(G) + \chi(\overline{G}) \leq n+1$$

$$n \leq \chi(G) \cdot \chi(\overline{G}) \leq \frac{(n+1)^2}{4}$$

անհավասարությունները:

**Ապացույց:** Դիտարկենք $G$ գրաֆի ճիշտ գագաթային $\chi(G)$-ներկումը: Պարզ է, որ

$$V(G) = V_1 \cup V_2 \cup \cdots \cup V_{\chi(G)} \text{ և } V_i \cap V_j = \emptyset, \text{ երբ } 1 \leq i \neq j \leq \chi(G),$$

որտեղ $V_i$-ն $i$-րդ գույնով ներկված գագաթների բազմությունն է ($1 \leq i \leq \chi(G)$): Այստեղից հետևում է, որ

$$n = |V(G)| = \sum_{i=1}^{\chi(G)} |V_i|,$$

ուստի

$$\max_{1 \leq i \leq \chi(G)} |V_i| \geq \frac{n}{\chi(G)}:$$

Քանի որ $V_i$-ն անկախ բազմություն է, ուստի $\overline{G}$-ում գոյություն ունի $|V_i|$ գագաթ պարունակող լրիվ ենթագրաֆ ($1 \leq i \leq \chi(G)$): Այստեղից և թեորեմ 8.1.1-ից հետևում է, որ

$$\chi(\overline{G}) \geq \omega(\overline{G}) \geq \max_{1 \leq i \leq \chi(G)} |V_i| \geq \frac{n}{\chi(G)}:$$

Այսպիսով, $\chi(G) \cdot \chi(\overline{G}) \geq n$:

Մյուս կողմից, քանի որ

$$\left(\chi(G) + \chi(\overline{G})\right)^2 \geq 4\chi(G) \cdot \chi(\overline{G}),$$

ուստի

$$\left(\chi(G) + \chi(\overline{G})\right)^2 \geq 4\chi(G) \cdot \chi(\overline{G}) \geq 4n:$$

Այստեղից հետևում է, որ $\chi(G) + \chi(\overline{G}) \geq 2\sqrt{n}$:

Այժմ ցույց տանք, որ $\chi(G) + \chi(\overline{G}) \leq n + 1$: Ապացույցը կատարենք մակածման եղանակով ըստ $n$-ի: Եթե $n = 1$, ապա պնդումն ակնհայտ է: Ենթադրենք, պնդումը ճիշտ է ցանկացած $G'$ գրաֆի համար, երբ $|V(G')| \leq n$: Դիտարկենք $n + 1$ գագաթ ունեցող որևէ $G$ գրաֆ: Դիցուք $v \in V(G)$: Դիտարկենք $H = G - v$ գրաֆը: Պարզ է, որ $\chi(G) \leq \chi(H) + 1$ և $\chi(\overline{G}) \leq \chi(\overline{H}) + 1$:

Եթե $\chi(G) = \chi(H)$ կամ $\chi(\overline{G}) = \chi(\overline{H})$, ապա, ըստ մակածման ենթադրության, կունենանք

$$\chi(G) + \chi(\overline{G}) \leq \chi(H) + \chi(\overline{H}) + 1 \leq |V(H)| + 2 = n + 2:$$

Այժմ ենթադրենք, որ $\chi(G) = \chi(H) + 1$ և $\chi(\overline{G}) = \chi(\overline{H}) + 1$: Քանի որ $H$ և $\overline{H}$ գրաֆներին $v$ գագաթն ավելացնելուց քրոմատիկ թիվը մեծանում է, ուստի $d_G(v) \geq \chi(H)$ և $d_{\overline{G}}(v) \geq \chi(\overline{H})$: Այստեղից հետևում է, որ

$$\chi(H) + \chi(\overline{H}) \leq d_G(v) + d_{\overline{G}}(v) = n:$$

Այսպիսով,

$$\chi(G) + \chi(\overline{G}) = \chi(H) + \chi(\overline{H}) + 2 \leq n + 2:$$

Մյուս կողմից, պարզ է, որ

$$\chi(G) \cdot \chi(\overline{G}) \leq \frac{\left(\chi(G)+\chi(\overline{G})\right)^2}{4} \leq \frac{(n+1)^2}{4}: \quad \blacksquare$$

Նշենք, որ թեորեմ 8.1.8-ում բերված ստորին և վերին գնահատականները հասանելի են: Իրոք, դիտարկենք $K_n$ լրիվ գրաֆը, չորս և հինգ երկարություն ունեցող $C_4$ և $C_5$ պարզ ցիկլերը. հեշտ է տեսնել, որ տեղի ունեն $\chi(K_n) + \chi(\overline{K_n}) = n + 1$, $\chi(K_n) \cdot \chi(\overline{K_n}) = n$ և $\chi(C_4) + \chi(\overline{C_4}) = 4$, $\chi(C_5) \cdot \chi(\overline{C_5}) = 9$ հավասարությունները:

Նշենք նաև առանց ապացույցի, որ հայտնի է համարյա բոլոր գրաֆների քրոմատիկ թվի արժեքը:

**Թեորեմ 8.1.9 (Ա. Կորշունով):** $n$ գագաթ ունեցող համարյա բոլոր $G$ գրաֆների համար տեղի ունի

$$\chi(G) \sim \frac{n}{2\log_2 n}$$

առնչությունը։

Այս պարագրաֆի վերջում անդրադառնանք նաև հարթ գրաֆների քրոմատիկ թվի գտնելու խնդրին։ Այդ խնդիրը հայտնի է *չորս գույների հիպոթեզ* անվամբ և ունի շուրջ 150 տարվա պատմություն։ Խնդրի առաջին հիշատակումը կապված է Ֆրենսիս Գուտրիի հետ, որը մոտավորապես 1850 թվականին այդ խնդիրը ձևակերպել է դե Մորգանին։ Այդ խնդիրը կայանում է հետևյալում.

**Չորս գույների հիպոթեզ**։ Հնարավոր է արդյոք ամեն մի աշխարհագրական քարտեզ ներկել չորս գույների միջոցով այնպես, որ յուրաքանչյուր երկրի տարածք ներկված լինի մեկ գույնով, իսկ ընդհանուր սահման ունեցող երկրները ներկված լինեն տարբեր գույներով։ Այստեղ հասկանում ենք, որ յուրաքանչյուր երկրի տարածքը կազմված է մեկ կապակցված տիրույթից և երկու երկիր համարում ենք հարևան, եթե նրանք ունեն ընդհանուր սահման գծի տեսքով, ոչ թե կետի։

Տանք այդ խնդրի մաթեմատիկական ձևակերպումը. եթե աշխարհագրական քարտեզին համապատասխանեցնենք գրաֆ, որի գագաթները քարտեզի երկրներն են, իսկ կողերը՝ ընդհանուր սահման ունեցող երկրների զույգերը, ապա ստացված գրաֆը կլինի հարթ և չորս գույների հիպոթեզը բերվում է հարթ գրաֆի ճիշտ գագաթային ներկման խնդրին չորսից ոչ ավելի գույների միջոցով։ Այսպիսով, «Չորս գույների հիպոթեզը» կունենա հետևյալ ձևակերպումը. եթե $G$-ն հարթ գրաֆ է, ապա $\chi(G) \leq 4$։ Նկատենք, որ գոյություն ունեն հարթ գրաֆներ, որոնց ճիշտ գագաթային ներկման համար անհրաժեշտ է չորս գույն. այդպիսի գրաֆի մի օրինակ է հանդիսանում $K_4$ լրիվ գրաֆը։

«Չորս գույների հիպոթեզով» զբաղվել են շատ մաթեմատիկոսներ և այդ հիպոթեզի առաջին, սխալ պարունակող ապացույցը, տրվել է Կեմպեի [22] կողմից։ Ապացույցում սխալը հայտնաբերել էր Հիվուդը [18], որը քիչ անց ցույց տվեց, որ եթե $G$-ն հարթ գրաֆ է, ապա $\chi(G) \leq 5$։ Այդ հիպոթեզը վերջնական լուծում գտավ 1976 թվականին Ապպելի և Հակենի աշխատանքներում [2,3]։ Նշենք, որ այս հիպոթեզի ապացույցի համար հեղինակները օգտվել են «լիցքերի բեռնաթափման» եղանակից, որի դեպքում անհրաժեշտ եղավ դիտարկել 1400-ից ավելի *բերվող կոնֆիգուրացիաներ* և որոնց ստուգումը կատարվել է համակարգչի օգնությամբ։ 1996 թվականին Ռոբերտսոնի,

Մանդերսի, Սեյմուրի և Տոմասի կողմից տրվեց այդ թեորեմի էապես ավելի պարզ ապացույց [31], որը ևս պահանջում է համակարգչի մասնակցություն, սակայն այս դեպքում դիտարկվել են 633 բերվող կոնֆիգուրացիաներ:

**Թեորեմ 8.1.10 (Ապպել, Հակեն):** Եթե $G$-ն հարթ գրաֆ է, ապա $\chi(G) \leq 4$:

Այժմ ապացուցենք Հիվուդի թեորեմը հինգ գույների մասին:

**Թեորեմ 8.1.11 (Հիվուդ):** Եթե $G$-ն հարթ գրաֆ է, ապա $\chi(G) \leq 5$:

**Ապացույց:** Դիցուք $|V(G)| = n$: Թեորեմի ապացույցը կատարենք մակածման եղանակով ըստ $n$-ի: Եթե $n \leq 5$, ապա թեորեմի պնդումն ակնհայտ է: Ենթադրենք, թեորեմը ստույգ է ցանկացած $G'$ գրաֆի համար, երբ $|V(G')| \leq n - 1$: Դիտարկենք $n$ գագաթ ունեցող որևէ $G$ հարթ գրաֆ: Ցույց տանք, որ $G$ գրաֆը ունի ճիշտ գագաթային 5-ներկում: Ենթադրենք, որ $G$ հարթ գրաֆը հարթության վրա պատկերված է այնպես, որ ցանկացած կող չունենա ինքնահատում և ցանկացած երկու կողեր չունենան ընդհանուր կետեր, բացի գագաթներից: Համաձայն հետևանք 7.1.3-ի, $G$ գրաֆում գոյություն ունի գագաթ, որի աստիճանը հինգից ավելի չէ: Դիցուք այդ գագաթը $u$-ն է: Դիտարկենք $H = G - u$ գրաֆը: Ըստ մակածման ենթադրության, $H$ գրաֆն ունի ճիշտ գագաթային 5-ներկում: Եթե $d_G(u) \leq 4$ կամ այդ ներկման դեպքում $u$ գագաթի հարևան գագաթների գույները կրկնվում են, ապա գոյություն կունենա $1, 2, 3, 4, 5$ գույներից մեկը, որը չի օգտագործվում $u$-ի հարևան գագաթների համար և, հետևաբար, ներկելով $u$ գագաթը այդ գույնով, մենք կստանանք $G$ գրաֆի ճիշտ գագաթային 5-ներկում: Դիցուք $N_G(u) = \{v_1, v_2, v_3, v_4, v_5\}$: Առանց ընդհանրությունը խախտելու կարող ենք ենթադրել, որ գոյություն ունի $H$ գրաֆի $\alpha$ ճիշտ գագաթային 5-ներկում, որի դեպքում $v_1, v_2, v_3, v_4, v_5$ գագաթների գույները $1, 2, 3, 4, 5$-ն են: Դիցուք $G_{ij} = H[S]$, որտեղ $S = \{w: w \in V(H)$ և $(\alpha(w) = i$ կամ $\alpha(w) = j)\}$ $(1 \leq i \neq j \leq 5)$:

Դիտարկենք երկու դեպք:

Դեպք 1: $G$ գրաֆի $v_1$ և $v_3$ գագաթները պատկանում են $G_{13}$ գրաֆի տարբեր կապակցված բաղադրիչների:

Իրոք, եթե $v_1$ և $v_3$ գագաթները պատկանում են $G_{13}$ գրաֆի տարբեր կապակցված բաղադրիչներին, ապա, վերաներկելով $v_1$ գագաթը պարունակող կապակցված բաղադրիչում $1$ գույնով ներկված գագաթները $3$ գույնով, իսկ $3$-ով ներկված գագաթները՝ $1$ գույնով, կստանանք $H$ գրաֆի $\alpha'$ ճիշտ գագաթային 5-ներկում, որի դեպքում $v_1$ և $v_3$ գագաթները ներկված են $3$ գույնով: Այդ դեպքում ներկենք $G$ գրաֆի $u$ գագաթը $1$ գույնով:

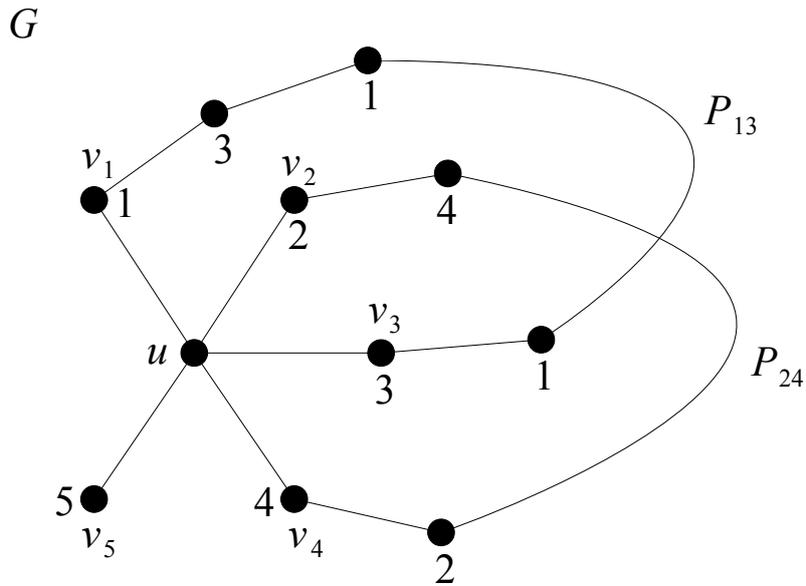

Նկ. 8.1.8

Դեպք 2: $G$ գրաֆի $v_1$ և $v_3$ գագաթները պատկանում են $G_{13}$ գրաֆի միևնույն կապակցված բաղադրիչին:

Հեշտ է տեսնել, որ այդ դեպքում $G$ գրաֆում գոյություն կունենա $v_1$ և $v_3$ գագաթները միացնող $P_{13}$ պարզ ճանապարհի, որի գագաթները ներկված են **1** և **3** գույներով: Քանի որ $G$ գրաֆը հարթ է, ուստի $G$ գրաֆի $v_2$ և $v_4$ գագաթները կպատկանեն $G_{24}$ գրաֆի տարբեր կապակցված բաղադրիչների: Իրոք, հակառակ դեպքում, $v_2$ և $v_4$ գագաթները միացնող $P_{24}$ պարզ ճանապարհը և $P_{13}$ պարզ ճանապարհը կհատվեն $G$ գրաֆի որևէ գագաթում, իսկ դա հակասում է նրան, որ այդ ճանապարհների գագաթները ներկված են տարբեր գույներով (նկ. 8.1.8): Վերաներկենք $v_2$ գագաթը պարունակող կապակցված բաղադրիչում **2** գույնով ներկված գագաթները **4** գույնով, իսկ **4**-ով ներկված գագաթները՝ **2** գույնով, կստանանք $H$ գրաֆի $\alpha''$ ճիշտ գագաթային **5**-ներկում, որի դեպքում $v_2$ և $v_4$ գագաթները ներկված են **4** գույնով: Այդ դեպքում ներկենք $G$ գրաֆի $u$ գագաթը **2** գույնով:

Երկու դեպքում էլ ստացանք, որ $G$ գրաֆը ունի ճիշտ գագաթային **5**-ներկում, ուստի $\chi(G) \leq 5$: ∎

Նշենք նաև առանց ապացույցի Գրոտշի թեորեմը եռանկյուն չպարունակող հարթ գրաֆների քրոմատիկ թվի մասին:

**Թեորեմ 8.1.12 (Գրոտշ):** Եթե $G$-ն եռանկյուն չպարունակող հարթ գրաֆ է, ապա $\chi(G) \leq 3$:

Հարթ գրաֆների քրոմատիկ թվի հետ է կապված նաև գրաֆների տեսության հայտնի և բարդ հիպոթեզներից մեկը, որը ձևակերպել է Հադվիգերը [16]:

**Հիպոթեզ 8.1.2 (Հադվիգեր):** Ցանկացած $k$-քրոմատիկ $G$ գրաֆ պարունակում է $K_k$ որպես մինոր:

Հայտնի է, որ այս հիպոթեզը ճիշտ է $k \leq 4$-ի դեպքում [13]: Եթե $k = 5$, ապա Հադվիգերի հիպոթեզի համաձայն, ցանկացած 5-քրոմատիկ $G$ գրաֆ պարունակում է $K_5$ որպես մինոր: Մյուս կողմից, թեորեմ 7.2.3-ից ստացվում է, որ այդ $G$ գրաֆը չի կարող լինել հարթ, ուստի Հադվիգերի հիպոթեզի $k = 5$ դեպքից հետևում է թեորեմ 8.1.10-ը: Վագները [37] ցույց է տվել, որ իրականում թեորեմ 8.1.10-ը համարժեք է Հադվիգերի հիպոթեզի $k = 5$ դեպքի: Այստեղից հետևում է, որ Հադվիգերի հիպոթեզը ճիշտ է նաև $k = 5$-ի դեպքում: 1993 թվականին Ռոբերտսոնի, Սեյմուրի և Տոմասի կողմից տրվեց Հադվիգերի հիպոթեզի ապացույցը $k = 6$-ի դեպքում [32]: Չլուծված են մնում Հադվիգերի հիպոթեզի $k \geq 7$ դեպքերը: Հադվիգերի հիպոթեզի մի ընդհանրացում դիտարկվել է Հայոշի կողմից և հայտնի է Հայոշի հիպոթեզ անվամբ:

**Հիպոթեզ 8.1.3 (Հայոշ):** Ցանկացած $k$-քրոմատիկ $G$ գրաֆ պարունակում է ենթագրաֆ, որը $K_k$-ի ենթատրոհում է:

Դիրակը ցույց է տվել, որ այս հիպոթեզը ճիշտ է $k \leq 4$-ի դեպքում [13]: Մյուս կողմից, Կատլինը հերքեց Հայոշի հիպոթեզը ցանկացած $k \geq 7$-ի դեպքում [10]: Այսպիսով, Հայոշի հիպոթեզը մնում է բաց $k = 5$ և $k = 6$ դեպքերում:

## § 8.2. Գրաֆների կողային ներկումներ

Դիցուք $G = (V, E)$-ն գրաֆ է:

**Սահմանում 8.2.1:** $G$ գրաֆի *կողային $k$-ներկում* կոչվում է $\alpha: E(G) \to \{1, \ldots, k\}$ արտապատկերումը, իսկ՝ $1, \ldots, k$ թվերը կոչվում են *գույներ*:

**Սահմանում 8.2.2:** $G$ գրաֆի $\alpha$ կողային $k$-ներկումը կոչվում է *ճիշտ կողային $k$-ներկում*, եթե ցանկացած $e, e' \in E(G)$ հարևան կողերի համար ստույգ է $\alpha(e) \neq \alpha(e')$ պայմանը: Այլ կերպ ասած, ճիշտ կողային ներկումն այնպիսի ներկում է, որի դեպքում հարևան կողերը ներկվում են տարբեր գույներով:

**Սահմանում 8.2.3:** $G$ գրաֆը կոչվում է *կողային $k$-ներկելի*, եթե գոյություն ունի $G$ գրաֆի ճիշտ կողային $k$-ներկում: Այն նվազագույն $k$-ն, որի դեպքում $G$-ն կողային $k$-ներկելի է կոչվում է $G$ գրաֆի *քրոմատիկ ինդեքս:* $\chi'(G)$-ով նշանակենք $G$ գրաֆի

քրոմատիկ ինդեքսը:

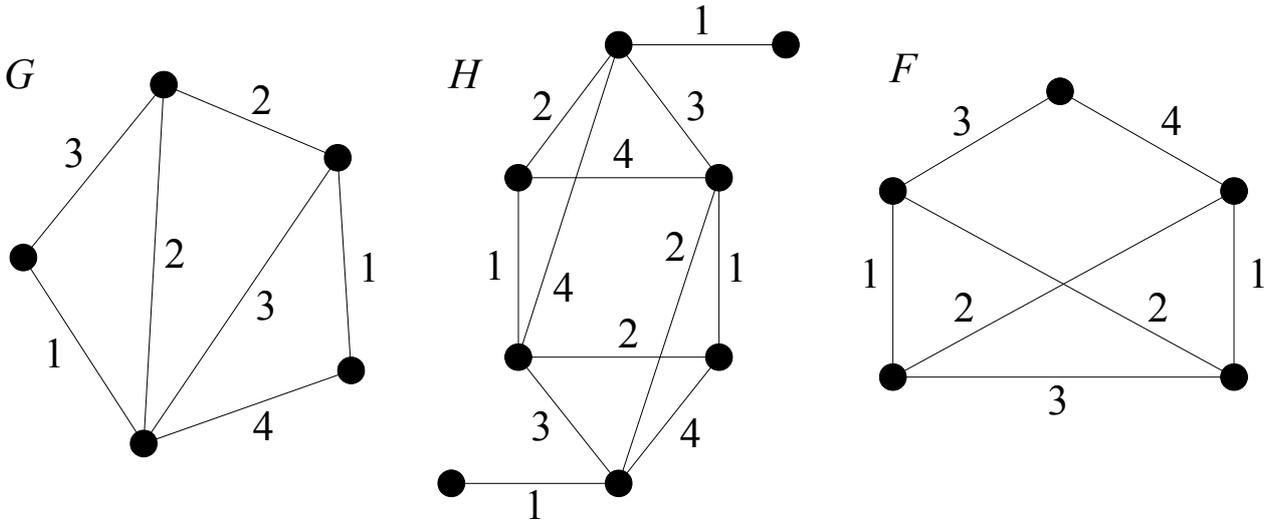

Նկ. 8.2.1

Դիտարկենք նկ. 8.2.1-ում պատկերված $G, H$ և $F$ գրաֆները: Հեշտ է տեսնել, որ այդ նկարում բերված է $G$ գրաֆի կողային **4**-ներկում, որը, սակայն, ճիշտ կողային **4**-ներկում չէ: Մյուս կողմից, հեշտ է տեսնել նաև, որ նկ. 8.2.1-ում պատկերված են $H$ և $F$ գրաֆների ճիշտ կողային **4**-ներկումներ: Ավելին, դժվար չէ համոզվել, որ $\chi'(H) = \chi'(F) = \mathbf{4}$:

Նկատենք, որ ցանկացած $G$ գրաֆի համար $\chi'(G) \geq \Delta(G)$: Իրոք, քանի որ $G$ գրաֆի ճիշտ կողային ներկման դեպքում ցանկացած գագաթին կից կողերը ներկված են զույգ առ զույգ տարբեր գույներով, ուստի այդ ներկման դեպքում օգտագործվող գույների քանակը չի կարող փոքր լինի $\Delta(G)$-ից: Մյուս կողմից, ցանկացած $G = (V, E)$ գրաֆի համար մենք կարող ենք դիտարկել § 1.5-ում սահմանված $\Omega(V, E)$ հատումների գրաֆը: Հիշեցնենք այդ կողային գրաֆի սահմանումը: $G$ գրաֆի կողային $L(G)$ գրաֆը սահմանվում է հետևյալ կերպ. $V(L(G)) = E(G)$ և $E(L(G)) = \{ee': e, e' \in E(G), e$ և $e' - ը$ *հարևան են*$\}$: Հեշտ է տեսնել, որ ցանկացած $G$ գրաֆի համար տեղի ունի $\chi'(G) = \chi(L(G))$ հավասարությունը: Այստեղից, հաշվի առնելով $\Delta(L(G)) \leq 2(\Delta(G) - 1)$ անհավասարությունը և համաձայն թեորեմ 8.1.1-ի, ստանում ենք հետևյալը.

$$\chi'(G) = \chi(L(G)) \leq \Delta(L(G)) + 1 \leq 2\Delta(G) - 1:$$

Այս վերին գնահատականը հասանելի է պարզ ցիկլերի դեպքում: Իրոք, քանի որ $L(C_n) \cong C_n$, երբ $n \geq 3$, ուստի $\chi'(C_n) = \chi(C_n)$: Այսպիսով, ստանում ենք, որ ցանկացած $n \geq 3$-ի համար տեղի ունի

$$\chi'(C_n) = \begin{cases} 2, & \text{եթե } n-ը\ \text{զույգ է}, \\ 3, & \text{եթե } n-ը\ \text{կենտ է}, \end{cases}$$

հավասարությունը։

Պարզվում է, քրոմատիկ ինդեքսի ճշգրիտ արժեքը հայտնի է նաև երկկողմանի և լրիվ գրաֆների դեպքում։

**Թեորեմ 8.2.1 (Դ. Քյոնիգ):** Եթե $G$-ն երկկողմանի գրաֆ է, ապա $\chi'(G) = \Delta(G)$։

**Ապացույց:** Ինչպես նշել ենք, տեղի ունի $\chi'(G) \geq \Delta(G)$ անհավասարությունը։ Ցույց տանք, որ երկկողմանի $G$ գրաֆը ունի ճիշտ կողային $\Delta(G)$-ներկում։ Նախ ապացուցենք, որ ցանկացած երկկողմանի $G$ գրաֆի համար գոյություն ունի այնպիսի $\Delta(G)$-համասեռ երկկողմանի $H$ գրաֆ, որ $G \subseteq H$։ Եթե $G$-ն համասեռ երկկողմանի գրաֆ է, ապա վերցնենք $H = G$։ Հակառակ դեպքում, սահմանենք $G_1$ գրաֆը հետևյալ կերպ. վերցնենք $G$ գրաֆի երկու օրինակ և միացնենք կողով առաջին $G$ գրաֆի յուրաքանչյուր $v$ գագաթ, որի համար $d_G(v) < \Delta(G)$, նույն գագաթի հետ երկրորդ $G$ գրաֆից։ Հեշտ է տեսնել, որ ստացված $G_1$ գրաֆը երկկողմանի է և $\delta(G_1) = \delta(G) + 1$։ Եթե $G_1$-ը համասեռ երկկողմանի գրաֆ է, ապա վերցնենք $H = G_1$։ Հակառակ դեպքում, մակածման եղանակով սահմանենք $G_{i+1}$ գրաֆը հետևյալ կերպ. վերցնենք $G_i$ գրաֆի երկու օրինակ և միացնենք կողով առաջին $G_i$ գրաֆի յուրաքանչյուր $v$ գագաթ, որի համար $d_{G_i}(v) < \Delta(G)$, նույն գագաթի հետ երկրորդ $G_i$ գրաֆից։ Այսպիսով, վերցնելով $H = G_{\Delta(G)-\delta(G)}$, մենք կստանանք անհրաժեշտ համասեռ երկկողմանի գրաֆը։

Քանի որ $H$-ը $\Delta(G)$-համասեռ երկկողմանի գրաֆ է, ուստի համաձայն թեորեմ 5.2.4-ի, այն պարունակում է $M_1$ կատարյալ զուգակցում։ Դիտարկենք $G - M_1$ գրաֆը։ Նկատենք, որ այն $(\Delta(G) - 1)$-համասեռ երկկողմանի գրաֆ է։ Համաձայն թեորեմ 5.2.4-ի, այն պարունակում է $M_2$ կատարյալ զուգակցում։ Դիտարկենք $G - M_1 - M_2$ գրաֆը։ Նկատենք, որ այն $(\Delta(G) - 2)$-համասեռ երկկողմանի գրաֆ է։ Նշված քայլերը կիրառելով $\Delta(G)$ անգամ, մենք կստանանք $H$ $\Delta(G)$-համասեռ երկկողմանի գրաֆի կողերի բաժանումը իրարից զույգ առ զույգ չհատվող կատարյալ զուգակցումների. $E(H) = M_1 \cup M_2 \cup \cdots \cup M_{\Delta(G)}$։ Այժմ $M_i$-րդ զուգակցման կողերը ներկենք $i$-րդ գույնով ($1 \leq i \leq \Delta(G)$)։ Այստեղից հետևում է, որ $H$-ը ունի ճիշտ կողային $\Delta(G)$-ներկում և, հետևաբար, նաև $G$ գրաֆը ունի ճիշտ կողային $\Delta(G)$-ներկում, քանի որ $G \subseteq H$։ Այսպիսով, $\chi'(G) \leq \chi'(H) = \Delta(G)$։ ∎

**Թեորեմ 8.2.2 (Վ. Վիզինգ):** Ցանկացած $n \geq 2$-ի համար տեղի ունի

$$\chi'(K_n) = \begin{cases} n - 1, & \text{եթե } n-\text{ը զույգ է,} \\ n, & \text{եթե } n-\text{ը կենտ է,} \end{cases}$$

հավասարությունը:

**Ապացույց:** Նախ դիտարկենք $n$-ի զույգ թիվ լինելու դեպքը: Դիցուք $n = 2l$ ($l \in \mathbb{N}$) և $V(K_{2l}) = \{v_0, v_1, \ldots, v_{2l-1}\}$: Դիտարկենք $K_{2l}$ գրաֆի կողերի $M_1, M_2, \ldots, M_{2l-1}$ բազմությունները, որտեղ

$$M_1 = \{v_0v_1, v_2v_{2l-1}, v_3v_{2l-2}, \ldots, v_lv_{l+1}\},$$
$$M_2 = \{v_0v_2, v_3v_1, v_4v_{2l-1}, \ldots, v_{l+1}v_{l+2}\},$$
$$\ldots \ldots \ldots \ldots \ldots \ldots \ldots$$
$$M_{2l-1} = \{v_0v_{2l-1}, v_1v_{2l-2}, v_2v_{2l-3}, \ldots, v_{l-1}v_l\}:$$

Այստեղ $M_{i+1}$-ի կողերի բազմությունը ստացվում է $M_i$-ի կողերից հետևյալ գործողության միջոցով. յուրաքանչյուր կողին կից ամեն մի գագաթի ինդեքսին, բացի $v_0$ գագաթից, գումարում ենք $1$ ըստ մոդուլ $(2l-1)$-ի: Հեշտ է տեսնել, որ կողերի $M_1, M_2, \ldots, M_{2l-1}$ բազմությունները զույգ առ զույգ չհատվող կատարյալ զուգակցումներ են: Այժմ $M_i$-րդ զուգակցման կողերը ներկենք $i$-րդ գույնով ($1 \leq i \leq 2l-1$): Այստեղից հետևում է, որ $K_{2l}$ գրաֆը ունի ճիշտ կողային $(2l-1)$-ներկում և, հետևաբար, $\chi'(K_{2l}) \leq \Delta(K_{2l}) = 2l-1$: Մյուս կողմից, քանի որ $\chi'(K_{2l}) \geq \Delta(K_{2l}) = 2l-1$, ուստի $\chi'(K_{2l}) = 2l-1$:

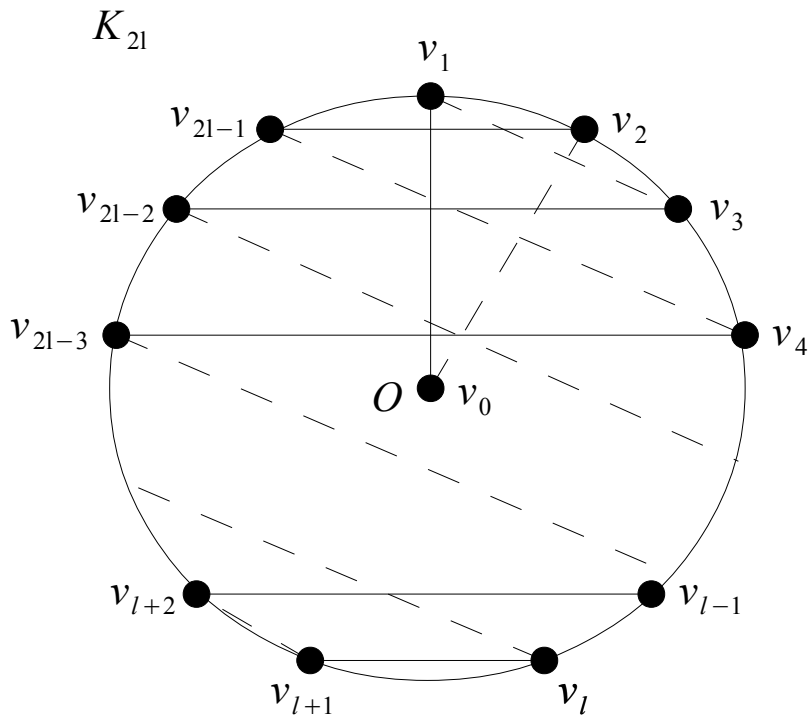

Նկ. 8.2.2

Ներկման այս եղանակն ունի հետևյալ երկրաչափական մեկնաբանումը: Որպես $K_{2l}$ գրաֆի գագաթներ վերցնենք շրջանագծին ներգծած կանոնավոր $(2l-1)$-անկյուն բազմանկյան գագաթները և $O$ կենտրոնը, իսկ որպես կողեր՝ համապատասխան գագաթները միացնող հատվածները։ Առաջին գույնով ներկվող կողերը պատկերված են նկ. 8.2.2-ում։ Յուրաքանչյուր հաջորդ գույնով ներկվող կողերը կստացվեն, եթե նկ. 8.2.2-ի պատկերը ժամացույցի սլաքի ուղղությամբ պտտենք բազմանկյան $O$ կենտրոնի շուրջը $\left(\frac{360}{2l-1}\right)^\circ$ չափով։ Քանի որ յուրաքանչյուր պտույտի ժամանակ հերթական գույնով ներկվող կողերի ուղղությունները փոխվում են, ուստի ոչ մի կող երկու անգամ չի ներկվի և պարզ է նաև, որ բոլոր կողերը կներկվեն $2l-1$ գույների միջոցով։

Այժմ դիտարկենք $n$-ի կենտ թիվ լինելու դեպքը։ Դիցուք $n = 2l+1$ ($l \in \mathbb{N}$)։ Քանի որ $K_{2l+1}$ գրաֆը ստացվում է $K_{2l+2}$ գրաֆից մեկ գագաթ հեռացնելով, ուստի $\chi'(K_{2l+1}) \leq \chi'(K_{2l+2}) = \Delta(K_{2l+2}) = 2l+1$։ Ցույց տանք, որ $\chi'(K_{2l+1}) \geq 2l+1$։ Դիտարկենք $K_{2l+1}$ գրաֆի ճիշտ կողային $\chi'(K_{2l+1})$-ներկում։ Պարզ է, որ

$$E(K_{2l+1}) = E_1 \cup E_2 \cup \cdots \cup E_{\chi'(K_{2l+1})} \text{ և } E_i \cap E_j = \emptyset, \text{ երբ } 1 \leq i \neq j \leq \chi'(K_{2l+1}),$$

որտեղ $E_i$-ն $i$-րդ գույնով ներկված կողերի բազմությունն է ($1 \leq i \leq \chi'(K_{2l+1})$)։ Նկատենք, որ յուրաքանչյուր $i$-ի համար $E_i$-ն զուգակցում է ($1 \leq i \leq \chi'(K_{2l+1})$)։ Հետևաբար, ցանկացած $i$-ի համար ($1 \leq i \leq \chi'(K_{2l+1})$) ստույգ է $|E_i| \leq \alpha'(K_{2l+1}) = l$ անհավասարությունը։ Այսպիսով,

$$l \cdot (2l+1) = \binom{2l+1}{2} = |E(K_{2l+1})| = \sum_{i=1}^{\chi'(K_{2l+1})} |E_i| \leq \chi'(K_{2l+1}) \cdot \alpha'(K_{2l+1}) = \chi'(K_{2l+1}) \cdot l,$$

ուստի $\chi'(K_{2l+1}) \geq 2l+1$։ ∎

Ստորև ձևակերպենք և ապացուցենք Վիզինգի թեորեմը գրաֆների քրոմատիկ ինդեքսի մասին։

**Թեորեմ 8.2.3 (Վ. Վիզինգ):** Կամայական $G$ գրաֆի համար տեղի ունի

$$\Delta(G) \leq \chi'(G) \leq \Delta(G) + 1$$

անհավասարությունը։

**Ապացույց:** Ինչպես նշել ենք, ցանկացած $G$ գրաֆի համար տեղի ունի $\chi'(G) \geq \Delta(G)$ անհավասարությունը։

Ցույց տանք, որ $G$ գրաֆը ունի ճիշտ կողային ($\Delta(G) + 1$)-ներկում։

Ենթադրենք հակառակը. գոյություն ունեն $H$ գրաֆներ, որոնք չունեն ճիշտ կողային ($\Delta(H) + 1$)-ներկում։ Ընտրենք այդ գրաֆներից նվազագույն քանակությամբ կողեր ունեցող $G = (V, E)$ գրաֆը։

Վերցնենք $G$ գրաֆի որևէ $ab$ կող և դիտարկենք $G' = G - ab$ գրաֆը։ $G$ գրաֆի ընտրությունից հետևում է, որ $G'$ գրաֆը ունի ճիշտ կողային $(\Delta(G) + 1)$-ներկում, քանի որ $\Delta(G') \leq \Delta(G)$։ Դիցուք այդ ճիշտ կողային $(\Delta(G) + 1)$-ներկումն $\alpha$-ն է։ Քանի որ $G'$ գրաֆի յուրաքանչյուր $v$ գագաթի համար $d_{G'}(v) \leq \Delta(G)$, ուստի $\alpha$ ներկման դեպքում յուրաքանչյուր $v$ գագաթում բացակայում է առնվազն մեկ գույն, որով այդ գագաթին կից կողերը չեն ներկվել։ $C_\alpha(v)$-ով նշանակենք այն գույների բազմությունը, որոնք բացակայում են $v$ գագաթում $\alpha$ ներկման դեպքում։

Դիտարկենք $G$ գրաֆի $ab$ կողը։ Եթե $C_\alpha(a) \cap C_\alpha(b) \neq \emptyset$, ապա $G$ գրաֆի $ab$ կողը ներկենք կամայական $i_0 \in C_\alpha(a) \cap C_\alpha(b)$ գույնով, որը $G'$ գրաֆի $\alpha$ ներկման հետ միասին կտա $G$ գրաֆի ճիշտ կողային $(\Delta(G) + 1)$-ներկում։

Այժմ ենթադրենք, $C_\alpha(a) \cap C_\alpha(b) = \emptyset$։ Ցույց տանք, որ $G'$ գրաֆի $\alpha$ ներկումից կարելի է անցնել նոր $\beta$ ճիշտ կողային $(\Delta(G) + 1)$-ներկման, որի դեպքում $C_\beta(a) \cap C_\beta(b) \neq \emptyset$։

Դիցուք $S \in C_\alpha(b)$ և $t \in C_\alpha(a)$։ Քանի որ $S \notin C_\alpha(a)$, ուստի գոյություն կունենա $v_1$ գագաթ այնպես, որ $av_1 \in E(G)$ և $\alpha(av_1) = S$։ Ընտրենք $v_1$ գագաթում որևէ $S_1 \in C_\alpha(v_1)$ գույն և նշենք $a$-ին հարևան այն $v_2$ գագաթը, որի համար $\alpha(av_2) = S_1$, այնուհետև ընտրենք $v_2$ գագաթում որևէ $S_2 \in C_\alpha(v_2)$ գույն և նշենք $a$-ին հարևան այն $v_3$ գագաթը, որի համար $\alpha(av_3) = S_2$ և այլն։ Շարունակելով այս պրոցեսը, կստանանք $a$-ին հարևան գագաթների հետևյալ հաջորդականությունը՝ $v_1, v_2, v_3, \ldots$ և գույների $S_1, S_2, S_3, \ldots$ հաջորդականությունը, որոնց համար $S_i \in C_\alpha(v_i)$ և $\alpha(av_i) = S_{i-1}$ (այստեղ $S_0 = S$) (նկ. 8.2.3)։

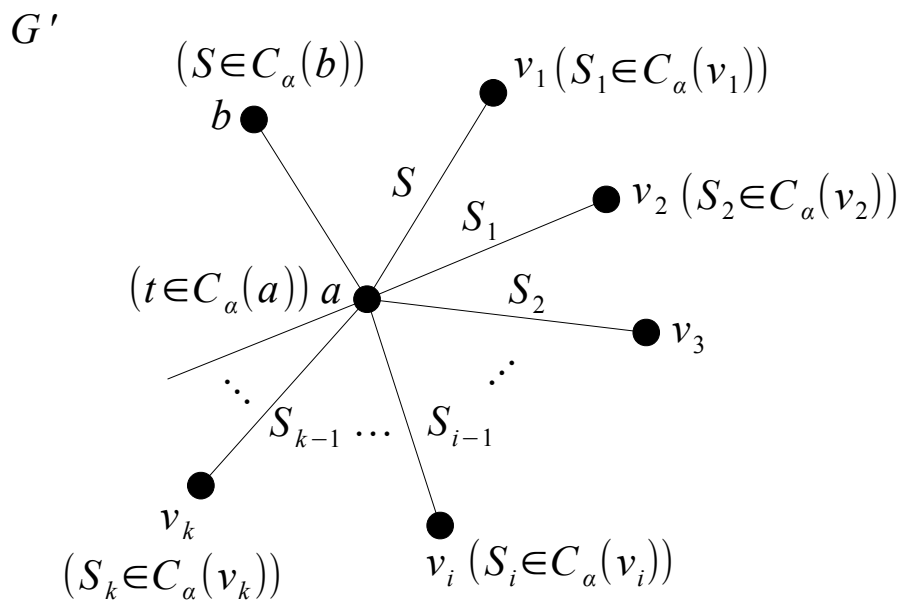

Նկ. 8.2.3

Տրամաբանորեն հնարավոր են հետևյալ դեպքերը։

Դեպք ա): Ստացել ենք գազաթների $v_1, v_2, \ldots, v_k$ հաջորդականությունը: $S_1, S_2, \ldots, S_k$ գույները միմյանցից տարբեր են և հնարավոր չէ նշել $a$-ին հարևան $v_{k+1}$ գազաթ, որի համար $\alpha(av_{k+1}) = S_k$:

Այդ դեպքում $G'$ գրաֆում կատարենք հետևյալ վերաներկումը. $av_k$ կողը ներկենք $S_k$ գույնով, $av_{k-1}$ կողը ներկենք $S_{k-1}$ գույնով, $av_{k-2}$ կողը՝ $S_{k-2}$ գույնով և այլն, $av_1$ կողը՝ $S_1$ գույնով: Պարզ է, որ արդյունքում կստանանք $G'$ գրաֆի $\alpha'$ ճիշտ կողային $(\Delta(G)+1)$-ներկում, որի դեպքում $S \in C_{\alpha'}(a) \cap C_{\alpha'}(b)$:

Դեպք բ): Ստացել ենք գազաթների $v_1, v_2, \ldots, v_k$ հաջորդականությունը: $S_1, S_2, \ldots, S_k$ գույները միմյանցից տարբեր են և $S = S_k$:

Դիտարկենք $G'$ գրաֆի $a, b$ և $v_k$ գազաթները: Պարզ է, որ $S \in C_\alpha(b)$, $t \in C_\alpha(a)$ և $S \in C_\alpha(v_k)$:

Դիտարկենք $G'$ գրաֆի $G_{St} = (V, E_{St})$ ենթագրաֆը, որտեղ $E_{St}$-ն $G'$ գրաֆի այն կողերի բազմությունն է, որոնք $\alpha$ ներկման ժամանակ ներկվել են $S$ կամ $t$ գույներով: Քանի որ $\Delta(G_{St}) \leq 2$, իսկ $d_{G_{St}}(a) = d_{G_{St}}(b) = d_{G_{St}}(v_k) = 1$, ուստի պարզ է, որ $a, b$ և $v_k$ գազաթները միևնույն կապակցված բաղադրիչին չեն պատկանում: Հնարավոր են հետևյալ երկու դեպքերը.

1. $a$ և $b$ գազաթները պատկանում են $G_{St}$ գրաֆի տարբեր կապակցված բաղադրիչներին:

   Այդ դեպքում, վերաներկելով $G_{St}$ գրաֆի $a$ գազաթը պարունակող կապակցված բաղադրիչում $S$ գույնով ներկված կողերը $t$ գույնով, իսկ $t$-ով ներկված կողերը՝ $S$ գույնով, կստանանք $G'$ գրաֆի $\alpha'$ ճիշտ կողային $(\Delta(G)+1)$-ներկում, որի դեպքում $S \in C_{\alpha'}(a) \cap C_{\alpha'}(b)$:

2. $a$ և $b$ գազաթները պատկանում են $G_{St}$ գրաֆի միևնույն կապակցված բաղադրիչին:

   Այդ դեպքում $v_k$ գազաթը չի պատկանի այդ բաղադրիչին և, վերաներկելով $G_{St}$ գրաֆի $v_k$ գազաթը պարունակող կապակցված բաղադրիչում $S$ գույնով ներկված կողերը $t$ գույնով, իսկ $t$-ով ներկված կողերը՝ $S$ գույնով, կստանանք $G'$ գրաֆի $\alpha'$ ճիշտ կողային $(\Delta(G)+1)$-ներկում, որի դեպքում $t \in C_{\alpha'}(v_k)$: Ստացվեց արդեն քննարկված ա) դեպքը:

Դեպք գ): Ստացել ենք գազաթների $v_1, v_2, \ldots, v_i, \ldots, v_k$ հաջորդականությունը և գույների $S_1, S_2, \ldots, S_i, \ldots, S_k$ հաջորդականությունը, որտեղ $S_i = S_k$ ($i + 1 < k$), իսկ մնացած բոլոր գույները միմյանցից տարբեր են:

Դիտարկենք $G'$ գրաֆի $a, v_i$ և $v_k$ գագաթները: Պարզ է, որ $t \in C_\alpha(a)$, $S_k \in C_\alpha(v_i)$ և $S_k \in C_\alpha(v_k)$: Հնարավոր են հետևյալ երկու դեպքերը.

1. $a$ և $v_i$ գագաթները պատկանում են $G_{S_k t}$ գրաֆի տարբեր կապակցված բաղադրիչներին:

    Այդ դեպքում, վերաներկելով $G_{S_k t}$ գրաֆի $a$ գագաթը պարունակող կապակցված բաղադրիչում $S_k$ գույնով ներկված կողերը $t$ գույնով, իսկ $t$-ով ներկված կողերը՝ $S_k$ գույնով, կստանանք $G'$ գրաֆի $\alpha'$ ճիշտ կողային $(\Delta(G) + 1)$-ներկում, որի դեպքում $\alpha'(av_{i+1}) = t$ և $S_k \in C_{\alpha'}(a) \cap C_{\alpha'}(v_i)$: Ստացվեց արդեն քննարկված ա) դեպքը:

2. $a$ և $v_i$ գագաթները պատկանում են $G_{S_k t}$ գրաֆի միևնույն կապակցված բաղադրիչին:

    Այդ դեպքում $v_k$ գագաթը չի պատկանի այդ բաղադրիչին և վերաներկելով $G_{S_k t}$ գրաֆի $v_k$ գագաթը պարունակող կապակցված բաղադրիչում $S_k$ գույնով ներկված կողերը $t$ գույնով, իսկ $t$-ով ներկված կողերը՝ $S$ գույնով, կստանանք $G'$ գրաֆի $\alpha'$ ճիշտ կողային $(\Delta(G) + 1)$-ներկում, որի դեպքում $t \in C_{\alpha'}(v_k)$: Նորից ստացվեց արդեն քննարկված ա) դեպքը:

Քանի որ $a$ գագաթին հարևան գագաթների քանակը վերջավոր է, ուստի այլ դեպքեր հնարավոր չեն: Այսպիսով, մենք ցույց տվեցինք, որ միշտ հնարավոր է $G'$ գրաֆի $\alpha$ ներկումից անցնել $\beta$ ճիշտ կողային $(\Delta(G) + 1)$-ներկմանը, որի դեպքում $S \in C_\beta(a) \cap C_\beta(b)$: Այնուհետև, ներկելով $G$ գրաֆի $ab$ կողը $S$ գույնով, կստանանք $G$ գրաֆի ճիշտ կողային $(\Delta(G) + 1)$-ներկում, որը հակասում է մեր սկզբնական ենթադրությանը: ∎

Թեորեմ 8.2.2-ը ցույց է տալիս, որ թեորեմ 8.2.3-ի գնահատականները հնարավոր չէ լավացնել:

Վիզինգի թեորեմը հնարավորություն է տալիս բոլոր գրաֆների բազմությունը տրոհել երկու ենթաբազմությունների:

**Սահմանում 8.2.4։** $G$ գրաֆը կոչվում է *առաջին դասի* գրաֆ, եթե $\chi'(G) = \Delta(G)$, հակառակ դեպքում՝ *երկրորդ դասի* գրաֆ:

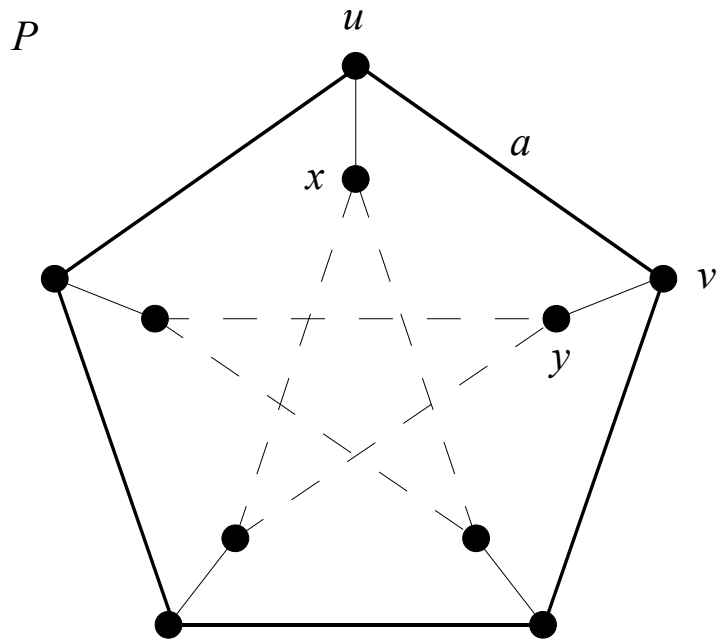

Նկ. 8.2.4

Այժմ դիտարկենք համասեռ գրաֆների քրոմատիկ ինդեքս գտնելու խնդիրը: Նկատենք, որ $G$ համասեռ գրաֆը ունի ճիշտ կողային $\Delta(G)$-ներկում այն և միայն այն դեպքում, երբ $G$-ն 1-ֆակտորիզացվող գրաֆ է: Այսպես, օրինակ, թեորեմ 8.2.2-ից հետևում է, որ լրիվ $K_{2n}$ գրաֆը 1-ֆակտորիզացվող է, իսկ $K_{2n+1}$-ը՝ 1-ֆակտորիզացվող չէ: Ինչպես նշել ենք, եթե $G$-ն կենտ երկարություն ունեցող պարզ ցիկլ է, ապա $G$-ն երկրորդ դասից է: Խորանարդ գրաֆների դեպքում այդպիսի օրինակ է հանդիսանում Պետերսենի $P$ գրաֆը (նկ. 8.2.4):

**Թեորեմ 8.2.4:** Եթե $P$-ն Պետերսենի գրաֆ է, ապա $\chi'(P) = 4$:

**Ապացույց:** Նախ նկատենք, որ Պետերսենի գրաֆը կարելի է պատկերել հինգ երկարություն ունեցող երկու պարզ ցիկլերի և այդ ցիկլերի գագաթները միացնող կատարյալ զուգակցման միջոցով: Նկ. 8.2.4-ում պատկերված հաստեցված գծերով նշված պարզ ցիկլը կանվանենք արտաքին ցիկլ, իսկ կետագծերով նշված պարզ ցիկլը՝ ներքին ցիկլ:

Ենթադրենք հակառակը. $P$-ն ունի ճիշտ կողային 3-ներկում: Դիցուք այդ ներկումն $\alpha$-ն է: Քանի որ արտաքին և ներքին պարզ ցիկլերը կենտ են, ուստի այդ ցիկլերի կողերը ներկելու համար անհրաժեշտ է երեք գույն: Դիտարկենք արտաքին ցիկլի $uv$ կողը: Դիցուք $\alpha(uv) = a$: Քանի որ խորանարդ գրաֆների ճիշտ կողային 3-ներկման դեպքում յուրաքանչյուր գույն ներկա է այդ գրաֆի ցանկացած գագաթում, ուստի $\alpha(ux) \neq a$ և $\alpha(vy) \neq a$ (նկ. 8.2.4): Մյուս կողմից, քանի որ $xy \notin E(G)$, ուստի $a$ գույնը ներկա է ներքին պարզ ցիկլի երկու տարբեր կողերի վրա, որոնք կից են $x$ և $y$ գագաթներին: Այստեղից, հաշվի առնելով, որ արտաքին ցիկլի վրա ներկա են երեք տարբեր գույներ, ստանում ենք, որ ներքին ցիկլի երկարությունը առնվազն վեց է, ինչը հակասություն է: Այսպիսով, $\chi'(P) \geq 4$: Մյուս կողմից, համաձայն թեորեմ 8.2.3-ի $\chi'(P) \leq 4$: ∎

Ստորև ցույց կտրվի, որ կենտ քանակությամբ գագաթներ ունեցող բոլոր համասեռ գրաֆները երկրորդ դասից են:

**Թեորեմ 8.2.5 (Վ. Վիզինգ):** Եթե $G$-ն $r$-համասեռ ($r \in \mathbb{N}$) գրաֆ է և $|V(G)|$-ն կենտ է, ապա $\chi'(G) = r + 1$:

**Ապացույց:** Դիցուք $|V(G)| = n$ և $n \geq 3$: Ցույց տանք, որ $\chi'(G) \geq r + 1$: Դիտարկենք $G$ գրաֆի ճիշտ կողային $\chi'(G)$-ներկում: Պարզ է, որ

$$E(G) = E_1 \cup E_2 \cup \cdots \cup E_{\chi'(G)} \text{ և } E_i \cap E_j = \emptyset, \text{ եթե } 1 \leq i \neq j \leq \chi'(G),$$

որտեղ $E_i$-ն $i$-րդ գույնով ներկված կողերի բազմությունն է ($1 \leq i \leq \chi'(G)$): Նկատենք, որ յուրաքանչյուր $i$-ի համար $E_i$-ն զուգակցում է ($1 \leq i \leq \chi'(G)$): Այստեղից, հաշվի առնելով, որ $n$-ը կենտ է, ստանում ենք, որ ցանկացած $i$-ի համար ($1 \leq i \leq \chi'(G)$) ստույգ է $|E_i| \leq \alpha'(G) \leq \frac{n-1}{2}$ անհավասարությունը: Այսպիսով,

$$\frac{r \cdot n}{2} = |E(G)| = \sum_{i=1}^{\chi'(G)} |E_i| \leq \chi'(G) \cdot \alpha'(G) \leq \chi'(G) \cdot \frac{n-1}{2},$$

ուստի $\chi'(G) \geq \frac{r \cdot n}{n-1} = r + \frac{r}{n-1} > r$: Մյուս կողմից, համաձայն թեորեմ 8.2.3-ի, $\chi'(G) \leq r + 1$: ∎

**Սահմանում 8.2.5:** $G$ գրաֆը կոչվում է *գերցված*, եթե $|E(G)| > \left\lfloor \frac{n}{2} \right\rfloor \cdot \Delta(G)$:

Նկատենք, որ գերցված գրաֆները ունեն կենտ քանակությամբ գագաթներ և երկրորդ դասից են: Գերցված գրաֆների հետ է կապված Չետվինդի և Հիլտոնի հանրահայտ հիպոթեզը, որը ձևակերպված է ստորև:

**Հիպոթեզ 8.2.1:** Դիցուք $G$-ն $n$ գագաթ պարունակող գրաֆ է, որի համար $\Delta(G) \geq \frac{n}{3}$: Այդ դեպքում $G$-ն երկրորդ դասից է այն և միայն այն դեպքում, երբ $G$-ն ունի $H$ գերցված ենթագրաֆ, որի համար $\Delta(H) = \Delta(G)$:

Հայտնի է, որ այս հիպոթեզից հետևում է Գլուխ 5-ում բերված 1-ֆակտորիզացիայի հիպոթեզը (հիպոթեզ 5.4.1):

Նշենք նաև առանց ապացույցի, որ հայտնի է համարյա բոլոր գրաֆների քրոմատիկ ինդեքսի արժեքը:

**Թեորեմ 8.2.6 (Պ. Էրդյոշ, Ռ. Վիլսոն):** Համարյա բոլոր գրաֆները առաջին դասից են:

Այժմ անդրադառնանք հարթ գրաֆների քրոմատիկ ինդեքս գտնելու խնդրին: Այդ խնդիրը հետազոտվել է Վիզինգի [44] կողմից, որը ցույց է տվել, որ բոլոր հարթ $G$ գրաֆները, որոնց համար $\Delta(G) \geq 10$, առաջին դասից են: Հետագայում նա ուժեղացրեց այդ արդյունքը և ապացուցեց, որ բոլոր հարթ $G$ գրաֆները, որոնց համար $\Delta(G) \geq 8$, առաջին դասից են: Մյուս կողմից, հեշտ է տեսնել, որ գոյություն ունեն հարթ $G$ գրաֆներ, որոնց համար $2 \leq \Delta(G) \leq 5$ և որոնք երկրորդ դասից են: Իրոք, դրա համար բավական է դիտարկել կենտ երկարություն ունեցող պարզ ցիկլը և այն գրաֆները, որոնք ստացվում են նկ. 7.1.6-ում պատկերված 3-, 4- և 5-համասեռ հարթ գրաֆներից ճիշտ մեկ կողը տրոհելով և համոզվել, որ բոլոր այդ գրաֆները գերլցված են: Այսպիսով, մենք գալիս ենք հարթ գրաֆների քրոմատիկ ինդեքսի մասին Վիզինգի հանրահայտ հիպոթեզի ձևակերպմանը:

**Հիպոթեզ 8.2.2 (Վ. Վիզինգ):** Եթե $G$-ն հարթ գրաֆ է, որի համար $6 \leq \Delta(G) \leq 7$, ապա $\chi'(G) = \Delta(G)$:

Այս հիպոթեզի $\Delta(G) = 7$ դեպքը հաստատվել է 2000 թվականին Ժանգի [40] կողմից, իսկ $\Delta(G) = 6$ դեպքը մնում է բաց:

Այժմ ապացուցենք հարթ գրաֆների քրոմատիկ ինդեքսի մասին Վիզինգի թեորեմներից մեկը:

**Թեորեմ 8.2.7 (Վ. Վիզինգ):** Եթե $G$-ն հարթ գրաֆ է, որի համար $\Delta(G) \geq 10$, ապա $\chi'(G) = \Delta(G)$:

**Ապացույց:** Ենթադրենք հակառակը. գոյություն ունեն $H$ հարթ գրաֆներ, որոնց համար $\Delta(H) \geq 10$ և որոնք չունեն ճիշտ կողային $\Delta(H)$-ներկում: Ընտրենք այդ գրաֆներից նվազագույն քանակությամբ կողեր ունեցող $G = (V, E)$ հարթ գրաֆը:

Դիցուք $U = \{u : u \in V(G) \text{ և } d_G(u) \leq 5\}$: Քանի որ $G$-ն հարթ գրաֆ է, ուստի, համաձայն հետևանք 7.1.3-ի, $U \neq \emptyset$: Քանի որ $G[V \backslash U]$-ը ևս հարթ գրաֆ է, ուստի, համաձայն հետևանք 7.1.3-ի, $G$ գրաֆում գոյություն ունի $a \in V \backslash U$ գագաթ, որի համար $d_{G[V \backslash U]}(a) \leq 5$: Մյուս կողմից, քանի որ $a \notin U$, ուստի $G$ գրաֆում գոյություն ունի այնպիսի $b \in U$, որ $ab \in E(G)$: Դիտարկենք $G' = G - ab$ գրաֆը: $G$ գրաֆի ընտրությունից հետևում է, որ $G'$ գրաֆը ունի ճիշտ կողային $\Delta(G)$-ներկում, քանի որ $\Delta(G') \leq \Delta(G)$: Դիցուք այդ ճիշտ կողային $\Delta(G)$-ներկումն $\alpha$-ն է: $C_\alpha(v)$-ով նշանակենք այն գույների բազմությունը, որոնք բացակայում են $v$ գագաթում $\alpha$ ներկման դեպքում: Պարզ է, որ $|C_\alpha(a)| \geq 1$ և $|C_\alpha(b)| \geq \Delta(G) - 4$: Եթե $C_\alpha(a) \cap C_\alpha(b) \neq \emptyset$, ապա $G$ գրաֆի $ab$ կողը ներկենք կամայական $i_0 \in C_\alpha(a) \cap C_\alpha(b)$ գույնով, որը $G'$ գրաֆի $\alpha$ ներկման հետ միասին կորոշի $G$ գրաֆի ճիշտ կողային $\Delta(G)$-ներկում:

Այժմ ենթադրենք, որ $C_\alpha(a) \cap C_\alpha(b) = \emptyset$: Այստեղից հետևում է, որ $a$ գագաթին կից են առնվազն $\Delta(G) - 4$ կողեր, որոնք ներկված են $b$ գագաթում բացակայող գույների միջոցով: Դիցուք այդ կողերը $av_1, av_2, \ldots, av_k$-ն են, որտեղ $k \geq \Delta(G) - 4$: Քանի որ $k \geq \Delta(G) - 4 \geq 6$, ուստի $v_1, v_2, \ldots, v_k$ գագաթներից առնվազն մեկը կպատկանի $U$ բազմությանը: Որոշակիության համար ենթադրենք, որ $v_1 \in U$: Դիցուք $\alpha(av_1) = S_1$ և $S_1 \in C_\alpha(b)$: Քանի որ $\Delta(G) \geq 10$, $|C_\alpha(v_1)| \geq \Delta(G) - 5$ և $|C_\alpha(b)| \geq \Delta(G) - 4$, ուստի $C_\alpha(v_1) \cap C_\alpha(b) \neq \emptyset$: Դիցուք $S \in C_\alpha(v_1) \cap C_\alpha(b)$ և $t \in C_\alpha(a)$:

Դիտարկենք $G'$ գրաֆի $a, b$ և $v_1$ գագաթները: Պարզ է, որ $S \in C_\alpha(b)$, $t \in C_\alpha(a)$, $S \in C_\alpha(v_1)$ և $S \neq t$ (քանի որ $C_\alpha(a) \cap C_\alpha(b) = \emptyset$):

Դիտարկենք $G'$ գրաֆի $G_{St} = (V, E_{St})$ ենթագրաֆը, որտեղ $E_{St}$-ն $G'$ գրաֆի այն կողերի բազմությունն է, որոնք $\alpha$ ներկման ժամանակ ներկվել են $S$ կամ $t$ գույներով: Քանի որ $\Delta(G_{St}) \leq 2$, իսկ $d_{G_{St}}(a) = d_{G_{St}}(b) = d_{G_{St}}(v_1) = 1$, ուստի պարզ է, որ $a, b$ և $v_1$ գագաթները միևնույն կապակցված բաղադրիչին չեն պատկանում: Դիտարկենք երկու դեպք:

Դեպք 1: $a$ և $b$ գագաթները պատկանում են $G_{St}$ գրաֆի տարբեր կապակցված բաղադրիչներին:

Այդ դեպքում, վերաներկելով $G_{St}$ գրաֆի $a$ գագաթը պարունակող կապակցված բաղադրիչում $S$ գույնով ներկված կողերը $t$ գույնով, իսկ $t$-ով ներկված կողերը՝ $S$ գույնով, կստանանք $G'$ գրաֆի $\alpha'$ ճիշտ կողային $\Delta(G)$-ներկում, որի դեպքում $S \in C_{\alpha'}(a) \cap C_{\alpha'}(b)$, իսկ դա հակասում է $\chi'(G) = \Delta(G) + 1$ պայմանին:

Դեպք 2: $a$ և $b$ գագաթները պատկանում են $G_{St}$ գրաֆի միևնույն կապակցված բաղադրիչին:

Այդ դեպքում $v_1$ գագաթը չի պատկանի այդ բաղադրիչին և, վերաներկելով $G_{St}$ գրաֆի $v_1$ գագաթը պարունակող կապակցված բաղադրիչում $S$ գույնով ներկված կողերը $t$ գույնով, իսկ $t$-ով ներկված կողերը՝ $S$ գույնով, կստանանք $G'$ գրաֆի $\alpha'$ ճիշտ կողային $\Delta(G)$-ներկում, որի դեպքում $t \in C_{\alpha'}(v_1)$: Այնուհետև, ներկենք $av_1$ կողը $t$ գույնով և $ab$ կողը՝ $S_1$ գույնով, արդյունքում կստանանք $G$ գրաֆի ճիշտ կողային $\Delta(G)$-ներկում, իսկ դա հակասում է $\chi'(G) = \Delta(G) + 1$ պայմանին: ∎

Այս պարագրաֆի վերջում նշենք նաև, որ հատուկ հետքրքրություն է ներկայացնում մուլտիգրաֆների քրոմատիկ ինդեքս գտնելու խնդիրը, որի մասին կարելի է մանրամասն ծանոթանալ [34] գրքում։

## § 8.3. Գրաֆների տոտալ ներկումներ

Դիցուք $G = (V, E)$-ն գրաֆ է։

**Սահմանում 8.3.1:** $G$ գրաֆի *տոտալ k-ներկում* կոչվում է $\alpha: V(G) \cup E(G) \to \{1, \ldots, k\}$ արտապատկերումը, իսկ՝ $1, \ldots, k$ թվերը կոչվում են *գույներ*։

**Սահմանում 8.3.2:** $G$ գրաֆի $\alpha$ տոտալ $k$-ներկումը կոչվում է *ճիշտ տոտալ k-ներկում*, եթե ցանկացած $uv \in E(G)$-ի համար ստույգ է $\alpha(u) \neq \alpha(v)$ պայմանը, ցանկացած $e, e' \in E(G)$-ի հարևան կողերի համար ստույգ է $\alpha(e) \neq \alpha(e')$ պայմանը և ցանկացած $v \in V(G)$-ի և նրան կից $e \in E(G)$-ի կողի համար ստույգ է $\alpha(v) \neq \alpha(e)$ պայմանը։ Այլ կերպ ասած, ճիշտ տոտալ ներկումն այնպիսի ներկում է, որի դեպքում հարևան գագաթները և կողերը ներկվում են տարբեր գույներով և ցանկացած գագաթ և նրան կից կող ևս ներկվում են տարբեր գույներով։

**Սահմանում 8.3.3:** $G$ գրաֆը կոչվում է *տոտալ k-ներկելի*, եթե գոյություն ունի $G$ գրաֆի ճիշտ տոտալ $k$-ներկում։ Այն նվազագույն $k$-ն, որի դեպքում $G$-ն տոտալ $k$-ներկելի է կոչվում է $G$ գրաֆի *տոտալ քրոմատիկ թիվ*։ $\chi''(G)$-ով նշանակենք $G$ գրաֆի տոտալ քրոմատիկ թիվը։

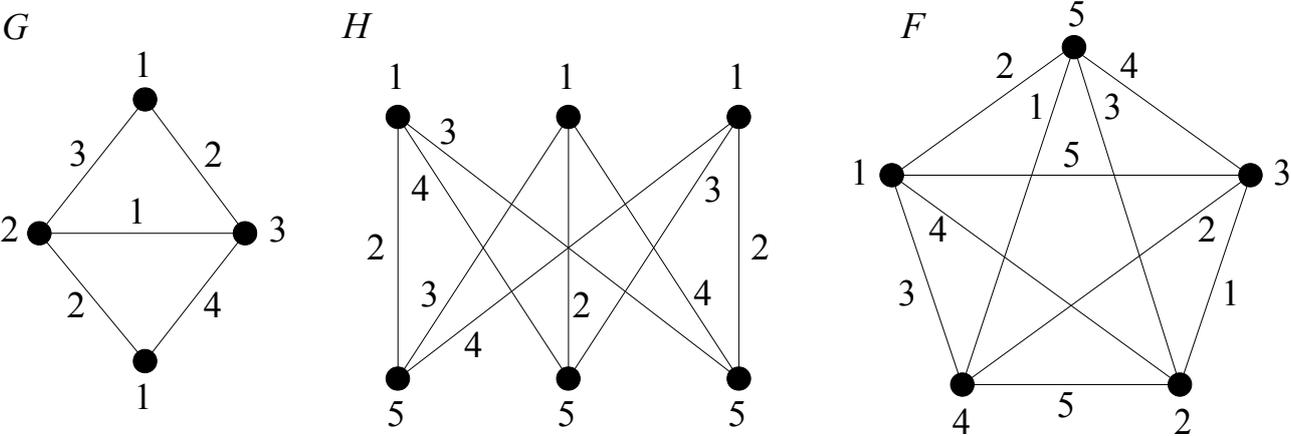

Նկ. 8.3.1

Դիտարկենք նկ. 8.3.1-ում պատկերված $G, H$ և $F$ գրաֆները։ Հեշտ է տեսնել, որ այդ նկարում բերված է $G$ գրաֆի տոտալ 4-ներկում, որը, սակայն, ճիշտ տոտալ 4-ներկում չէ։

Մյուս կողմից, հեշտ է տեսնել նաև, որ նկ. 8.3.1-ում պատկերված են $H$ և $F$ գրաֆների ճիշտ տոտալ 5-ներկումներ։ Ավելին, դժվար չէ համոզվել, որ $\chi''(H) = \chi''(F) = 5$։

Նկատենք, որ ցանկացած $G$ գրաֆի համար $\chi''(G) \geq \Delta(G) + 1$։ Իրոք, քանի որ $G$ գրաֆի ճիշտ տոտալ ներկման դեպքում ցանկացած գագաթին կից կողերը և այդ գագաթը ներկված են զույգ առ զույգ տարբեր գույներով, ուստի այդ ներկման դեպքում օգտագործվող գույների քանակը չի կարող փոքր լինել $(\Delta(G) + 1)$-ից։ Մյուս կողմից, ցանկացած $G = (V, E)$ գրաֆի համար, որտեղ $V = \{v_1, v_2, \ldots, v_n\}$, մենք կարող ենք դիտարկել § 1.5-ում սահմանված $\Omega(V, \{\{v_1\} \cup \partial_G(v_1), \{v_2\} \cup \partial_G(v_2), \ldots, \{v_n\} \cup \partial_G(v_n)\} \cup E)$ հատումների գրաֆը։ Տանք այդ տոտալ գրաֆի սահմանումը։ $G$ գրաֆի տոտալ $T(G)$ գրաֆը սահմանվում է հետևյալ կերպ. $V(T(G)) = V(G) \cup E(G)$ և $E(T(G)) = E(G) \cup E(L(G)) \cup \{ue : u \in V(G), e \in E(G)$ և $u-$ն կից է $e-$ին$\}$։ Նկ. 8.3.2-ում պատկերված է $G$ գրաֆը և նրա տոտալ $T(G)$ գրաֆի օրինակը․

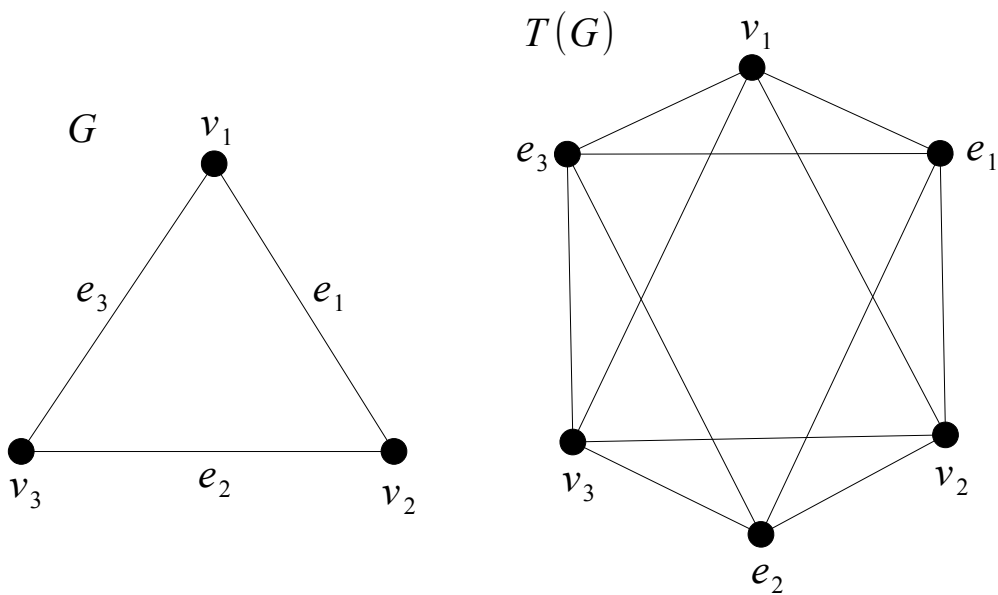

Նկ. 8.3.2

Հեշտ է տեսնել, որ ցանկացած $G$ գրաֆի համար տեղի ունի $\chi''(G) = \chi(T(G))$ հավասարությունը։ Այստեղից, հաշվի առնելով $\Delta(T(G)) \leq 2\Delta(G)$ անհավասարությունը, և համաձայն թեորեմ 8.1.1-ի, ստանում ենք հետևյալը․

$$\chi''(G) = \chi(T(G)) \leq \Delta(T(G)) + 1 \leq 2\Delta(G) + 1:$$

Նշենք գրաֆների տոտալ քրոմատիկ թվի համար ևս մի պարզ վերին գնահատական. ցանկացած $G$ գրաֆի համար $\chi''(G) \leq \chi(G) + \chi'(G)$։ Իրոք, ներկենք սկզբից $G$ գրաֆի գագաթները $1, \ldots, \chi(G)$ գույներով, այնպես որ հարևան գագաթները ներկվեն տարբեր գույներով, այնուհետև ներկենք այդ գրաֆի կողերը $\chi(G) + 1, \ldots, \chi(G) +$ $\chi'(G)$ գույներով, այնպես որ հարևան կողերը ներկվեն տարբեր գույներով։ Պարզ է, որ արդյունքում կստանանք $G$ գրաֆի ճիշտ տոտալ ներկում, որի դեպքում օգտագործվող գույների քանակը կլինի հավասար $(\chi(G) + \chi'(G))$-ին։ Այդ վերին գնահատականը հասանելի է, օրինակ, լրիվ երկկողմանի $K_{3,3}$ գրաֆի համար, քանի որ $\chi(K_{3,3}) = 2$ և $\chi'(K_{3,3}) = 3$, իսկ $\chi''(K_{3,3}) = 5$։ Պարզվում է, ավելի ընդհանուր փաստ տեղի ունի։

**Թեորեմ 8.3.1 (Բեհզադ, Չարտրանդ, Կուպեր):** Եթե $G$-ն ունի առնվազն երկու գագաթ և $\chi''(G) = \chi(G) + \chi'(G)$, ապա $G$-ն երկկողմանի գրաֆ է։

**Ապացույց:** Նախ նկատենք, որ եթե $G$-ն չունի կող, ապա $\chi'(G) = 0$ և $\chi''(G) = \chi(G) = 1$, իսկ $G$-ն ակնհայտաբար, երկկողմանի է։ Ենթադրենք, որ $G$-ն ունի առնվազն մեկ կող և $\chi''(G) = \chi(G) + \chi'(G)$, բայց $G$-ն երկկողմանի չէ։ Այդ դեպքում պարզ է, որ $\chi(G) \geq 3$։

Դիտարկենք $G$ գրաֆի գագաթների և կողերի բազմությունների հետևյալ տրոհումները.

$$V(G) = V_1 \cup V_2 \cup \cdots \cup V_{\chi(G)} \text{ և } V_i \cap V_j = \emptyset, \text{ երբ } 1 \leq i \neq j \leq \chi(G), \text{ և}$$

$$E(G) = E_1 \cup E_2 \cup \cdots \cup E_{\chi'(G)} \text{ և } E_s \cap E_t = \emptyset, \text{ երբ } 1 \leq s \neq t \leq \chi'(G),$$

որտեղ $V_i$-ն անկախ բազմություն է ($1 \leq i \leq \chi(G)$), իսկ $E_j$-ն՝ զուգակցում է ($1 \leq j \leq \chi'(G)$)։ Այստեղից հետևում է, որ $\left(\cup_{i=1}^{\chi(G)} V_i\right) \cup \left(\cup_{j=1}^{\chi'(G)} E_j\right)$-ն հանդիսանում է $V(G) \cup E(G)$ բազմության տրոհում անկախ գագաթների և կողերի բազմությունների։ Քանի որ $\chi(G) \geq 3$, ուստի $E_1$ զուգակցման յուրաքանչյուր կողի համար միշտ կգտնվի այնպիսի $V_{i_1}$ անկախ բազմություն, որի գագաթները կից չեն այդ կողին։ Այստեղից հետևում է, որ ավելացնելով $E_1$ զուգակցման կողերը նրանց համապատասխան $V_{i_1}$ անկախ բազմություններին, մենք կստանանք նոր գագաթների և կողերի $V_1', V_2', \ldots, V_{\chi(G)}'$ բազմություններ, որտեղ $V_i'$-ի գագաթները զույգ առ զույգ հարևան չեն, կողերը զույգ առ զույգ հարևան չեն, իսկ գագաթները և կողերը՝ կից չեն ($1 \leq i \leq \chi(G)$)։ Այսպիսով, դեն նետելով $E_1$ զուգակցումը $G$ գրաֆի կողերի բազմության տրոհումից, մենք կստանանք այդ գրաֆի ճիշտ տոտալ $(\chi(G) + \chi'(G) - 1)$-ներկում հետևյալ եղանակով. $V_i'$-րդ բազմության գագաթները և կողերը ներկենք $i$-րդ գույնով ($1 \leq i \leq \chi(G)$), իսկ $E_j$-րդ զուգակցման կողերը ներկենք $(\chi(G) + j - 1)$-րդ գույնով ($2 \leq j \leq \chi'(G)$)։ Այստեղից հետևում է, որ $\chi''(G) < \chi(G) + \chi'(G)$, իսկ դա հակասում է թեորեմի պայմանին։ ∎

Այժմ ցույց տանք, որ ցանկացած երկկողմանի $G$ գրաֆի համար ստույգ է $\chi''(G) \leq \Delta(G) + 2$ անհավասարությունը։

**Թեորեմ 8.3.2:** Եթե $G$-ն երկկողմանի գրաֆ է, ապա $\chi''(G) \leq \Delta(G) + 2$:

**Ապացույց:** Դիցուք $V(G) = V_1 \cup V_2$ երկկողմանի $G$ գրաֆի գագաթների բազմության համապատասխան տրոհումն է: Համաձայն թեորեմ 8.2.1-ի, երկկողմանի $G$ գրաֆն ունի ճիշտ կողային $\Delta(G)$-ներկում: Դիցուք այդ ճիշտ կողային $\Delta(G)$-ներկումն $\alpha$-ն է: Սահմանենք $G$ գրաֆի տոտալ $\beta$ ներկումը հետևյալ կերպ.

1. ցանկացած $u \in V_1$-ի համար $\beta(u) = \Delta(G) + 1$,
2. ցանկացած $v \in V_2$-ի համար $\beta(v) = \Delta(G) + 2$,
3. ցանկացած $e \in E(G)$-ի համար $\beta(e) = \alpha(e)$:

Հեշտ է տեսնել, որ $\beta$-ն հանդիսանում է $G$ գրաֆի ճիշտ տոտալ $(\Delta(G) + 2)$-ներկում, ուստի $\chi''(G) \leq \Delta(G) + 2$: ∎

Պարզվում է, տոտալ քրոմատիկ թվի ճշգրիտ արժեքը հայտնի է լրիվ և լրիվ երկկողմանի գրաֆների դեպքում:

**Թեորեմ 8.3.3 (Բեհզադ, Չարտրանդ, Կուպեր):** Ցանկացած $n \in \mathbb{N}$-ի համար տեղի ունի

$$\chi''(K_n) = \begin{cases} n + 1, & \text{եթե } n-ը \text{ զույգ է}, \\ n, & \text{եթե } n-ը \text{ կենտ է}, \end{cases}$$

հավասարությունը:

**Ապացույց:** Ինչպես նշել ենք, տեղի ունի $\chi''(K_n) \geq \Delta(K_n) + 1 = n$ անհավասարությունը:

Նախ դիտարկենք $n$-ի կենտ թիվ լինելու դեպքը: Համաձայն թեորեմ 8.2.2-ի, $K_n$ գրաֆն ունի ճիշտ կողային $n$-ներկում: Դիցուք այդ ճիշտ կողային $n$-ներկումն $\alpha$-ն է: Նկատենք, որ $\alpha$ ներկման դեպքում ցանկացած գագաթում բացակայում է $1, 2, \ldots, n$ գույներից ճիշտ մեկը, ընդ որում տարբեր գագաթներում բացակայող գույները զույգ առ զույգ տարբեր են: Ներկենք $K_n$ գրաֆի յուրաքանչյուր գագաթ այն գույնով, որը բացակայում է այդ գագաթում $\alpha$ ներկման դեպքում: Հեշտ է տեսնել, որ այդ ճիշտ գագաթային $n$-ներկումը $\alpha$ ներկման հետ միասին կորոշի $K_n$ գրաֆի ճիշտ տոտալ $n$-ներկում: Այսպիսով, $\chi''(K_n) \leq n$ և, հետևաբար, $\chi''(K_n) = n$:

Այժմ դիտարկենք $n$-ի զույգ թիվ լինելու դեպքը: Պարզ է, որ $|V(K_n)| + |E(K_n)| = n + \binom{n}{2} = \frac{n(n+1)}{2}$: Քանի որ $K_n$ գրաֆի ցանկացած ճիշտ տոտալ ներկման դեպքում յուրաքանչյուր գույնով կարող է ներկված լինել ամենաշատը մեկ գագաթ, ուստի միևնույն գույնով ներկված գագաթների և կողերի քանակը չի գերազանցում $\frac{n}{2}$-ը: Այստեղից և $|V(K_n)| + |E(K_n)| = \frac{n(n+1)}{2}$ հավասարությունից հետևում է, որ $\chi''(K_n) \geq n + 1$: Դիտարկենք $K_{n+1}$ գրաֆը: Քանի որ $K_n$ գրաֆը ստացվում է $K_{n+1}$ գրաֆից մեկ գագաթ հեռացնելով, ուստի $\chi''(K_n) \leq \chi''(K_{n+1}) = n + 1$: Հետևաբար $\chi''(K_n) = n + 1$: ∎

**Թեորեմ 8.3.4 (Բեհզադ, Չարտրանդ, Կուպեր):** Ցանկացած $m, n \in \mathbb{N}$-ի համար տեղի ունի

$$\chi''(K_{m,n}) = \begin{cases} \max\{m, n\} + 1, & \text{եթե } m \neq n, \\ n + 2, & \text{եթե } m = n, \end{cases}$$

հավասարությունը:

**Ապացույց:** Ինչպես նշել ենք, տեղի ունի $\chi''(K_{m,n}) \geq \Delta(K_{m,n}) + 1 = \max\{m, n\} + 1$ անհավասարությունը:

Նախ դիտարկենք $m \neq n$ դեպքը: Առանց ընդհանրությունը խախտելու կարող ենք ենթադրել, որ $m < n$: Դիցուք $V(K_{m,n}) = \{u_1, u_2, \ldots, u_m, v_1, v_2, \ldots, v_n\}$ և $E(K_{m,n}) = \{u_i v_j : 1 \leq i \leq m, 1 \leq j \leq n\}$:

Սահմանենք $K_{m,n}$ գրաֆի կողային $\alpha$ ներկումը հետևյալ կերպ.

$$\alpha(u_i v_j) = \begin{cases} (i + j - 1) \bmod n, & \text{եթե } i + j \neq n + 1, \\ n, & \text{եթե } i + j = n + 1, \end{cases}$$

որտեղ $1 \leq i \leq m, 1 \leq j \leq n$:

Դժվար չէ համոզվել, որ $\alpha$-ն հանդիսանում է $K_{m,n}$ գրաֆի ճիշտ կողային $n$-ներկում: $C_\alpha(v)$-ով նշանակենք այն գույների բազմությունը, որոնք բացակայում են $v$ գագաթում $\alpha$ ներկման դեպքում: Քանի որ $m < n$, ուստի յուրաքանչյուր $j$-ի համար գոյություն կունենա այնպիսի $c_j$ գույն, որ $c_j \in C_\alpha(v_j)$ ($1 \leq j \leq n$): Այժմ սահմանենք $K_{m,n}$ գրաֆի տոտալ $\beta$ ներկումը հետևյալ կերպ.

1. ցանկացած $i$-ի համար $\beta(u_i) = n + 1$, որտեղ $1 \leq i \leq m$,
2. ցանկացած $j$-ի համար $\beta(v_j) = c_j$, որտեղ $1 \leq j \leq n$,
3. ցանկացած $e \in E(G)$-ի համար $\beta(e) = \alpha(e)$:

Հեշտ է տեսնել, որ $\beta$-ն հանդիսանում է $G$ գրաֆի ճիշտ տոտալ $(n+1)$-ներկում, ուստի $\chi''(K_{m,n}) \leq n + 1$: Հետևաբար, $\chi''(K_{m,n}) = n + 1$:

Այժմ դիտարկենք $m = n$ դեպքը: Պարզ է, որ $|V(K_{n,n})| + |E(K_{n,n})| = 2n + n^2 = n(n + 2)$: Հեշտ է տեսնել, որ $K_{n,n}$ գրաֆի ցանկացած ճիշտ տոտալ ներկման դեպքում միևնույն գույնով ներկված գագաթների և կողերի քանակը չի գերազանցում $n$-ը: Այստեղից և $|V(K_{n,n})| + |E(K_{n,n})| = n(n + 2)$ հավասարությունից հետևում է, որ $\chi''(K_{n,n}) \geq n+2$: Մյուս կողմից, համաձայն թեորեմ 8.3.2-ի, $\chi''(K_{n,n}) \leq n+2$, ուստի $\chi''(K_{n,n}) = n+2$: ∎

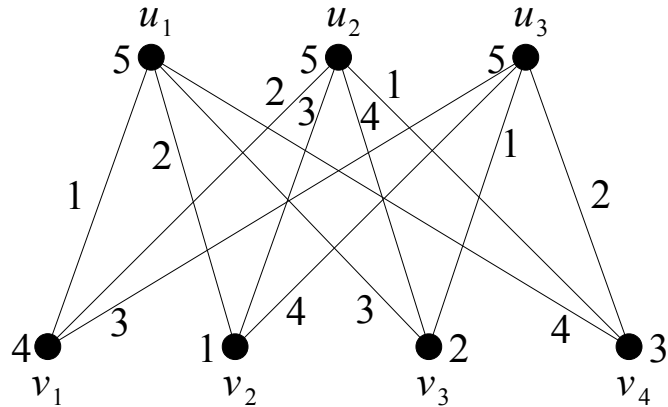

Նկ. 8.3.3

Նկ. 8.3.3-ում պատկերված է թեորեմ 8.3.4-ի ապացույցում բերված $K_{3,4}$ գրաֆի $\beta$ ճիշտ տոտալ 5-ներկումը:

Ինչպես նշել ենք, ցանկացած $G$ գրաֆի տոտալ քրոմատիկ թիվը բավարարում է $\chi''(G) \leq 2\Delta(G)+1$ անհավասարությանը: Մյուս կողմից, մեր բոլոր ապացուցված թեորեմներում այդ թիվը չէր գերազանցում գրաֆի առավելագույն աստիճան գումարած երկուս, ավելին, այդպիսի գրաֆի օրինակ, որի տոտալ քրոմատիկ թիվը մեծ է նրա առավելագույն աստիճան գումարած երկուսից, մինչև այժմ հայտնի չէ: Հաշվի առնելով այդ փաստը, Բեհզադը և Վիզինգը 1965 թվականին ձևակերպեցին նրանց հանրահայտ հիպոթեզը:

**Հիպոթեզ 8.3.1 (Բեհզադ, Վիզինգ):** Կամայական $G$ գրաֆի համար տեղի ունի
$$\Delta(G) + 1 \leq \chi''(G) \leq \Delta(G) + 2$$
անհավասարությունը:

Հայտնի է, որ այս հիպոթեզը, բաց մնալով ընդհանուր դեպքում, ճիշտ է մի շարք գրաֆների դասերի համար: Նշենք դրանցից մի քանիսը:

**Թեորեմ 8.3.5 (Ռոզենֆելդ, Վիժայադիտյա):** Եթե $G$ գրաֆի համար տեղի ունի $\Delta(G) \leq 3$ պայմանը, ապա $\chi''(G) \leq 5$:

**Թեորեմ 8.3.6 (Կոստոչկա):** Եթե $G$ գրաֆի համար տեղի ունի $\Delta(G) \leq 4$ պայմանը, ապա $\chi''(G) \leq 6$:

**Թեորեմ 8.3.7 (Կոստոչկա):** Եթե $G$ գրաֆի համար տեղի ունի $\Delta(G) \leq 5$ պայմանը,

ապա $\chi''(G) \leq 7$:

**Թեորեմ 8.3.8 (Յափ, Վանգ, Ժանգ):** Եթե $n$ գագաթ ունեցող $G$ գրաֆի համար տեղի ունի $\Delta(G) \geq n - 5$ պայմանը, ապա $\chi''(G) \leq \Delta(G) + 2$:

**Թեորեմ 8.3.9 (Հիլտոն, Հինդ):** Եթե $n$ գագաթ ունեցող $G$ գրաֆի համար տեղի ունի $\Delta(G) \geq \frac{3n}{4}$ պայմանը, ապա $\chi''(G) \leq \Delta(G) + 2$:

Հայտնի են նաև այդ հիպոթեզի ապացույցի ուղղությամբ որոշ արդյունքներ, որոնք լավացնում են տոտալ քրոմատիկ թվի վերին գնահատականները: Նշենք դրանցից երկուսը:

**Թեորեմ 8.3.10 (Կոստոչկա):** Եթե $G$ գրաֆի համար տեղի ունի $\Delta(G) \geq 4$ պայմանը, ապա $\chi''(G) \leq \left\lfloor \frac{3}{2}\Delta(G) \right\rfloor$:

**Թեորեմ 8.3.11 (Մոլլոյ, Ռիդ):** Կամայական $G$ գրաֆի համար տեղի ունի

$$\chi''(G) \leq \Delta(G) + 10^{26}$$

անհավասարությունը:

Այս պարագրաֆի վերջում անդրադառնանք նաև հարթ գրաֆների տոտալ քրոմատիկ թվի գտնելու խնդրին: Հայտնի է, որ եթե $G$ հարթ գրաֆի համար տեղի ունի $\Delta(G) \geq 8$ պայմանը, ապա $\chi''(G) \leq \Delta(G) + 2$: 1999 թվականին Սանդերսի և Ժաոյի [33] կողմից ապացուցվեց հետևյալ թեորեմը:

**Թեորեմ 8.3.12 (Սանդերս, Ժաո):** Եթե $G$ հարթ գրաֆի համար տեղի ունի $\Delta(G) = 7$ պայմանը, ապա $\chi''(G) \leq 9$:

Մյուս կողմից, թեորեմ 8.3.7-ից հետևում է, որ եթե $G$ հարթ գրաֆում $\Delta(G) \leq 5$, ապա $\chi''(G) \leq 7$: Այսպիսով, հարթ գրաֆների դեպքում հիպոթեզ 8.3.1-ը բաց է մնում միայն այն $G$ հարթ գրաֆների համար, որոնցում $\Delta(G) = 6$:

Գրաֆների տոտալ ներկումներին կարելի է ավելի մանրամասն ծանոթանալ [39] գրքում:

# ԳՐԱԿԱՆՈՒԹՅՈՒՆ

1. J. Akiyama, M. Kano, Factors and Factorizations of Graphs, (Proof Techniques in Factor Theory), Springer-Verlag Berlin Heidelbelg, 2011.
2. K. Appel, W. Haken, Every planar map is four colorable, Part I, Discharging, Illinois Journal of Mathematics 21, 1977, pp. 429-490.
3. K. Appel, W. Haken, J. Koch, Every planar map is four colorable, Part II, Reducibility, Illinois Journal of Mathematics 21, 1977, pp. 491-567.
4. A.S. Asratian, T.M.J. Denley, R. Haggkvist, Bipartite Graphs and their Applications, Cambridge University Press, Cambridge, 1998.
5. M. Behzad, G. Chartrand, J.K. Cooper Jr., The colour numbers of complete graphs, J. London Math. Soc. 42, 1967, pp. 226-228.
6. C. Berge, Graphs and Hypergraphs, North Holland, 1973.
7. B. Bollobas, Extremal Graph Theory, London Mathematical Society Monographs, Academic Press, London, 1978.
8. B. Bollobas, Modern Graph Theory, Springer, 1998.
9. J.A. Bondy, U.S.R. Murty, Graph Theory, Springer, 2008.
10. P.A. Catlin, Hajós's graph-colouring conjecture: variations and counterexamples, Journal of Combinatorial Theory B 26, 1979, pp. 268-274.
11. G. Chartrand, P. Zhang, Chromatic Graph Theory, Discrete Mathematics and Its Applications, CRC Press, 2009.
12. B. Chen, M. Matsumoto, J. Wang, Z. Zhang, J. Zhang, A short proof of Nash-Williams' theorem for arboricity of a graph, Graphs and Combinatorics 10, 1994, pp. 27-28.
13. G.A. Dirac, A property of **4**-chromatic graphs and some remarks on critical graphs, J. London Math. Soc. 27, 1952, pp. 85-92.
14. H. Fleischner, Eulerian Graphs and Related Topics, Part 1, Volume 1, Annals of Discrete Mathematic 45, Elsevier, North-Holland, Amsterdam, 1990.
15. M.R. Garey, D.S. Johnson, Crossing number is $NP$-complete, SIAM J. Alg. Discr. Meth. 4 (3), 1983, pp. 312-316.
16. H. Hadwiger, Über eine Klassifikation der Streckenkomplexe, Vierteljschr. Naturforsch. Ges. Zürich 88, 1943, pp. 133-143.
17. R. Hammack, W. Imrich, S. Klavzar, Handbook of Product Graphs, Second Edition, CRC Press, 2011.
18. P.J. Heawood, Map-colour theorem, Quarterly Journal of Mathematics, Oxford 24, 1890, pp. 332-338.
19. P. Hell, J. Nešetřil, Graph and Homomorphisms, Oxford University Press, New York, 2004.
20. T.R. Jensen, B. Toft, Graph Coloring Problems, Wiley Interscience Series in Discrete Mathematics and Optimization, 1995.
21. P.J. Kelly, A congruence theorem for trees, Pacific J. Math. 7, 1957, pp. 961-968.

22. A.B. Kempe, On the geographical problem of four colors, Amer. J. Math. 2, 1879, pp. 193-200.
23. D.J. Kleitman, The crossing number of $K_{5,n}$, Journal of Combinatorial Theory 9, 1971, pp. 315-323.
24. L. Lovasz, A note on the line reconstruction problem, Journal of Combinatorial Theory B 13, 1972, pp. 309-310.
25. L. Lovasz, Three short proofs in graph theory, Journal of Combinatorial Theory B 19, 1975, pp. 111-113.
26. L. Lovasz, M.D. Plummer, Matching Theory, Annals of Discrete Mathematic 29, North-Holland Publishing, 1986.
27. L.S. Melnikov, V.G. Vizing, New proof of Brooks theorem, Journal of Combinatorial Theory 7, 1969, pp. 289-290.
28. B. Mohar, C. Thomassen, Graphs on Surfaces, The Johns Hopkins University Press, 2001.
29. V. Muller, The edge-reconstruction hypothesis is true for graphs with more than $n \cdot \log_2 n$ edges, Journal of Combinatorial Theory B 22, 1977, pp. 281-283.
30. R. Naserasr, R. Škrekovski, The Petersen graph is not 3-edge-colorable - a new proof, Discrete Mathematics 268, 2003, pp. 325-326.
31. N. Robertson, D. P. Sanders, P. D. Seymour, R. Thomas, The four colour theorem, Journal of Combinatorial Theory B 70, 1997, pp. 2-44.
32. N. Robertson, P.D. Seymour, R. Thomas, Hadwiger's conjecture for $K_6$-free graphs, Combinatorica 13 (3), 1993, pp. 279-361.
33. D.P. Sanders, Y. Zhao, On total $9$-coloring planar graphs of maximum degree seven, Journal of Graph Theory 31, 1999, pp. 67-73.
34. M. Stiebitz, B. Toft, D. Scheide, L.M. Favrholdt, Graph edge colouring: Vizing's theorem and Goldberg's conjecture, Wiley Interscience, 2012.
35. C. Thomassen, The graph genus problem is $NP$-complete, J. of Algorithms 10 (4), 1989, pp. 568-576.
36. S. M. Ulam, A Collection of Mathematical Problems, Wiley, New York, 1960.
37. K. Wagner, Über eine Eigenschaft der ebenen Komplexe, Mathematische Annalen 114, 1937, pp. 570-590.
38. D.B. West, Introduction to Graph Theory, Prentice-Hall, New Jersey, 2001.
39. H.P. Yap, Total Colorings of Graphs, Lecture Notes in Mathematics 1623, Springer-Verlag, 1996.
40. L. Zhang, Every planar graph with maximum degree $7$ is of class $1$, Graphs and Combinatorics 16, 2000, pp. 467-495.
41. М. Айгнер, Комбинаторная теория, Пер. с англ.-М.: Мир, 1982.
42. К. Берж, Теория графов и ее применения, Пер. с франц.-М.: ИЛ, 1962.
43. В.Г. Визинг, Об оценке хроматического класса $p$-графа, Дискретный анализ 3, 1964, стр. 25-30.
44. В.Г. Визинг, Хроматический класс мультиграфа, Кибернетика 3, 1965, стр. 29-39.

45. В.А. Емеличев, О.И. Мельников, В.И. Сарванов, Р.И. Тышкевич, Лекции по теории графов, М.: Наука, 1990.
46. А.А. Зыков, Основы теории графов, М.: Наука, 1987.
47. Ф.А. Новиков, Дискретная математика для программистов, 3-е изд., СПб.: Питер, 2008.
48. О. Оре, Теория графов, Пер. с англ.-М.: Наука, 1980.
49. М. Свами, К. Тхуласираман, Графы, сети и алгоритмы, Пер. с англ.-М.: Мир, 1984.
50. Р. Уилсон, Введение в теорию графов, Пер. с англ.-М.: Мир, 1977.
51. Ф. Харари, Теория графов, Пер. с англ.-М.: Мир, 1973.
52. Ф. Харари, Э. Палмер, Перечисление графов, Пер. с англ.-М.: Мир, 1977.
53. Ի.Ա. Կարապետյան, Գրաֆների տեսություն (մեթոդական ցուցումներ), Երևան, ՀՊՃՀ, 2006.
54. Հ.Ց. Հակոբյան, Գրաֆների տեսության ներածություն (մեթոդական ցուցումներ), Երևան, ԵՊՃ, 1982.
55. Հ.Ց. Հակոբյան, Ա.Ս. Հասրաթյան, Գրաֆների տեսության խնդիրների ժողովածու, Երևան, ԵՊՃ, 1985.
56. Ժ.Գ. Նիկողոսյան, Դիսկրետ մաթեմատիկա, Գյումրիի տեղեկատվական տեխնոլոգիաների կենտրոն, Գյումրի, 2007.
57. Ռ.Ն. Տոնոյան, Դիսկրետ մաթեմատիկայի տարրերը, Երևան, ԵՊՃ, 1982.
58. Ռ.Ն. Տոնոյան, Դիսկրետ մաթեմատիկայի դասընթաց 1 (դասախոսություններ և առաջադրանքներ), Երևան, ԵՊՃ, 1997.
59. Ռ.Ն. Տոնոյան, Դիսկրետ մաթեմատիկայի դասընթաց, Երևան, ԵՊՃ, 1999.

# ԳՐԱՖՆԵՐԻ ՏԵՍՈՒԹՅՈՒՆ

ՈՒՍՈՒՄՆԱՄԵԹՈԴԱԿԱՆ ՁԵՌՆԱՐԿ